**МИНИСТЕРСТВО ОБРАЗОВАНИЯ И НАУКИ УКРАИНЫ**

**ЗАПОРОЖСКИЙ НАЦИОНАЛЬНЫЙ УНИВЕРСИТЕТ**

**С.В. Курапов**
**М.В. Давидовский**

# АЛГОРИТМИЧЕСКИЕ МЕТОДЫ КОНЕЧНЫХ ДИСКРЕТНЫХ СТРУКТУР

# ТОПОЛОГИЧЕСКИЙ РИСУНОК ГРАФА

## часть 1

**(монография)**

**Рецензенты:**
Доктор физико-математических наук, профессор
***В.А.Перепелица***
Доктор физико-математических наук, профессор
***Козин И.В.***

Утверждено Ученым советом ЗНУ (протокол № 2 от 30.10.2020 г.)

**Курапов С.В.**

**Давидовский М.В.**

Современные методы теории графов описывают граф с точностью до изоморфизма, что затрудняет создавать математические модели для визуализации рисунка графа на плоскости. Топологический рисунок плоской части графа позволяет описывать процесс планаризации алгебраическими методами, не производя никаких геометрических построений на плоскости. Получение вращения вершин графа сразу решает две важнейшие задачи теории графов: задачу проверки графа на планарность и задачу построения топологического рисунка плоского графа. Показано, что задачу построения рисунка непланарного графа можно свести к задаче построения рисунка плоского графа с учетом введения дополнительных вершин характеризующих пересечение рёбер. Естественно, что создание такой математической структуры позволит в будущем решать следующие задачи теории графов: проверка планарности графа, выделение максимально плоского суграфа, определение толщины графа, получение графа с минимальным количеством пересечений и т.д.

Для научных работников, преподавателей, студентов и аспирантов высших учебных заведений специализирующихся в области прикладной математике и информатики.

# Содержание

## Введение

Задачи структурного анализа и синтеза возникают при разработке практически любых объектов или систем на всех этапах, начиная с эскизного проектирования и заканчивая выпуском конструкторской документации. Успешное решение задач синтеза и анализа структуры сложных систем невозможно без их формализации. При решении таких задач в качестве аппарата формализации объектов проектирования широко применяется теория графов.

Очевидно, что решение задачи анализа и синтеза структуры систем тесно связано со смежными проблемами теории графов, таких как исследование планарности неориентированных графов, выделение максимальной планарной части, укладка графов на плоскости, определение изоморфизма графов, решение задачи автоморфизма, раскраски графов и других актуальных задач. Естественно, что при этом приходится так-то производить изображение графов. Таким образом, возникает необходимость в математическом описании самого рисунка графа в целях построения и передачи информации о структуре графа.

Что такое рисунок графа? Как его изобразить на плоскости? Ведь для заданной матрицы смежностей графа можно построить множество рисунков. Например, граф $K_5$ можно представить в виде следующего множества рисунков:

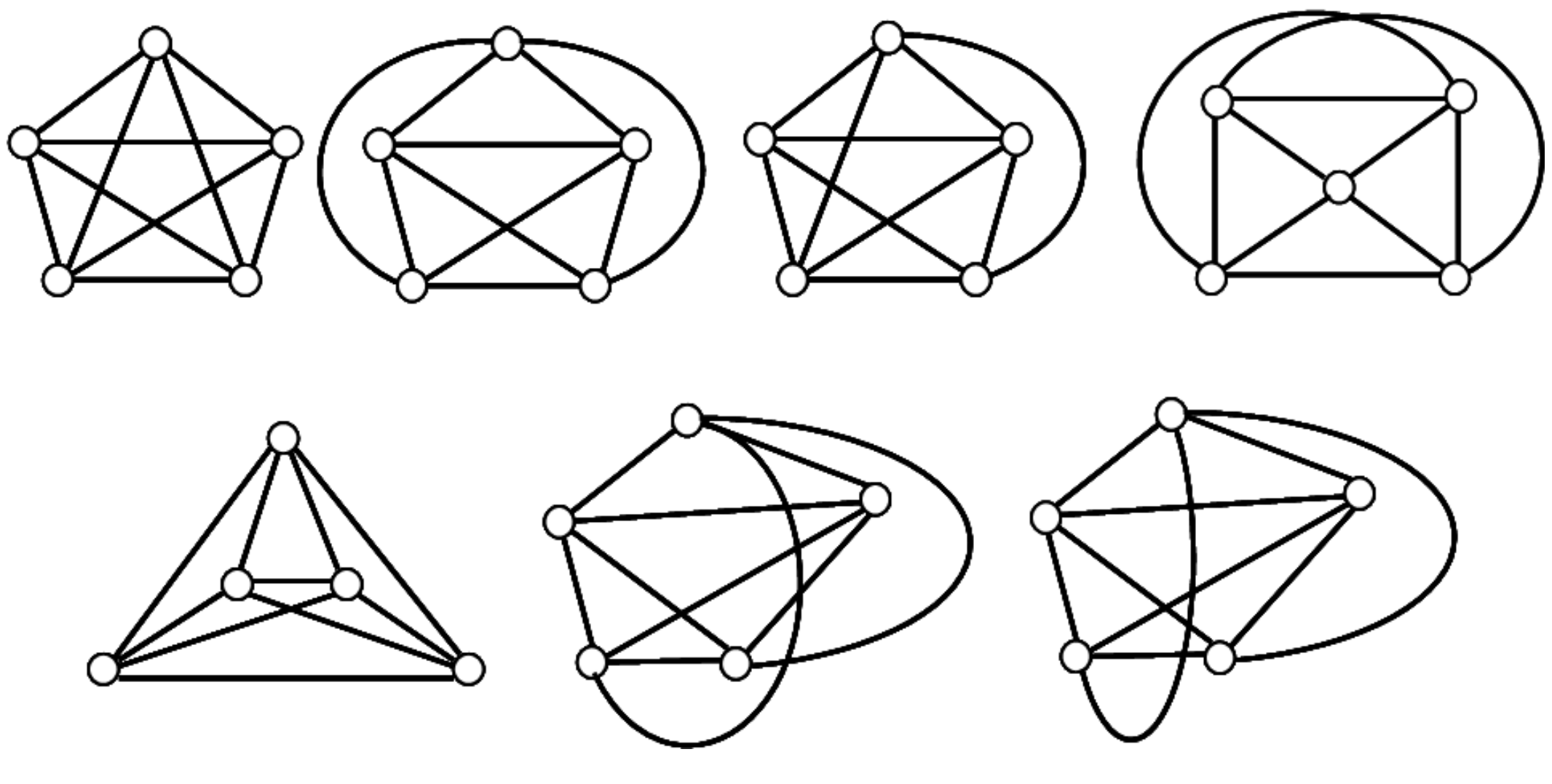

Рис. 1. Различные рисунки графа $K_5$.

Первые, кто столкнулся с задачей построения рисунка графа на плоскости, были разработчики автоматизированных систем проектирования плоских конструктивов, таких как печатные платы, микросборки, интегральные микросхемы и т.д. Конструктор располагал контакты элементов, в виде вершин графа задавая их геометрические координаты, и после их расположения начиналась решаться задача проведения соединений. Но если характеризовать соединения прямыми линиями, то они могут пересекаться, а это недопустимо для электрических соединений (см. рис. 2,а). Очевидно, что для того чтобы избежать пересечения, лучше характеризовать соединения ломаными линиями (см. рис.2,б). Но как проводить ломаные линии? Решение было найдено с появлением волнового алгоритма Ли [5,37]. Здесь вся плоскость проведения соединений разбивалась на клетки определенного размера, и последовательно проводились соединения путем распространения волны (см. рис. 2,в). Данный способ являлся метрическим, так как рассматривал проведение соединений в некотором пространстве $R^2$ с заданной эвклидовой метрикой. Основными недостатками волнового алгоритма являлись его большая вычислительная сложность и последовательный характер проведения соединений, так как проведение последующих соединений напрямую зависит от проведения предыдущих. И такой способ не гарантировал оптимального конструкторского решения.

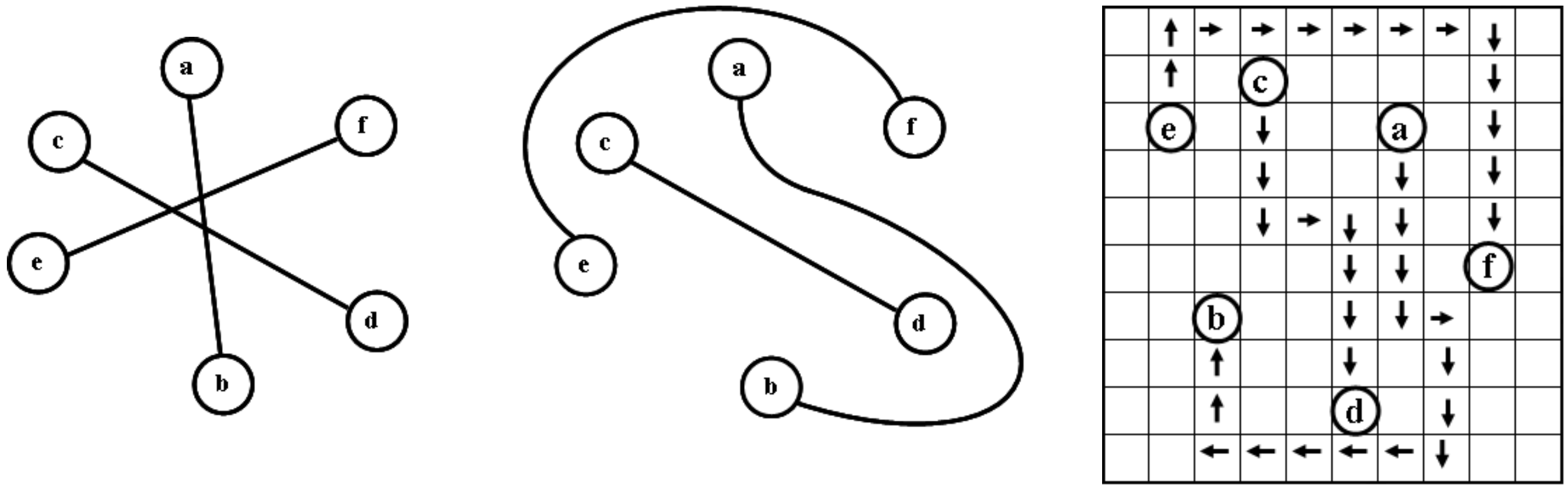

а) б) в)

Рис. 2. а) соединения пересекаются, б) соединения не пересекаются, в) проведение соединений алгоритмом Ли.

Для устранения недостатков волнового алгоритма Ли было создано множество модификаций: лучевые алгоритмы, канальные, магистральные и т.д., что позволило значительно сократить вычислительную сложность. На базе волнового алгоритма Ли и его модификаций созданы мощные системы автоматизации проектирования плоских конструктивов такие как PCAD, OrCAD, TOPOR и множество других.

Последовательный характер проведения соединений и, как следствие, разбиение единого процесса проведения соединений на два: размещение элементов и трассировка соединений, с потерей оптимальности решения, не устраивало разработчиков автоматизированных систем проектирования. И это вызвало целую серию новых математических подходов к решению поставленной проблемы, которые можно условно назвать топологическими методами использующими теоретико-множественное представление [9,12,28,32,40].

Одной из первых отечественных публикаций по представлению топологического описания рисунка графа с определением количества пересечений соединений только для двух циклических фрагментов была работа [43].

Автору удалось преобразовать данную задачу к решению известной задачи теории графов - представлению двудольного графа с минимальным числом пересечения ребер (см. рис. 3). Данный подход привел к развитию методов цепочного построения рисунка графа для определения количества пересечений [5]. Однако вопрос поиска оптимальной цепочки соединений для получения рисунка графа с минимальным числом пересечения ребер остался открытым.

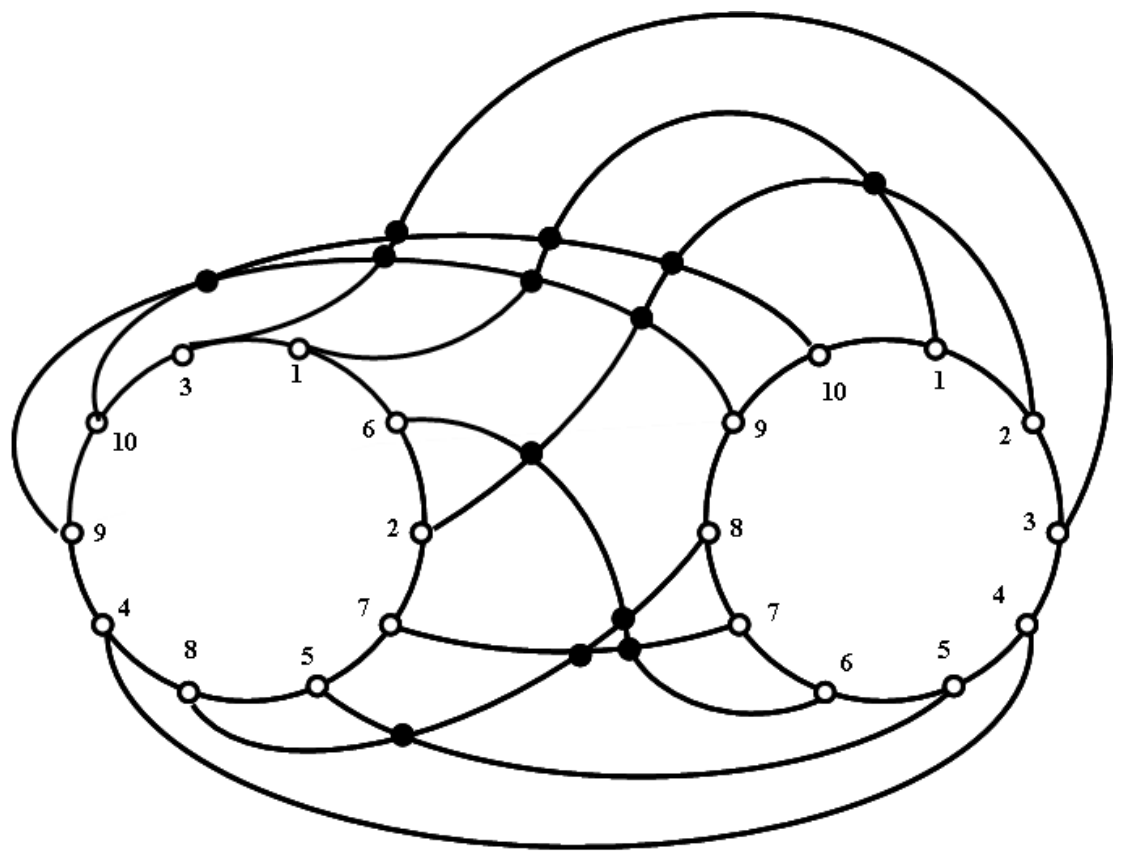

Рис. 3. Пересечения соединений двух фрагментов.

Дальнейшее развитие методов определения пересечений соединений представлено в работах Раппопорта Л.И., Мороговского Б.Н., Поливцева С.А. [34-36]. Они первые высказали мысль, что при построении рисунка графа анализ отношения пересечения ребер можно производить в топологическом пространстве, в котором метрические свойства не определены. Они разработали основы векторной алгебры пересечений, где для полного и непротиворечивого описания рисунка графа авторы вводят понятие координатно-базисной системы (КБС), и относительно её устанавливают проекции всех соединений с целью определения пересечения ребер по их проекциям.

Как правило, такая координатно-базисная система может быть построена относительно любого дерева графа. Несмотря на революционность подхода, ограничения, накладываемые на построение координатно-базисной системы запрещающие соединениям пересекать базовые вектора КБС, несколько сократили применимость метода. Так как различные деревья порождают различные КБС с различным количеством пересечений соединений, и тогда возникает проблема оптимального, с точки зрения минимального количества пересечения соединений, выбора дерева графа. Выбор такого дерева является трудноерешаемой задачей требующей перебора всех деревьев графа и построения рисунков на их основе.

Трудности, связанные с топологическим способом определения пересечений соединений (ребер графа), привели к более тщательному исследованию процесса планаризации и созданию алгоритмов для плоской укладки графа.

Имеется несколько критериев планарности графа, данные Л.С.Понтрягиным, К.Куратовским, К.Вайнером, С.Маклейном [9-11,26-28,39,41]. Критерии планарности таковы, что если даже удалось установить планарность графа, то нет информации о том, как строить его укладку на плоскости. В тоже время, при решении практических задач недостаточно знать, что граф планарен, а необходимо построить его плоское изображение. Практическое использование перечисленных критериев затруднено из-за необходимости полного перебора простых циклов для содержательного рассмотрения графа. Поэтому были разработаны эвристические методы определения планарности. Условно их можно разделить на три класса: циклические, матричные и комбинированные. В работе [28] приведены ссылки на опубликованные работы по созданию алгоритмов планаризации. Подробно рассмотрены все применяемые типы алгоритмов на тот временной период и указаны их основные недостатки для практического применения.

Один из лучших алгоритмов проверки графов на планарность разработан в 1970 г. Хопкрофтом и Тарьяном [42]. Они нашли алгоритм, требующий $O(|N|\log|N|)$ единиц времени, который они в конечном счете улучшили до $O(|N|)$. Данный алгоритм проверяет граф на планарность и, если он планарен, производит его плоскую укладку. Однако, данный алгоритм лишь косвенно содержит информацию о рисунке графа и его применение затруднительно для дальнейшего решения задачи построения и выделения максимально плоского суграфа из непланарного графа, а также для решения задачи построения рисунка графа с минимальным числом пересечений для непланарного графа.

Опыт построения рисунка графа для электрических соединений был перенесен на другие области деятельности человека связанные с сетевыми структурами, такие как картография и биоинформатика, биологические науки, искусственный интеллект, анализ финансовой информации, анализ социальных сетей и так далее.

В настоящее время сформировалось новое научное направление, связанное с рисованием графов и визуализацией сетей на компьютере. Основные вопросы, решаемые данным научным направлением, можно сформулировать следующим образом: получить графическое представление укладки графа направленное на удобное отображение некоторых свойств моделируемого объекта представленного в виде графа.

В обзорных работах [3,4,12,47] рассказывается об эволюции методов визуализации графов в связи с потребностями практики.

Следует заметить, что понятия «рисование графов» (Graph Drawing) и «визуализация графов» (Graph Visualization), появившиеся первоначально в англоязычной литературе, отличаются по смыслу. Термин Graph Drawing чаще употребляется в контексте теоретических работ по рисованию графов на плоскости, а термин визуализация графов является самостоятельным подмножеством в быстро развивающемся направлении «Визуализация информации». Понятие «визуализация графов» является более широким по отношению к понятию «рисование графов». С алгоритмической точки зрения, визуализация графов рассматривает не только проблемы рисования графов, но и вопросы интерактивного взаимодействия с графами, навигацию по очень большим графам, вопросы удобства использования тех или иных изображений графа.

Исторически, выделение направления Graph Drawing в качестве самостоятельного научного направления, принято связывать с первой международной конференцией, которая состоялась в 1992 году в Риме. С тех пор издательство Springer-Verlag ежегодно публикует в сериях LNCS труды этих конференций, которые содержат новые алгоритмы рисования и визуализации графов, теоретические результаты, касающиеся эффективности этих алгоритмов, а также демонстрации новых систем. Недавно начал издаваться электронный журнал «Journal of Graph Algorithms and Applications», посвященный конструированию и анализу алгоритмов визуализации графов, а также экспериментам и приложениям.

Среди книг, посвященных рисованию графов, следует отметить книгу основоположников этого направления (G. di Battista P. Eades R. Tamassia I.G. Tollis). Книги по визуализации информации, использующие методы визуализации графов, гораздо более многочисленны (S. K. Card, J. D. MacKinlay, B. Shneiderman и другие) [44,45].

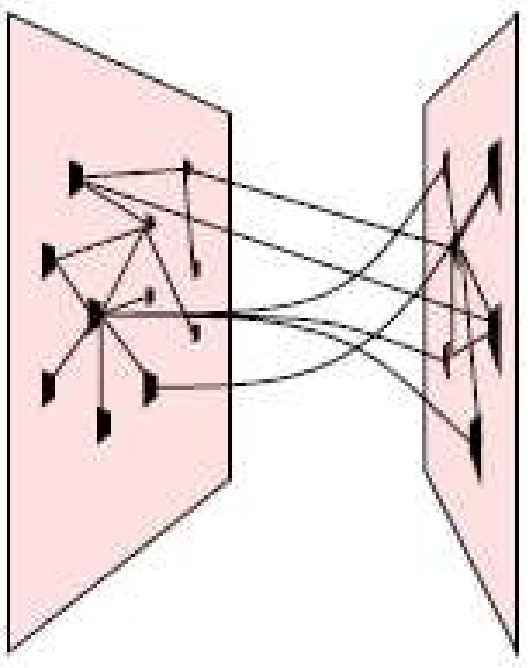

а)

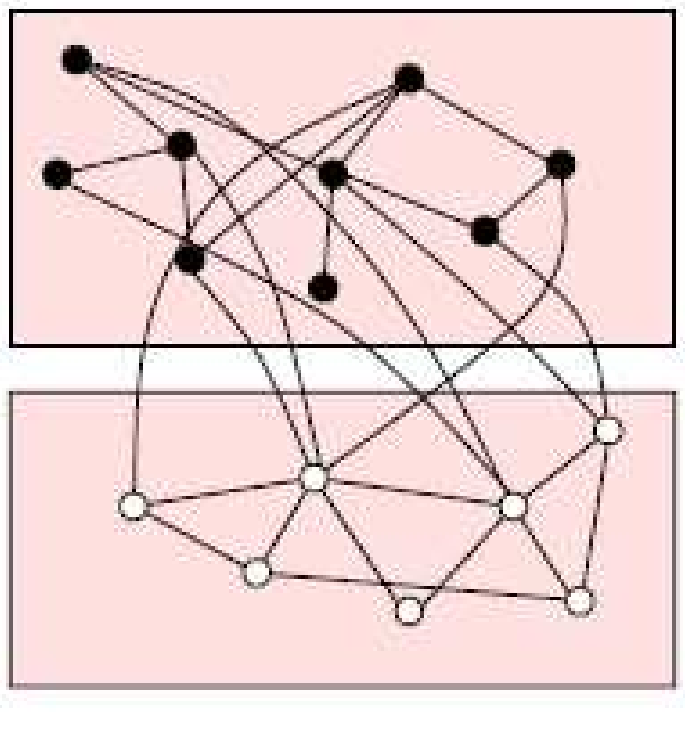

б)

Рис. 4. а) рисунок графа в пространстве, б) изображение на плоскости.

В настоящее время имеется большое количество статей, предлагающих различные способы отображения и рассматривающие вопрос о том, что такое хорошее изображение графа. И для характеристики можно придумать огромное количество разных способов классификации изображений (см. рис. 4-5). Во-первых, изображения графов можно классифицировать в соответствии с *типом изображаемого графа*. Во-вторых, можно разбить графы на группы. Графы, которые можно нарисовать без пересечений ребер, образуют группу планарных графов. Также принято рассматривать отдельно алгоритмы визуализации деревьев, ориентированных и неориентированных графов и т.д.

Много ограничений выражается также в терминах эстетических критериев, налагаемых на финальное изображение графа. Эстетические критерии применяются для улучшения читаемости (восприятия) изображения при помощи выделения целевых функций и их оптимизации.

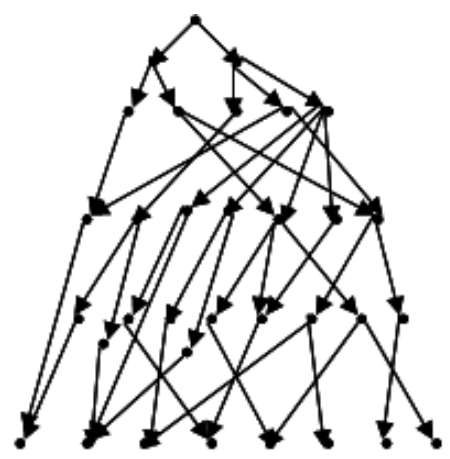

а) поуровневое изображение

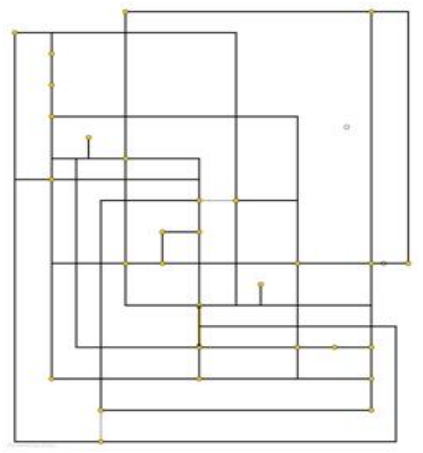

б) ортогональное изображение

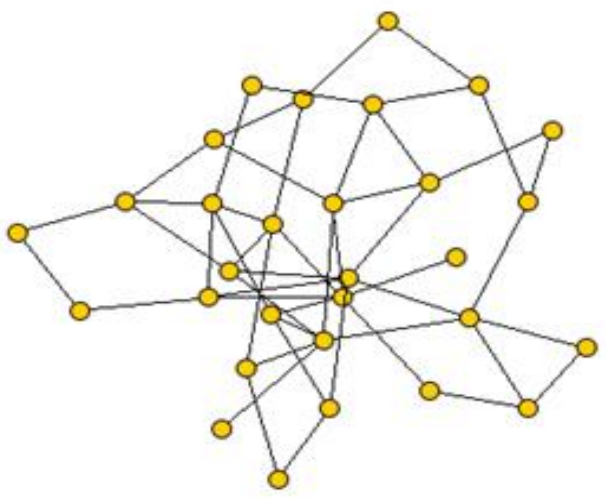

с) изображение, полученное на основе физических аналогий

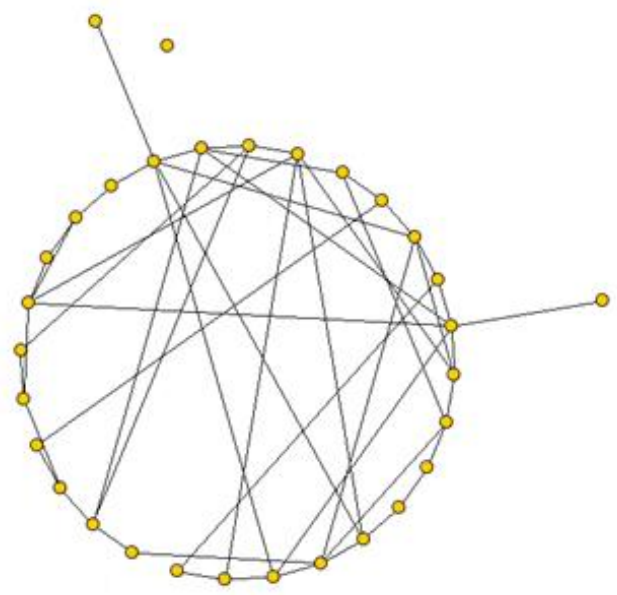

д) круговое изображение.

Рис. 5. Общепринятые статические изображения графов.

Обзор методов и способов представления рисунка графа на плоскости (визуализация графа) показывает, что в основном все имеющиеся на сегодняшний день способы представления рисунков графов сводятся к геометрическому заданию координат вершин тем или другим способом (см. рис.4-5). Заведомо расположив вершины в точках координатной сетки, или расположив их в иерархическом порядке по уровням, или по окружности (эллипсу) и представив соединения (ребра) прямыми, ломаными линиями или криволинейными отрезками описывают рисунок графа.

При таком подходе не существует строгого и недвусмысленного математического понятия плоского рисунка графа на плоскости для планарных графов. Отсутствует и понятие рисунка с минимальным числом пересечений ребер для непланарных графов. Отсутствие таких базовых понятий привело к тому, что рисунок графа строится метрическими методами при полном игнорировании его топологических свойств. И поэтому возникает вопрос, а возможно ли в принципе дать топологическое описание рисунка графа на плоскости? Покажем, что это в принципе возможно.

Целью данной работы является рассмотрение иллюстративно-технологической цепочки применения математических методов и приемов для решения поставленной задачи.

Данная работа рассчитана на прикладных математиков и может служить учебным пособием для дисциплин связанных с вычислительными методами теории графов. В работе приведены необходимые теоретические обоснования метода и подробное описание процессов вычислений. Для более углубленного понимания методов приведено необходимое количество решений примеров и задач.

Авторы выражают искреннюю благодарность сотрудникам математического факультета Запорожского национального университета: к.т.н. Борю С.Ю., д.ф.-м.н. Козину И.В., к.ф.-м.н. Стеганцевой П.Г., к.ф.-м.н. Зиновееву И.В., к.ф.-м.н. Кудину А.В., за высказанные замечаний и ценные советы при обсуждении данной работы.

# Глава 1. ИЗОМЕТРИЧЕСКИЕ ЦИКЛЫ ГРАФА

## 1.1. Фундаментальные циклы и разрезы

Как известно, в пространстве суграфов можно выделить два подпространства называемых подпространством разрезов S(G) и подпространством циклов C(G) [77].

Обычно в теории графов для поиска базисов применяется фундаментальная система циклов и разрезов. Данная система образуется в результате выделения случайного дерева графа (ациклического суграфа), тем самым разделяя ребра графа на ветви дерева и хорды. Ребра, принадлежащие дереву, называются ветвями, а не принадлежащие дереву – хордами. Каждый фундаментальный цикл образуется как объединение одной хорды и ветвей дерева. Рассматривая все хорды для выделенного дерева графа, строим матрицу фундаментальных циклов. Если представить (0,1)-матрицу фундаментальных циклов в виде единичной матрицы хорд и блочной матрицы $\pi$ состоящей из ветвей дерева, то можно получить матрицу фундаментальных разрезов графа в виде единичной матрицы ветвей дерева и блочной матрицы $\rho = \pi^t$ состоящей из хорд.

Например, для графа G (рис. 1.1) относительно дерева T = $\{e_1,e_3,e_7,e_{10},e_{11}\}$ матрица фундаментальных циклов имеет вид [80]:

$C_ф$=

| | $e_2$ | $e_4$ | $e_5$ | $e_6$ | $e_8$ | $e_9$ | $e_1$ | $e_3$ | $e_7$ | $e_{10}$ | $e_{11}$ |
|---|---|---|---|---|---|---|---|---|---|---|---|
| $e_2$ | 1 | | | | | | | 1 | | 1 | 1 |
| $e_4$ | | 1 | | | | | 1 | 1 | 1 | 1 | 1 |
| $e_5$ | | | 1 | | | | 1 | 1 | | 1 | 1 |
| $e_6$ | | | | 1 | | | 1 | 1 | | | |
| $e_8$ | | | | | 1 | | | | 1 | 1 | |
| $e_9$ | | | | | | 1 | | | 1 | 1 | 1 |

Единичная блочная подматрица | Блочная подматрица $\pi$

Матрица фундаментальных разрезов имеет вид

$S_ф$=

| | $e_1$ | $e_3$ | $e_7$ | $e_{10}$ | $e_{11}$ | $e_2$ | $e_4$ | $e_5$ | $e_6$ | $e_8$ | $e_9$ |
|---|---|---|---|---|---|---|---|---|---|---|---|
| $e_1$ | 1 | | | | | | 1 | 1 | 1 | | |
| $e_3$ | | 1 | | | | 1 | 1 | 1 | 1 | | |
| $e_7$ | | | 1 | | | | 1 | | | 1 | 1 |
| $e_{10}$ | | | | 1 | | 1 | 1 | 1 | | 1 | 1 |
| $e_{11}$ | | | | | 1 | 1 | 1 | | | | 1 |

Единичная блочная подматрица | Транспонированная блочная подматрица $\rho = \pi^t$

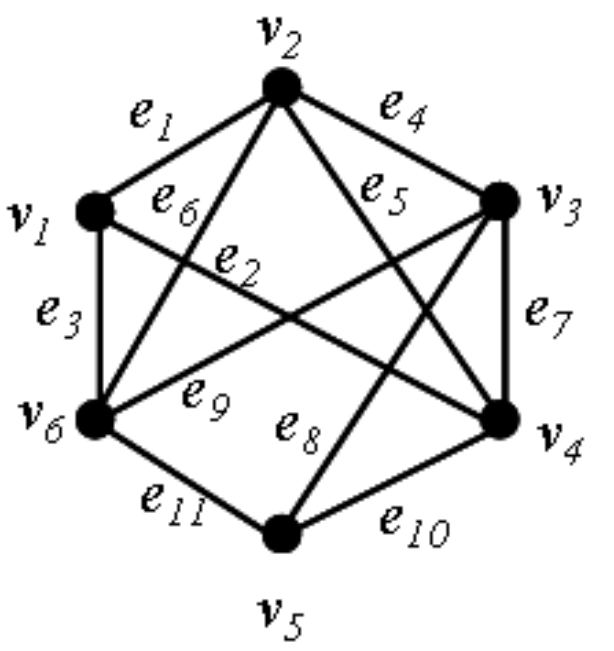

Рис. 1.1. Граф G.

Таким образом, относительно ациклического суграфа (дерева) строится система фундаментальных циклов и фундаментальных разрезов графа. Фундаментальная система циклов и фундаментальная система разрезов, в свою очередь, служат для формирования базиса подпространства циклов и базиса подпространства разрезов.

Для определения базиса подпространства циклов и базиса подпространства разрезов в теории графов применяется фундаментальная система циклов и разрезов. Количество фундаментальных циклов определяется цикломатическим числом графа $\nu(G) = m–n+1$, а количество фундаментальных разрезов определяется рангом графа $\rho(G) = n – 1$ [84].

Любой суграф, принадлежащий подпространству разрезов S графа G**,** в общем случае, является квалиразрезом. Количество ребер составляющих квалиразрез называется *длиной квалиразреза.* Например, фундаментальные разрезы графа G относительно дерева T = $\{e_1,e_3,e_7,e_{10},e_{11}\}$ (рис. 1.2):

$s_1 = \{\mathbf{e}_1,e_4,e_5,e_6\}$; $s_2 = \{e_2,\mathbf{e}_3,e_4,e_5,e_6\}$; $s_3 = \{e_4,\mathbf{e}_7,e_8,e_9\}$;

$s_4 = \{e_2,e_4,e_5,e_8,e_9,\mathbf{e}_{10}\}$; $s_5 = \{e_2,e_4,e_5,e_9,\mathbf{e}_{11}\}$

имеют суммарную длину $l = l_1 + l_2 + l_3 + l_4 + l_5 = 4 + 5 + 4 + 6 + 5 = 24$.

А следующие базисные разрезы графа G полученные путем линейной комбинации фундаментальных разрезов

$s_6 = s_1 \oplus s_2 = \{e_1,e_4,e_5,e_6\} \oplus \{e_2,e_3,e_4,e_5,e_6\} = \{\mathbf{e}_1,e_2,\mathbf{e}_3\}$;
$s_7 = s_1 = \{\mathbf{e}_1,e_4,e_5,e_6\}$;
$s_8 = s_3 = \{e_4,\mathbf{e}_7,e_8,e_9\}$;
$s_9 = s_3 \oplus s_4 = \{e_4,e_7,e_8,e_9\} \oplus \{e_2,e_4,e_5,e_8,e_9,e_{10}\} = \{e_2,e_5,\mathbf{e}_7,\mathbf{e}_{10}\}$;
$s_{10} = s_4 \oplus s_5 = \{e_2,e_4,e_5,e_8,e_9,e_{10}\} \oplus \{e_2,e_4,e_5,e_9,e_{11}\} = \{e_8,\mathbf{e}_{10},\mathbf{e}_{11}\}$,

имеют меньшую суммарную длину $l = l_6 + l_7 + l_8 + l_9 + l_{10} = 3+4+4+4+3 = 18$. Такие разрезы, имеющие минимально возможную длину, будем называть центральными разрезами графа G [21,22].

С другой стороны, длина центрального разреза ставит в соответствие локальную степень вершины и наоборот.

Рассмотрим фундаментальные циклы для выбранного дерева графа G (рис. 1.5):

$c_1 = \{e_2, \mathbf{e}_3, \mathbf{e}_{10}, \mathbf{e}_{11}\}$; $c_2 = \{\mathbf{e}_1, \mathbf{e}_3, e_4, \mathbf{e}_7, \mathbf{e}_{10}, \mathbf{e}_{11}\}$; $c_3 = \{\mathbf{e}_1, \mathbf{e}_3, e_5, \mathbf{e}_{10}, \mathbf{e}_{11}\}$;
$c_4 = \{\mathbf{e}_1, \mathbf{e}_3, e_6\}$; $c_5 = \{\mathbf{e}_7, e_8, \mathbf{e}_{10}\}$; $c_5 = \{\mathbf{e}_7, e_9, \mathbf{e}_{10}, \mathbf{e}_{11}\}$

имеют суммарную длину $l = l_1 + l_2 + l_3 + l_4 + l_5 + l_6 = 4 + 6 + 5 + 3 + 3 + 4 = 25$.

А следующие базисные циклы графа G, полученные путем линейной комбинации фундаментальных циклов

$c_7 = c_1 = \{e_2, \mathbf{e}_3, \mathbf{e}_{10}, \mathbf{e}_{11}\}$;
$c_8 = c_4 = \{\mathbf{e}_1, \mathbf{e}_3, e_6\}$;
$c_9 = c_5 = \{\mathbf{e}_7, e_8, \mathbf{e}_{10}\}$;
$c_{10} = c_5 \oplus c_6 = \{e_7, e_8, e_{10}\} \oplus \{e_7, e_9, e_{10}, e_{11}\} = \{e_8, e_9, \mathbf{e}_{11}\}$;
$c_{11} = c_2 \oplus c_4 \oplus c_6 = \{e_1, e_3, e_4, e_7, e_{10}, e_{11}\} \oplus \{e_7, e_9, e_{10}, e_{11}\} \oplus \{e_1, e_3, e_6\} = \{e_4, e_6, e_9\}$;
$c_{12} = c_2 \oplus c_3 = \{e_1, e_3, e_4, e_7, e_{10}, e_{11}\} \oplus \{e_1, e_3, e_5, e_{10}, e_{11}\} = \{e_4, e_5, \mathbf{e}_7\}$

имеют меньшую суммарную длину $l = l_7 + l_8 + l_9 + l_{10} + l_{11} + l_{12} = 4 + 3 + 3 + 3 + 3 + 3 = 19$.

Будем рассматривать такие базисные циклы, суммарная длина которых минимальна.

## 1.2. Метрика графов. Расстояние в графе

Введем понятие расстояния между двумя вершинами графа[21,22].

*Расстоянием* $\rho$(x,y) графе **G** между вершинами x и y графа G = (V,E) называется длина кратчайшего из маршрутов (а значит – кратчайшей из простых цепей), соединяющих эти вершины; если x и y отделены в G, то $\rho$ (x,y) = $+\infty$. Функция $\rho = \rho(x,y)$ определенная на множестве всех пар вершин графа G и принимающая целые неотрицательные значения (к числу которых относится и бесконечное), заслуживает названия метрики графа, поскольку она удовлетворяет трем аксиомам Фреше [50]:

$$\forall x, y \in X[\rho(x, y) = 0 \Leftrightarrow x = y], \tag{1.1}$$

$$\forall x, y \in X[\rho(x, y) = \rho(y, x_1)], \tag{1.2}$$

$$\forall x, y, z \in X[\rho(x, y) + \rho(y, z) \geq \rho(x, z)]. \tag{1.3}$$

Выполнение первых двух аксиом тривиально, проверим третью (неравенство треугольника).

Если вершины x, y или вершины y, z отделены, то по крайней мере одна из двух величин $\rho$(x, y) и $\rho$(y, z) есть $\infty$. Если же ни x и y, ни z и y не отделены, то пусть

$$x\ e_1\ x_1\ e_2\ x_2 \ldots x_{\rho(x, y)-1}\ e_{\rho(x, y)}\ y;$$

и

$$y\ u_1\ y_1\ u_2\ y_2 \ldots y_{\rho(y, z)-1}\ u_{\rho(y, z)}\ z$$

– какие-либо из кратчайших цепей, соединяющих эти пары вершин.

Маршрут

$$x\ e_1\ x_1\ e_2\ x_2 \ldots x_{\rho(x, y)-1}\ e_{\rho(x, y)}\ y\ u_1\ y_1\ u_2\ y_2 \ldots y_{\rho(y, z)-1}\ u_{\rho(y, z)}\ z$$

обладает длиной $\rho$(x,y) + $\rho$(y,z), значит длина $\rho$(x,z) кратчайшей цепи между x и z не превышает $\rho$(x,y) $+\rho$(y,z). Таким образом, в обоих случаях неравенство треугольника выполнено.

Введем следующее понятие, связанное с метрикой графа

**Определение 1.1**[40]**.** *Изометрический подграф* – подграф $G^*$ графа G, у которого все расстояния внутри $G^*$ те же самые, что и в G.

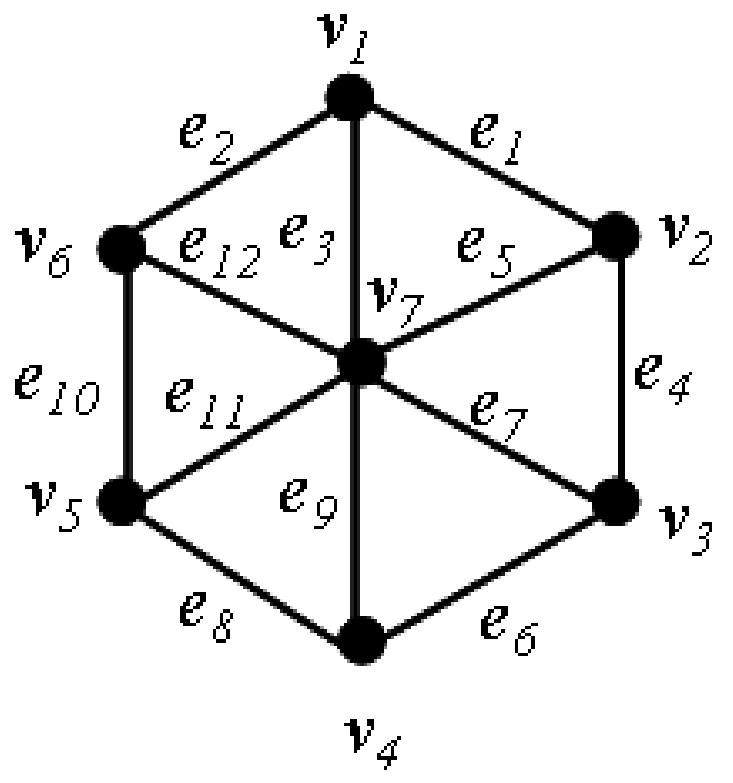

Рис. 1.2. Граф G.

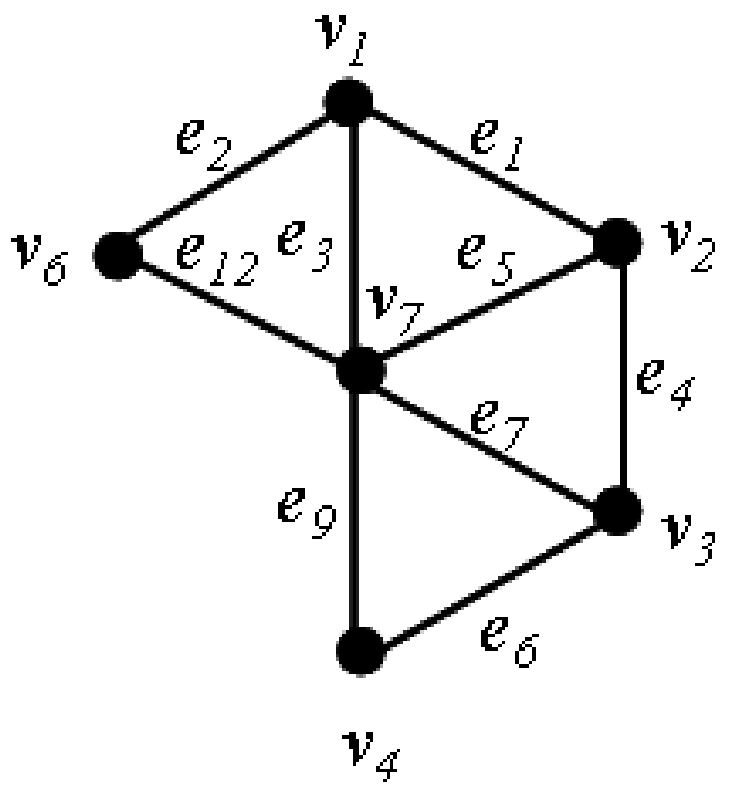

Рис. 1.3. Изометрический подграф графа G.

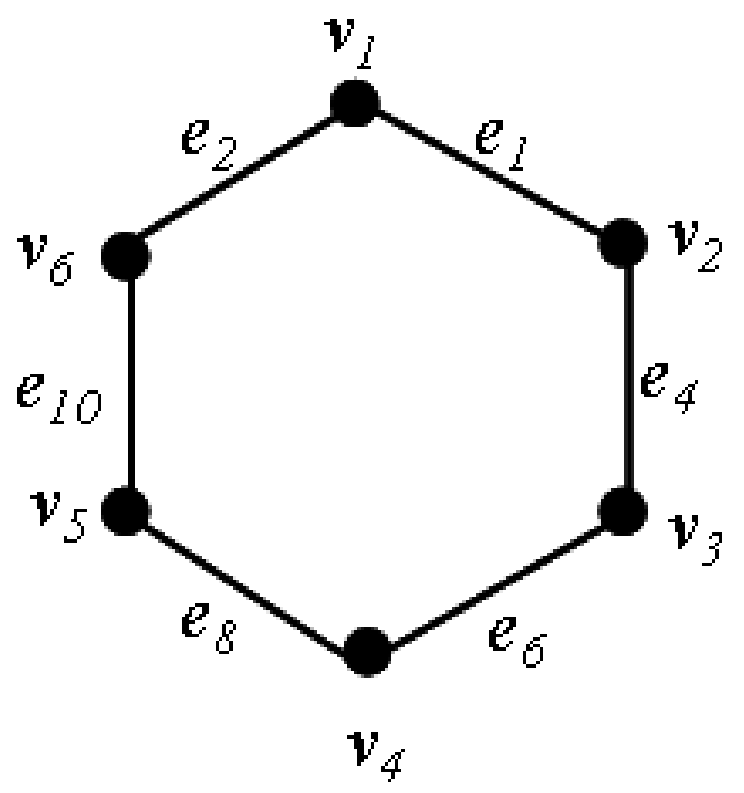

Рис. 1.4. Неизометрический подграф графа G.

На рис. 1.2. представлен граф G, имеющий следующую матрицу расстояний:

$\rho =$

| | $v_1$ | $v_2$ | $v_3$ | $v_4$ | $v_5$ | $v_6$ | $v_7$ |
|---|---|---|---|---|---|---|---|
| $v_1$ | | 1 | 2 | 2 | 2 | 1 | 1 |
| $v_2$ | 1 | | 1 | 2 | 2 | 2 | 1 |
| $v_3$ | 2 | 1 | | 1 | 2 | 2 | 1 |
| $v_4$ | 2 | 2 | 1 | | 1 | 2 | 1 |
| $v_5$ | 2 | 2 | 2 | 1 | | 1 | 1 |
| $v_6$ | 1 | 2 | 2 | 2 | 1 | | 1 |
| $v_7$ | 1 | 1 | 1 | 1 | 1 | 1 | |

Подматрица расстояний для изометрического подграфа представленного на рис. 1.3.

$\rho =$

| | $v_1$ | $v_2$ | $v_3$ | $v_4$ | $v_6$ | $v_7$ |
|---|---|---|---|---|---|---|
| $v_1$ | | 1 | 2 | 2 | 1 | 1 |
| $v_2$ | 1 | | 1 | 2 | 2 | 1 |
| $v_3$ | 2 | 1 | | 1 | 2 | 1 |
| $v_4$ | 2 | 2 | 1 | | 2 | 1 |
| $v_6$ | 1 | 2 | 2 | 2 | | 1 |
| $v_7$ | 1 | 1 | 1 | 1 | 1 | |

Подматрица расстояний для неизометрического подграфа представленного на рис. 1.4.

$\rho =$

| | $v_1$ | $v_2$ | $v_3$ | $v_4$ | $v_5$ | $v_6$ |
|---|---|---|---|---|---|---|
| $v_1$ | | 1 | 2 | 3 | 2 | 1 |
| $v_2$ | 1 | | 1 | 2 | 3 | 2 |
| $v_3$ | 2 | 1 | | 1 | 2 | 3 |
| $v_4$ | 3 | 2 | 1 | | 1 | 2 |
| $v_5$ | 2 | 3 | 2 | 1 | | 1 |
| $v_6$ | 1 | 2 | 3 | 2 | 1 | |

Как видим, матрица расстояний для подграфа на рис. 1.3 является подматрицей расстояний графа G, а матрица расстояний для подграфа на рис. 1.4 не является подматрицей расстояний графа G.

### 1.3. Множества изометрических циклов и центральных разрезов

Что касается базисной изометрической системы циклов и центральных разрезов, принцип их выбора отличается от выбора фундаментальной системы циклов и разрезов, так как понятие выделенного дерева здесь не несет полезной информации.

Центральные разрезы образуют множество размерностью равное количеству вершин графа G. Обозначим данное множество символом $\mathbf{S_e}$.

**Свойство 1.1.** Множество центральных разрезов графа обладает следующим свойством: кольцевая сумма всех центральных разрезов для графа G с *n* вершинами есть пустое множество:

$$\sum_{i=1}^{n} s_i = \varnothing. \qquad (1.4)$$

Количество центральных разрезов равное рангу графа **G** определяет базис подпространства разрезов. Алгоритм выделения множества центральных разрезов довольно прост, он может быть построен перечислением инцидентных ребер для *n-1* вершин графа **G** за линейное время.

Для графа, представленного на рис. 1.5, множество центральных разрезов $\mathbf{S_e}$ = $\{s_1,s_2,s_3,s_4,s_5,s_6\}$, где:

$s_1 = \{e_1,e_2,e_3\}$; $s_2 = \{e_1,e_4,e_5,e_6\}$; $s_3 = \{e_4,e_7,e_8,e_9\}$; $s_4 = \{e_2,e_5,e_7,e_9\}$;
$s_5 = \{e_8,e_{10},e_{11}\}$; $s_6 = \{e_3,e_6,e_{10},e_{11}\}$.

Любой суграф, принадлежащий подпространству циклов C, в общем случае является квазициклом. *Простые циклы* – это *квазициклы*, у которых локальная степень вершин в точности равна двум [21].

Мощность подмножества простых циклов в графе меньше мощности множества квазициклов. Подмножество простых циклов обозначим $C_R$:

$$\text{card } C_R \leq \text{card } C \qquad (1.5)$$

Однако существует подмножество с мощностью еще меньшей, чем подмножество простых циклов, обладающее определенными характерными свойствами.

**Определение 1.2.** *Изометрическим циклом* в графе называется простой цикл, для которого кратчайший путь между любыми двумя его вершинами состоит из рёбер этого цикла. Изометрический цикл – частный случай изометрического подграфа [40].

Или, другими словами, изометрическим циклом в графе называется подграф $\mathbf{G'}$ в виде простого цикла, если между двумя любыми несмежными вершинами данного подграфа в соответствуюшем графе G не существует маршрутов меньшей длины, чем маршруты, принадлежащие данному циклу.

Подмножество, состоящее из изометрических циклов, назовем подмножеством изометрических циклов и обозначим $C_\tau$. Сказанное поясним на примерах. Рассмотрим суграф, состоящий из ребер $\{e_1,e_3,e_{13},e_{15}\}$ графа $G_a$, представленного на рис. 1.5,а. Как видно, это простой цикл. Но в то же время, это не изометрический цикл, так как между вершинами $v_7$ и $v_8$ в графе существуют маршрут меньшей длины, проходящий по ребру $e_{14}$.

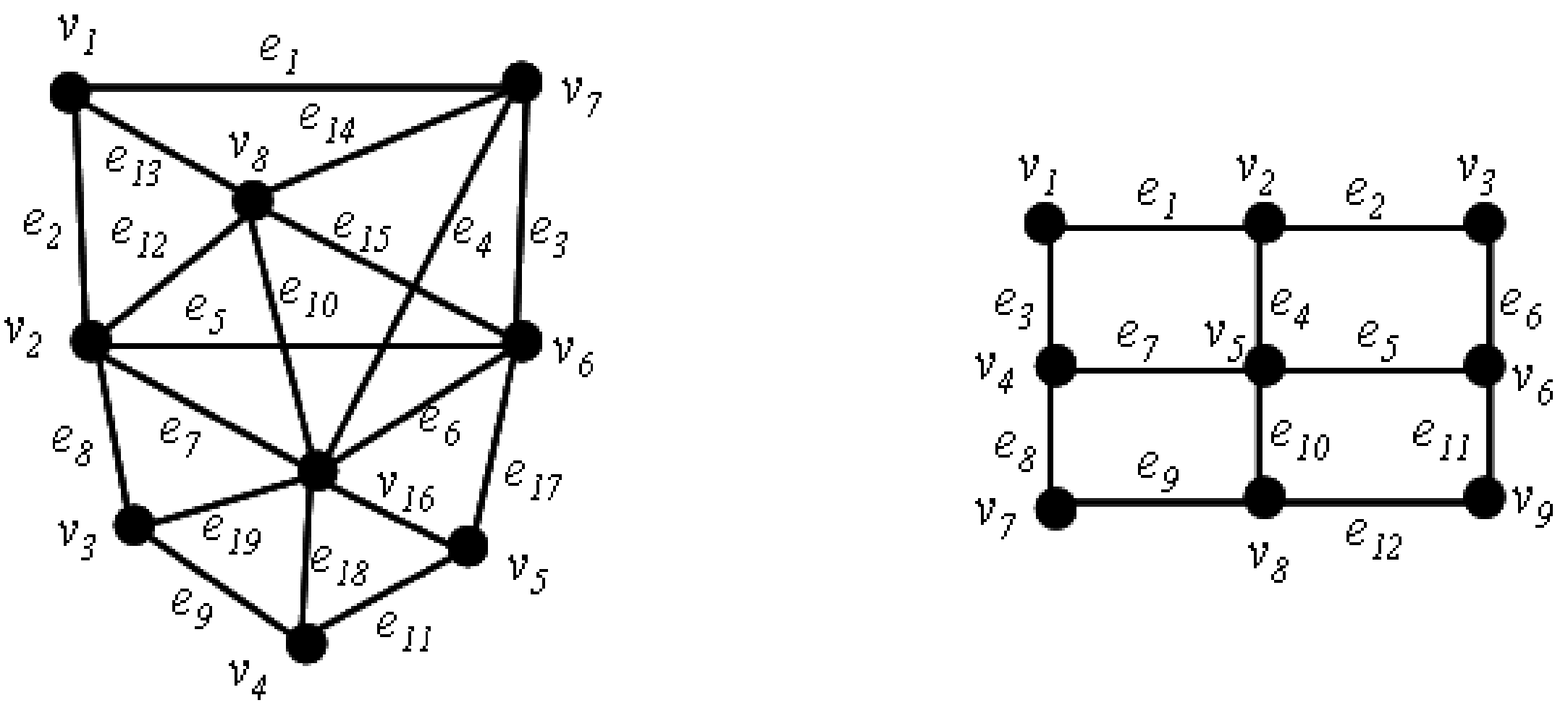

а) Граф $G_a$ б) Граф $G_b$

Рис. 1.5. Графы $G_a$ и $G_b$.

Рассмотрим граф $G_b$, представленный на рис. 1.5,б. Пусть цикл состоит из ребер $e_1,e_2,e_3,e_6,e_8,e_9,e_{11},e_{12}$. Данный суграф есть простой цикл. Однако этот суграф не может быть изометрическим циклом, так как в соответствующем графе между вершинами $v_2$ и $v_8$ имеется маршрут меньшей длины (а именно, маршрут, проходящий по ребрам $e_4$ и $e_{10}$), чем маршруты, принадлежащие этому суграфу (например, маршрут, проходящий по ребрам $e_1,e_3,e_8,e_9$ или $e_2,e_6,e_{11},e_{12}$).

Следует заметить, что в полных графах множество изометрических циклов совпадает с множеством циклов минимальной длины. В целях сокращения записи иногда будем обозначать вершины и ребра целыми числами.

Для изучения свойств изометрических циклов нам понадобится следующая теорема.

**Теорема 1.1.** Для любого несепарабельного графа линейное подпространство квазициклов имеет независимое подмножество изометрических циклов с мощностью равной цикломатическому числу графа [37].

*Доказательство*. Рассмотрим множество деревьев M(T) графа G.

Рассмотрим систему фундаментальных циклов порожденных деревом **T**. Пусть вершины $A_1,A_2,...,A_p$ образуют фундаментальный цикл содержащий хорду $(A_p,A_1)$. Если между несмежными вершинами этого цикла в графе не существует путей меньшей или равной длины, чем пути, принадлежащие циклу – то это изометрический цикл. Если в цикле существуют две несмежные вершины графа $A_i$ и $A_j$ ($i<j$), а в графе существует путь меньшей длины для выбранных вершин, чем путь по циклу $A_i,B_1,B_2,...,B_r,A_j$, то образуются циклы, кольцевая сумма которых есть исходный цикл. Оставляем цикл, содержащий хорду. Это и есть изометрический цикл. Так как количество фундаментальных циклов независимо и определяется цикломатическим числом, то количество изометрических циклов, полученных описанным выше способом, также независимо и равно цикломатическому числу графа включая все хорды.

Теорема доказана.

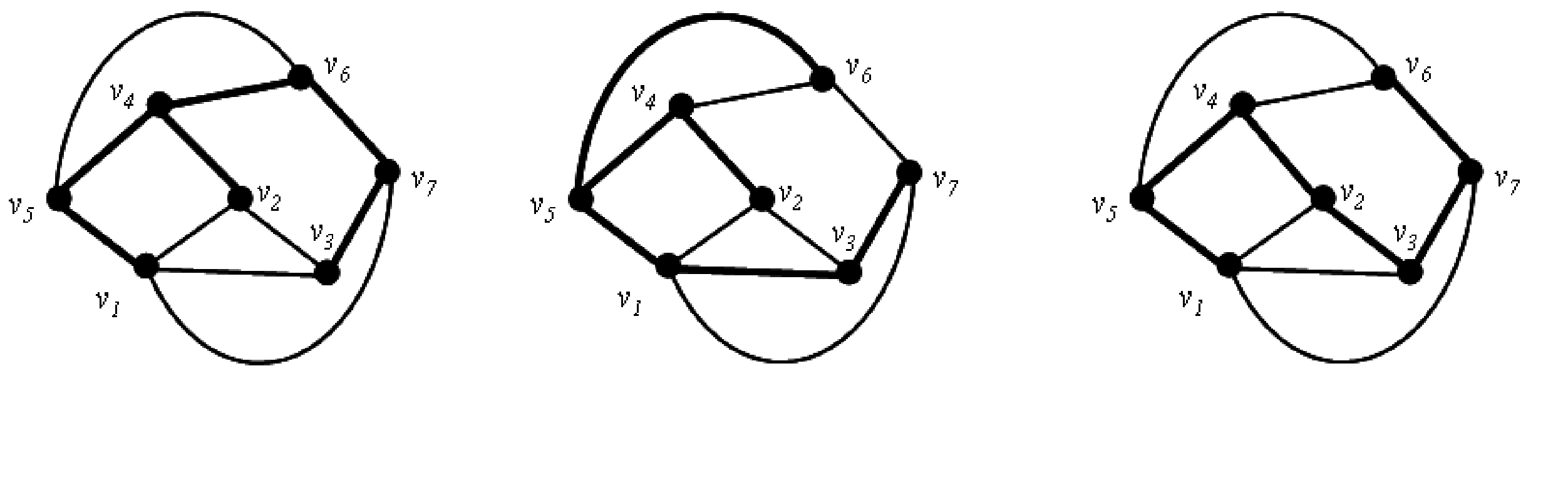

Рис. 1.6. К теореме 1.1.

## 1.4. Методы выделения множества изометрических циклов графа

Понятие изометрического цикла графа G тесно связано с минимальными (s-t) маршрутами графа [80].

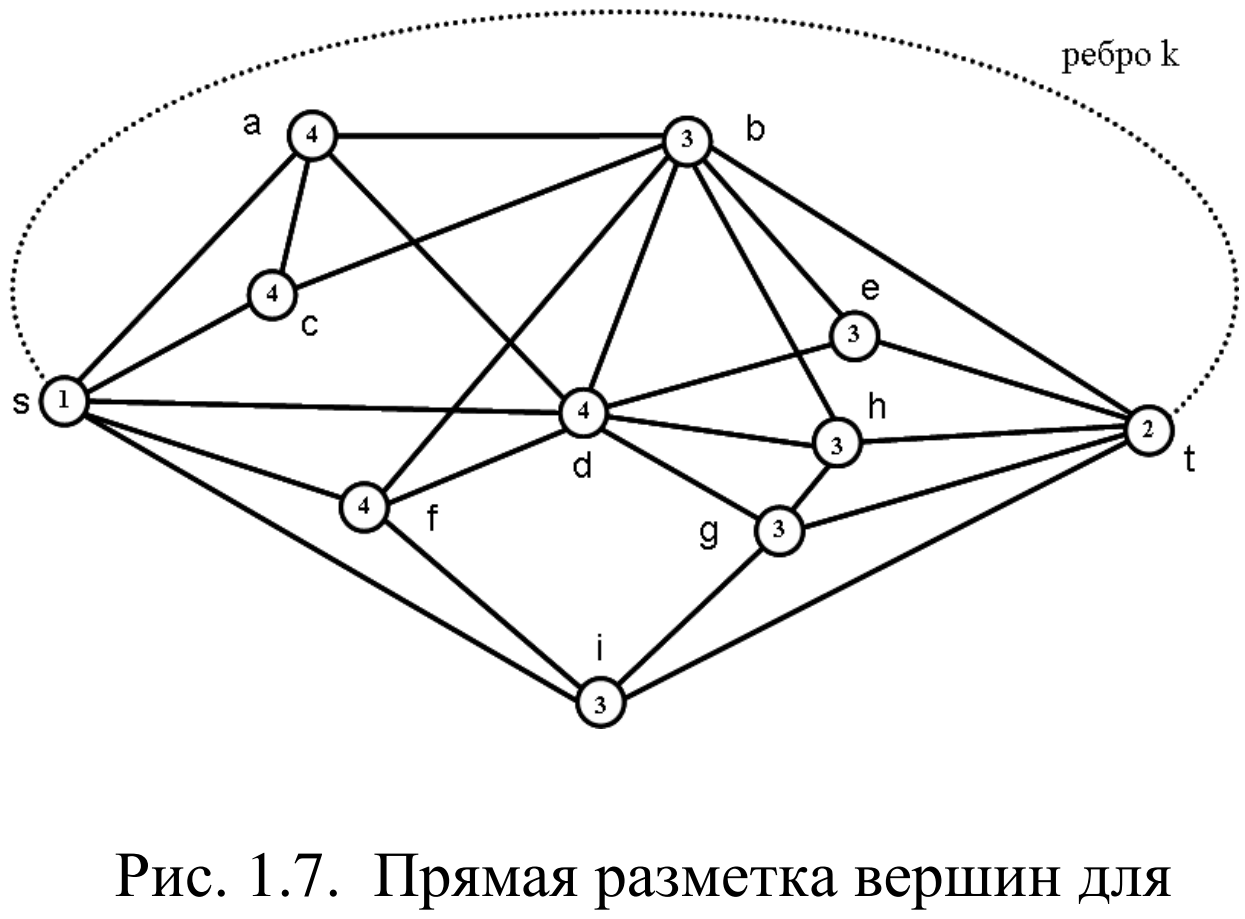

Рис. 1.7. Прямая разметка вершин для ребра k.

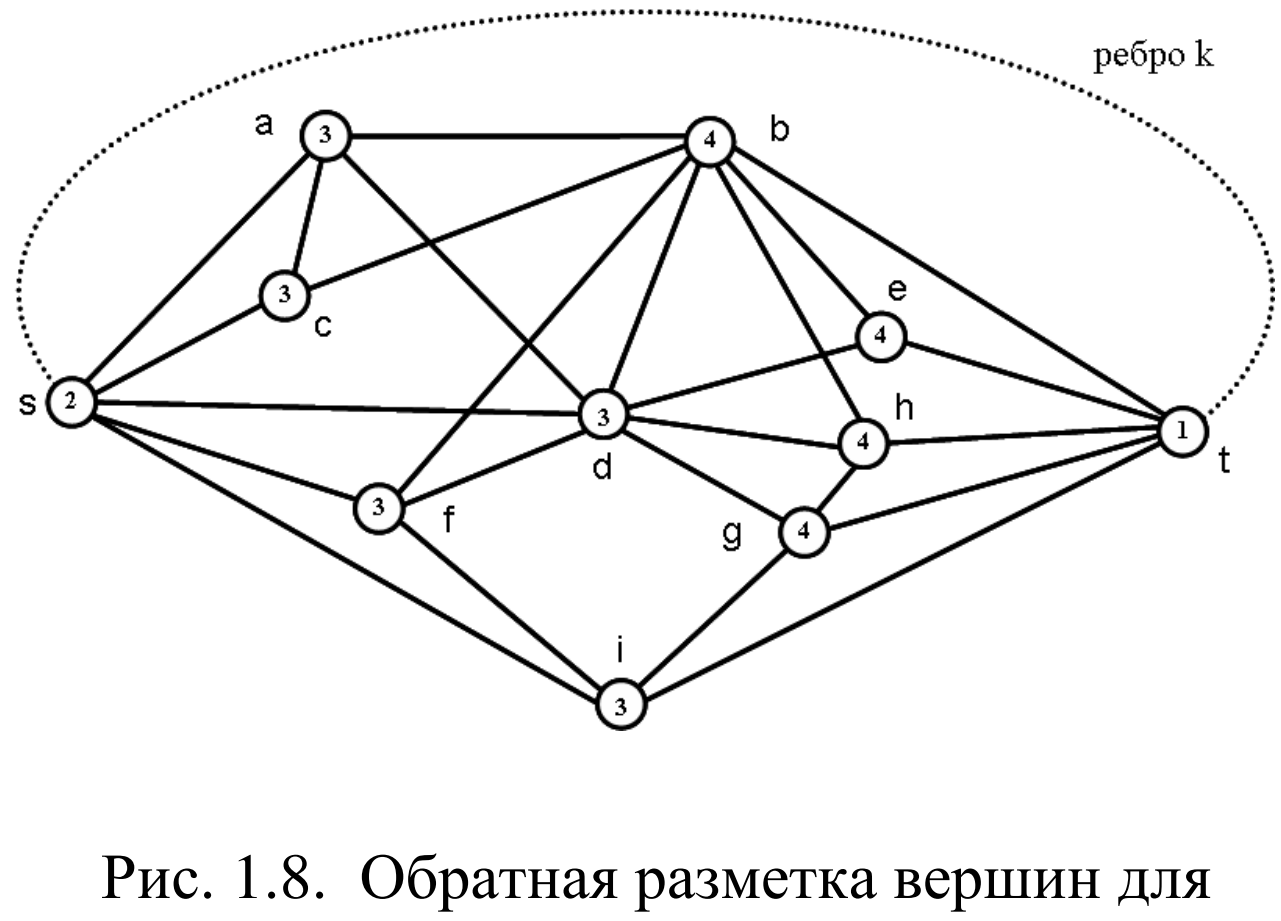

Рис. 1.8. Обратная разметка вершин для ребра k.

С этой целью рассмотрим изометрические циклы, проходящие по k-му ребру соединяющему вершины s и t графа G (рис. 1.7). Удалим из графа ребро k. Получим граф **G**-k, где вершины s и t теперь несмежны. Применим алгоритм поиска в ширину. Вершине s поставим в соответствие фронт волны 1, а вершине t – фронт волны 2. Тогда вершины, смежные с вершиной 2 и еще не помеченные пометим цифрой 3 и так далее. Другими словами, применим алгоритм поиска в ширину относительно вершины t. Выделим все простые цепи минимальной длины, образованные алгоритмом поиска в ширину, осуществляя проход от вершины с большим номером к вершине с меньшим номером (рис. 1.7). Сформируем множество $C_{st}$ циклов, где элементами множеств являются вершины:

$C_{st}$ = {{s,a,b,t},{s,c,b,t},{s,f,b,t},{s,d,b,t},{s,d,e,t},{s,d,h,t},{s,d,g,t},{s,f,b,t}, {s,f,i,t},{s,i,t}}.

Теперь вершине t поставим в соответствие фронт волны 1, а вершине s – фронт волны 2. Тогда вершины, смежные с вершиной 2 и еще не помеченные пометим цифрой 3, и так далее Другими словами, применим алгоритм поиска в ширину относительно вершины s. Выделим все простые цепи минимальной длины, образованные алгоритмом поиска в ширину, осуществляя проход от вершины с большим номером к вершине с меньшим номером (рис. 1.8). Сформируем множество $C_{ts}$ циклов:

$C_{ts}$ = {{s,a,b,t},{s,c,b,t},{s,f,b,t},{s,d,b,t},{s,d,e,t},{s,d,h,t},{s,d,g,t},{s,f,b,t}, {s,i,g,t},{s,i,t}}.

Циклы {s,a,b,t},{s,c,b,t},{s,f,b,t},{s,d,b,t},{s,d,e,t},{s,d,h,t},{s,d,g,t},{s,f,b,t},{s,i,t} принадлежащие одновременно и множеству $M_{st}$ и множеству $M_{ts}$ суть изометрические циклы. Циклы {s,f,i,t} и {s,i,g,t} принадлежат только одному из множеств и поэтому не являются изометрическими циклами.

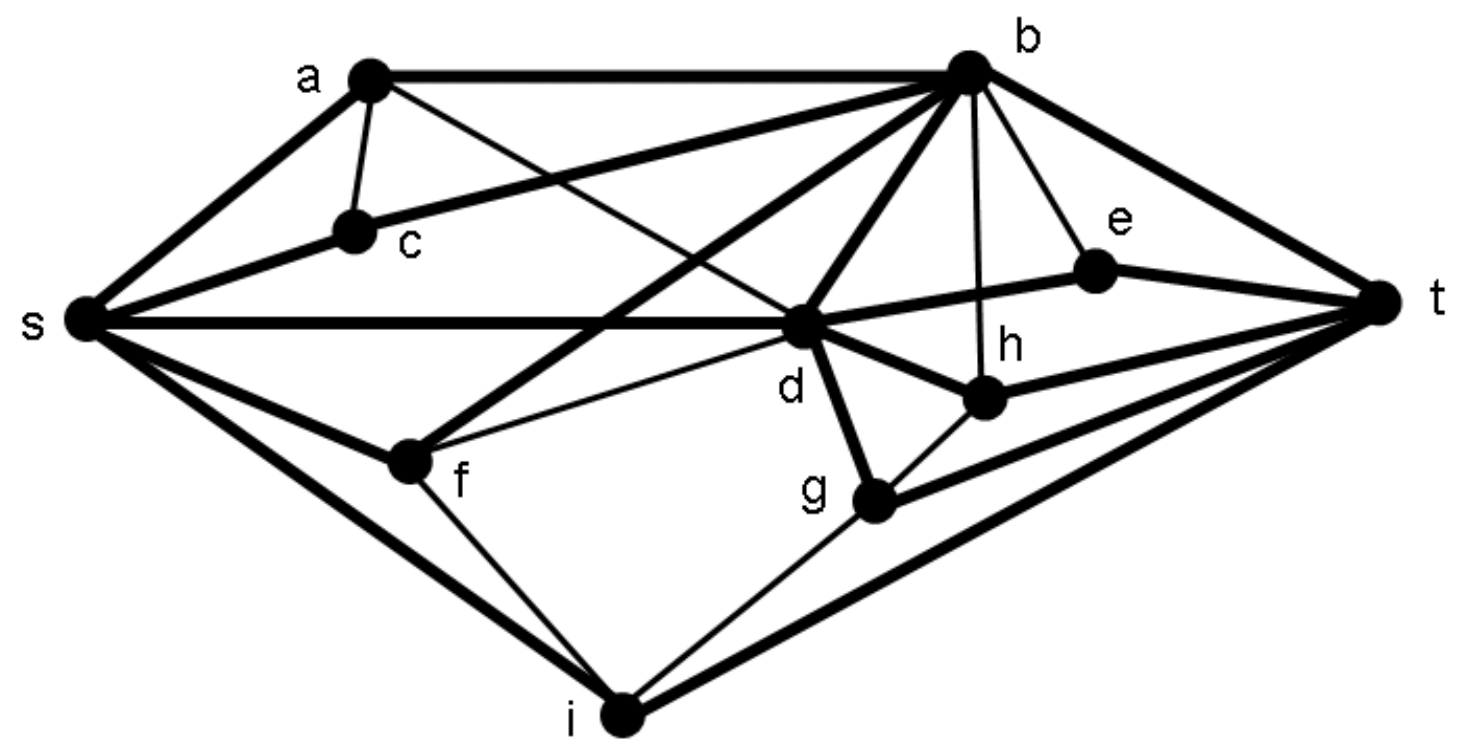

Рис. 1.9. Минимальные s-t цепи

Данные рассуждения можно применить ко всем ребрам графа G и сформировать множество изометрических циклов графа.

Вычислительная сложность такого алгоритма [] определится громоздкостью операции сравнения циклов выделенных при прямой и обратной разметке вершин. Здесь максимальное количество циклов определяется как (n-2) для каждой вершины, а количество сравнений тогда будет равно $(n-2)^2$. Полученное выражение нужно умножить на количество ребер = n(n-

1)/2. Окончательно получим, что вычислительная сложность равна

$$O(n^4) = n(n-1)(n-2)^2/2. \qquad (1.6)$$

Количество изометрических циклов в полном графе определяется по формуле

$$card(C_e) = n(n-1)(n-2)/6. \qquad (1.7)$$

Рассмотрим другой способ построения изометрических циклов. Выделим в полном графе множество циклов длиной три. Это, очевидно, будет множество изометрических циклов для полного графа.

Будем последовательно удалять ребра из полного графа $K_n$. Естественно, что тогда будут удалены и изометрические циклы, содержащие данное ребро, или два удаляемых цикла образуют новый изометрический цикл равный их кольцевой сумме в случае, если вновь образованный цикл не образован оставшимися изометрическими циклами. Таким образом, процесс удаления ребер из полного графа приводит к уменьшению количества изометрических циклов в графе, и, естественно, что их количество в произвольном графе не может превышать величины n(n-1)(n-2)/6.

Сказанное рассмотрим на примере графа $K_5$. Множество изометрических циклов $\mathbf{C}_\tau$ для графа $K_5$ (рис. 1.10):

$c_1 = \{e_1,e_2,e_5\}$; $c_2 = \{e_1,e_3,e_6\}$; $c_3 = \{e_1,e_4,e_7\}$; $c_4 = \{e_2,e_3,e_8\}$; $c_5 = \{e_2,e_4,e_9\}$;
$c_6 = \{e_3,e_4,e_{10}\}$; $c_7 = \{e_5,e_6,e_8\}$; $c_8 = \{e_5,e_7,e_9\}$; $c_9 = \{e_6,e_7,e_{10}\}$;
$c_{10} = \{e_8,e_9,e_{10}\}$.

Удалим из графа ребро $e_{10}$. Тогда из множества изометрических циклов $C_\tau$ удаляются все циклы включающие 10-ое ребро:

$c_6 = \{e_3,e_4,e_{10}\}$; $c_9 = \{e_6,e_7,e_{10}\}$; $c_{10} = \{e_8,e_9,e_{10}\}$.

Остаются изометрические циклы:

$c_1 = \{e_1,e_2,e_5\}$; $c_2 = \{e_1,e_3,e_6\}$; $c_3 = \{e_1,e_4,e_7\}$; $c_4 = \{e_2,e_3,e_8\}$; $c_5 = \{e_2,e_4,e_9\}$;
$c_7 = \{e_5,e_6,e_8\}$; $c_8 = \{e_5,e_7,e_9\}$.

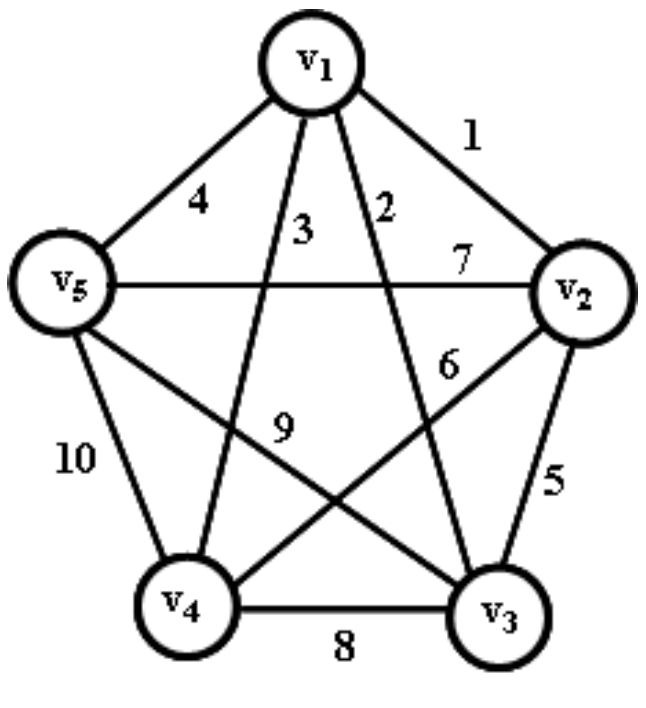

Рис. 1.10. Граф $K_5$ с пронумерованными ребрами.

В перспективе должны образоваться новые изометрические циклы длиной четыре, образованные из удаленных циклов:

$c_6 \oplus c_9 = \{e_3,e_4,e_{10}\} \oplus \{e_6,e_7,e_{10}\} = \{e_3,e_4,e_6,e_7\}$;
$c_6 \oplus c_{10} = \{e_3,e_4,e_{10}\} \oplus \{e_8,e_9,e_{10}\} = \{e_3,e_4,e_8,e_9\}$;

$c_9 \oplus c_{10} = \{e_6,e_7,e_{10}\} \oplus \{e_8,e_9,e_{10}\} = \{e_6,e_7,e_8,e_9\}$.

Однако их включение во множество оставшихся изометрических циклов невозможно, так как они могут быть образованы как результат кольцевого суммирования из оставшихся изометрических циклов:

$c_2 \oplus c_3 = \{e_1,e_3,e_6\} \oplus \{e_1,e_4,e_7\} = \{e_3,e_4,e_6,e_7\}$;
$c_4 \oplus c_5 = \{e_2,e_3,e_8\} \oplus \{e_2,e_4,e_9\} = \{e_3,e_4,e_8,e_9\}$;
$c_7 \oplus c_8 = \{e_5,e_6,e_8\} \oplus \{e_5,e_7,e_9\} = \{e_6,e_7,e_8,e_9\}$.

Если мы продолжим удаление ребра $e_2$ из графа, то из множества изометрических циклов $C_\tau = \{c_1,c_2,c_3,c_4,c_5,c_7,c_8\}$ удаляются все циклы включающие 2-ое ребро:

$c_1 = \{e_1,e_2,e_5\}$; $c_4 = \{e_2,e_3,e_8\}$; $c_5 = \{e_2,e_4,e_9\}$.

Остаются изометрические циклы:

$c_2 = \{e_1,e_3,e_6\}$; $c_3 = \{e_1,e_4,e_7\}$; $c_7 = \{e_5,e_6,e_8\}$; $c_8 = \{e_5,e_7,e_9\}$.

В перспективе должны образоваться новые изометрические циклы длиной четыре:

$c_1 \oplus c_4 = \{e_1,e_2,e_5\} \oplus \{e_2,e_3,e_8\} = \{e_1,e_3,e_5,e_8\}$;
$c_1 \oplus c_5 = \{e_1,e_2,e_5\} \oplus \{e_2,e_4,e_9\} = \{e_1,e_4,e_5,e_9\}$;
$c_4 \oplus c_5 = \{e_2,e_3,e_8\} \oplus \{e_2,e_4,e_9\} = \{e_3,e_4,e_8,e_9\}$.

Однако включение во множество оставшихся изометрических циклов двух первых невозможно, так как они могут быть образованы как результат кольцевого суммирования из оставшихся изометрических циклов:

$c_1 \oplus c_4 = \{e_1,e_3,e_6\} \oplus \{e_5,e_6,e_8\} = \{e_1,e_3,e_5,e_8\}$;
$c_1 \oplus c_5 = \{e_1,e_4,e_7\} \oplus \{e_5,e_7,e_9\} = \{e_1,e_4,e_5,e_9\}$.

Вновь образованный изометрический цикл $c_{4,5}$ включается во множество оставшихся изометрических циклов: $c_4 \oplus c_5 = \{e_3,e_4,e_8,e_9\}$.

Таким образом, множество изометрических циклов $C_\tau$ для графа полученного путем удаления 10-го и 2-го ребер из графа $K_5$ состоит из следующих изометрических циклов:

$c_2 = \{e_1,e_3,e_6\}$; $c_3 = \{e_1,e_4,e_7\}$; $c_7 = \{e_5,e_6,e_8\}$; $c_8 = \{e_5,e_7,e_9\}$;
$c_{4,5} = \{e_3,e_4,e_8,e_9\}$.

Характерная и особая роль изометрических циклов в теории графов определяется тем, что в плоских графах они являются границами граней. В свою очередь, характерная особенность центральных разрезов проявляется в том, что их длина определяет локальную степень вершин.

Ввиду важности вопроса выделения конечного множества изометрических циклов из множества квазициклов, предлагается алгоритм выделения множества изометрических циклов в графе.

Построение алгоритма начинается с выделения всех ребер в графе G. Выберем очередное ребро графа. Одну из вершин такого выбранного ребра пометим индексом 1, другую – индексом 2. Вершины графа смежные с вершиной, имеющей индекс 2, и ещё не помеченные, пометим индексом 3. Вершины графа смежные с вершиной, имеющей индекс 3, и ещё не помеченные, пометим индексом 4 и т. д. Число, выражающее индекс последней помеченной вершины (вершин) графа, называется глубиной проникновения волны, относительно выбранного ребра. Данный процесс представляет собой разметку вершин графа, относительно выбранного ребра волновым алгоритмом (алгоритмом поиска в ширину).

Построим простые циклы, проходящие по выбранному ребру, относительно первоначальной ориентации. С этой целью выберем все вершины графа G смежные с вершиной, помеченной индексом 1. Будем идти от любой выбранной вершины, имеющей глубину проникновения *d*, к вершинам, имеющим глубину проникновения (*d*-1), проходя при этом по ребрам графа, затем от вершины (*d*-1) к вершинам (*d*-2) и т.д. Остановим этот процесс тогда, когда подойдем к вершине, имеющей индекс 2. Пройдя по всем таким образом построенным маршрутам, построим систему циклов, проходящих по выбранному ребру *j*. Обозначим такое множество циклов через $S_j^1$. Переориентируем направление разметки, т.е. вершина, имеющая индекс 1, будет иметь индекс 2, а вершина, имеющая индекс 2, будет иметь индекс 1. И вновь построим разметку вершин. Описанным выше методом выделим систему циклов. Изометрические циклы, проходящие по выбранному ребру *j*, будут образованы как:

$$C_j = C_j^1 \cap C_j^2. \qquad (1.8)$$

Множество изометрических циклов графа G будет образовано как объединение всех циклов, проходящих по всем ребрам графа:

$$C_\tau = \bigcup_{i=1}^{m} C_i, \quad i = 1,2,...,m. \qquad (1.9)$$

### 1.5. Алгоритм выделения множества изометрических циклов графа

**Алгоритм 1.1. [Выделение множества изометрических циклов методом поиска в ширину]** [53,74].

**Шаг 1.** [**Выбор ребра**]. Выбираем ребро, идем на шаг 2. Если количество ребер исчерпано, то конец работы алгоритма.

**Шаг 2**. [**Прямая разметка вершин относительно вершины s ребра**]. Алгоритмом поиска в ширину производим прямую разметку вершин относительно вершины s выбранного ребра. Идем на шаг 3.

**Шаг 3**. [**Формирование множества циклов $C_s$ при прямой разметке вершин**].

Производим формирование множества циклов при прямой разметке вершин. Идем на шаг 4.

**Шаг 4**. [**Обратная разметка вершин относительно вершины t ребра**]. Алгоритмом поиска в ширину производим обратную разметку вершин относительно вершины t выбранного ребра. Идем на шаг 5.

**Шаг 5**. [**Формирование множества циклов $C_t$ при обратной разметке вершин**]. Производим формирование множества циклов при обратной разметке вершин. Идем на шаг 6.

**Шаг 6.** [**Проверка циклов $C_s$ и $C_t$ на совпадение**]**.** Приверяем циклы $C_s$ и $C_t$ для выбранного ребра на совпадение. Несовпадающие циклы исключаем из рассмотрения. Идем на шаг 7.

**Шаг 7.** [**Запись циклов во множество изометрических циклов**]. Проверяем сформированные на предыдущем шаге циклы с ранее записанными циклами во множестве изометрических циклов, и в случае их отсутствия добавляем их во множество изометрических циклов. Идем на шаг 1.

***Пример 1.1.*** В качестве примера рассмотрим граф G (рис. 1.11).

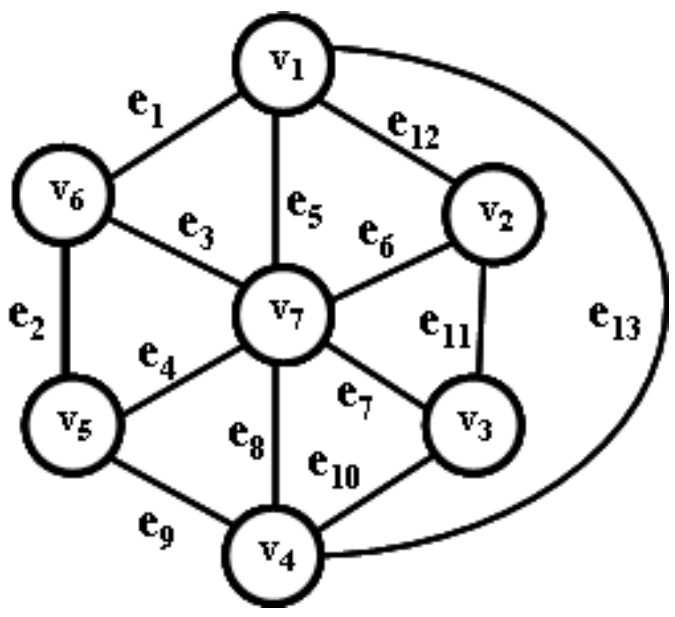

Рис. 1.11. Граф **G**.

Если в качестве выбранного ребра взять ребро $e_{13}$, то процесс разметки вершин имеет вид, представленный на рис. 1.12. Система циклов, проходящих по ребру $e_{13}$, для разметки, показанной на рис. 1.12,*а* (номера разметки представлены внутри вершины):

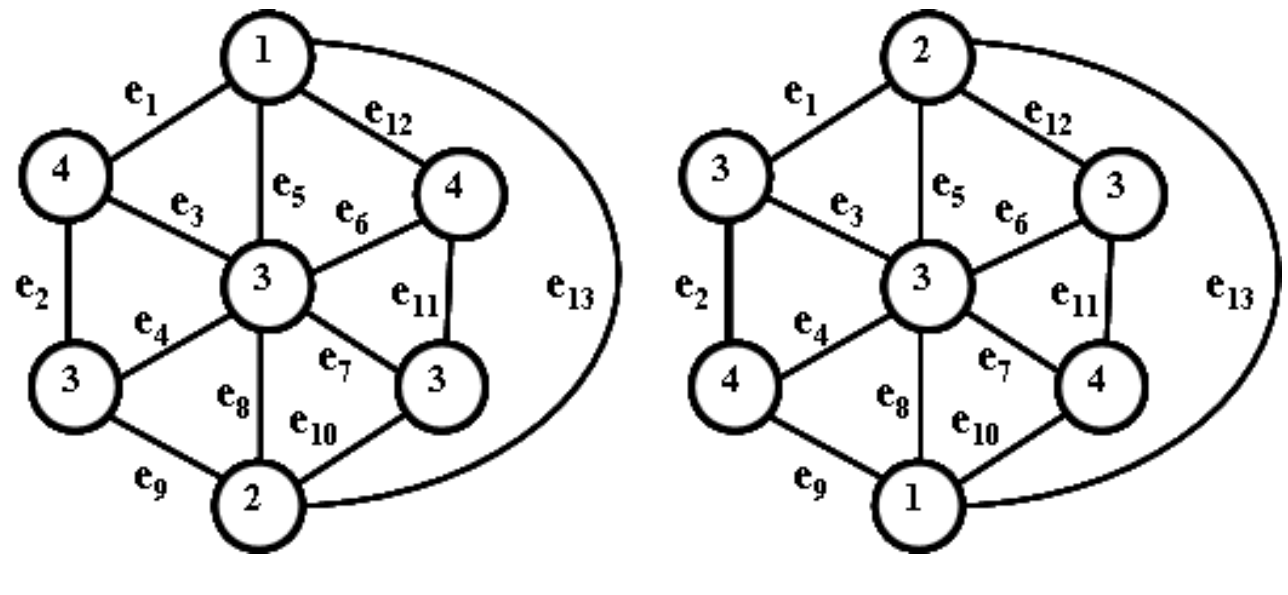

а) б)

Рис. 1.12. Прямой и обратный процесс разметки вершин
(числа в вершинах фронт волны).

$C_{13}^1 = \{\{e_5,e_8,e_{13}\}, \{e_1,e_3,e_8,e_{13}\}, \{e_1,e_2,e_9,e_{13}\}, \{e_6,e_8,e_{12},e_{13}\},\{e_{10},e_{11},e_{12},e_{13}\}\}.$

Система циклов, проходящих по ребру $e_{13}$, для разметки, представленной на рис. 1.11,б (номера разметки представлены внутри вершины):

$C_{13}^2 = \{\{e_5,e_8,e_{13}\},\{e_4,e_5,e_9,e_{13}\},\{e_1,e_2,e_9,e_{13}\},\{e_{10},e_{11},e_{12},e_{13}\},\{e_5,e_7,e_{10},e_{13}\}\}.$

Пересечение множеств $C_{13}^1$ и $C_{13}^2$:

$C_{13} = C_{13}^1 \cap C_{13}^2 = \{\{e_5,e_8,e_{13}\},\{e_1,e_2,e_9,e_{13}\},\{e_{10},e_{11},e_{12},e_{13}\}\}.$

Каждому ребру принадлежат следующие изометрические циклы:

$C_1 = \{\{e_1,e_3,e_5\},\{e_1,e_2,e_9,e_{13}\}\};$
$C_2 = \{\{e_2,e_3,e_4\},\{e_1,e_2,e_9,e_{13}\}\};$
$C_3 = \{\{e_2,e_3,e_4\},\{e_1,e_3,e_5\}\};$
$C_4 = \{\{e_4,e_8,e_9\},\{e_2,e_3,e_4\}\};$
$C_5 = \{\{e_1,e_3,e_5\},\{e_5,e_6,e_{12}\},\{e_5,e_8,e_{13}\}\};$
$C_6 = \{\{e_5,e_6,e_{12}\},\{e_6,e_7,e_{11}\}\};$
$C_7 = \{\{e_6,e_7,e_{11}\},\{e_7,e_8,e_{10}\}\};$
$C_8 = \{\{e_5,e_8,e_{13}\},\{e_7,e_8,e_{10}\},\{e_3,e_8,e_9\}\};$
$C_9 = \{\{e_4,e_8,e_9\},\{e_1,e_2,e_9,e_{13}\}\};$
$C_{10} = \{\{e_7,e_8,e_{10}\},\{e_{10},e_{11},e_{12},e_{13}\}\};$
$C_{11} = \{\{e_6,e_7,e_{11}\},\{e_{10},e_{11},e_{12},e_{13}\}\};$
$C_{12} = \{\{e_5,e_6,e_{12}\},\{e_{10},e_{11},e_{12},e_{13}\}\};$
$C_{13} = \{\{e_5,e_8,e_{13}\},\{e_1,e_2,e_9,e_{13}\},\{e_{10},e_{11},e_{12},e_{13}\}\}.$

Множество изометрических циклов получим как объединение:

$C_\tau = C_1 \cup C_2 \cup C_3 \cup C_4 \cup C_5 \cup C_6 \cup C_7 \cup C_8 \cup C_9 \cup C_{10} \cup C_{11} \cup C_{12} \cup C_{13} =$
$= \{\{e_1,e_3,e_5\}, \{e_2,e_3,e_4\}, \{e_4,e_8,e_9\}, \{e_5,e_6,e_{12}\}, \{e_5,e_8,e_{13}\}, \{e_6,e_7,e_{11}\},$
$\{e_7,e_8,e_{10}\},\{e_1,e_2,e_9,e_{13}\}, \{e_{10},e_{11},e_{12},e_{13}\}\}.$

Таким образом, множество изометрических циклов состоит из 9-ти элементов. Цикломатическое число графа G равно 7. Следовательно, для построения базиса нужно удалить два изометрических цикла. Очевидно, что для любого трехсвязного и более графа G, множество изометрических циклов имеет мощность меньшую, чем мощность множества простых циклов, но большую или равную, цикломатическому числу графа:

$$\nu(G) \le \text{card } C_\tau \le \text{card } C_R \le \text{card } C. \qquad (1.10)$$

Теперь покажем, что построение множества изометрических циклов должно производиться относительно всего множества ребер графа.

***Пример 1.2.*** Следующий пример демонстрирует невозможность получения полного множества изометрических циклов, если построение производится только относительно хорд для выбранного дерева графа. Рассмотрим граф, представленный на рис. 1.13.

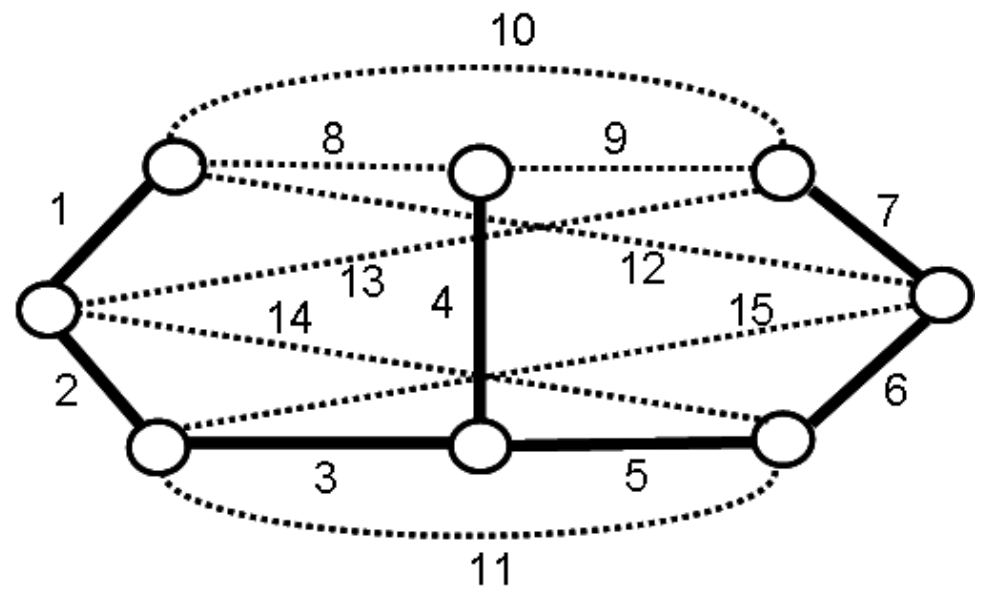

Рис. 1.13. Граф G и его дерево.

Рассмотрим изометрические циклы относительно 8-ой хорды:

| | |
|---|---|
| $c_1 = \{e_8,e_9,e_{10}\}$; | $c_1' = \{e_8,e_9,e_{10}\}$; |
| $c_2 = \{e_1,e_8,e_9,e_{13}\}$; | $c_2' = \{e_4,e_5,e_6,e_8,e_{12}\}$; |
| $c_3 = \{e_7,e_8,e_9,e_{12}\}$; | $c_3' = \{e_1,e_2,e_3,e_4,e_8\}$; |
| | $c_4' = \{e_4,e_5,e_8,e_{12},e_{15}\}$; |
| | $c_5' = \{e_4,e_5,e_9,e_{13},e_{14}\}$. |

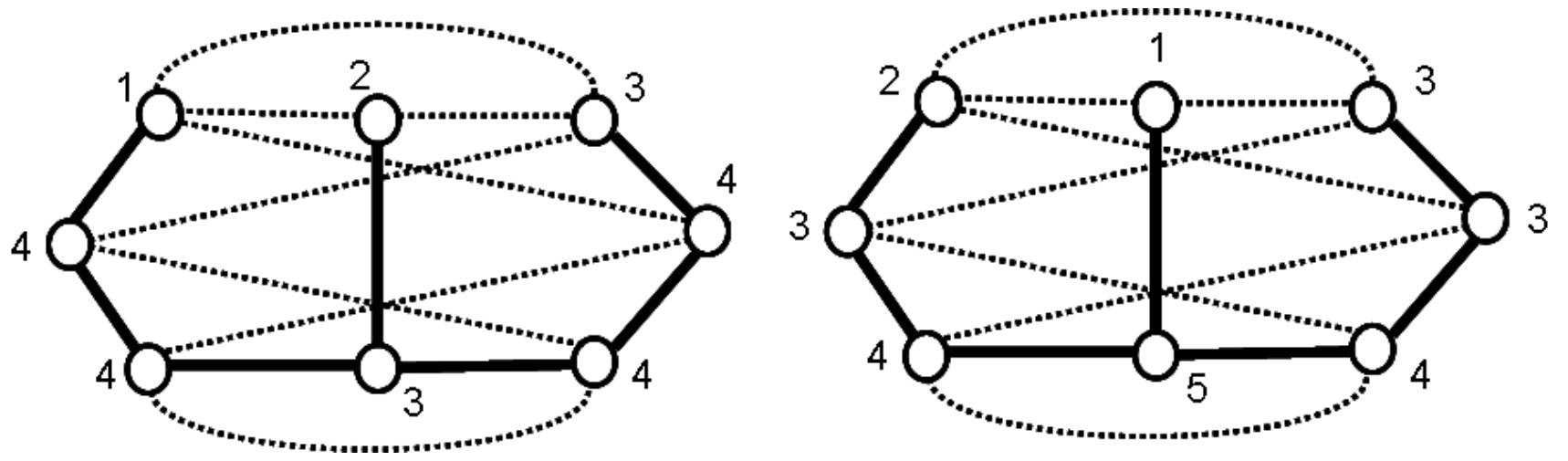

Рис. 1.14. Изометрические циклы относительно 8-ой хорды.

Изометрический цикл относительно 8-ой хорды: $\{e_8,e_9,e_{10}\}$.

Рассмотрим изометрические циклы относительно 9-ой хорды:

| | |
|---|---|
| $c_1 = \{e_8,e_9,e_{10}\}$; | $c_1' = \{e_8,e_9,e_{10}\}$; |
| $c_2 = \{e_4,e_5,e_6,e_8,e_{12}\}$; | $c_2' = \{e_1,e_8,e_9,e_{13}\}$; |
| $c_3 = \{e_2,e_3,e_4,e_9,e_{13}\}$; | $c_3' = \{e_7,e_8,e_9,e_{12}\}$. |
| $c_4 = \{e_3,e_4,e_7,e_9,e_{15}\}$; | |
| $c_5 = \{e_3,e_4,e_9,e_{12},e_{15}\}$; | |

Изометрический цикл относительно 9-ой хорды: $\{e_8,e_9,e_{10}\}$.

Рассмотрим изометрические циклы относительно 10-ой хорды:

| | |
|---|---|
| $c_1 = \{e_8,e_9,e_{10}\}$; | $c_1' = \{e_8,e_9,e_{10}\}$; |
| $c_2 = \{e_7,e_{10},e_{12}\}$; | $c_2' = \{e_7,e_{10},e_{12}\}$; |
| $c_3 = \{e_1,e_{10},e_{13}\}$; | $c_3' = \{e_1,e_{10},e_{13}\}$. |

Изометрические циклы относительно 10-ой хорды: $\{e_8,e_9,e_{10}\}$, $\{e_1,e_{10},e_{13}\}$, $\{e_7,e_{10},e_{12}\}$.

Рассмотрим изометрические циклы относительно 11-ой хорды:

| | |
|---|---|
| $c_1 = \{e_3,e_5,e_{11}\}$; | $c_1' = \{e_3,e_5,e_{11}\}$; |
| $c_2 = \{e_2,e_{11},e_{14}\}$; | $c_2' = \{e_2,e_{11},e_{14}\}$; |
| $c_3 = \{e_6,e_{11},e_{15}\}$; | $c_3' = \{e_6,e_{11},e_{15}\}$. |

Изометрические циклы относительно 11-ой хорды: $\{e_3,e_5,e_{11}\}$, $\{e_2,e_{11},e_{14}\}$, $\{e_6,e_{11},e_{15}\}$.

Рассмотрим изометрические циклы относительно 12-ой хорды:

| | |
|---|---|
| $c_1 = \{e_7,e_{10},e_{12}\}$; | $c_1' = \{e_7,e_{10},e_{12}\}$; |
| $c_2 = \{e_7,e_8,e_9,e_{12}\}$; | $c_2' = \{e_1,e_2,e_{12},e_{15}\}$; |
| $c_3 = \{e_1,e_2,e_{12},e_{15}\}$; | $c_3' = \{e_1,e_6,e_{12},e_{14}\}$. |
| $c_4 = \{e_1,e_7,e_{12},e_{13}\}$; | |
| $c_5 = \{e_1,e_6,e_{12},e_{14}\}$; | |

Изометрические циклы относительно 12-ой хорды: $\{e_7,e_{10},e_{12}\}$, $\{e_1,e_2,e_{12},e_{15}\}$, $\{e_1,e_6,e_{12},e_{14}\}$.

Рассмотрим изометрические циклы относительно 13-ой хорды:

| | |
|---|---|
| $c_1 = \{e_1,e_{10},e_{13}\}$; | $c_1' = \{e_1,e_{10},e_{13}\}$; |
| $c_2 = \{e_6,e_7,e_{13},e_{14}\}$; | $c_2' = \{e_1,e_7,e_{12},e_{13}\}$; |
| $c_3 = \{e_2,e_7,e_{13},e_{15}\}$; | $c_3' = \{e_6,e_7,e_{13},e_{14}\}$; |
| | $c_4' = \{e_2,e_7,e_{13},e_{15}\}$; |
| | $c_5' = \{e_1,e_8,e_9,e_{13}\}$. |

Изометрические циклы относительно 13-ой хорды: $\{e_1,e_{10},e_{13}\}$, $\{e_6,e_7,e_{13},e_{14}\}$, $\{e_2,e_7,e_{13},e_{15}\}$.

Рассмотрим изометрические циклы относительно 14-ой хорды:

| | |
|---|---|
| $c_1 = \{e_2,e_{11},e_{14}\}$; | $c_1' = \{e_2,e_{11},e_{14}\}$; |
| $c_2 = \{e_6,e_7,e_{13},e_{14}\}$; | $c_2' = \{e_2,e_6,e_{14},e_{15}\}$; |
| $c_3 = \{e_1,e_6,e_{12},e_{14}\}$; | $c_3' = \{e_2,e_3,e_5,e_{14}\}$; |
| | $c_4' = \{e_6,e_7,e_{13},e_{14}\}$; |
| | $c_5' = \{e_1,e_6,e_{12},e_{14}\}$. |

Изометрические циклы относительно 14-ой хорды: $\{e_2,e_{11},e_{14}\}$, $\{e_6,e_7,e_{13},e_{14}\}$, $\{e_1,e_6,e_{12},e_{14}\}$.

Рассмотрим изометрические циклы относительно 15-ой хорды:

| | |
|---|---|
| $c_1 = \{e_6,e_{11},e_{15}\}$; | $c_1' = \{e_6,e_{11},e_{15}\}$; |
| $c_2 = \{e_3,e_5,e_6,e_{15}\}$; | $c_2' = \{e_2,e_7,e_{13},e_{15}\}$; |
| $c_3 = \{e_2,e_6,e_{14},e_{15}\}$; | $c_3' = \{e_1,e_2,e_{12},e_{15}\}$. |
| $c_4 = \{e_2,e_7,e_{13},e_{15}\}$; | |
| $c_5 = \{e_1,e_2,e_{12},e_{15}\}$; | |

Изометрические циклы относительно 15-ой хорды: $\{e_6,e_{11},e_{15}\}$, $\{e_2,e_7,e_{13},e_{15}\}$, $\{e_1,e_2,e_{12},e_{15}\}$.

Рассмотрим изометрические циклы относительно 4-ого ребра:

| | |
|---|---|
| $c_1 = \{e_1,e_2,e_3,e_4,e_8\}$; | $c_1' = \{e_1,e_2,e_3,e_4,e_8\}$; |
| $c_2 = \{e_3,e_4,e_8,e_{12},e_{15}\}$; | $c_2' = \{e_3,e_4,e_8,e_{12},e_{15}\}$; |
| $c_3 = \{e_4,e_5,e_6,e_8,e_{12}\}$; | $c_3' = \{e_4,e_5,e_6,e_8,e_{12}\}$; |
| $c_4 = \{e_4,e_5,e_6,e_7,e_9\}$; | $c_4' = \{e_4,e_5,e_6,e_7,e_9\}$; |
| $c_5 = \{e_2,e_3,e_4,e_9,e_{13}\}$; | $c_5' = \{e_2,e_3,e_4,e_9,e_{13}\}$; |
| $c_6 = \{e_4,e_5,e_9,e_{13},e_{14}\}$; | $c_6' = \{e_4,e_5,e_9,e_{13},e_{14}\}$. |

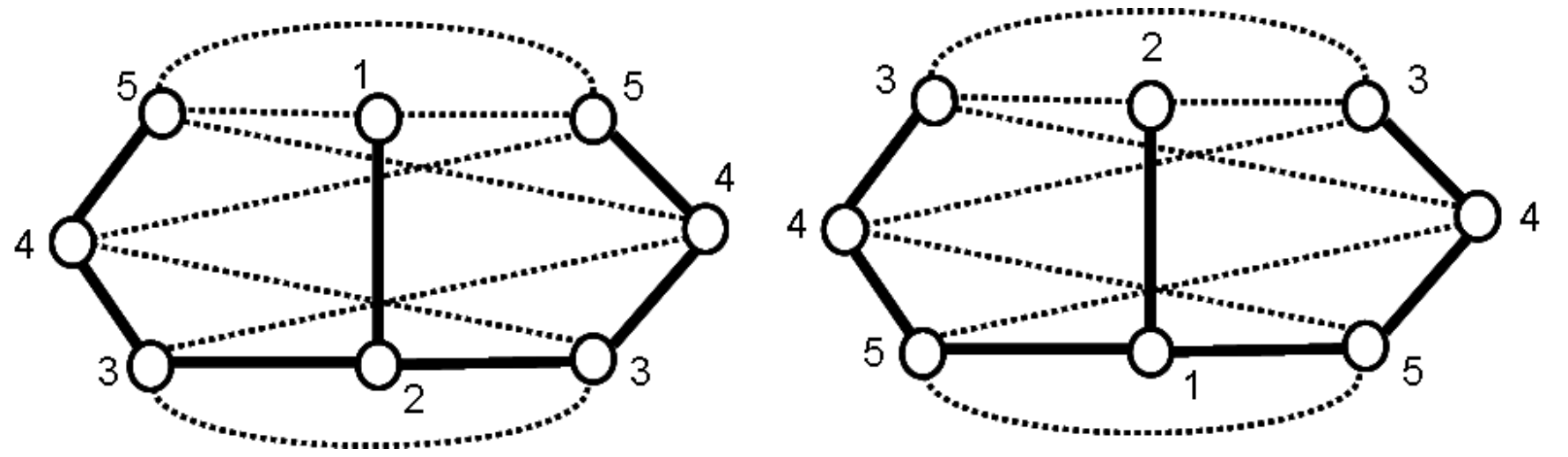

Рис. 1.15. Изометрические циклы относительно 4-го ребра дерева.

Изометрические циклы относительно 4-го ребра: $\{e_1,e_2,e_3,e_4,e_8\}$, $\{e_3,e_4,e_8,e_{12},e_{15}\}$, $\{e_4,e_5,e_6,e_8,e_{12}\}$, $\{e_4,e_5,e_6,e_7,e_9\}$, $\{e_2,e_3,e_4,e_9,e_{13}\}$, $\{e_4,e_5,e_9,e_{13},e_{14}\}$.

Как видно из данного примера, если построение производится только относительно хорд для выбранного дерева графа, то множество изометрических циклов будет не полно. В данном примере, во множество изометрических циклов не вошли изометрические циклы, проходящие по четвертому ребру.

### 1.6. Свойства изометрических циклов и центральных разрезов

Рассмотрим более подробно свойства изометрических циклов и центральных разрезов на примере графа $K_5$.

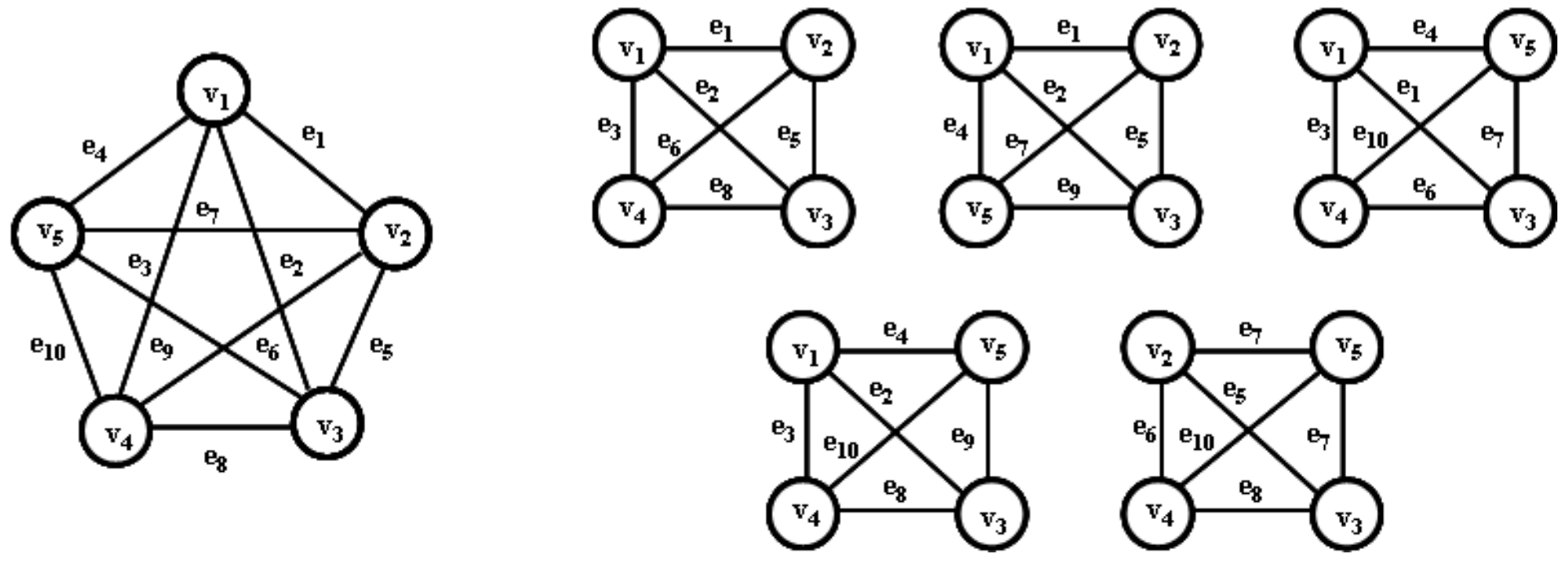

Рис. 1.16. Граф $K_5$ и его пять четырехвершинных подграфов.

Множество центральных разрезов **S** для графа $K_5$:

$s_1 = \{e_1,e_2,e_3,e_4\}$; $s_2 = \{e_1,e_5,e_6,e_7\}$; $s_3 = \{e_2,e_5,e_8,e_9\}$; $s_4 = \{e_3,e_6,e_8,e_{10}\}$;

$s_5 = \{e_4,e_7,e_9,e_{10}\}$.

Множество изометрических циклов $C_\tau$ для графа $K_5$:

$c_1 = \{e_1,e_2,e_5\}$; $c_2 = \{e_1,e_3,e_6\}$; $c_3 = \{e_1,e_4,e_7\}$; $c_4 = \{e_2,e_3,e_8\}$; $c_5 = \{e_2,e_4,e_9\}$;

$c_6 = \{e_3,e_4,e_{10}\}$; $c_7 = \{e_5,e_6,e_8\}$; $c_8 = \{e_5,e_7,e_9\}$; $c_9 = \{e_6,e_7,e_{10}\}$; $c_{10} = \{e_8,e_9,e_{10}\}$.

Как видно, все изометрические циклы принадлежат четырехвершинным подграфам (см. рис. 1.16).

Любой граф G можно получить из полного графа путем удаления соответствующих ребер.

Проиллюстрируем данный процесс на примере графа $K_5$. Удалим из графа ребро $e_{10}$ (см. рис. 1.17).

Из множества изометрических циклов C удаляются все циклы, включающие 10-е ребро:

$c_6 = \{e_3,e_4,e_{10}\}$; $c_9 = \{e_6,e_7,e_{10}\}$; $c_{10} = \{e_8,e_9,e_{10}\}$.

Остаются изометрические циклы:

$c_1 = \{e_1,e_2,e_5\}$; $c_2 = \{e_1,e_3,e_6\}$; $c_3 = \{e_1,e_4,e_7\}$; $c_4 = \{e_2,e_3,e_8\}$; $c_5 = \{e_2,e_4,e_9\}$;
$c_7 = \{e_5,e_6,e_8\}$; $c_8 = \{e_5,e_7,e_9\}$.

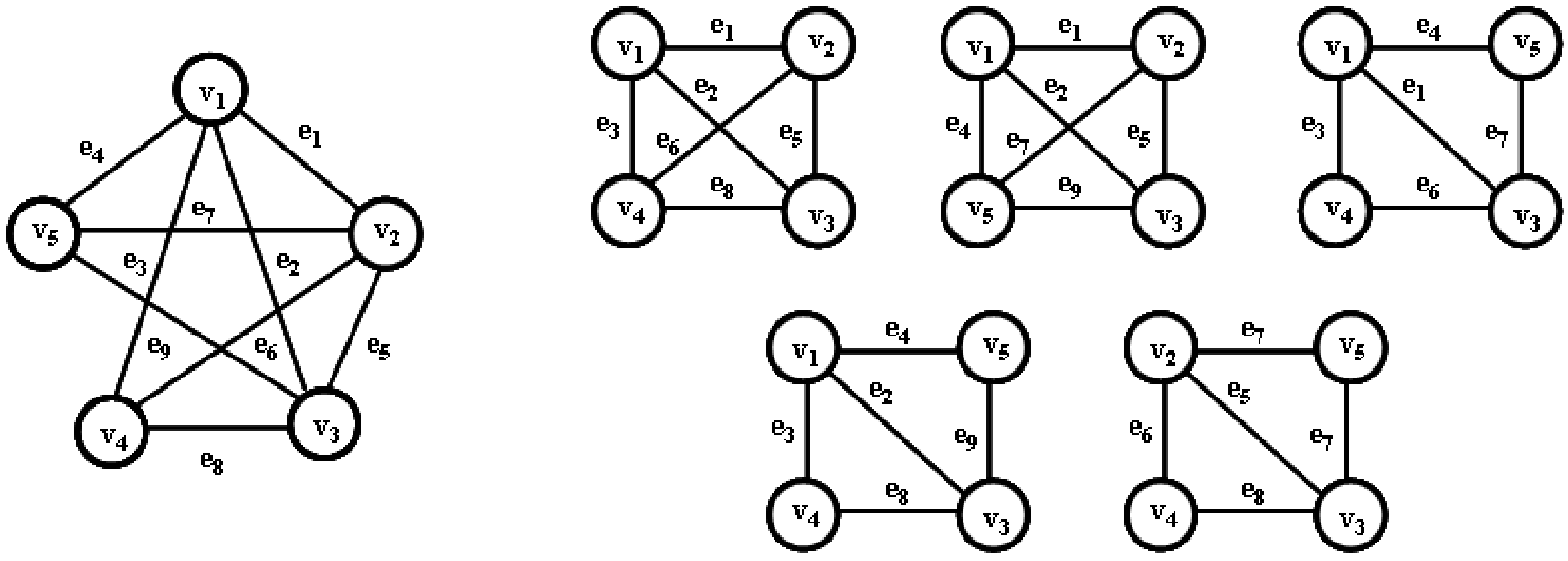

Рис. 1.17. Граф $K_5$ с удаленным ребром $e_{10}$ и его пять четырехвершинных подграфов.

И должны, в перспективе, образоваться новые изометрические циклы длиной четыре:

$c_6 \oplus c_9 = \{e_3,e_4,e_{10}\} \oplus \{e_6,e_7,e_{10}\} = \{e_3,e_4,e_6,e_7\}$;
$c_6 \oplus c_{10} = \{e_3,e_4,e_{10}\} \oplus \{e_8,e_9,e_{10}\} = \{e_3,e_4,e_8,e_9\}$;
$c_9 \oplus c_{10} = \{e_6,e_7,e_{10}\} \oplus \{e_8,e_9,e_{10}\} = \{e_6,e_7,e_8,e_9\}$.

Однако, их включение во множество оставшихся изометрических циклов невозможно, так как они образованы как результат попарного кольцевого суммирования оставшихся изометрических циклов (см. рис. 1.17):

$c_2 \oplus c_3 = \{e_1,e_3,e_6\} \oplus \{e_1,e_4,e_7\} = \{e_3,e_4,e_6,e_7\}$;
$c_4 \oplus c_5 = \{e_2,e_3,e_8\} \oplus \{e_2,e_4,e_9\} = \{e_3,e_4,e_8,e_9\}$;
$c_7 \oplus c_8 = \{e_5,e_6,e_8\} \oplus \{e_5,e_7,e_9\} = \{e_6,e_7,e_8,e_9\}$.

Таким образом, множество изометрических циклов $C_\tau$ для графа $G_1$ (см. рис. 1.17), полученного путем удаления 10-го ребра из графа $K_5$, состоит из следующих изометрических циклов:

$c_1 = \{e_1,e_2,e_5\}$; $c_2 = \{e_1,e_3,e_6\}$; $c_3 = \{e_1,e_4,e_7\}$; $c_4 = \{e_2,e_3,e_8\}$; $c_5 = \{e_2,e_4,e_9\}$;
$c_7 = \{e_5,e_6,e_8\}$; $c_8 = \{e_5,e_7,e_9\}$.

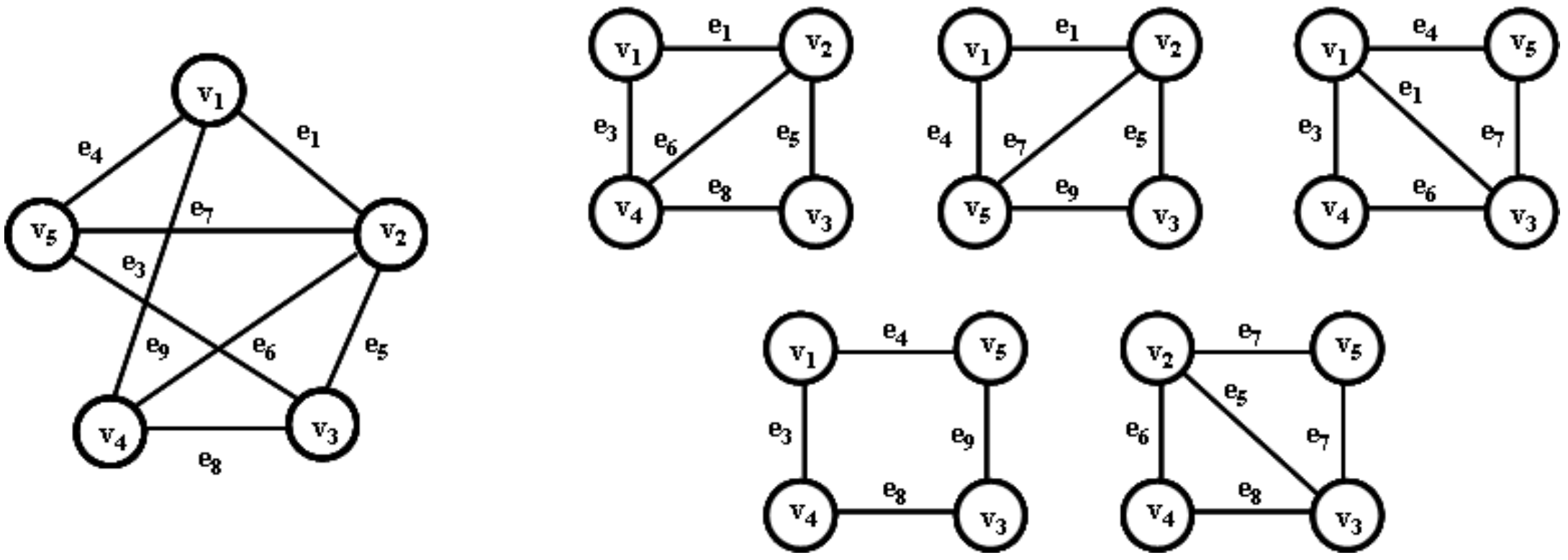

Рис. 1.18. Граф $K_5$ с удаленными $e_{10}$ и $e_2$ ребрами и его пять четырехвершинных подграфов.

Продолжаем удалять ребра из графа $K_5$. Удалим ребро $e_2$ (см. рис. 1.18).

Из множества изометрических циклов $C_\tau$ = $\{c_1,c_2,c_3,c_4,c_5,c_7,c_8\}$ удаляются все циклы, включающие 2-е ребро:

$c_1 = \{e_1,e_2,e_5\}$; $c_4 = \{e_2,e_3,e_8\}$; $c_5 = \{e_2,e_4,e_9\}$.

Остаются изометрические циклы:

$c_2 = \{e_1,e_3,e_6\}$; $c_3 = \{e_1,e_4,e_7\}$; $c_7 = \{e_5,e_6,e_8\}$; $c_8 = \{e_5,e_7,e_9\}$.

И должны, в перспективе, образоваться новые изометрические циклы длиной четыре:

$c_1 \oplus c_4 = \{e_1,e_2,e_5\} \oplus \{e_2,e_3,e_8\} = \{e_1,e_3,e_5,e_8\}$;
$c_1 \oplus c_5 = \{e_1,e_2,e_5\} \oplus \{e_2,e_4,e_9\} = \{e_1,e_4,e_5,e_9\}$;
$c_4 \oplus c_5 = \{e_2,e_3,e_8\} \oplus \{e_2,e_4,e_9\} = \{e_3,e_4,e_8,e_9\}$.

Однако, включение во множество $C_\tau$ оставшихся изометрических циклов двух первых невозможно, так как они образованы как результат попарного кольцевого суммирования оставшихся изометрических циклов (см. рис. 1.18):

$c_1 \oplus c_4 = \{e_1,e_3,e_6\} \oplus \{e_5,e_6,e_8\} = \{e_1,e_3,e_5,e_8\}$;
$c_1 \oplus c_5 = \{e_1,e_4,e_7\} \oplus \{e_5,e_7,e_9\} = \{e_1,e_4,e_5,e_9\}$.

А вот вновь образованный изометрический цикл $c_{4,5}$ может быть включен во множество оставшихся изометрических циклов: $c_4 \oplus c_5 = \{e_3,e_4,e_8,e_9\}$, так как для него не существует попарного кольцевого суммирования оставшихся изометрических циклов.

Таким образом, множество изометрических циклов $C_\tau$ для графа $G_2$ (см. рис. 1.19), полученного путем удаления 10-го и 2-го ребер из графа $K_5$, состоит из следующих изометрических циклов:

$c_2 = \{e_1,e_3,e_6\}$; $c_3 = \{e_1,e_4,e_7\}$; $c_7 = \{e_5,e_6,e_8\}$; $c_8 = \{e_5,e_7,e_9\}$; $c_{4,5} = \{e_3,e_4,e_8,e_9\}$.

Продолжаем удалять ребра из графа $K_5$. Удалим ребро $e_7$ (см. рис. 1.19).

Из множества изометрических циклов $C_2$ = $\{c_2,c_3,c_7,c_8\ c_{4,5}\}$ удаляются все циклы, включающие 7-ое ребро: $c_3 = \{e_1,e_4,e_7\}$; $c_8 = \{e_5,e_7,e_9\}$.

Остаются изометрические циклы:

$c_2 = \{e_1,e_3,e_6\}$; $c_7 = \{e_5,e_6,e_8\}$, $c_{4,5} = \{e_3,e_4,e_8,e_9\}$.

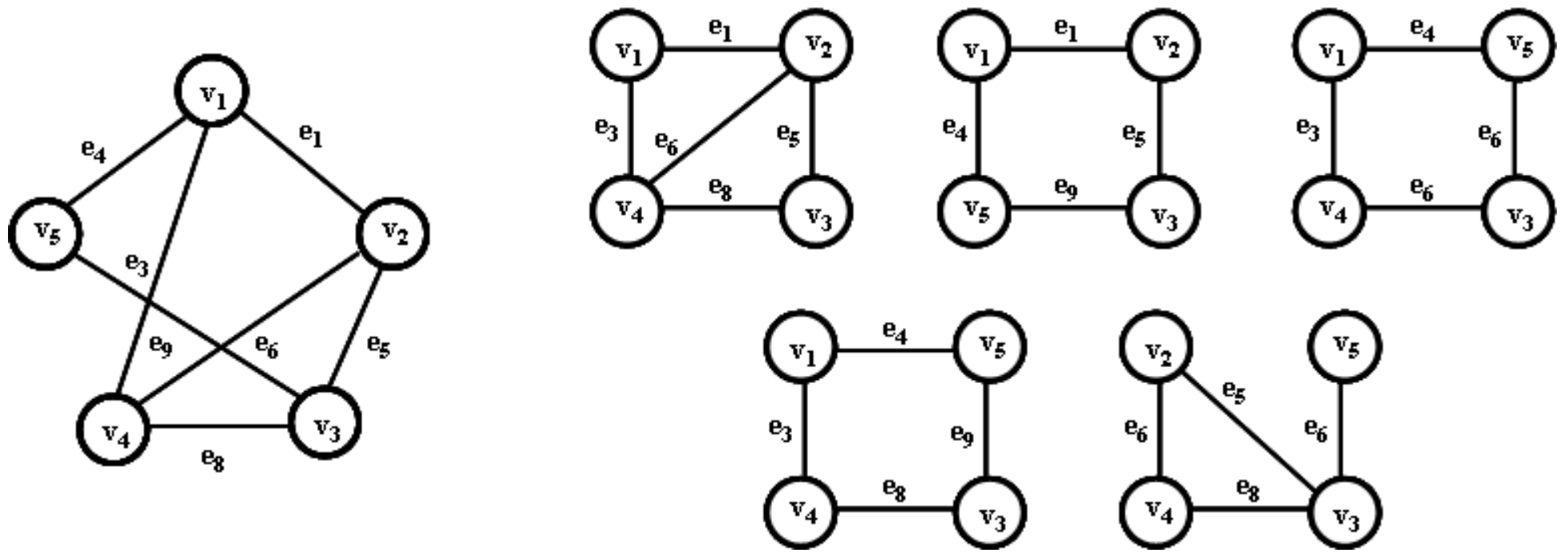

Рис. 1.19. Граф $K_5$ с удаленными $e_{10},e_2,e_7$ ребрами и его пять четырехвершинных подграфов.

И должны, в перспективе, образоваться новые изометрические циклы:

$c_3 \oplus c_8 = \{e_1,e_4,e_7\} \oplus \{e_5,e_7,e_9\} = \{e_1,e_4,e_5,e_9\}$.

Данный цикл $c_{3,8}$ может быть включен во множество оставшихся изометрических циклов, так как он не может быть образован как результат попарного кольцевого суммирования из оставшихся изометрических циклов (см. рис. 1.19):

Таким образом, множество изометрических циклов $C_\tau$ для графа $G_3$, полученного путем удаления 10-го, 2-го и 7-го ребер из графа $K_5$, состоит из следующих изометрических циклов:

$c_2 = \{e_1,e_3,e_6\}$; $c_7 = \{e_5,e_6,e_8\}$, $c_{4,5} = \{e_3,e_4,e_8,e_9\}$; $c_{3,8} = \{e_1,e_4,e_5,e_9\}$.

Продолжаем удалять ребра из графа $K_5$. Удалим ребро $e_1$ (см. рис. 1.20).

Из множества изометрических циклов $C_\tau$ = $\{c_2,c_7,c_{4,5},c_{3,8}\}$ удаляются все циклы, включающие 1-е ребро:

$c_2 = \{e_1,e_3,e_6\}$; $c_{3,8} = \{e_1,e_4,e_5,e_9\}$.

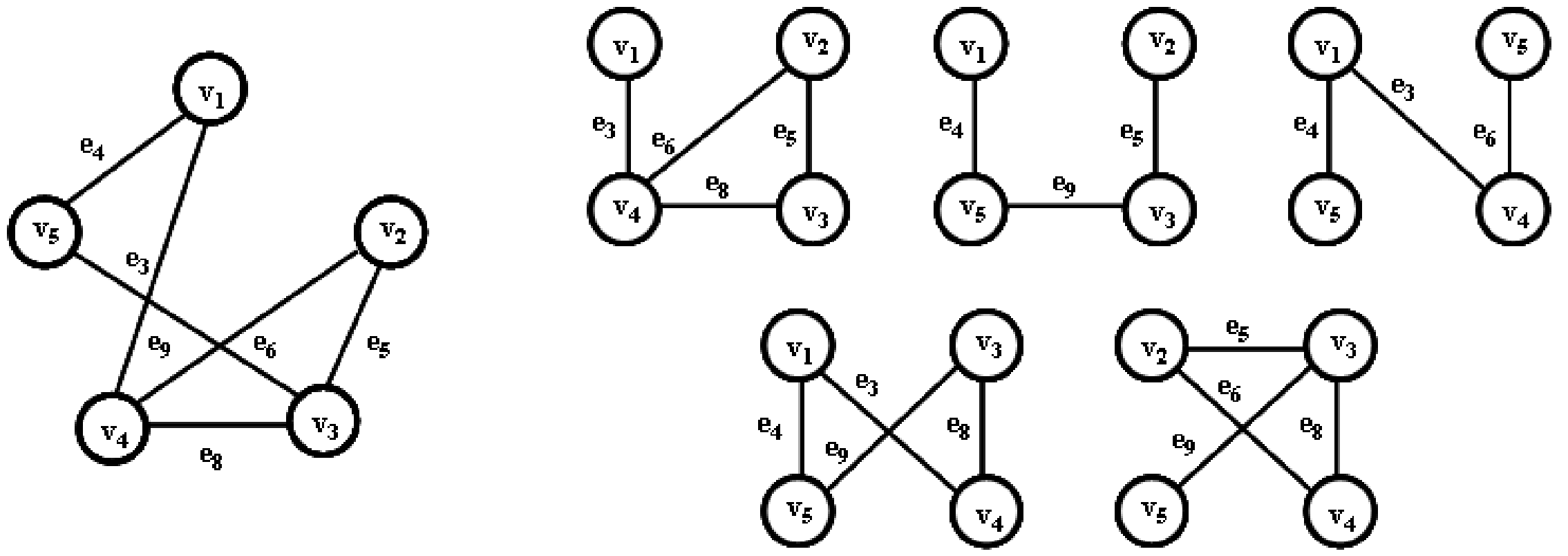

Рис. 1.20. Граф $K_5$ с удаленными $e_{10},e_2,e_7,e_1$ ребрами и его пять четырехвершинных подграфов.

Остаются изометрические циклы: $c_7 = \{e_5,e_6,e_8\}$, $c_{4,5} = \{e_3,e_4,e_8,e_9\}$.

И должны, в перспективе, образоваться новые изометрические циклы:

$c_2 \oplus c_{3,8} = \{e_1,e_3,e_6\} \oplus \{e_1,e_4,e_5,e_9\} = \{e_3,e_4,e_5,e_6,e_9\}$.

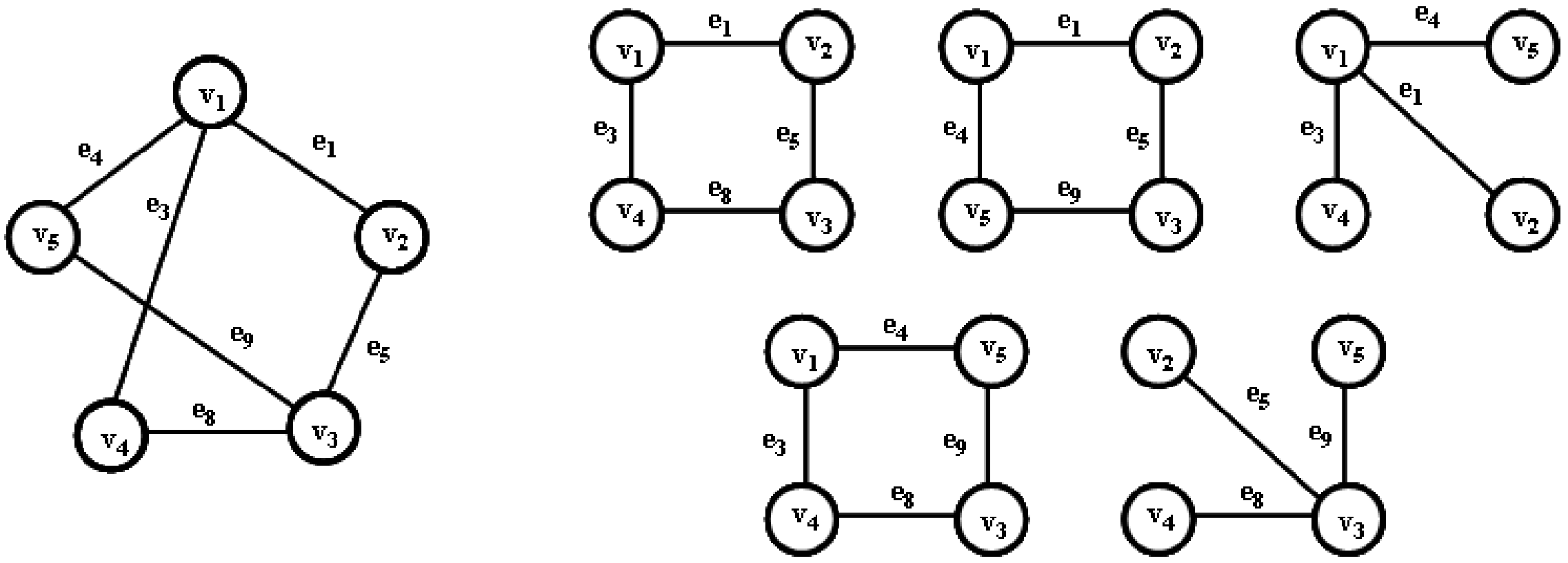

Рис. 1.21. Граф $K_5$ с удаленными $e_{10},e_2,e_7,e_6$ ребрами и его пять четырехвершинных подграфов.

Дальнейшее образование путем удаления ребер можно приостановить, так как полученный цикл равен кольцевой сумме оставшихся циклов:

$c_7 \oplus c_{4,5} = \{e_5,e_6,e_8\} \oplus \{e_3,e_4,e_8,e_9\} = \{e_3,e_4,e_5,e_6,e_9\}$.

Таким образом, множество изометрических циклов для графа $G_5$ (см. рис. 1.20), полученного путем удаления 10-го, 2-го, 7-го и 1-го ребер из графа $K_5$ состоит из следующих изометрических циклов: $c_7 = \{e_5,e_6,e_8\}$ и $c_{4,5} = \{e_3,e_4,e_8,e_9\}$.

Если удалить 6-ое ребро вместо 1-го ребра (см. рис. 1.21), то из множества изометрических циклов $C_\tau = \{c_2,c_7, c_{4,5},c_{3,8}\}$ удаляются все циклы, включающие 6-е ребро:

$c_2 = \{e_1,e_3,e_6\}$; $c_7 = \{e_5,e_6,e_8\}$.

Остаются изометрические циклы:

$c_{4,5} = \{e_3,e_4,e_8,e_9\}$; $c_{3,8} = \{e_1,e_4,e_5,e_9\}$.

И должны, в перспективе, образоваться новые изометрические циклы:

$c_2 \oplus c_7 = \{e_1,e_4,e_7\} \oplus \{e_5,e_7,e_9\} = \{e_1,e_3,e_5,e_8\}$.

Дальнейшее образование путем удаления ребер можно приостановить, так как полученный цикл $c_{2,}c_7$ равен кольцевой сумме оставшихся циклов:

$c_{4,5} \oplus c_{3,8} = \{e_3,e_4,e_8,e_9\} \oplus \{e_1,e_4,e_5,e_9\} = \{e_1,e_3,e_5,e_8\}$.

Таким образом, множество изометрических циклов для графа $G_5$, полученного путем удаления 10-го, 2-го, 7-го и 6-го ребер из графа $K_5$, состоит из следующих изометрических циклов:

$c_{4,5} = \{e_3,e_4,e_8,e_9\}$ и $c_{3,8} = \{e_1,e_4,e_5,e_9\}$.

Аналогичные рассуждения можно провести для любого графа G с *n* вершинами, удаляя соответствующие ребра из полного графа $K_n$.

Опишем алгоритм формирования множества изометрических циклов из полного графа

методом удаления ребер

**[Инициализация].** Существует множество изометрических циклов полного графа $K_n$. Изометрические циклы записанны в виде подмножества ребер $M_e$ и в виде подмножество вершин $M_v$, хранящихся в виде связного списка. Имеется массив номеров исключаемых ребер из полного графа.

**Шаг 1.** [**Перебор исключаемых ребер**]. Последовательно выбираем текущее ребро для исключения. Если список ребер не исчерпан, идем на шаг 2. Иначе конец работы алгоритма.

**Шаг 2.** [**Выбор исключенного ребра графа**]. Выбираем очередное исключаемое ребро графа. Формируем множество новых простых циклов как кольцевую сумму всех попарно пересекающихся по данному ребру изометрических циклов в записи по ребрам. Одновременно формируем такое же связное множество циклов в записи по вершинам, как объединение выбранных циклов по вершинам. Идем на шаг 2.

**Шаг 3.** [**Удаление изометрических циклов**]. Связно удаляем изометрические циклы, имеющие исключаемое ребро. Идем на шаг 4.

**Шаг 4.** [**Проверка на включение**]. Проверяем множество новых простых циклов на включение с множеством изометрических циклов по записям в виде вершин. Если включение имеется, то такие простые циклы не включаются во множество изометрических циклов. Если включения нет, то производится запись в массив изометрических циклов. Идем на шаг 1.

Сказанное рассмотрим на примере графа $K_5$ (см. рис. 1.16).

Пусть задано множество ребер $\{e_2,e_6,e_7\}$ для удаления из полного графа.

Связное множество изометрических циклов запишем в виде связного списка:

| **Цикл** | **Множество $M_e$** | | **Множество $M_v$** |
|---|---|---|---|
| $c_1$ | $\{e_1,e_2,e_5\}$ | $\rightarrow$ | $\{v_1,v_2,v_3\}$ |
| $c_2$ | $\{e_1,e_3,e_6\}$ | $\rightarrow$ | $\{v_1,v_2,v_4\}$ |
| $c_3$ | $\{e_1,e_4,e_7\}$ | $\rightarrow$ | $\{v_1,v_2,v_5\}$ |
| $c_4$ | $\{e_2,e_3,e_8\}$ | $\rightarrow$ | $\{v_1,v_3,v_4\}$ |
| $c_5$ | $\{e_2,e_4,e_9\}$ | $\rightarrow$ | $\{v_1,v_3,v_5\}$ |
| $c_6$ | $\{e_3,e_4,e_{10}\}$ | $\rightarrow$ | $\{v_1,v_4,v_5\}$ |
| $c_7$ | $\{e_5,e_6,e_8\}$ | $\rightarrow$ | $\{v_2,v_3,v_4\}$ |
| $c_8$ | $\{e_5,e_7,e_9\}$ | $\rightarrow$ | $\{v_2,v_3,v_5\}$ |
| $c_9$ | $\{e_6,e_7,e_{10}\}$; | $\rightarrow$ | $\{v_2,v_4,v_5\}$ |
| $c_{10}$ | $\{e_8,e_9,e_{10}\}$ | $\rightarrow$ | $\{v_3,v_4,v_5\}$ |

Исключаем второе ребро $e_2$. Формируем новое множество простых циклов:

| **Цикл** | | **Новое $M_e$** | | | **Новое $M_v$** |
|---|---|---|---|---|---|
| $c_1 \oplus c_4 =$ | $\{e_1,e_2,e_5\} \oplus \{e_2,e_3,e_8\}$ | $=\{e_1,e_3,e_5,e_8\}$ | $\rightarrow$ | $\{v_1,v_2,v_3\} \cup \{v_1,v_3,v_4\}$ | $=\{v_1,v_2,v_3,v_4\}$ |
| $c_1 \oplus c_5 =$ | $\{e_1,e_2,e_5\} \oplus \{e_2,e_4,e_9\}$ | $=\{e_1,e_4,e_5,e_9\}$ | $\rightarrow$ | $\{v_1,v_2,v_3\} \cup \{v_1,v_3,v_5\}$ | $=\{v_1,v_2,v_3,v_5\}$ |
| $c_4 \oplus c_5 =$ | $\{e_2,e_3,e_8\} \oplus \{e_2,e_4,e_9\}$ | $=\{e_1,e_3,e_4,e_5\}$ | $\rightarrow$ | $\{v_1,v_3,v_4\} \cup \{v_1,v_3,v_5\}$ | $=\{v_1,v_3,v_4,v_5\}$ |

Как видно, после исключения изометрических циклов $c_1,c_4,c_5$, имеющих ребро $e_2$, новые простые циклы включают следующие изометрические циклы: $c_2,c_7,c_9,c_6,c_3,c_8,c_{10}$. Например, новый цикл $(c_1 \oplus c_4) \rightarrow \{v_1,v_2,v_3,v_4\}$ включает $\{v_1,v_2,v_4\}$ и $\{v_2,v_3,v_4\}$. Новый цикл $(c_1 \oplus c_5) \rightarrow \{v_1,v_2,v_3,v_5\}$ включает $\{v_1,v_2,v_5\}$ и $\{v_2,v_3,v_5\}$. Новый цикл $(c_4 \oplus c_5) \rightarrow \{v_1,v_3,v_4,v_5\}$ включает $\{v_1,v_4,v_5\}$ и $\{v_3,v_4,v_5\}$. Поэтому, во множестве новых циклов отсутствуют новые изометрические циклы для включения во множество изометрических циклов.

| **Цикл** | **Множество $M_e$** | | **Множество $M_v$** |
|---|---|---|---|
| $c_2$ | $\{e_1,e_3,e_6\}$ | $\rightarrow$ | $\{v_1,v_2,v_4\}$ |
| $c_3$ | $\{e_1,e_4,e_7\}$ | $\rightarrow$ | $\{v_1,v_2,v_5\}$ |
| $c_6$ | $\{e_3,e_4,e_{10}\}$ | $\rightarrow$ | $\{v_1,v_4,v_5\}$ |
| $c_7$ | $\{e_5,e_6,e_8\}$ | $\rightarrow$ | $\{v_2,v_3,v_4\}$ |
| $c_8$ | $\{e_5,e_7,e_9\}$ | $\rightarrow$ | $\{v_2,v_3,v_5\}$ |
| $c_9$ | $\{e_6,e_7,e_{10}\}$; | $\rightarrow$ | $\{v_2,v_4,v_5\}$ |
| $c_{10}$ | $\{e_8,e_9,e_{10}\}$ | $\rightarrow$ | $\{v_3,v_4,v_5\}$ |

Исключаем ребро $e_6$. Формируем новое множество простых циклов:

| **Цикл** | | **Новое $M_e$** | | | **Новое $M_v$** |
|---|---|---|---|---|---|
| $c_2 \oplus c_7=$ | $\{e_1,e_3,e_6\}\oplus\{e_5,e_6,e_8\}$ | $=\{e_1,e_3,e_5,e_8\}$ | $\rightarrow$ | $\{v_1,v_2,v_4\}\cup\{v_2,v_3,v_4\}$ | $=\{v_1,v_2,v_3,v_4\}$ |
| $c_2 \oplus c_9=$ | $\{e_1,e_3,e_6\}\oplus\{e_6,e_7,e_{10}\}$ | $=\{e_1,e_3,e_7,e_{10}\}$ | $\rightarrow$ | $\{v_1,v_2,v_4\}\cup\{v_2,v_4,v_5\}$ | $=\{v_1,v_2,v_4,v_5\}$ |
| $c_7 \oplus c_9=$ | $\{e_5,e_6,e_8\}\oplus\{e_6,e_7,e_{10}\}$ | $=\{e_5,e_7,e_8,e_{10}\}$ | $\rightarrow$ | $\{v_2,v_3,v_4\}\cup\{v_2,v_4,v_5\}$ | $=\{v_2,v_3,v_4,v_5\}$ |

После исключения изометрических циклов $c_2,c_7,c_7$, имеющих ребро $e_6$, некоторые новые простые циклы включают в себя оставшиеся изометрические циклы, а некоторые - нет.

Например, новый цикл $(c_2 \oplus c_7) \rightarrow \{v_1,v_2,v_3,v_4\}$ не включает в себя оставшиеся изометрические циклы и поэтому может быть включен во множество изометрических циклов. Новый цикл $(c_2 \oplus c_9) \rightarrow \{v_1,v_2,v_4,v_5\}$ включает $\{v_1,v_2,v_5\}$ и $\{v_1,v_4,v_5\}$. Новый цикл $(c_7 \oplus c_9) \rightarrow \{v_2,v_3,v_4,v_5\}$ включает $\{v_2,v_3,v_5\}$ и $\{v_3,v_4,v_5\}$. Поэтому циклы $(c_2 \oplus c_9)$ и $(c_7 \oplus c_9)$ не могут быть включенными во множество изометрических циклов.

| **Цикл** | **Множество $M_e$** | | **Множество $M_v$** |
|---|---|---|---|
| $c_3$ | $\{e_1,e_4,e_7\}$ | $\rightarrow$ | $\{v_1,v_2,v_5\}$ |
| $c_6$ | $\{e_3,e_4,e_{10}\}$ | $\rightarrow$ | $\{v_1,v_4,v_5\}$ |
| $c_8$ | $\{e_5,e_7,e_9\}$ | $\rightarrow$ | $\{v_2,v_3,v_5\}$ |
| $c_{10}$ | $\{e_8,e_9,e_{10}\}$ | $\rightarrow$ | $\{v_3,v_4,v_5\}$ |

Множество изометрических циклов изменится.

| **Цикл** | **Множество $M_e$** | | **Множество $M_v$** |
|---|---|---|---|
| $c_3$ | $\{e_1,e_4,e_7\}$ | $\rightarrow$ | $\{v_1,v_2,v_5\}$ |
| $c_6$ | $\{e_3,e_4,e_{10}\}$ | $\rightarrow$ | $\{v_1,v_4,v_5\}$ |
| $c_8$ | $\{e_5,e_7,e_9\}$ | $\rightarrow$ | $\{v_2,v_3,v_5\}$ |
| $c_{10}$ | $\{e_8,e_9,e_{10}\}$ | $\rightarrow$ | $\{v_3,v_4,v_5\}$ |
| $c_2 \oplus c_7$ | $\{e_1,e_3,e_5,e_8\}$ | $\rightarrow$ | $\{v_1,v_2,v_3,v_4\}$ |

Исключаем ребро $e_7$. Формируем новое множество простых циклов:

| Цикл | | Новое $M_e$ | | | Новое $M_v$ |
|---|---|---|---|---|---|
| $c_3 \oplus c_8$ | $\{e_1,e_4,e_7\} \oplus \{e_5,e_7,e_9\}$ | $=\{e_1,e_4,e_5,e_9\}$ | $\rightarrow$ | $\{v_1,v_2,v_5\} \cup \{v_2,v_3,v_5\}$ | $=\{v_1,v_2,v_3,v_5\}$ |

После исключения изометрических циклов $c_3$,$c_8$, имеющих ребро $e_7$, новый простой цикл не включает в себя оставшиеся изометрические циклы. Поэтому он может быть включен во множество изометрических циклов.

| Цикл | Множество $M_e$ | | Множество $M_v$ |
|---|---|---|---|
| $c_6$ | $\{e_3,e_4,e_{10}\}$ | $\rightarrow$ | $\{v_1,v_4,v_5\}$ |
| $c_{10}$ | $\{e_8,e_9,e_{10}\}$ | $\rightarrow$ | $\{v_3,v_4,v_5\}$ |
| $c_2 \oplus c_7$ | $\{e_1,e_3,e_5,e_8\}$ | $\rightarrow$ | $\{v_1,v_2,v_3,v_4\}$ |

Множество изометрических циклов изменится.

| Цикл | Множество $M_e$ | | Множество $M_v$ |
|---|---|---|---|
| $c_6$ | $\{e_3,e_4,e_{10}\}$ | $\rightarrow$ | $\{v_1,v_4,v_5\}$ |
| $c_{10}$ | $\{e_8,e_9,e_{10}\}$ | $\rightarrow$ | $\{v_3,v_4,v_5\}$ |
| $c_2 \oplus c_7$ | $\{e_1,e_3,e_5,e_8\}$ | $\rightarrow$ | $\{v_1,v_2,v_3,v_4\}$ |
| $c_3 \oplus c_8$ | $\{e_1,e_4,e_5,e_9\}$ | $\rightarrow$ | $\{v_1,v_2,v_3,v_5\}$ |

Таким образом, для графа с удаленными ребрами $e_2$,$e_6$,$e_7$ из полного графа, оставшиеся изометрические циклы образуют множество изометрических циклов усеченного графа.

Можно сказать, что формирование множества изометрических циклов можно осуществить двумя способами:

- первый способ основан на выделении и сравнении циклов, проходящих по каждому ребру алгоритмом поиска в ширину с учетом минимальных s,t –маршрутов графа;
- второй способ основан на методе удаления ребер из полного графа с соответствующим удалением изометрических циклов полного графа и, в случае необходимости, включения в оставшееся множество изометрических циклов, которые образуются как кольцевая сумма из удаленных изометрических циклов.

Произведя сравнительный анализ, можно сказать следующее:

- количество изометрических циклов в графе является постоянной величиной, равной или большей цикломатического числа графа, и не зависит от способа их выделения, в то время как множество фундаментальных циклов в точности равно цикломатическому числу графа и зависит от выбора дерева;
- длина центральных разрезов графа определяет локальные степени вершин.

Рассмотрим основные свойства множества изометрических циклов графа. Введем фундаментальное понятие 0-подмножества изометрических циклов и рассмотрим их основные свойства.

По поводу изометрических циклов следует сказать следующее. Так как каждый суграф графа G представляет собой вектор из пространства суграфов $L_G$ размерностью $m$, то система точек $x_0,x_1,x_2,\dots,x_k$ $m$ – мерного линейного пространства $L_G$ называется независимой, если система векторов:

$$(x_1 - x_0), (x_2 - x_0),\dots, (x_k - x_0) \quad (1.11)$$

линейно независима. Очевидно, что независимость возможна и при $k < m$. Система (1.8) линейно независима тогда и только тогда, когда из соотношений

$$\lambda_0 x_0 + \lambda_1 x_1 + \lambda_2 x_2 + \dots + \lambda_k x_k = 0 \quad (1.12)$$

$$\lambda_0 + \lambda_1 + \lambda_2 + \dots + \lambda_k = 0 \quad (1.13)$$

вытекает:

$$\lambda_0 = \lambda_1 = \lambda_2 = \dots = \lambda_k = 0 \quad (1.14)$$

Здесь $\lambda_0,\lambda_1,\lambda_2,\dots\lambda_k$ – действительные числа. Таким образом, свойство системы $x_0,x_1,x_2,\dots,x_k$ быть независимой не зависит от порядка нумерации точек. Сверх того, ясно, что если система точек независима, то всякая её подсистема также независима.

Покажем, что если система векторов (1.11) линейно независима, то из соотношений (1.12) и (1.10) вытекает (1.14). В силу (1.13) соотношение (1.12) переписывается в форме:

$(\lambda_1 + \lambda_2 + \dots + \lambda_k)\, x_0 + \lambda_1 x_1 + \lambda_2 x_2 + \dots + \lambda_k x_k = 0,$

или иначе: $\lambda_1(x_1 - x_0) + \lambda_2(x_2 - x_0) + \dots + \lambda_k(x_k - x_0) = 0.$

Но, так как, система (1.11) линейно независима, то из последнего вытекает $\lambda_0 = \lambda_1 = \lambda_2 = \dots = \lambda_k = 0$, а отсюда ввиду (1.13) следует и $\lambda_0 = 0$. Покажем теперь, что если из соотношения (1.12) и (1.13) вытекает (1.14), то система (1.11) линейно независима.

Пусть:

$$\lambda_1(x_1 - x_0) + \lambda_2(x_2 - x_0) + \dots + \lambda_k(x_k - x_0) = 0. \quad (1.15)$$

Полагая $\lambda_0 = -(\lambda_1 + \lambda_2 + \dots + \lambda_k)$, мы можем переписать соотношение (1.16) в виде:

$$\lambda_0 x_0 + \lambda_1 x_1 + \lambda_2 x_2 + \dots + \lambda_k x_k = 0,$$

причем для наших чисел $\lambda_0,\lambda_1,\lambda_2,\dots\lambda_k$ выполнено условие (1.13). Таким образом, в силу предположения имеем $\lambda_0 = \lambda_1 = \lambda_2 = \dots = \lambda_k = 0$, т.е. из (3.16) вытекает $\lambda_1 = \lambda_2 = \dots = \lambda_k = 0$, а это и означает линейную независимость системы (1.12).

**Определение 1.5.** Кольцевую сумму всех изометрических циклов будем называть ***ободом графа.***

Теперь можно ввести понятие 0-подмножества графа как зависимой системы изометрических циклов. Рассмотрим основные свойства.

**Определение 1.6.** ***0-подмножество изометрических циклов*** – это множество изометрических циклов, кольцевая сумма которых есть пустое множество.

Например, пусть заданы изометрические циклы {a,b,c}, {c,d,e}, {b,g,e}, {a,d,g}. Здесь a,b,c,d,e,g – ребра графа (см. рис. 1.20). Их кольцевая сумма есть пустое множество.

Кольцевая сумма изометрических циклов полного двудольного графа и обода есть пустое множество.

Кольцевая сумма полного графа с четным количеством вершин есть пустое множество. Кольцевая сумма полного графа с нечетным количеством и его обода есть пустое иножество.

**Определение 1.7.** ***Дубль – цикл*** – это простой цикл, который допускает, по крайней мере, два различных нетривиальных разложений в сумму изометрических циклов.

В нашем случае может быть образован следующий дубль-цикл:

$\{a,b,c\} \oplus \{c,d,e\} = \{b,g,e\} \oplus \{a,d,g\} = \{a,b,d,e\}$,

или можно организовать следующие дубль-циклы:

$\{a,b,c\} \oplus \{b,g,e\} = \{c,d,e\} \oplus \{a,d,g\} = \{a,c,e,g\}$,
$\{a,b,c\} \oplus \{c,d,e\} \oplus \{b,g,e\} = \{a,d,g\}$.

По сути дела, множество дубль-циклов – это подмножество простых циклов. Поэтому дубль-циклы, обладают всеми свойствами простых циклов. Изометрические циклы, простые циклы и дубль-циклы, являются частным случаем квазициклов.

***Пример 1.3***. Рассмотрим следующий граф G (см. рис. 1.22). Выделим изометрические циклы в этом графе:

$C^{\tau} = \{\{e_1,e_4,e_5\},\{e_2,e_3,e_5\},\{e_1,e_2,e_6\},\{e_3,e_4,e_6\}\}$.

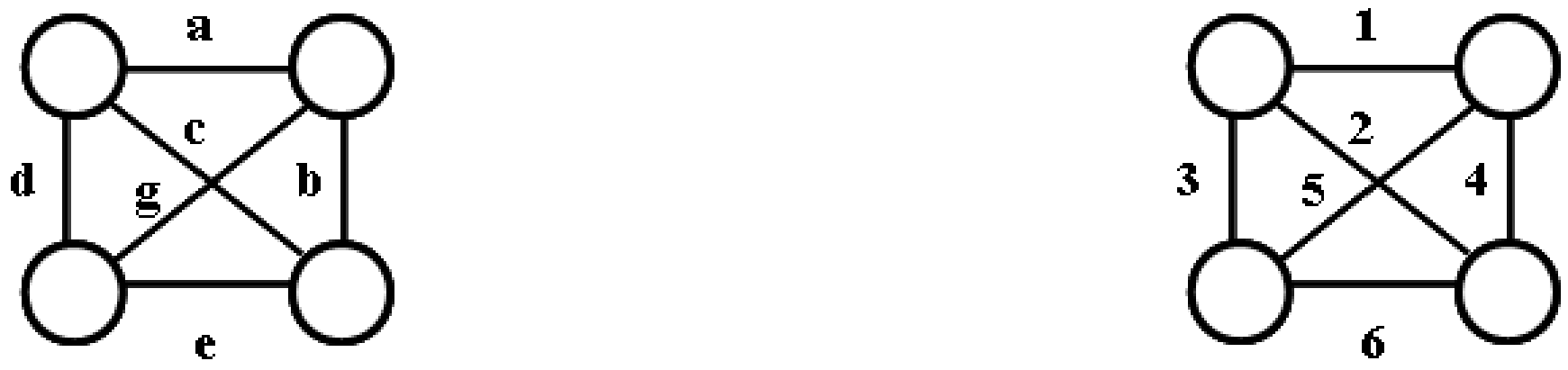

Рис. 1.22. Граф $K_4$.

Кольцевая сумма всех циклов есть пустое множество $\{e_1,e_2,e_4\} \oplus \{e_2,e_3,e_6\} \oplus \{e_1,e_3,e_5\} \oplus \{e_4,e_5,e_6\} = \varnothing$.

Теперь выделим дубль-циклы длины четыре:

$C_1^d = \{e_1,e_2,e_4\} \oplus \{e_2,e_3,e_6\} = \{e_1,e_3,e_5\} \oplus \{e_4,e_5,e_6\} = \{e_1,e_3,e_4,e_6\}$;
$C_2^d = \{e_1,e_2,e_4\} \oplus \{e_1,e_3,e_5\} = \{e_2,e_3,e_6\} \oplus \{e_4,e_5,e_6\} = \{e_2,e_3,e_4,e_5\}$;
$C_3^d = \{e_1,e_2,e_4\} \oplus \{e_4,e_5,e_6\} = \{e_2,e_3,e_6\} \oplus \{e_1,e_3,e_5\} = \{e_1,e_2,e_5,e_6\}$.

## 1.7 Инварианты, построенные на множестве изометрических циклов и центральных разрезов графа

Как мы успели убедиться, запись замкнутого маршрута может быть осуществлена через подмножество ребер или через подмножество вершин графа. Запись цикла через подмножество ребер будем называть реберной записью цикла. Соответственно, вершинная запись цикла состоит из подмножества вершин, принадлежащих рассматриваемому циклу. В основном для записи суграфов и операций над ними применяется реберная запись. Вершинная запись применяется несколько реже и характеризует несколько иные свойства циклов.

Имея множество изометрических циклов графа можно построить вектор количества изометрических циклов проходящих по ребру (впредь будем его называть – вектор циклов по ребрам) . Например, для графа G представленного на рис. 1.11 множество изометрических циклов в реберной записи имеет вид:

$c_1 = \{e_1,e_3,e_5\}$; $c_2 = \{e_2,e_3,e_4\}$; $c_3 = \{e_4,e_8,e_9\}$; $c_4 = \{e_5,e_6,e_{12}\}$;
$c_5 = \{e_5,e_8,e_{13}\}$; $c_6 = \{e_6,e_7,e_{11}\}$; $c_7 = \{e_7,e_8,e_{10}\}$; $c_8 = \{e_1,e_2,e_9,e_{13}\}$;
$c_9 = \{e_{10},e_{11},e_{12},e_{13}\}$.

И тогда вектор циклов по ребрам можно записать в виде кортежа:

$P_e$ =

| $e_1$ | $e_2$ | $e_3$ | $e_4$ | $e_5$ | $e_6$ | $e_7$ | $e_8$ | $e_9$ | $e_{10}$ | $e_{11}$ | $e_{12}$ | $e_{13}$ |
|---|---|---|---|---|---|---|---|---|---|---|---|---|
| 2 | 2 | 2 | 2 | 3 | 2 | 2 | 3 | 2 | 2 | 2 | 2 | 3 |

Или в виде $P_e$ = <2,2,2,2,3,2,2,3,2,2,2,2,3>.

Если записать изометрические циклы через вершины:

$c_1 = \{v_1,v_6,v_7\}$; $c_2 = \{v_5,v_6,v_7\}$; $c_3 = \{v_4,v_5,v_7\}$; $c_4 = \{v_1,v_2,v_7\}$;
$c_5 = \{v_1,v_4,v_7\}$; $c_6 = \{v_2,v_3,v_7\}$; $c_7 = \{v_3,v_4,v_7\}$; $c_8 = \{v_1,v_4,v_5,v_6\}$;
$c_9 = \{v_1,v_2,v_3,v_4\}$,

то вектор количества изометрических циклов проходящих по вершинам графа (впредь будем его называть – вектор циклов по вершинам) запишется в виде кортежа:

$P_v$

| $v_1$ | $v_2$ | $v_3$ | $v_4$ | $v_5$ | $v_6$ | $v_7$ |
|---|---|---|---|---|---|---|
| 5 | 3 | 3 | 5 | 3 | 3 | 7 |

Или в виде $P_v$ = <5,3,3,5,3,3,7>.

Но существует запись циклов в виде замкнутого ориентированного маршрута, так как любое неориентированное ребро может быть представлено двумя разнонаправленными ориентированными ребрами. Такая запись циклов (замкнутых маршрутов) характерна только для плоских графов с учетом заданного направления обхода. Например, для графа, представленного на рис. 1.10 базисная система изометрических циклов и обод, характеризующие плоский граф, могут быть записана в векторном виде:

$c_1 = (v_1,v_7) + (v_7,v_6) + (v_6,v_1)$;
$c_2 = (v_6,v_7) + (v_7,v_5) + (v_5,v_6)$;
$c_3 = (v_5,v_7) + (v_7,v_4) + (v_4,v_5)$;
$c_4 = (v_1,v_2) + (v_2,v_7) + (v_7,v_1)$;
$c_6 = (v_2,v_3) + (v_3,v_7) + (v_7,v_2)$;
$c_7 = (v_7,v_3) + (v_3,v_4) + (v_4,v_7)$;
$c_9 = (v_1,v_4) + (v_4,v_3) + (v_3,v_2) + (v_2,v_1)$;
$c_0 = c_8 = (v_1,v_6) + (v_6,v_5) + (v_5,v_5) + (v_4,v_1)$.

Инвариант графа – это число (функция) графа G, которое принимает одно и то же значение на любом графе, изоморфном G. Пусть f – функция, соотносящая каждому графу G некоторый элемент f(G) из множества M произвольной природы (элементами множества M чаще всего служат числа и системы чисел, векторы, многочлены, матрицы). Эту функцию будем называть инвариантом, если на изоморфных графах её значения совпадают, т.е. для любых G и G':

$$G \cong G' \Rightarrow f(G) = f(G'). \tag{1.16}$$

Подпространство разрезов и подпространство циклов являются нормированными пространствами, так как любому элементу подпространства можно поставить в соответствие неотрицательное вещественное число $\|l\|$, называемое *нормой*. В данном случае, для подпространства разрезов – это длина разреза, а для подпространства циклов – это длина цикла. Причем, введенное понятие удовлетворяет следующим условиям:

- $\|l\| > 0$ при $l \neq 0$, $\|0\| = 0$,
- $\|l_1 + l_2\| \leq \|l_1\| + \|l_2\|$ для любых $l_1 \in$ R, $l_2 \in$ R,
- $\|\alpha l\| = |\alpha| \|l\|$ для любого $l \in$ R и вещественного числа $\alpha$.

Множеству центральных разрезов можно поставить в соответствие, так называемый вектор локальных степеней вида, который также будет инвариантом графа:

$$P_s = (p_1 \times l_1, p_2 \times l_2, ...) \tag{1.17}$$

где $p_1$ – количество центральных разрезов во множестве **S** длиной $l_1$;

$p_2$ – количество центральных разрезов длиной $l_2$ во множестве **S** и т.д.

Причем $l_1 < l_2 < l_3 < ….$, то есть длина циклов расставлена в порядке неубывания.

Множеству изометрических циклов можно также поставить в соответствие вектор вида (1.16), который также будет инвариантом графа:

$$P_c = (p_1 \times l_1, p_2 \times l_2, ...) \tag{1.18}$$

где $p_1$ – количество изометрических циклов во множестве $C_\tau$ длиной $l_1$;

$p_2$ – количество изометрических циклов длиной $l_2$ во множестве $C_\tau$ и т.д.

Причем $l_1 < l_2 < l_3 < ….$, то есть длина циклов, расставленая в порядке неубывания.

***Пример 1.4.*** Выделить множество изометрических циклов в графе $\mathbf{G_1}$.

Количество вершин графа = 10.
Количество ребер графа = 23.
Количество изометрических циклов графа = 31.

Матрица смежностей графа представленная в виде списка смежных вершин:

вершина $v_1$: $\{v_2,v_6,v_7,v_{10}\}$;
вершина $v_2$: $\{v_1,v_3,v_5,v_7\}$;
вершина $v_3$: $\{v_2,v_4,v_5,v_9\}$;
вершина $v_4$: $\{v_3,v_5,v_6,v_7,v_9\}$;
вершина $v_5$: $\{v_2,v_3,v_4,v_6,v_8,v_{10}\}$;
вершина $v_6$: $\{v_1,v_4,v_5,v_7,v_9\}$;

вершина $v_7$: $\{v_1,v_2,v_4,v_6,v_8\}$;
вершина $v_8$: $\{v_5,v_7,v_9,v_{10}\}$;
вершина $v_9$: $\{v_3,v_4,v_6,v_8,v_{10}\}$;
вершина $v_{10}$: $\{v_1,v_5,v_8,v_9\}$.

Элементы матрицы инциденций:

ребро $e_1$: $(v_1,v_2)$ или $(v_2,v_1)$;
ребро $e_2$: $(v_1,v_6)$ или $(v_6,v_1)$;
ребро $e_3$: $(v_1,v_7)$ или $(v_7,v_1)$;
ребро $e_4$: $(v_1,v_{10})$ или $(v_{10},v_1)$;
ребро $e_5$: $(v_2,v_3)$ или $(v_3,v_2)$;
ребро $e_6$: $(v_2,v_5)$ или $(v_5,v_2)$;
ребро $e_7$: $(v_2,v_7)$ или $(v_7,v_2)$;
ребро $e_8$: $(v_3,v_4)$ или $(v_4,v_3)$;
ребро $e_9$: $(v_3,v_5)$ или $(v_5,v_3)$;
ребро $e_{10}$: $(v_3,v_9)$ или $(v_9,v_3)$;
ребро $e_{11}$: $(v_4,v_5)$ или $(v_5,v_4)$;
ребро $e_{12}$: $(v_4,v_6)$ или $(v_6,v_4)$;
ребро $e_{13}$: $(v_4,v_7)$ или $(v_7,v_4)$;
ребро $e_{14}$: $(v_4,v_9)$ или $(v_9,v_4)$;
ребро $e_{15}$: $(v_5,v_6)$ или $(v_6,v_5)$;
ребро $e_{16}$: $(v_5,v_8)$ или $(v_8,v_5)$;
ребро $e_{17}$: $(v_5,v_{10})$ или $(v_{10},v_5)$;
ребро $e_{18}$: $(v_6,v_7)$ или $(v_7,v_6)$;
ребро $e_{19}$: $(v_6,v_9)$ или $(v_9,v_6)$;
ребро $e_{20}$: $(v_7,v_8)$ или $(v_8,v_7)$;
ребро $e_{21}$: $(v_8,v_9)$ или $(v_9,v_8)$;
ребро $e_{22}$: $(v_8,v_{10})$ или $(v_{10},v_8)$;
ребро $e_{23}$: $(v_9,v_{10})$ или $(v_{10},v_9)$.

Множество изометрических циклов графа

| **Циклы** | **Множество изометрических циклов графа в виде ребер:** | **Множество изометрических циклов графа в виде вершин:** |
|---|---|---|
| цикл $c_1$ | $\{e_1,e_2,e_6,e_{15}\}$ | $\{v_1,v_2,v_5,v_6\}$; |
| цикл $c_2$ | $\{e_1,e_3,e_7\}$ | $\{v_1,v_2,v_7\}$; |
| цикл $c_3$ | $\{e_1,e_4,e_6,e_{17}\}$ | $\{v_1,v_2,v_5,v_{10}\}$; |
| цикл $c_4$ | $\{e_2,e_3,e_{18}\}$ | $\{v_1,v_6,v_7\}$; |
| цикл $c_5$ | $\{e_2,e_4,e_{15},e_{17}\}$ | $\{v_1,v_5,v_6,v_{10}\}$; |
| цикл $c_6$ | $\{e_2,e_4,e_{19},e_{23}\}$ | $\{v_1,v_6,v_9,v_{10}\}$; |
| цикл $c_7$ | $\{e_3,e_4,e_{20},e_{22}\}$ | $\{v_1,v_7,v_8,v_{10}\}$; |
| цикл $c_8$ | $\{e_1,e_2,e_5,e_{10},e_{19}\}$ | $\{v_1,v_2,v_3,v_6,v_9\}$; |
| цикл $c_9$ | $\{e_1,e_4,e_5,e_{10},e_{23}\}$ | $\{v_1,v_2,v_3,v_9,v_{10}\}$; |
| цикл $c_{10}$ | $\{e_5,e_6,e_9\}$ | $\{v_2,v_3,v_5\}$; |
| цикл $c_{11}$ | $\{e_5,e_7,e_8,e_{13}\}$ | $\{v_2,v_3,v_4,v_7\}$; |
| цикл $c_{12}$ | $\{e_6,e_7,e_{11},e_{13}\}$ | $\{v_2,v_4,v_5,v_7\}$; |
| цикл $c_{13}$ | $\{e_6,e_7,e_{15},e_{18}\}$ | $\{v_2,v_5,v_6,v_7\}$; |
| цикл $c_{14}$ | $\{e_6,e_7,e_{16},e_{20}\}$ | $\{v_2,v_5,v_7,v_8\}$; |
| цикл $c_{15}$ | $\{e_8,e_9,e_{11}\}$ | $\{v_3,v_4,v_5\}$; |
| цикл $c_{16}$ | $\{e_8,e_{10},e_{14}\}$ | $\{v_3,v_4,v_9\}$; |
| цикл $c_{17}$ | $\{e_9,e_{10},e_{15},e_{19}\}$ | $\{v_3,v_5,v_6,v_9\}$; |
| цикл $c_{18}$ | $\{e_9,e_{10},e_{16},e_{21}\}$ | $\{v_3,v_5,v_8,v_9\}$; |
| цикл $c_{19}$ | $\{e_9,e_{10},e_{17},e_{23}\}$ | $\{v_3,v_5,v_9,v_{10}\}$; |
| цикл $c_{20}$ | $\{e_{11},e_{12},e_{15}\}$ | $\{v_4,v_5,v_6\}$; |
| цикл $c_{21}$ | $\{e_{11},e_{13},e_{16},e_{20}\}$ | $\{v_4,v_5,v_7,v_8\}$; |
| цикл $c_{22}$ | $\{e_{11},e_{14},e_{16},e_{21}\}$ | $\{v_4,v_5,v_8,v_9\}$; |
| цикл $c_{23}$ | $\{e_{11},e_{14},e_{17},e_{23}\}$ | $\{v_4,v_5,v_9,v_{10}\}$; |
| цикл $c_{24}$ | $\{e_{12},e_{13},e_{18}\}$ | $\{v_4,v_6,v_7\}$; |
| цикл $c_{25}$ | $\{e_{12},e_{14},e_{19}\}$ | $\{v_4,v_6,v_9\}$; |
| цикл $c_{26}$ | $\{e_{13},e_{14},e_{20},e_{21}\}$ | $\{v_4,v_7,v_8,v_9\}$; |
| цикл $c_{27}$ | $\{e_{15},e_{16},e_{18},e_{20}\}$ | $\{v_5,v_6,v_7,v_8\}$; |
| цикл $c_{28}$ | $\{e_{15},e_{16},e_{19},e_{21}\}$ | $\{v_5,v_6,v_8,v_9\}$; |
| цикл $c_{29}$ | $\{e_{15},e_{17},e_{19},e_{23}\}$ | $\{v_5,v_6,v_9,v_{10}\}$; |
| цикл $c_{30}$ | $\{e_{16},e_{17},e_{22}\}$ | $\{v_5,v_8,v_{10}\}$; |
| цикл $c_{31}$ | $\{e_{18},e_{19},e_{20},e_{21}\}$ | $\{v_6,v_7,v_8,v_9\}$; |
| цикл $c_{32}$ | $\{e_{21},e_{22},e_{23}\}$ | $\{v_8,v_9,v_{10}\}$. |

***Пример 1.5.*** Выделить множество изометрических циклов в графе $G_2$.

Количество вершин графа = 12
Количество ребер графа = 35
Количество изометрических циклов графа = 56

Матрица смежностей графа, представленная в виде списка смежных вершин:

вершина $v_1$: $\{v_2,v_3,v_4,v_5,v_6,v_7,v_8,v_9,v_{10},v_{12}\}$;
вершина $v_2$: $\{v_1,v_3,v_4,v_7,v_8,v_{12}\}$;
вершина $v_3$: $\{v_1,v_2,v_4,v_6,v_7,v_8,v_{10}\}$;
вершина $v_4$: $\{v_1,v_2,v_3,v_6,v_{11}\}$;
вершина $v_5$: $\{v_1,v_8,v_9\}$;
вершина $v_6$: $\{v_1,v_3,v_4,v_7,v_9,v_{12}\}$;
вершина $v_7$: $\{v_1,v_2,v_3,v_6,v_8,v_9,v_{11},v_{12}\}$;
вершина $v_8$: $\{v_1,v_2,v_3,v_5,v_7,v_{10},v_{11}\}$;
вершина $v_9$: $\{v_1,v_5,v_6,v_7,v_{10}\}$;
вершина $v_{10}$: $\{v_1,v_3,v_8,v_9,v_{12}\}$;
вершина $v_{11}$: $\{v_4,v_7,v_8\}$;
вершина $v_{12}$: $\{v_1,v_2,v_6,v_7,v_{10}\}$.

Элементы матрицы инциденций:

ребро $e_1$: $(v_1,v_2)$ или $(v_2,v_1)$;
ребро $e_2$: $(v_1,v_3)$ или $(v_3,v_1)$;
ребро $e_3$: $(v_1,v_4)$ или $(v_4,v_1)$;
ребро $e_4$: $(v_1,v_5)$ или $(v_5,v_1)$;
ребро $e_5$: $(v_1,v_6)$ или $(v_6,v_1)$;
ребро $e_6$: $(v_1,v_7)$ или $(v_7,v_1)$;
ребро $e_7$: $(v_1,v_8)$ или $(v_8,v_1)$;
ребро $e_8$: $(v_1,v_9)$ или $(v_9,v_1)$;
ребро $e_9$: $(v_1,v_{10})$ или $(v_{10},v_1)$;
ребро $e_{10}$: $(v_1,v_{12})$ или $(v_{12},v_1)$;
ребро $e_{11}$: $(v_2,v_3)$ или $(v_3,v_2)$;
ребро $e_{12}$: $(v_2,v_4)$ или $(v_4,v_2)$;
ребро $e_{13}$: $(v_2,v_7)$ или $(v_7,v_2)$;
ребро $e_{14}$: $(v_2,v_8)$ или $(v_8,v_2)$;
ребро $e_{15}$: $(v_2,v_{12})$ или $(v_{12},v_2)$;
ребро $e_{16}$: $(v_3,v_4)$ или $(v_4,v_3)$;
ребро $e_{17}$: $(v_3,v_6)$ или $(v_6,v_3)$;
ребро $e_{18}$: $(v_3,v_7)$ или $(v_7,v_3)$;
ребро $e_{19}$: $(v_3,v_8)$ или $(v_8,v_3)$;
ребро $e_{20}$: $(v_3,v_{10})$ или $(v_{10},v_3)$;
ребро $e_{21}$: $(v_4,v_6)$ или $(v_6,v_4)$;
ребро $e_{22}$: $(v_4,v_{11})$ или $(v_{11},v_4)$;
ребро $e_{23}$: $(v_5,v_8)$ или $(v_8,v_5)$;
ребро $e_{24}$: $(v_5,v_9)$ или $(v_9,v_5)$;
ребро $e_{25}$: $(v_6,v_7)$ или $(v_7,v_6)$;
ребро $e_{26}$: $(v_6,v_9)$ или $(v_9,v_6)$;
ребро $e_{27}$: $(v_6,v_{12})$ или $(v_{12},v_6)$;
ребро $e_{28}$: $(v_7,v_8)$ или $(v_8,v_7)$;
ребро $e_{29}$: $(v_7,v_9)$ или $(v_9,v_7)$;
ребро $e_{30}$: $(v_7,v_{11})$ или $(v_{11},v_7)$;
ребро $e_{31}$: $(v_7,v_{12})$ или $(v_{12},v_7)$;
ребро $e_{32}$: $(v_8,v_{10})$ или $(v_{10},v_8)$;
ребро $e_{33}$: $(v_8,v_{11})$ или $(v_{11},v_8)$;
ребро $e_{34}$: $(v_9,v_{10})$ или $(v_{10},v_9)$;
ребро $e_{35}$: $(v_{10},v_{12})$ или $(v_{12},v_{10})$.

Множество изометрических циклов графа:

| Циклы | Множество изометрических циклов графа в виде ребер: | Множество изометрических циклов графа в виде вершин: |
|---|---|---|
| цикл $c_1$ | $\{e_1,e_2,e_{11}\}$; | $\{v_1,v_2,v_3\}$; |
| цикл $c_2$ | $\{e_1,e_3,e_{12}\}$; | $\{v_1,v_2,v_4\}$; |
| цикл $c_3$ | $\{e_1,e_6,e_{13}\}$; | $\{v_1,v_2,v_7\}$; |
| цикл $c_4$ | $\{e_1,e_7,e_{14}\}$; | $\{v_1,v_2,v_8\}$; |
| цикл $c_5$ | $\{e_1,e_{10},e_{15}\}$; | $\{v_1,v_2,v_{12}\}$; |
| цикл $c_6$ | $\{e_2,e_3,e_{16}\}$; | $\{v_1,v_3,v_4\}$; |
| цикл $c_7$ | $\{e_2,e_5,e_{17}\}$; | $\{v_1,v_3,v_6\}$; |
| цикл $c_8$ | $\{e_2,e_6,e_{18}\}$; | $\{v_1,v_3,v_7\}$; |
| цикл $c_9$ | $\{e_2,e_7,e_{19}\}$; | $\{v_1,v_3,v_8\}$; |
| цикл $c_{10}$ | $\{e_2,e_9,e_{20}\}$; | $\{v_1,v_3,v_{10}\}$; |
| цикл $c_{11}$ | $\{e_3,e_5,e_{21}\}$; | $\{v_1,v_4,v_6\}$; |
| цикл $c_{12}$ | $\{e_3,e_6,e_{22},e_{30}\}$; | $\{v_1,v_4,v_7,v_{11}\}$; |
| цикл $c_{13}$ | $\{e_3,e_7,e_{22},e_{33}\}$; | $\{v_1,v_4,v_8,v_{11}\}$; |
| цикл $c_{14}$ | $\{e_4,e_7,e_{23}\}$; | $\{v_1,v_5,v_8\}$; |
| цикл $c_{15}$ | $\{e_4,e_8,e_{24}\}$; | $\{v_1,v_5,v_9\}$; |
| цикл $c_{16}$ | $\{e_5,e_6,e_{25}\}$; | $\{v_1,v_6,v_7\}$; |
| цикл $c_{17}$ | $\{e_5,e_8,e_{26}\}$; | $\{v_1,v_6,v_9\}$; |
| цикл $c_{18}$ | $\{e_5,e_{10},e_{27}\}$; | $\{v_1,v_6,v_{12}\}$; |
| цикл $c_{19}$ | $\{e_6,e_7,e_{28}\}$; | $\{v_1,v_7,v_8\}$; |
| цикл $c_{20}$ | $\{e_6,e_8,e_{29}\}$; | $\{v_1,v_7,v_9\}$; |
| цикл $c_{21}$ | $\{e_6,e_{10},e_{31}\}$; | $\{v_1,v_7,v_{12}\}$; |
| цикл $c_{22}$ | $\{e_7,e_9,e_{32}\}$; | $\{v_1,v_8,v_{10}\}$; |
| цикл $c_{23}$ | $\{e_8,e_9,e_{34}\}$; | $\{v_1,v_9,v_{10}\}$; |
| цикл $c_{24}$ | $\{e_9,e_{10},e_{35}\}$; | $\{v_1,v_{10},v_{12}\}$; |
| цикл $c_{25}$ | $\{e_{11},e_{12},e_{16}\}$; | $\{v_2,v_3,v_4\}$; |
| цикл $c_{26}$ | $\{e_{11},e_{13},e_{18}\}$; | $\{v_2,v_3,v_7\}$; |
| цикл $c_{27}$ | $\{e_{11},e_{14},e_{19}\}$; | $\{v_2,v_3,v_8\}$; |
| цикл $c_{28}$ | $\{e_{11},e_{15},e_{17},e_{27}\}$; | $\{v_2,v_3,v_6,v_{12}\}$; |
| цикл $c_{29}$ | $\{e_{11},e_{15},e_{20},e_{35}\}$; | $\{v_2,v_3,v_{10},v_{12}\}$; |
| цикл $c_{30}$ | $\{e_{12},e_{13},e_{21},e_{25}\}$; | $\{v_2,v_4,v_6,v_7\}$; |
| цикл $c_{31}$ | $\{e_{12},e_{13},e_{22},e_{30}\}$; | $\{v_2,v_4,v_7,v_{11}\}$; |
| цикл $c_{32}$ | $\{e_{12},e_{14},e_{22},e_{33}\}$; | $\{v_2,v_4,v_8,v_{11}\}$; |
| цикл $c_{33}$ | $\{e_{12},e_{15},e_{21},e_{27}\}$; | $\{v_2,v_4,v_6,v_{12}\}$; |
| цикл $c_{34}$ | $\{e_{13},e_{14},e_{28}\}$; | $\{v_2,v_7,v_8\}$; |
| цикл $c_{35}$ | $\{e_{13},e_{15},e_{31}\}$; | $\{v_2,v_7,v_{12}\}$; |
| цикл $c_{36}$ | $\{e_{14},e_{15},e_{32},e_{35}\}$; | $\{v_2,v_8,v_{10},v_{12}\}$; |
| цикл $c_{37}$ | $\{e_{16},e_{17},e_{21}\}$; | $\{v_3,v_4,v_6\}$; |
| цикл $c_{38}$ | $\{e_{16},e_{18},e_{22},e_{30}\}$; | $\{v_3,v_4,v_7,v_{11}\}$; |
| цикл $c_{39}$ | $\{e_{16},e_{19},e_{22},e_{33}\}$; | $\{v_3,v_4,v_8,v_{11}\}$; |
| цикл $c_{40}$ | $\{e_{17},e_{18},e_{25}\}$; | $\{v_3,v_6,v_7\}$; |
| цикл $c_{41}$ | $\{e_{17},e_{20},e_{26},e_{34}\}$; | $\{v_3,v_6,v_9,v_{10}\}$; |
| цикл $c_{42}$ | $\{e_{17},e_{20},e_{27},e_{35}\}$; | $\{v_3,v_6,v_{10},v_{12}\}$; |
| цикл $c_{43}$ | $\{e_{18},e_{19},e_{28}\}$; | $\{v_3,v_7,v_8\}$; |
| цикл $c_{44}$ | $\{e_{18},e_{20},e_{29},e_{34}\}$; | $\{v_3,v_7,v_9,v_{10}\}$; |
| цикл $c_{45}$ | $\{e_{18},e_{20},e_{31},e_{35}\}$; | $\{v_3,v_7,v_{10},v_{12}\}$; |
| цикл $c_{46}$ | $\{e_{19},e_{20},e_{32}\}$; | $\{v_3,v_8,v_{10}\}$; |
| цикл $c_{47}$ | $\{e_{21},e_{22},e_{25},e_{30}\}$; | $\{v_4,v_6,v_7,v_{11}\}$; |
| цикл $c_{48}$ | $\{e_{23},e_{24},e_{28},e_{29}\}$; | $\{v_5,v_7,v_8,v_9\}$; |
| цикл $c_{49}$ | $\{e_{23},e_{24},e_{32},e_{34}\}$; | $\{v_5,v_8,v_9,v_{10}\}$; |
| цикл $c_{50}$ | $\{e_{25},e_{26},e_{29}\}$; | $\{v_6,v_7,v_9\}$; |
| цикл $c_{51}$ | $\{e_{25},e_{27},e_{31}\}$; | $\{v_6,v_7,v_{12}\}$; |
| цикл $c_{52}$ | $\{e_{26},e_{27},e_{34},e_{35}\}$; | $\{v_6,v_9,v_{10},v_{12}\}$; |
| цикл $c_{53}$ | $\{e_{28},e_{29},e_{32},e_{34}\}$; | $\{v_7,v_8,v_9,v_{10}\}$; |
| цикл $c_{54}$ | $\{e_{28},e_{30},e_{33}\}$; | $\{v_7,v_8,v_{11}\}$; |
| цикл $c_{55}$ | $\{e_{28},e_{31},e_{32},e_{35}\}$; | $\{v_7,v_8,v_{10},v_{12}\}$; |
| цикл $c_{56}$ | $\{e_{29},e_{31},e_{34},e_{35}\}$. | $\{v_7,v_9,v_{10},v_{12}\}$. |

## 1.8. Текст программы Raschet1 определения изометрических циклов как множество ребер

```
program Raschet1;

type
    TMasy = array[1..1000] of integer;
    TMass = array[1..4000] of integer;
var
    F1,F2 : text;
    i,ii,j,jj,K,K1,Np,Nv,Kzikl,M,MakLin : integer;
    Ziklo,Nr,KKK,AB,K9 : integer;
    Masy: TMasy;
    Mass: TMass;
    Masi: TMass;
    MasyT: TMasy;
    MassT : TMass;
    MasMdop : TMasy;
    MasMy1 : TMasy;
    MasMs1 : TMass;
    MasMy2 : TMasy;
    MasMs2 : TMass;
    MasMy3 : TMasy;
```

```
        MasMs3 : TMass;
        MasMcg : TMasy;
        MasMcg1 : TMasy;
        MasKol : TMasy;
        Mass1 : TMass;
{*****************************************************************}
procedure FormVolna(var Nv,Nv1,Nv2 : integer;
                    var My : TMasy;
                    var Ms : TMass;
                    var Mdop : TMasy);
{ Nv -  количество вершин в графе;                                   }
{ Nv1 -  номер первой вершины;                                       }
{ Nv2 -  номер второй вершины;                                       }
{ My -  массив указателей для маттрицы смежностей графа;             }
{ Ms -  массив элементов матрицы смежностей;                         }
{ Mdop -  массив глубины распространения волны.                      }
{                                                                    }
{   Данная процедура формирует массив распространения                }
{   волны, здесь номер уровень волны.                                }
 var I,Im,J,Kum,Istart,Istop : integer;
 label 1,2,3,4;
 begin
       for I:=1 to Nv do Mdop[I]:=0;
       Mdop[Nv1]:=1;
       Mdop[Nv2]:=2;
       Im:=2;
      1:  Im:=Im+1;
       for I:=1 to Nv do
        if Mdop[I]=0 then goto 2;
       goto 3;
      2:  for J:=1 to Nv do
       begin
          if Mdop[J]<>Im-1 then goto 4;
          Istart:=My[J];
          Istop:=My[J+1]-1;
          for I:=Istart to Istop do
          begin
             Kum:=Ms[I];
             if Mdop[Kum]=0 then Mdop[Kum]:=Im;
          end;
      4:  end;
       goto 1;
      3:;
      end; {FormVolna}
{*******************************************************}
 procedure FormKpris(var Kzikl : integer;
                     var Myy : TMasy;
                     var Mss : TMass;
                     var My : TMasy;
                     var Ms : TMass;
                     var Ms2 : TMass);
{ Kzikl  - количество t-циклов в графе;                  }
{ My  - массив указателей для матрицы смежностей;        }
{ Ms  - массив элементов матрицы смежностей;             }
{ Ms2 -  массив элементов матрицы инциденций;            }
{ Myy - массив указателей для матрицы t-циклов;          }
{ Mss - массив элементов матрицы t-циклов;               }
{                                                        }
{   Процедура переводит запись циклов в виде вершин      }
{   в запись в виде ребер.                               }
{                                                        }
   var I,J,JJ,Ip,Ip1,JJJ,Npn,KK : integer;
   label 1,2;
```

```
   begin
     for J:= 1 to Kzikl do
     begin
      for JJ:=Myy[J] to Myy[J+1]-2 do
      begin
        Ip:= Mss[JJ];
        Ip1:=Mss[JJ+1];
        for JJJ:=My[Ip] to My[IP+1]-1 do
         begin
           if Ms[JJJ]<>Ip1 then goto 1;
           Mss[JJ]:=Ms2[JJJ];
1:          end;
     end;
     Mss[Myy[J+1]-1]:=0;
     end;
     Npn:=Myy[Kzikl+1]-1;
     KK:=0;
     for I:=1 to Npn do
     begin
       if Mss[I]=0 then goto 2;
       KK:=KK+1;
       Mss[KK]:=Mss[I];
2:      end;
     for I:=1 to Kzikl do Myy[I+1]:=Myy[I+1]-I;
 end; {FormKpris}
{***********************************************************}
 procedure FormDiz(var Mm1,Mm2,Mm4 : integer;
                  var M1 : TMasy;
                  var M2 : TMasy;
                  var M4 : TMasy);
 { Mm1 - количество элементов в первом цикле;                    }
 { Mm2 - количество элементов во втором цикле;                   }
 { Mm4 - количество элементов в их пересечении;                  }
 { M1 - массив элементов первого цикла;                          }
 { M2 - массив элементов второго цикла;                          }
 { M4 - массив элементов пересечения.                            }
 {                                                               }
 {   Процедура определения пересечения двух циклов               }
 {                                                               }
  var J,I : integer;
  label 1,2;
     begin
      Mm4:=0;
      for I:=1 to Mm1 do
      begin
         for J:=1 to Mm2 do
         begin
            if M1[I]<>M2[J] then goto 1;
            Mm4:=Mm4+1;
            M4[Mm4]:=M2[J];
            goto 2;
     1:     end;
     2: end;
     end; {FormDiz}
{***********************************************************}
 procedure FormSwigug(var Kzikl1,Kzikl2 : integer;
                   var My1 : TMasy;
                   var Ms1 : TMass;
                   var My2 : TMasy;
                   var Ms2 : TMass;
                   var Mcg : TMasy;
                   var Mcg1 : TMasy;
                   var Kol : Tmasy);
```

```
{ Kzikl1 - количество первых t-циклов в графе;                    }
{ Kzikl2 - количество вторых t-циклов в графе;                    }
{ My1 - массив указателей для первых t-циклов в графе;            }
{ Ms1 - массив для первых t-циклов в графе;                       }
{ My2 - массив указателей для вторых t-циклов в графе;            }
{ Ms2 - массив для вторых t-циклов в графе;                       }
{ Mcg -  вспомогательный массив;                                  }
{ Mcg1 -  вспомогательный массив;                                 }
{ Kol -  вспомогательный массив.                                  }
{                                                                 }
{    Процедура формирования и сравнения t-циклов                  }
{                                                                 }
    label 1;
    var I,J,JJ,II,NN1,NN2,NN3 : integer;
    begin
     for I:=1 to Kzikl1 do
     begin
       NN1:=My1[I+1]-My1[I];
       for II:=1 to NN1 do Mcg[II]:=Ms1[My1[I]-1+II];
       for J:=1 to Kzikl2 do
       begin
         NN2:=My2[J+1]-My2[J];
         for JJ:=1 to NN2 do Mcg1[JJ]:=Ms2[My2[J]-1+JJ];
         FormDiz(NN1,NN2,NN3,Mcg,Mcg1,Kol);
         if NN3=NN1 then goto 1;
       end;
       for II:=1 to NN1 do Ms1[My1[I]-1+II]:=0;
    1:  end;
    end; {FormSwigug}
{*************************************************************}
procedure FormDozas(var Kzikl,Kzikl1 : integer;
                    var My1 : TMasy;
                    var Ms1 : TMass;
                    var MyT : TMasy;
                    var MsT : TMass;
                    var Mcg : TMasy;
                    var Mcg1 : TMasy;
                    var Kol : Tmasy);
{ Kzikl  - количество t-циклов в графе;                           }
{ Ms1  - массив указателей для матрицы первых t-циклов;           }
{ Ms1  - массив элементов матрицы первых t-циклов;                }
{ MyT - массив указателей для матрицы t-циклов;                   }
{ MsT - массив элементов матрицы t-циклов;                        }
{ Mcg -  вспомогательный массив;                                  }
{ Mcg1 -  вспомогательный массив;                                 }
{ Kol -  вспомогательный массив.                                  }
{                                                                 }
{    Процедура формирования t-циклов                              }
{                                                                 }
label 1,2,3,4;
var I,J,JJ,II,III,Kp,Ks,NN1,NN2,NN3 : integer;
    begin
     if Kzikl=0  then goto 1;
     Kp:=Kzikl;
     Ks:=MyT[Kzikl+1]-1;
     for I:=1 to Kzikl1 do
     begin
       NN1:=My1[I+1]-My1[I];
       for II:=1 to NN1 do
       begin
         if Ms1[My1[I]-1+II]=0 then goto 2;
         Mcg[II]:=Ms1[My1[I]-1+II];
       end;
```

```
      for J:=1 to Kzikl do
      begin
        NN2:=MyT[J+1]-MyT[J];
        for JJ:=1 to NN2 do Mcg1[JJ]:=MsT[MyT[J]-1+JJ];
        FormDiz(NN1,NN2,NN3,Mcg,Mcg1,Kol);
        if NN3=NN1 then goto 2;
      end;
      Kp:=Kp+1;
      MyT[Kp+1]:=MyT[Kp]+NN1;
      for III:=1 to NN1 do
      begin
        Ks:=Ks+1;
        MsT[Ks]:=Mcg[III];
      end;
   2:  end;
    Kzikl:=Kp;
    goto 3;
   1:  Ks:=0;
    MyT[1]:=1;
    for I:=1 to Kzikl1 do
    begin
      NN1:=My1[I+1]-My1[I];
      for II:=1 to NN1 do
      begin
        if Ms1[My1[I]-1+II]=0 then goto 4;
        Ks:=Ks+1;
        MsT[Ks]:=Ms1[My1[I]-1+II];
      end;
      Kzikl:=Kzikl+1;
      MyT[Kzikl+1]:=MyT[Kzikl]+NN1;
   4:  end;
   3:;
   end; {FormDozas}
{*************************************************************}
 procedure FormSoasda(var Nv1,L,I2 : integer;
                      var Key1 : Boolean;
                      var My1 : TMasy;
                      var Ms1 : TMass;
                      var My : TMasy;
                      var Ms : TMass;
                      var Mcg : TMasy;
                      var Mcg1 : TMasy;
                      var Kol : TMasy);
{ Nv1 - номер вершины в графе;                                }
{ Key1  - признак;                                            }
{ I2 - признак;                                               }
{ L - длина цикла;                                            }
{ My  - массив указателей для матрицы смежностей;             }
{ Ms - массив для элементов мматрицы смежностей;              }
{ My1 - массив указателей для матрицы t-циклов;               }
{ Ms1 - массив элементов матрицы t-циклов;                    }
{ Mcg -  вспомогательный массив;                              }
{ Mcg1 -  вспомогательный массив;                             }
{ Kol -  вспомогательный массив.                              }
{                                                             }
{    Процедура формирования цикла заданной длины              }
{                                                             }
   label 2,3,14,11,5,6,12;
   var I,J,Ip1,Isu,Ip2 : integer;
   begin
    Key1:=false;
    if I2>0 then goto 2;
    for I:=1 to L do
```

```
        begin
          Mcg[I]:=1;
          Kol[I]:=My1[I+1]-My1[I];
        end;
        Mcg[L]:=0;
       2:  Mcg[L]:=Mcg[L]+1;
       14:;
        for I:=1 to L do
        begin
          if Mcg[I]<=Kol[I] then goto 3;
          if Mcg[1]>Kol[1] then goto 11;
          Mcg[I-1]:=Mcg[I-1]+1;
          for J:=I to L do Mcg[J]:=1;
          goto 14;
       3:  end;
        Mcg1[1]:=Nv1;
        Mcg1[L+1]:=Nv1;
        for I:=2 to L do
        begin
          Ip1:=My1[I]+Mcg[I]-1;
          Isu:=Ms1[Ip1];
          Mcg1[I]:=Isu;
          Ip2:=Mcg1[I-1];
          for J:=My[Ip2] to My[Ip2+1]-1 do
          begin
            if Ms[J]<>Isu then goto 6;
            goto 5;
       6:     end;
          Mcg[I]:=Mcg[I]+1;
          goto 14;
       5:  end;
        Key1:=true;
        I2:=1;
        goto 12;
       11: Key1:=false;
       12:;
       end; {FormSoasda}
{*************************************************************}
 procedure FormWegin(var Nv,Nv1,Nv2,Kzikl : integer;
                     var My : TMasy;
                     var Ms : TMass;
                     var My1 : TMasy;
                     var Ms1 : TMass;
                     var MyT : TMasy;
                     var MsT : TMass;
                     var Mdop : TMasy;
                     var Mcg : TMasy;
                     var Mcg1 : TMasy;
                     var Kol : Tmasy);
 { Nv - количество вершин в графе;                              }
 { Nv1 -  номер первой вершины;                                 }
 { Nv2 -  номер второй вершины;                                 }
 { Kzikl  - количество t-циклов в графе;                        }
 { My  - массив указателей для матрицы смежностей;              }
 { Ms  - массив элементов матрицы смежностей;                   }
 { My1 - массив для формирования указателей t-циклов;           }
 { Ms1 -  массив формирования t-циклов;                         }
 { MyT - массив указателей для матрицы t-циклов;                }
 { MsT - массив элементов матрицы t-циклов;                     }
 { Mdop -  вспомогательный массив;                              }
 { Mcg -  вспомогательный массив;                               }
 { Mcg1 -  вспомогательный массив;                              }
 { Kol -  вспомогательный массив.                               }
```

```
{ Ms2 -  текущая начальная строка элементов;                         }
{                                                                     }
{    Процедура построения единичных циклов                            }
{                                                                     }
label 1,2,3,22,5,7,8,10,11,13,9;
var I,J,JJ,II,Kot,Istart,Istop,JJJ,KKK,JI,Imum : integer;
var Ji1,Kot1,Kot2,Isu,Irr,I2 : integer;
var Key1 : Boolean;
  begin
     Kot:=0;
     MyT[1]:=1;
     Kzikl:=0;
     for I:=1 to Nv do
     begin
       if Mdop[I]<=2 then goto 1;
       for J:=My[I] to My[I+1]-1 do
       begin
         if Ms[J]<>Nv1 then goto 2;
         My1[1]:=1;
         Isu:=Mdop[I];
         My1[2]:=2;
         Ms1[1]:=Nv1;
         My1[3]:=3;
         Ms1[2]:=I;
         Irr:=Isu;
         Imum:=2;
     3:     Irr:=Irr-1;
         if Irr=1 then goto 22;
         Istart:=My1[Isu-Irr+1];
         Istop:=My1[Isu-Irr+2]-1;
         for II:=1 to Nv do
         begin
           if Mdop[II]<>Irr then goto 5;
           for JJ:=My[II] to My[II+1]-1 do
           begin
             for JJJ:=Istart to Istop do
             begin
               if Ms[JJ]<>Ms1[JJJ] then goto 7;
               Imum:=Imum+1;
               Ms1[Imum]:=II;
     7:         end;
           end;
     5:    end;
        Kot1:=Istop+1-Imum;
        if Kot1=0 then goto 9;
        for II:=Istop+1 to Imum-1 do
        begin
          if Ms1[II]=0 then goto 8;
          for JJ:=II+1 to Imum do
           if Ms1[JJ]=Ms1[II] then Ms1[JJ]:=0;
     8:    end;
     9:    KKK:=Istop;
        for II:=Istop+1 to Imum do
        begin
          if Ms1[II]=0 then goto 10;
          KKK:=KKK+1;
          Ms1[KKK]:=Ms1[II];
     10:    end;
        Imum:=KKK;
        My1[Isu-Irr+3]:=Imum+1;
        Kot2:=Isu-Irr+2;
        goto 3;
     2:   end;
```

```
          goto 1;
        22:   I2:=0;
        11:   Kzikl:=Kzikl+1;
          FormSoasda(Nv1,Isu,I2,Key1,My1,Ms1,My,Ms,Mcg,Mcg1,Kol);
          if Key1=false then goto 13;
          MyT[Kzikl+1]:=MyT[Kzikl]+Isu+1;
          for Ji:=1 to Isu+1 do
          begin
            Kot:=Kot+1;
            MsT[Kot]:=Mcg1[JI];
          end;
          goto 11;
        13:   Kzikl:=Kzikl-1;
        1:  end;
     end; {FormWegin}
{*************************************************************}
procedure FormIncide(var Nv : integer;
                     var My : TMasy;
                     var Ms : TMass;
                     var Ms3 : TMass);
{ Nv - количество вершин в графе;                              }
{ My  - массив указателей для матрицы смежностей;              }
{ Ms  - массив элементов матрицы смежностей;                   }
{ Ms3 - массив элементов матрицы инциденций.                   }
{    Формируется матрица инциденций графа в массиве            }
{    Ms3.                                                      }
{                                                              }
        var I,J,K,NNN,P,M,L : integer;
        begin
{       инициализация                          }
        NNN:=My[Nv+1]-1;
        K:=0;
        for J:= 1 to NNN do Ms3[J]:= 0;
{   определение номера элемента                 }
        for I:= 1 to Nv do
         for M:= My[I] to My[I+1]-1 do
          if Ms3[M]=0 then
          begin
            P:=Ms[M];
            K:=K+1;
            Ms3[M]:=K;
            for L:=My[P] to My[P+1]-1 do
             if Ms[L]=I then Ms3[L]:=K;
          end;
        end; {FormIncide}
{***********************************************************}
 procedure EinZikle(var Nv,Kzikl,M : integer;
                    var Masy : TMasy;
                    var Mass : TMass;
                    var Mass1 : TMass;
                    var Masi : TMass;
                    var MasyT : TMasy;
                    var MassT : TMass;
                    var MasMy1 : TMasy;
                    var MasMs1 : TMass;
                    var MasMy2 : TMasy;
                    var MasMs2 : TMass;
                    var MasMy3 : TMasy;
                    var MasMs3 : TMass;
                    var MasMdop : Tmasy;
                    var MasMcg : Tmasy;
                    var MasMcg1 : Tmasy;
                    var MasKol : Tmasy);
```

```
{Процедура создания множества единичных циклов графа                }
{                                                                    }
{ Nv - количество вершин в графе;                                    }
{ Kzikl  - количество единичных циклов в графе;                      }
{ Masy  - массив указателей для матрицы смежностей;                  }
{ Mass  - массив элементов матрицы смежностей;                       }
{ Mass1 - массив для несмежных элементов строки.                     }
{ Masi -  массив элементов матрицы инциденций;                       }
{ MasyT - массив указателей для матрицы единичных циклов;            }
{ MassT - массив элементов матрицы единичных циклов;                 }
{ MasMy1 -  вспомогательный массив;                                  }
{ MasMs1 -  вспомогательный массив;                                  }
{ MasMy2 -  вспомогательный массив;                                  }
{ MasMs2 -  вспомогательный массив;                                  }
{ MasMy3 -  вспомогательный массив;                                  }
{ MasMs3 -  вспомогательный массив;                                  }
{ MasMdop -  вспомогательный массив для хранения уровней;            }
{ MasMcg -  вспомогательный массив;                                  }
{ MasMcg1 -  вспомогательный массив;                                 }
{ MasKol -  вспомогательный массив.                                  }
label 4;
var I,J,JJ,Nv1,Nv2,Pr,Kzikl1,Kzikl2,Ip1: integer;
begin
        Pr:= Masy[Nv+1]-1;
        M:= Pr div 2;
        FormIncide(Nv,Masy,Mass,Masi);
        Kzikl:= 0;
        MasyT[1]:= 1;
        for I:= 1 to M do
        begin {1}
        Ip1:= I;
          for J:= 1 to Nv do
          begin {2}
            for JJ:=Masy[J] to Masy[J+1]-1 do
            begin {3}
              if Masi[JJ] = Ip1 then
              begin {4}
                Nv1:= J; {определение первой вершины ребра}
                Nv2:= Mass[JJ]; {определение второй вершины ребра}
                goto 4; {концевые вершины ребра определены}
              end;  {4}
            end;  {3}
          end;  {2}
4:        FormVolna(Nv,Nv1,Nv2,Masy,Mass,MasMdop); {алгоритм поиска в ширину}
           {для ориентированного ребра (Nv1,Nv2)}
          FormWegin(Nv,Nv1,Nv2,Kzikl1,Masy,Mass,MasMy3,MasMs3,
          MasMy1,MasMs1,MasMdop,MasMcg,MasMcg1,MasKol);
{построение кандидатов в единичные циклы проходящих по ребру (Nv1,Nv2)}
          FormKpris(Kzikl1,MasMy1,MasMs1,Masy,Mass,Masi);
          FormVolna(Nv,Nv2,Nv1,Masy,Mass,MasMdop);
{алгоритм поиска в ширину для ориентированного ребра (Nv2,Nv1)}
          FormWegin(Nv,Nv2,Nv1,Kzikl2,Masy,Mass,MasMy3,MasMs3,
          MasMy2,MasMs2,MasMdop,MasMcg,MasMcg1,MasKol);
{построение кандидатов в единичные циклы проходящих по ребру (Nv2,Nv1)}
          FormKpris(Kzikl2,MasMy2,MasMs2,Masy,Mass,Masi);
          FormSwigug(Kzikl1,Kzikl2,MasMy1,MasMs1,MasMy2,MasMs2,
          MasMcg,MasMcg1,MasKol);
          FormDozas(Kzikl,Kzikl1,MasMy1,MasMs1,MasyT,MassT,
          MasMcg,MasMcg1,MasKol);
        end;  {1}
end;{EinZikle}
{*********************************************************}
procedure  Shell(var N : integer;
```

```
                var A : TMasy);
{       процедура Шелла для упорядочивания элементов                          }
{                                                                             }
{     N - количество элементов в массиве;                                     }
{     A - сортируемый массив;                                                 }
var D,Nd,I,J,L,X : integer;
label 1,2,3,4,5;
begin
  D:=1;
1:D:=2*D;
  if D<=N then goto 1;
2:D:=D-1;
  D:=D div 2;
  if D=0 then goto 5;
  Nd:=N-D;
  for I:=1 to Nd do
  begin
    J:=I;
3:  L:=J+D;
    if A[L]>=A[J] then goto 4;
    X:=A[J];
    A[J]:=A[L];
    A[L]:=X;
    J:=J-D;
    if J>0 then goto 3;
4:end;
  goto 2;
5:end;{Shell}
{*************************************************************}
procedure  ProzYpor(var N : integer;
           var Masy: TMasy;
           var Mass: TMass;
           var A : TMasy);
{Расположение элементов массива в порядке возрастания                  }
var i,j,K: integer;
begin
    for i:= 1 to N do
    begin
      K:=0;
      for j:= MasY[i] to MasY[i+1]-1 do
      begin
        K:=K+1;
        A[K]:= Mass[j];
      end;
      Shell(K,A);
      for j:= 1 to K do Mass[MasY[i]-1+j]:=A[j];
    end;
end;{ProzYpor}
{*******************************************************************}
begin
         assign(F1,'D:\Isomorf\GRF\Петерсен.grf');
         reset(F1);
         readln(F1,Nv);
         for I:=1 to Nv+1 do
         begin
           if I<>Nv+1 then read(F1,Masy[I]);
           if I=Nv+1 then readln(F1,Masy[I]);
         end;
         Np:=Masy[Nv+1]-1;
         for I:=1 to Np do
         begin
           if I<>Np then read(F1,Mass[I]);
           if I=Np then read(F1,Mass[I]);
```

```
    end;
    close (F1);
    { Создаём новый файл и открываем его в режиме "для чтения и записи"}
    Assign(F2,'D:\Isomorf\EZI\Петерсен.ezi');
    Rewrite(F2);
    EinZikle(Nv,Kzikl,M,Masy,Mass,Mass1,Masi,MasyT,MassT,
            MasMy1,MasMs1,MasMy2,MasMs2,MasMy3,MasMs3,MasMdop,
            MasMcg,MasMcg1,MasKol);
    ProzYpor(Kzikl,MasyT,MassT,MasKol);
    writeln(F2,Nv);
    writeln(F2,M);
    writeln(F2,Kzikl);
    for I:=1 to Nv+1 do
    begin
     if i<>Nv+1 then write(F2,Masy[i],' ');
     if i=Nv+1 then writeln(F2,Masy[i]);
    end;
    for I:=1 to Nv do
    for j:=Masy[i] to Masy[i+1]-1 do
    begin
     if j<>Masy[i+1]-1 then write(F2,Mass[j],' ');
     if j=Masy[i+1]-1 then writeln(F2,Mass[j]);
    end;
    for I:=1 to Nv do
    for j:=Masy[i] to Masy[i+1]-1 do
    begin
     if j<>Masy[i+1]-1 then write(F2,Masi[j],' ');
     if j=Masy[i+1]-1 then writeln(F2,Masi[j]);
    end;
    for i:=1 to Kzikl+1 do
    begin
      if i<> Kzikl+1 then write(F2, MasyT [i],' ');
      if i= Kzikl+1 then writeln(F2, MasyT [i]);
    end;
    for I:=1 to Kzikl do
    begin
     for j:=MasyT[i] to MasyT[i+1]-1 do
     begin
       if j<>MasyT[i+1]-1 then write(F2,MassT[j],' ');
       if j=MasyT[i+1]-1 then writeln(F2,MassT[j]);
     end;
    end;
    close (F2);
    writeln('Конец расчета!');
end.
```

## 1.9. Входные и выходные файлы процедуры Raschet1

В качестве примера рассмотрим граф Петерсена (см. рис. 1.23).

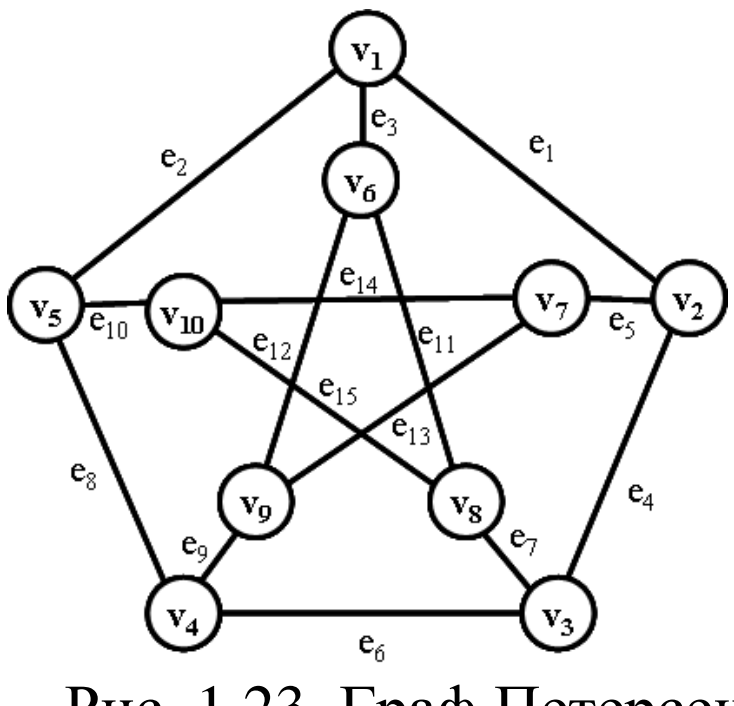

Рис .1.23. Граф Петерсена.

Входной файл Петерсен.grf

| | |
|---|---|
| 10 | {количество вершин} |
| 1 4 7 10 13 16 19 22 25 28 31 | {массив указателей} |
| 2 5 6 | {смежность вершины 1} |
| 1 3 7 | {смежность вершины 2} |
| 2 4 8 | {смежность вершины 3} |
| 3 5 9 | {смежность вершины 4} |
| 1 4 10 | {смежность вершины 5} |
| 1 8 9 | {смежность вершины 6} |
| 2 9 10 | {смежность вершины 7} |
| 3 6 10 | {смежность вершины 8} |
| 4 6 7 | {смежность вершины 9} |
| 5 7 8 | {смежность вершины 10} |

Выходной файл Петерсен.ezi

| | |
|---|---|
| 10 | {количество вершин} |
| 15 | {количество ребер} |
| 12 | {количество изометрических циклов} |
| 1 4 7 10 13 16 19 22 25 28 31 | {массив указателей для матрицы смежностей} |
| 2 5 6 | {смежность вершины 1} |
| 1 3 7 | {смежность вершины 2} |
| 2 4 8 | {смежность вершины 3} |
| 3 5 9 | {смежность вершины 4} |
| 1 4 10 | {смежность вершины 5} |
| 1 8 9 | {смежность вершины 6} |
| 2 9 10 | {смежность вершины 7} |
| 3 6 10 | {смежность вершины 8} |
| 4 6 7 | {смежность вершины 9} |
| 5 7 8 | {смежность вершины 10} |
| 1 2 3 | {совместимость ребер для вершины 1} |
| 1 4 5 | {совместимость ребер для вершины 2} |
| 4 6 7 | {совместимость ребер для вершины 3} |
| 6 8 9 | {совместимость ребер для вершины 4} |
| 2 8 10 | {совместимость ребер для вершины 5} |
| 3 11 12 | {совместимость ребер для вершины 6} |
| 5 13 14 | {совместимость ребер для вершины 7} |
| 7 11 15 | {совместимость ребер для вершины 8} |
| 9 12 13 | {совместимость ребер для вершины 9} |
| 10 14 15 | {совместимость ребер для вершины 10} |
| 1 6 11 16 21 26 31 36 41 46 51 56 61 | {массив указателей для циклов} |
| 1 2 4 6 8 | {изометрический цикл 1 в ребрах} |
| 1 2 5 10 14 | {изометрический цикл 2 в ребрах } |
| 1 3 4 7 11 | {изометрический цикл 3 в ребрах } |
| 1 3 5 12 13 | {изометрический цикл 4 в ребрах } |
| 2 3 10 11 15 | {изометрический цикл 5 в ребрах } |
| 2 3 8 9 12 | {изометрический цикл 6 в ребрах } |
| 4 5 6 9 13 | {изометрический цикл 7 в ребрах } |
| 4 5 7 14 15 | {изометрический цикл 8 в ребрах } |
| 6 7 9 11 12 | {изометрический цикл 9 в ребрах } |
| 6 7 8 10 15 | {изометрический цикл 10 в ребрах } |

8 9 10 13 14 {изометрический цикл 11 в ребрах }
11 12 13 14 15 {изометрический цикл 12 в ребрах }

**Примечания**

**примечание 1**. Массив указателей формируется как сложение локальных степеней вершин с указанием местоположения в списке смежных вершин

Например: 1 4 7 10 13 16 19 22 25 28 31 {MY - массив указателей}

Нужно определить смежность для вершины 4.

Элементы для 4-ой вершины начинаются с MY(4)=10 и заканчивается MY(5)-1=13-1=12.

Следовательно, вершины смежные к вершине 4 находятся с 10 по 12 местах в списке вершин

2 5 6 1 3 7 2 4 8 **3 5 9** 1 4 10 1 8 9 2 9 10 3 6 10 4 6 7 5 7 8

**примечание 2.** Цифрами, в машинной информации, обозначаются либо вершины, либо ребра графа в зависимости от их принадлежности.

### 1.10. Текст программы PrintZikNabc для вычисления изометрических циклов как множество вершин

```
type
      TMasy = array[1..1000] of integer;
      TMass = array[1..4000] of integer;
var
      F1,F2 : text;
      i,ii,j,jj,K,K1,M,Nv,Kzikl,Nr : integer;
      Siz,Np,Nk,M1,KCF,KCF1,Ln : integer;
      Masy: TMasy;
      Mass: TMass;
      Masi: TMass;
      MasMy1: TMasy;
      MasMs1: TMass;
      MasMy2: TMasy;
      MasMs2: TMass;
      MasMy5: TMasy;
      MasMs5: TMass;
      MasMy4: TMasy;
      MasMs4: TMass;
      MasMdop: TMasy;
{************************************************************}
procedure FormVerBasis1(var Nv,Ziklo : integer;
                  var My : TMasy;
                  var Ms : TMass;
                  var My1 : TMasy;
                  var Ms1 : TMass;
                  var Mass1 : TMass;
                  var Mdop : TMasy;
                  var Masi : TMass);
{Процедура записи базиса циклов через вершины                }
{                                                            }
{ Nv - количество вершин в графе;                            }
{ Ziklo - цикломатическое число графа;                       }
{ My - массив указателей для матрицы смежностей;             }
```

```
{ Ms - массив элементов матрицы смежностей;                                    }
{ My1 - массив указателей для матрицы базисных циклов                          }
{ Ms1 : массив элементов матрицы базисных циклов;                              }
{ Masi - массив элементов матрицы инциденций;                                  }
{ Mass1 - элементов матрицы базисных циклов через вершины;                     }
{ Mdop - вспомогательный массив;                                               }
{**********************************************************************}
label 1,2,3;
var i,j,ii,jj,iii,jjj,Nr,KKK,K1 :integer;
begin
  for i:= 1 to Ziklo do
  begin{1}
   {writeln(F2,'FormVerBasis: i = ',i);}
    KKK:=0;
    j:=My1[i]-1;
    2:
    j:=j+1;
    if j = My1[i+1] then goto 3 else
    begin {2}
     Nr:= Ms1[j];
     {writeln(F2,'FormVerBasis: Nr = ',Nr);}
     for ii:=1 to Nv do
     begin {3}
      for jj:= My[ii] to My[ii+1]-1 do
      begin {4}
       if Masi[jj]= Nr  then
       begin {5}
        KKK:=KKK+1;
        Mdop[KKK]:=ii;
        KKK:=KKK+1;
        Mdop[KKK]:=Ms[jj];
        {for jjj:= 1 to KKK do
        begin
         if jjj <> KKK then write(F2, Mdop[jjj],' ');
         if jjj = KKK then writeln(F2, Mdop[jjj],' ');
        end;}
        goto 2;
       end; {5}
      end;{4}
     end;{3}
    end; {2}
    3: {writeln(F2,'Запись [FormVerBasis]');}
    for jjj:=1 to KKK-1 do if Mdop[jjj]<>0 then
    begin {7}
     for iii:=jjj+1 to KKK do if Mdop[iii]<>0 then
     begin {8}
      if Mdop[iii]=Mdop[jjj] then Mdop[iii]:=0;
      begin {9}
       K1:=0;
       for jj:=1 to KKK do if Mdop[jj]<>0 then
       begin {6}
        K1:=K1+1;
        Mdop[K1]:=Mdop[jj];
       end; {6}
      end; {9}
     end; {8}
    end; {7}
    for jjj:=1 to K1 do Mass1[My1[i]-1 +jjj]:=Mdop[jjj];
 end; {1}
 {for i:= 1 to Ziklo do
   begin
    for j:= My1[i] to My1[i+1]-1 do
    begin
```

```
       if j<> My1[i+1]-1 then write(F2, Mass1[j],' ');
       if j= My1[i+1]-1 then writeln(F2, Mass1[j],' ');
      end;
    end;}
end;{FormVerBasis1}
{***********************************************************************}
procedure  Shell(var N : integer;
                 var A : TMasy);
{       процедура Шелла для упорядочивания элементов        }
{                                                   }
{     N - количество элементов в массиве;                   }
{     A - сортируемый массив;                          }
var D,Nd,I,J,L,X : integer;
label 1,2,3,4,5;
begin
  D:=1;
1:D:=2*D;
  if D<=N then goto 1;
2:D:=D-1;
  D:=D div 2;
  if D=0 then goto 5;
  Nd:=N-D;
  for I:=1 to Nd do
  begin
    J:=I;
3:  L:=J+D;
    if A[L]>=A[J] then goto 4;
    X:=A[J];
    A[J]:=A[L];
    A[L]:=X;
    J:=J-D;
    if J>0 then goto 3;
4:end;
  goto 2;
5:end;{Shell}
{***************************************************************}
label 1;
begin

         assign(F1,'D:\Isomorf\EZI\Петерсен.ezi');
         reset(F1);
         readln(F1,Nv);
         readln(F1,Nr);
         readln(F1,Kzikl);
         {Вводим матрицу смежностей неориентированного графа }
         for I:=1 to Nv+1 do
         begin {1}
           if i<>Nv+1 then read(F1,Masy[i]);
           if i=Nv+1 then readln(F1,Masy[i]);
         end;  {1}
         for I:=1 to Nv do
         for j:=Masy[i] to Masy[i+1]-1 do
         begin {2}
           if j<>Masy[i+1]-1 then read(F1,Mass[j]);
           if j=Masy[i+1]-1 then readln(F1,Mass[j]);
         end;  {2}
         {Вводим матрицу инциденций неориентированного графа }
         for I:=1 to Nv do
         for j:=Masy[i] to Masy[i+1]-1 do
         begin
           if j<>Masy[i+1]-1 then read(F1,Masi[j]);
           if j=Masy[i+1]-1 then readln(F1,Masi[j]);
         end;
```

```
{Вводим единичные циклы }
for I:=1 to Kzikl+1 do
begin
  if i<>Kzikl+1 then read(F1,MasMy4[i]);
  if i=Kzikl+1 then readln(F1,MasMy4[i]);
end;
for I:=1 to Kzikl do
begin
  for j:=MasMy4[i] to MasMy4[i+1]-1 do
  begin
    if j<>MasMy4[i+1]-1 then read(F1,MasMs4[j]);
    if j=MasMy4[i+1]-1 then readln(F1,MasMs4[j]);
  end;
end;
close (F1);
FormVerBasis1(Nv,Kzikl,Masy,Mass,MasMy4,MasMs4,MasMs2,MasMy5,Masi);
{ Создаём новый файл и открываем его в режиме "для чтения и записи"}
Assign(F2,'D:\Isomorf\MY4\Петерсен.my4');
Rewrite(F2);
writeln(F2,'  Количество вершин графа = ',Nv);
writeln(F2,'  Количество ребер графа = ',Nr);
writeln(F2,'  Количество единичных циклов = ',Kzikl);
writeln(F2,' ');
writeln(F2,'  Матрица смежностей графа: ');
writeln(F2,' ');
for I:=1 to Nv do
begin
    write(F2,' вершина ',I:3,':');
    Siz:=0;
    for K:=Masy[I] to Masy[I+1]-1 do
    begin
        Siz:=Siz+1;
        Np:=Siz mod 10;
        Nk:=Siz div 10;
        if K<>Masy[I+1]-1 then
        begin
         if (Np=1) and (Nk>0) then
         write(F2,' ',Mass[K]:3);
         if (Np=1) and (Nk=0) then
         write(F2,' ',Mass[K]:3);
         if (Np>1) and (Np<=9) then
         write(F2,' ',Mass[K]:3);
         if Np=0 then
          writeln(F2,' ',Mass[K]:3);
       end;
       if K=Masy[I+1]-1 then
       begin
         if (Np=1) and (Nk>0) then
          writeln(F2,' ',Mass[K]:3)
         else
        writeln(F2,' ',Mass[K]:3);
       end;
    end;
end;
writeln(F2,' ');
writeln(F2,'  Элементы матрицы инциденций: ');
writeln(F2,' ');
M1:=Masy[Nv+1]-1;
M:=M1 div 2;
for I:=1 to M do
begin
  for  K:=1 to Nv do
  begin
```

```
      for  M1:=Masy[K] to Masy[K+1]-1 do
      begin
        if Masi[M1]=I then
        begin
          write(F2,' ребро ',I:3,':');
          Np:=Mass[M1];
          write(F2,'  ( ',K:3);
          write(F2,' ',Np:3,' )  или  ');
          write(F2,' ( ',Np:3);
          writeln(F2,'  ',K:3,' )');
          goto 1;
        end;
      end;
     end;
1:    end;
    writeln(F2,' ');
    writeln(F2,'  Множество единичных циклов графа: ');
    writeln(F2,' ');
    for  I:=1 to Kzikl do
    begin
      write(F2,' цикл ',I:3,':');
      Siz:=0;
      for K:=MasMy4[I] to MasMy4[I+1]-1 do
      begin
        Siz:=Siz+1;
        Np:=Siz mod 10;
        Nk:=Siz div 10;
        if K<>MasMy4[I+1]-1 then
        begin
         if (Np=1) and (Nk>0) then
           write(F2,' ',MasMs4[K]:3);
         if (Np=1) and (Nk=0) then
          write(F2,' ',MasMs4[K]:3);
         if (Np>1) and (Np<=9) then
          write(F2,' ',MasMs4[K]:3);
         if Np=0 then
          writeln(F2,' ',MasMs4[K]:3);
        end;
        if K=MasMy4[I+1]-1 then
        begin
         if (Np=1) and (Nk>0) then
           writeln(F2,' ',MasMs4[K]:3)
         else
          writeln(F2,' ',MasMs4[K]:3);
        end;
      end;
    end;
    writeln(F2,' ');
    for I:=1 to Kzikl do
    begin
     Ln:= MasMy4[i+1] -MasMy4[I];
     ii:=0;
     for j:=MasMy4[i] to MasMy4[I+1]-1 do
     begin
      ii:=ii+1;
      MasMdop[ii]:=MasMs2[j];
     end;
     Shell(Ln,MasMdop);
     ii:=0;
     for j:=MasMy4[i] to MasMy4[I+1]-1 do
     begin
      ii := ii+1;
      MasMs2[j]:=MasMdop[ii];
```

```
      end;
    end;
    writeln(F2,'  Множество единичных циклов графа в записи вершин: ');
    writeln(F2,' ');
    for I:=1 to Kzikl do
    begin
      write(F2,' цикл ',I:3,':');
      Siz:=0;
      for K:=MasMy4[I] to MasMy4[I+1]-1 do
      begin
        Siz:=Siz+1;
        Np:=Siz mod 10;
        Nk:=Siz div 10;
        if K<>MasMy4[I+1]-1 then
        begin
         if (Np=1) and (Nk>0) then
           write(F2,' ',MasMs2[K]:3);
         if (Np=1) and (Nk=0) then
          write(F2,' ',MasMs2[K]:3);
         if (Np>1) and (Np<=9) then
          write(F2,' ',MasMs2[K]:3);
         if Np=0 then
          writeln(F2,' ',MasMs2[K]:3);
       end;
       if K=MasMy4[I+1]-1 then
       begin
        if (Np=1) and (Nk>0) then
          writeln(F2,' ',MasMs2[K]:3)
        else
         writeln(F2,' ',MasMs2[K]:3);
       end;
      end;
    end;
    close (F2);
    writeln('Конец расчета!');
end.
```

## 1.11. Входные и выходные файлы программы PrintZikNabc

Выходной файл Петерсен.my4

Количество вершин графа = 10
Количество ребер графа = 15
Количество единичных циклов = 12

Матрица смежностей графа:

```
вершина   1:   2   5   6
вершина   2:   1   3   7
вершина   3:   2   4   8
вершина   4:   3   5   9
вершина   5:   1   4  10
вершина   6:   1   8   9
вершина   7:   2   9  10
вершина   8:   3   6  10
вершина   9:   4   6   7
вершина  10:   5   7   8
```

Элементы матрицы инциденций:

ребро 1: ( 1 2 ) или ( 2 1 )
ребро 2: ( 1 5 ) или ( 5 1 )
ребро 3: ( 1 6 ) или ( 6 1 )
ребро 4: ( 2 3 ) или ( 3 2 )
ребро 5: ( 2 7 ) или ( 7 2 )
ребро 6: ( 3 4 ) или ( 4 3 )
ребро 7: ( 3 8 ) или ( 8 3 )
ребро 8: ( 4 5 ) или ( 5 4 )
ребро 9: ( 4 9 ) или ( 9 4 )
ребро 10: ( 5 10 ) или ( 10 5 )
ребро 11: ( 6 8 ) или ( 8 6 )
ребро 12: ( 6 9 ) или ( 9 6 )
ребро 13: ( 7 9 ) или ( 9 7 )
ребро 14: ( 7 10 ) или ( 10 7 )
ребро 15: ( 8 10 ) или ( 10 8 )

Множество единичных циклов графа в записи ребер:

цикл 1: 1 2 4 6 8
цикл 2: 1 2 5 10 14
цикл 3: 1 3 4 7 11
цикл 4: 1 3 5 12 13
цикл 5: 2 3 10 11 15
цикл 6: 2 3 8 9 12
цикл 7: 4 5 6 9 13
цикл 8: 4 5 7 14 15
цикл 9: 6 7 9 11 12
цикл 10: 6 7 8 10 15
цикл 11: 8 9 10 13 14
цикл 12: 11 12 13 14 15

Множество единичных циклов графа в записи вершин:

цикл 1: 1 2 3 4 5
цикл 2: 1 2 5 7 10
цикл 3: 1 2 3 6 8
цикл 4: 1 2 6 7 9
цикл 5: 1 5 6 8 10
цикл 6: 1 4 5 6 9
цикл 7: 2 3 4 7 9
цикл 8: 2 3 7 8 10
цикл 9: 3 4 6 8 9
цикл 10: 3 4 5 8 10
цикл 11: 4 5 7 9 10
цикл 12: 6 7 8 9 10

## Комментарии

Свойство минимальности расстояний между парой вершин графа позволяет выделить подмножество изометрических циклов в отдельный вид квазициклов с валентностью вершин равной двум. На основе метрических свойств графа, определено понятие изометрического цикла графа. Причем мощность подмножества изометрических циклов значительно меньше мощности множества простых циклов. И это свойство определяет меньшую вычислительную сложность алгоритмов использующих для вычислений подмножество изометрических циклов.

Показано отличие изометрических циклов от фундаментальных циклов графа. Доказана теорема о существовании базиса во множестве изометрических циклов несепарабельного графа. Показана связь маршрутов определяемых теоремой Менгера и маршрутами в изометрических циклах. Представлен алгоритм построения множества изометрических циклов для неориентированного несепарабельного графа. Рассматривая пример построения множества изометрических циклов, сделан вывод о необходимости проведения маршрутов минимальной длины относительно всех ребер графа.

Рассмотрены свойства изометрических циклов, для подграфов построенных методом удаления ребер из полного графа. Рассмотрены вопросы построения инвариантов графа на основе понятия изометрического цикла. Приведены примеры выделения множество изометрических циклов графа.

Вычислительную сложность алгоритма выделения множества изометрических циклов графа можно определить как $O(n^4)$.

# Глава 2. ЭЛЕМЕНТЫ АЛГЕБРЫ СТРУКТУРНЫХ ЧИСЕЛ И МАТРОИДЫ

## 2.1. Определение структурного числа

Пусть $\aleph$ - подмножество абстрактного пространства $\Re$. Элементы множества $\aleph$ обозначим

$\alpha_i, \beta_i, \gamma_i, ... \in \aleph$

Рассмотрим систему элементов в виде таблицы

$$\mathbf{A} = \begin{vmatrix} a_{11} & a_{12} & \dots & a_{1n} \\ a_{21} & a_{22} & \dots & a_{2n} \\ \dots & \dots & \dots & \dots \\ a_{m1} & a_{m2} & \dots & a_{mn} \end{vmatrix} \qquad (2.1)$$

Будем рассматривать эту систему как совокупность столбцов $\alpha_k$, т.е.

$$A = \{a_{\square}, a_2, ..., a_n\}, \quad a_i \neq a_j \quad (i \neq j). \qquad (2.2)$$

Столбцы $a_k$ в свою очередь представляют собой неупорядоченные множества элементов $a_{ik}$

$$a_k = \{a_{1k}, a_{2k, \dots,} a_{mk}\} \quad a_{ik} \neq a_{jk} \qquad (2.3)$$

Столбцы будем считать равными, если они содержат одинаковые элементы. Положим по определению, что система (2.1) не содержит одинаковых столбцов. Систему типа (2.1) будем рассматривать как элемент новой алгебры - алгебры структурных чисел. Согласно определениям абстрактной алгебры, алгебру структурных чисел можно отнести к категории операторных алгебр, т.е. ее можно характеризовать упорядоченной тройкой $<E,\Omega,e>$ [7], где $E$ - носитель алгебры (в нашем случае семейство множеств); $\Omega$ - двухэлементное множество операторов $w_1, w_2$ определяющих сумму и произведение; е - результат, т.е. функция, которая выражению $AwB$ ставит в соответствие элемент $C \in E$, являющийся результатом действия [7].

Введем вспомогательное понятие, которое используем при определении структурного числа.

Рассмотрим последовательность элементов $z_i$, необязательно различных:

$$<z_1, z_2, ..., z_i, ..., z_n> \qquad (2.4)$$

Обозначим через $r(z_k)$ - число одинаковых элементов последовательности (6.4).

Структурным числом называется система элементов $\alpha_{ik}$ вида (2.1) [с учетом (2.2) и

(2.3)], удовлетворяющая следующим определениям .

**Определение 2.1.** Два структурных числа считаются равными (A=B) тогда и только тогда, когда ($a \in$ A) $\Leftrightarrow$ ($a \in$ B) или

$$A = B \Leftrightarrow \underset{a}{\wedge} (a \in A \Leftrightarrow a \in B) \qquad (2.5)$$

**Определение 2.2.** Суммой структурных чисел A и B называется структурное число

$$C = \{z \mid (z \in A) \wedge (z \in B), z \notin A \cup B\} = A\square B; \qquad (2.6)$$

В этом случае можно написать C = A + B.

Выражение A□B в формуле (2.6) означает симметричную разность множеств A и B.

**Определение 2.3.** Произведением структурных чисел A и B называется структурное число

$$C = \{a \cup b \mid a \cap b = \varnothing, r(a \cup b) \in \{1,3,...,\}, a \in A, b \in B\}, \qquad (2.7)$$

которое записывается в виде C = AB.

В соответствии с определением суммы при сложении структурных чисел опускаются столбцы, одновременно присутствующие в обоих числах A и B, а в соответствии с определением произведения при умножении структурных чисел A и B опускаются те столбцы $a \cup b$, в которых какой-либо элемент повторяется, т.е. для которых $a \cap b \neq \varnothing$, а также опускается четное число идентичных столбцов.

Можно заметить, что равенство структурных чисел представляет собой отношение эквивалентности, т.е. является рефлексивным, симметричным и транзитивным.

***Пример 2.1.*** Равенство структурных чисел

$$\begin{vmatrix} 1 & 3 & 2 \\ 3 & 2 & 5 \end{vmatrix} = \begin{vmatrix} 5 & 1 & 3 \\ 2 & 3 & 2 \end{vmatrix} = \begin{vmatrix} 2 & 2 & 3 \\ 3 & 5 & 1 \end{vmatrix}$$

Два структурных числа равны, если содержат идентичные столбцы, независимо от порядка элементов в столбцах и порядка столбцов.

***Пример 2.2.*** Сложение структурных чисел

$$\begin{vmatrix} 2 & 3 & 4 \\ 7 & 5 & 7 \end{vmatrix} + \begin{vmatrix} 3 & 2 & 5 \\ 2 & 7 & 4 \end{vmatrix} = \begin{vmatrix} 2 & 3 & 4 & 3 & 2 & 5 \\ 7 & 5 & 7 & 2 & 7 & 4 \end{vmatrix} = \begin{vmatrix} 3 & 4 & 3 & 5 \\ 5 & 7 & 2 & 4 \end{vmatrix}$$

***Пример 2.3.*** Умножение структурных чисел

$$\begin{vmatrix} 1 & 3 & 1 \\ 2 & 4 & 3 \end{vmatrix} \times \begin{vmatrix} 3 & 1 & 4 \\ 4 & 2 & 5 \end{vmatrix} = \begin{vmatrix} 1 & 1 & 1 & 3 & 3 & 3 & 1 & 1 & 1 \\ 2 & 2 & 2 & 4 & 4 & 4 & 3 & 3 & 3 \\ 3 & 1 & 4 & 3 & 1 & 4 & 3 & 1 & 4 \\ 4 & 2 & 5 & 4 & 2 & 5 & 4 & 2 & 5 \end{vmatrix} = \begin{vmatrix} 1 & 1 & 3 & 1 \\ 2 & 2 & 4 & 3 \\ 3 & 4 & 1 & 4 \\ 4 & 5 & 2 & 5 \end{vmatrix} = \begin{vmatrix} 1 & 1 \\ 2 & 3 \\ 4 & 4 \\ 5 & 5 \end{vmatrix}$$

Произведением двух структурных чисел A и B называется структурное число, столбцы которого представляют собой суммы (согласно понятиям теории множеств) всех возможных комбинаций столбцов А и В, за исключением наибольшего четного числа идентичных столбцов и таких столбцов, в которых какой-либо элемент повторяется (произведение других столбцов не содерит).

Из определения суммы и произведения структурных чисел следует, что эти операции всегда можно выполнить на множестве этих чисел. Из тех же определений можно сделать вывод, что сложение и умножение структурных чисел коммутативны и ассоциативны, а умножение дистрибутивно относительно сложения.

Для трех произвольных структурных чисел имеют место следующие соотношения, подобные тем, которые справедливы для элементарной алгебры:

$$A + B = B + A,$$

$$AB = BA,$$

$$A(BC) = (AB)C,$$

$$A(B + C) = AB + AC. \qquad (2.8)$$

Следует различать структурное число [∅], содержащее один столбец, который является пустым множеством ∅, и структурное число [7], не содержащее ни одного столбца.

## 2.2. Геометрическое изображение структурного числа

Попробуем дать геометрическую интерпретацию структурного числа. Следует отметить, что геометрическая интерпретация встречается также и в других случаях, например в случае комплексных чисел, которым ставится в соответствие некоторые точки плоскости.

**Определение 2.4.** Если столбцы структурного числа А взаимно однозначно соответствуют остовам графа **G** так, что каждый столбец представляет собой множество значений описывающей функции соответствующего остова, то граф **G** называется геометрическим изображением числа А и записывается в виде

$$G = ob(A). \qquad (2.9)$$

Следовательно геометрическим изображением структурного числа А служит любой детерминированный граф, удовлетворяющий условию (2.9), или класс графов подобных структур. Из принятого определения следует, что геометрическое изображение структурного числа - не однозначное понятие, так как структурному числу может соответствовать многоэлементное семейство графов, составляющих класс с подобной структурой. Однако это в известном смысле является достоинством метода, так как становится возможным, например, в задачах синтеза электрических цепей, нахождение не одного, а множества вариантов цепи, удовлетворяющей заданным условиям.

**Теорема 2.1**[7]. Структурное число A с одинаковым числом элементов в строках, геометрическим изображением которого служит связный граф с вершинами $v_1, v_2, \ldots, v_n$, равно произведению $n$-1 простых однострочных сомножителей

$$A = V_1 V_2 \ldots V_{n-1} \qquad (2.10)$$

причем сомножители состоят из значений описывающей функции ребер, инцидентных произвольно выбранной вершине $v_i$ ($v_i \neq v_j$ , если $i \neq j$) графа **G**.

**Теорема 2.2**[7]. Необходимое и достаточное условие существования геометрического изображения структурного числа в виде связного графа состоят в том, чтобы структурное число A имело разложение на простые однострочные сомножители

$$A = V_1 V_2 \ldots V_p \qquad (2.11)$$

причем произвольный элемент $a_{ij}$ должен встречаться самое большее в двух простых числах $V_i$, $V_j$. Доказательство этих двух элегантных теорем приведено в [7]. Из теоремы 2.2 следует, что структурное число, у которого число элементов в строках различно, не имеет связного геометрического изображения.

## 2.3. Дополнительное структурное число и геометрическое обратное изображение

**Определение 2.5.** Дополнительным структурным числом для данного структурного числа A называется структурное число $A^d$, столбцы которого представляют собой дополнения столбцов числа A до множества элементов $a_{ij}$ из которых состоит структурное число A.

Если обозначить множество элементов $a_{ij}$, из которых состоит число A, через L, то столбцы $C_i^d$ числа $A^d$ определим как разность (в смысле понятий алгебры множеств)

$$C_1^d = L - C_1,\ C_2^d = L - C_2, \ldots, C_p^d == L - C_p \qquad (2.12)$$

где $C_1, C_2, \ldots, C_p$ - столбцы числа A.

Следует отметить справедливость такого свойства

$$(A + B)^d = A^d + B^d,\ L = L_A \cup L_B \qquad (2.13)$$

Способ получения дополнительного структурного числа иллюстрирует следующий пример.

***Пример 2.4.*** Определить структурное число $A^d$ по отношению к структурному числу

$$A = \begin{vmatrix} 2 & 6 & 9 \\ 3 & 2 & 8 \\ 5 & 8 & 3 \end{vmatrix}.$$

Множество элементов числа L таково:

L = {2, 3, 5, 6, 8, 9}.

Дополнительное структурное число равно

$$A^d = \begin{vmatrix} 6 & 3 & 2 \\ 8 & 5 & 5 \\ 9 & 9 & 6 \end{vmatrix}.$$

Оказывается, что для структурного числа удобно иметь дуальное геометрическое изображение, поэтому введем понятие обратного изображения геометрического структурного числа.

**Определение 2.6**. Граф **G** называется обратным изображением структурного числа A, если столбцы числа A взаимно однозначно соответствуют дополнениям остовов графа **G** так, что столбец числа A представляет собой множество значений описывающей функции соответствующего дополнения остова. Тогда напишем

$$G = \mathrm{cob}\ (A) \qquad (2.14)$$

Нетрудно заметить, что обратное изображение дополнительного числа $A^d$ одновременно служит изображением числа A и наоборот.

***Пример 2.5***. Граф G представленный на рис. 2.1 характеризуется структурным числом A

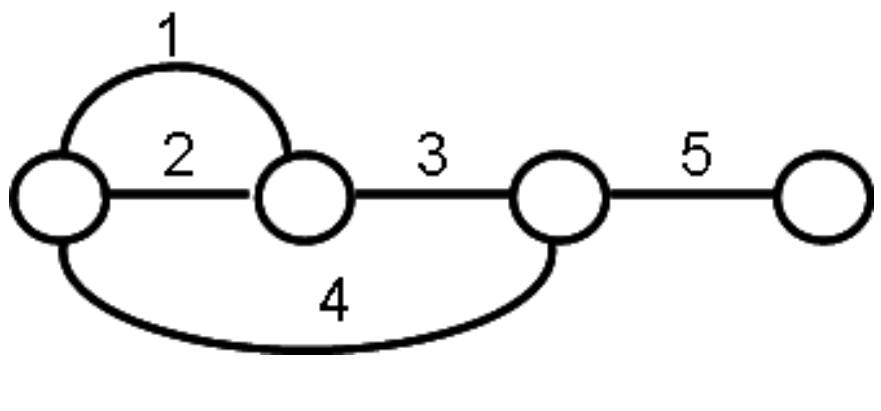

Рис. 2.1. Граф G.

Дополнительное структурное число по отношению к структурному числу A

$$A^d = \begin{vmatrix} 2 & 2 & 1 & 1 & 1 \\ 4 & 3 & 4 & 3 & 2 \end{vmatrix}.$$

Обратное геометрическое изображение структурного числа A представлено на рис. 2.2

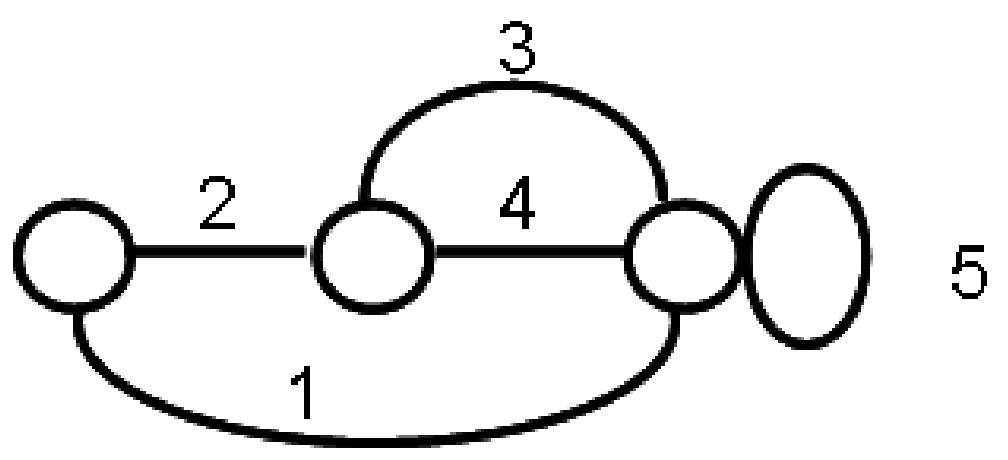

Рис. 2.2. Обратное геометрическое изображение структурного числа А.

### 2.4. Алгебраическая производная и обратная производная структурного числа

На множестве структурных чисел можно определить различные операции; одна из них - операция алгебраической производной.

**Определение 2.7**. Алгебраической производной структурного числа называется число $\partial A / \partial a$, определенное как

$$\frac{\partial A}{\partial a} = A|\text{столбцы, не содержащие элемент } a \text{ исключены.} \quad (2.15)$$

Легко доказать правильность следующих зависимостей, аналогичных «обычной» производной:

$$\frac{\partial A}{\partial a}(A_1 + A_2) = \frac{\partial A_1}{\partial a} + \frac{\partial A_2}{\partial a},$$

$$\frac{\partial A}{\partial a}(A_1 A_2) = \frac{\partial A_1}{\partial a} A_1 + \frac{\partial A_2}{\partial a} A_2. \quad (2.16)$$

Алгебраическую производную обозначим как $A_a$, т.е.

$$\frac{\partial A}{\partial a} = A_a \quad (2.17)$$

Следует заметить, что для одноэлементного структурного числа

$$\frac{\partial}{\partial a}[a] = 1. \quad (2.18)$$

***Пример 2.6.*** Нахождение алгебраической производной структурного числа:

$$A = \begin{vmatrix} 1 & 1 & 5 & 4 \\ 2 & 3 & 3 & 3 \\ 4 & 2 & 7 & 1 \end{vmatrix}, \quad \frac{\partial A}{\partial 1} = \begin{vmatrix} 2 & 3 & 4 \\ 4 & 2 & 3 \end{vmatrix}, \quad \frac{\partial A}{\partial 4} = \begin{vmatrix} 1 & 3 \\ 2 & 1 \end{vmatrix}.$$

По аналогии с математическим анализом нахождение производной будем называть дифференцированием.

Дифференцирование структурного числа имеет весьма простую геометрическую интерпретацию.

**Свойство 2.1.** Геометрическое изображение структурного .числа $\partial A / \partial a$ представляет собой геометрическое изображение структурного числа А с замкнутым ребром *a*.

Кроме алгебраической производной, сформулируем для структурных чисел еще одно понятие (в известном смысле дуальное по отношению к производной) - понятие обратной алгебраической производной.

Алгебраической обратной производной структурного числа называется структурное число $\delta \mathrm{A} / \delta a$, равное

$$\frac{\delta \mathrm{A}}{\delta a} = \mathrm{A} \Big| \text{столбцы, содержащие элемент } a \text{ опущены.} \tag{2.19}$$

Для обратной алгебраической производной имеют место соотношения

$$\frac{\delta}{\delta a}(\mathrm{A}_1 + \mathrm{A}_2) = \frac{\delta \mathrm{A}_1}{\delta a} + \frac{\delta \mathrm{A}_2}{\delta a},$$

$$\frac{\delta}{\delta a}(\mathrm{A}_1 \mathrm{A}_2) = \frac{\delta \mathrm{A}_1}{\delta a}\mathrm{A}_2 + \frac{\delta \mathrm{A}_2}{\delta a}\mathrm{A}_1 + \mathrm{A}_1\mathrm{A}_2, \tag{2.20}$$

справедливые для произвольных чисел $\mathrm{A}_1$ и $\mathrm{A}_2$.

Кроме того,

$$\frac{\delta}{\delta a}(\mathrm{A}_1 \mathrm{A}_2) = \frac{\delta \mathrm{A}_1}{\delta a}\frac{\delta \mathrm{A}_2}{\delta a}. \tag{2.21}$$

Для одноэлементного структурного числа имеем

$$\frac{\delta}{\delta a}[a] = 0. \tag{2.22}$$

Соотношение алгебраических производной и обратной производной можно записать следующим образом:

$$\frac{\partial \mathrm{A}}{\partial a}(\mathrm{A}[a]) = \frac{\delta \mathrm{A}}{\delta a}. \tag{2.23}$$

Алгебраическую обратную производную будем обозначать как

$$\frac{\delta \mathrm{A}}{\delta a} = \mathrm{A}^a. \tag{2.24}$$

***Пример 2.7.*** Расчет алгебраической обратной производной:

$$\mathrm{A} = \begin{vmatrix} 1 & 2 & 1 & 5 \\ 3 & 4 & 2 & 4 \\ 4 & 7 & 3 & 8 \end{vmatrix}, \quad \frac{\delta \mathrm{A}}{\delta 1} = \begin{vmatrix} 2 & 5 \\ 4 & 4 \\ 7 & 8 \end{vmatrix}, \quad \frac{\delta \mathrm{A}}{\delta 2} = \begin{vmatrix} 1 & 5 \\ 3 & 4 \\ 5 & 8 \end{vmatrix}.$$

Алгебраическая обратная производная имеет простую геометрическую интерпретацию.

**Свойство 2.2**. Геометрическое изображение структурного числа $\delta \mathrm{A} / \delta a$ представляет собой геометрическое изображение структурного числа А, в котором ребро *a* отключено в одной вершине и замкнуто в петлю.

Обратное геометрическое изображение структурного числа $\delta \mathrm{A} / \delta a$ представляет собой обратное изображение геометрического числа A с замкнутым числом $a$.

Отметим, что для структурного числа A всегда имеет место соотношение

$$\mathrm{A} = \frac{\delta \mathrm{A}}{\delta a} + [a] \frac{\partial \mathrm{A}}{\partial a}, \tag{2.25}$$

где $a$ - элемент числа A.

***Пример 2.8.*** Задача перечисления всех деревьев графа

Рассмотрим задачу перечисления всех деревьев графа G с применением методов алгебры структурных чисел. В качестве примера рассмотрим граф представленный на рис. 2.3.

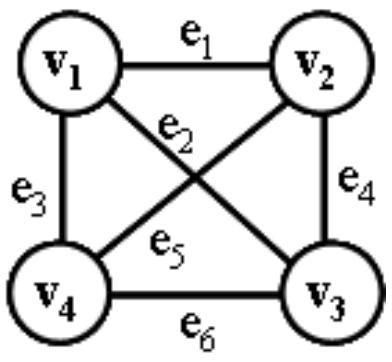

Рис. 2.3. Граф **G**

Представим множество центральных разрезов, записанных как однострочные структурные числа, имеет вид (здесь ребра представлены цифрами):

$s_1 \rightarrow [1\ 2\ 3]$;

$s_2 \rightarrow [1\ 4\ 5]$;

$s_3 \rightarrow [2\ 4\ 6]$;

$s_4 \rightarrow [3\ 5\ 6]$.

Предположим, что множество деревьев

$$\mathrm{T} = [1\ \ 2\ \ 3] \times [1\ \ 4\ \ 5] \times [2\ \ 4\ \ 6] =$$

$$= \begin{vmatrix} 1 & 1 & 2 & 2 & 2 & 3 & 3 & 3 \\ 4 & 5 & 1 & 4 & 5 & 1 & 4 & 5 \end{vmatrix} \times \begin{vmatrix} 2 & 4 & 6 \end{vmatrix} =$$

$$= \begin{vmatrix} 1 & 1 & 1 & 1 & 1 & 2 & 2 & 2 & 2 & 2 & 3 & 3 & 3 & 3 & 3 & 3 & 3 & 3 \\ 4 & 4 & 5 & 5 & 5 & 1 & 1 & 4 & 5 & 5 & 1 & 1 & 1 & 4 & 4 & 5 & 5 & 5 \\ 2 & 6 & 2 & 4 & 6 & 4 & 6 & 6 & 4 & 6 & 2 & 4 & 6 & 2 & 6 & 2 & 4 & 6 \end{vmatrix} =$$

$$= \begin{vmatrix} 1 & 1 & 1 & 1 & 2 & 2 & 2 & 2 & 3 & 3 & 3 & 3 & 3 & 3 & 3 & 3 \\ 4 & 5 & 5 & 5 & 1 & 4 & 5 & 5 & 1 & 1 & 1 & 4 & 4 & 5 & 5 & 5 \\ 6 & 2 & 4 & 6 & 6 & 6 & 4 & 6 & 2 & 4 & 6 & 2 & 6 & 2 & 4 & 6 \end{vmatrix}.$$

$$\mathrm{T}^{d} = \begin{vmatrix} 2 & 3 & 2 & 2 & 3 & 1 & 1 & 1 & 4 & 2 & 2 & 1 & 1 & 1 & 1 & 1 \\ 3 & 4 & 3 & 3 & 4 & 3 & 3 & 3 & 5 & 5 & 4 & 5 & 2 & 4 & 2 & 2 \\ 5 & 6 & 6 & 4 & 5 & 5 & 6 & 4 & 6 & 6 & 5 & 4 & 5 & 6 & 6 & 4 \end{vmatrix}.$$

если $T=\begin{bmatrix}1 & 2 & 3\end{bmatrix}\times\begin{bmatrix}1 & 4 & 5\end{bmatrix}\times\begin{bmatrix}3 & 5 & 6\end{bmatrix}=$

$$=\begin{vmatrix}1 & 1 & 2 & 2 & 2 & 3 & 3 & 3\\ 4 & 5 & 1 & 4 & 5 & 1 & 4 & 5\end{vmatrix}\times\begin{vmatrix}3 & 5 & 6\end{vmatrix}=$$

$$=\begin{vmatrix}1 & 1 & 1 & 1 & 1 & 2 & 2 & 2 & 2 & 2 & 2 & 2 & 2 & 3 & 3 & 3 & 3 & 3\\ 4 & 4 & 4 & 5 & 5 & 1 & 1 & 1 & 4 & 4 & 4 & 5 & 5 & 1 & 1 & 4 & 4 & 5\\ 3 & 5 & 6 & 3 & 6 & 3 & 5 & 6 & 3 & 5 & 6 & 3 & 6 & 5 & 6 & 5 & 6 & 6\end{vmatrix}=$$

$$=\begin{vmatrix}1 & 1 & 1 & 1 & 2 & 2 & 2 & 2 & 2 & 2 & 2 & 2 & 3 & 3 & 3 & 3\\ 4 & 4 & 4 & 5 & 1 & 1 & 1 & 4 & 4 & 4 & 5 & 5 & 1 & 4 & 4 & 5\\ 3 & 5 & 6 & 6 & 3 & 5 & 6 & 3 & 5 & 6 & 3 & 6 & 6 & 5 & 6 & 6\end{vmatrix}.$$

если $T=\begin{bmatrix}1 & 4 & 5\end{bmatrix}\times\begin{bmatrix}2 & 4 & 6\end{bmatrix}\times\begin{bmatrix}3 & 5 & 6\end{bmatrix}=$

$$=\begin{vmatrix}1 & 1 & 1 & 4 & 4 & 5 & 5 & 5\\ 2 & 4 & 6 & 2 & 6 & 2 & 4 & 6\end{vmatrix}\times\begin{vmatrix}3 & 5 & 6\end{vmatrix}=$$

$$=\begin{vmatrix}1 & 1 & 1 & 1 & 1 & 1 & 1 & 1 & 4 & 4 & 4 & 4 & 4 & 5 & 5 & 5 & 5 & 5\\ 2 & 2 & 2 & 4 & 4 & 4 & 6 & 6 & 2 & 2 & 2 & 6 & 6 & 2 & 2 & 4 & 4 & 6\\ 3 & 5 & 6 & 3 & 5 & 6 & 3 & 5 & 3 & 5 & 6 & 3 & 5 & 3 & 6 & 3 & 6 & 3\end{vmatrix}=$$

$$=\begin{vmatrix}1 & 1 & 1 & 1 & 1 & 1 & 1 & 1 & 4 & 4 & 4 & 4 & 5 & 5 & 5 & 5\\ 2 & 2 & 2 & 4 & 4 & 4 & 6 & 6 & 2 & 2 & 2 & 6 & 2 & 2 & 4 & 6\\ 3 & 5 & 6 & 3 & 5 & 6 & 3 & 5 & 3 & 5 & 6 & 3 & 3 & 6 & 3 & 3\end{vmatrix}.$$

если $T=\begin{bmatrix}1 & 2 & 3\end{bmatrix}\times\begin{bmatrix}2 & 4 & 6\end{bmatrix}\times\begin{bmatrix}3 & 5 & 6\end{bmatrix}=$

$$=\begin{vmatrix}1 & 1 & 1 & 2 & 2 & 3 & 3 & 3\\ 2 & 4 & 6 & 4 & 6 & 2 & 4 & 6\end{vmatrix}\times\begin{vmatrix}3 & 5 & 6\end{vmatrix}=$$

$$=\begin{vmatrix}1 & 1 & 1 & 1 & 1 & 1 & 1 & 1 & 2 & 2 & 2 & 2 & 2 & 3 & 3 & 3 & 3 & 3\\ 2 & 2 & 2 & 4 & 4 & 4 & 6 & 6 & 4 & 4 & 4 & 6 & 6 & 2 & 2 & 4 & 4 & 6\\ 3 & 5 & 6 & 3 & 5 & 6 & 3 & 5 & 3 & 5 & 6 & 3 & 5 & 5 & 6 & 5 & 6 & 5\end{vmatrix}=$$

$$=\begin{vmatrix}1 & 1 & 1 & 1 & 1 & 1 & 1 & 1 & 2 & 2 & 2 & 2 & 3 & 3 & 3 & 3\\ 2 & 2 & 2 & 4 & 4 & 4 & 6 & 6 & 4 & 4 & 4 & 6 & 2 & 4 & 4 & 6\\ 3 & 5 & 6 & 3 & 5 & 6 & 3 & 5 & 3 & 5 & 6 & 5 & 5 & 5 & 6 & 5\end{vmatrix}.$$

таким образом

T = [1 2 3]×[1 4 5] ×[2 4 6]=

= [1 2 3] ×[1 4 5] ×[3 5 6]=

= [1 4 5 ] ×[2 4 6] ×[3 5 6] =

= [1 2 3 ] ×[2 4 6] ×[3 5 6].

Что и требовалось показать.

Множество дополнений дерева - кодеревья (или множество хорд), можно получить, рассматривая произведения однострочных структурных чисел, где каждое однострочное структурное число представляет собой базисный изометрический цикл [].

$T^d$ = [1 2 4] ×[1 3 5] ×[4 5 6]=

Все деревья и дополнения дерева данного графа представлены на рис. 2.4 и рис. 2.5.

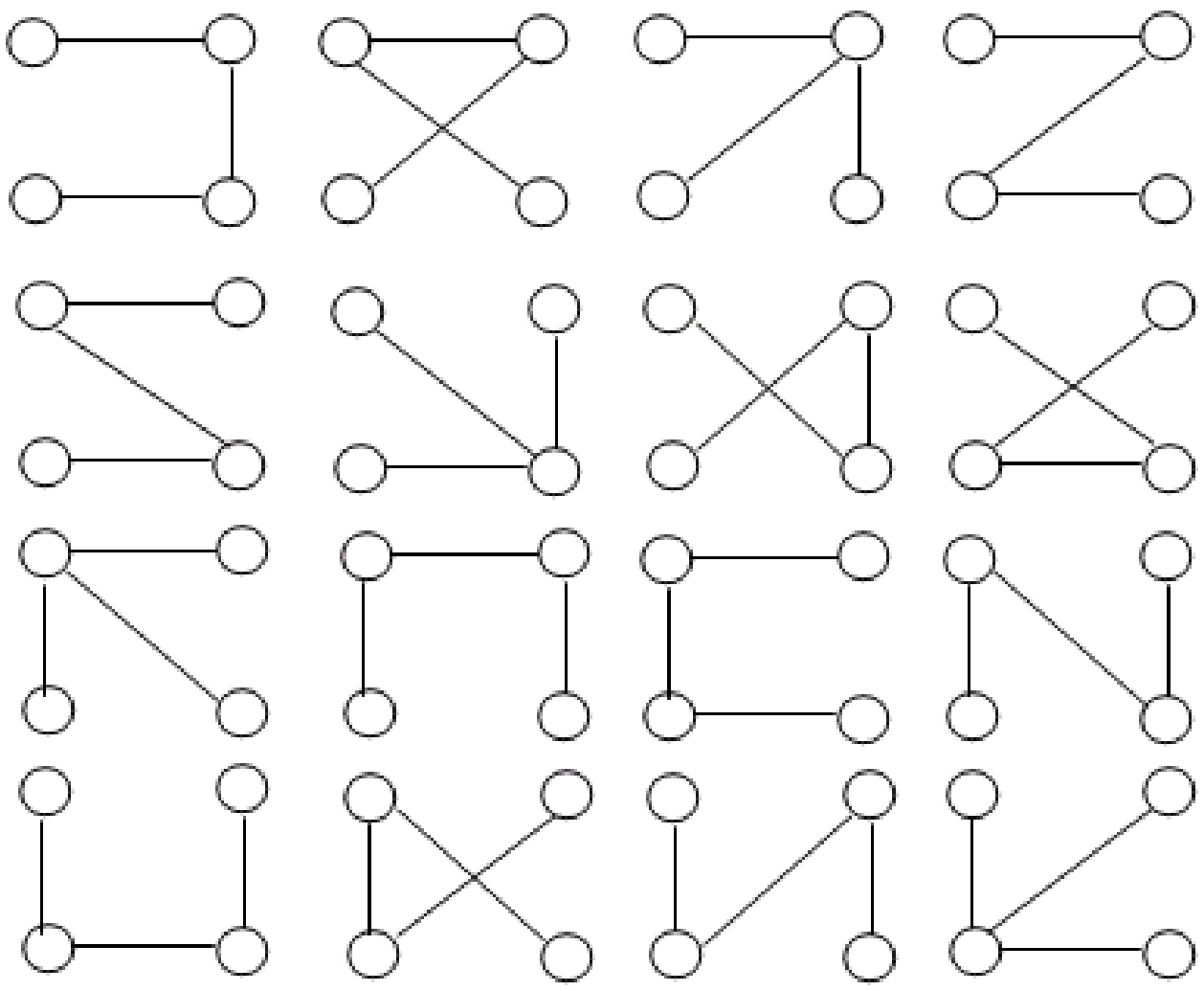

Рис. 2.4. Множество деревьев графа G.

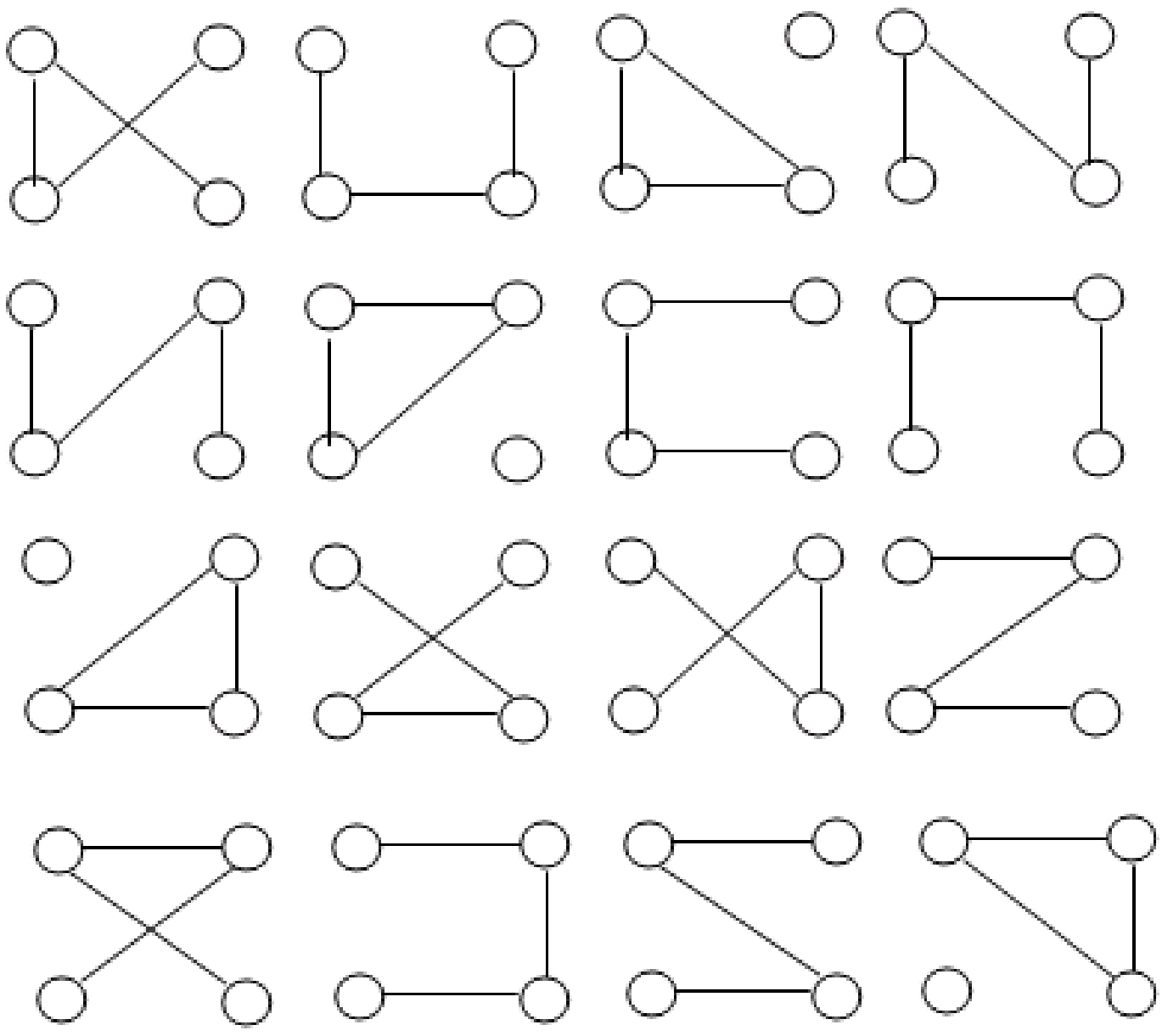

Рис. 2.5. Множество кодеревьев графа G.

## 2.5. Структурные числа и изометрические циклы

Согласно теореме 2.1, в неориентированном несепарабельном графе существуют независимые подмножества, состоящие из изометрических циклов с мощностью равной цикломатическому числу графа. Для описания таких независимых множеств состоящих из изометрических циклов можно применить методы алгебры структурных чисел. Множество таких подмножеств будем обозначать буквой $C_\tau$. С другой стороны, элемент такого множества, состоящего из изометрических циклов, можно представлять и описывать как элемент структурного числа $C_\tau$ [40].

Каждое подмножество, состоящее из изометрических циклов (в том числе и базисы), можно описывать характеристической функцией, называемой вектором количества циклов, проходящих по ребрам

$$P_e = \langle a_1,a_2,a_3,\ldots,a_m\rangle, \qquad (2.26)$$

где $a_i$ – количество циклов проходящих по ребру i.

Или характеристической функцией, называемой вектором количества циклов, проходящих по вершинам

$$P_v = \langle b_1,b_2,b_3,\ldots,b_n\rangle. \qquad (2.27)$$

где $b_j$ – количество циклов проходящих по вершине j.

Элемент множества базисов, состоящий из изометрических циклов, можно представлять и описывать как элемент $c \in C_\tau$ структурного числа $C_\tau$ [40].

Структурное число $C_\tau$ определим как произведение однострочных структурных чисел, построенных на хордах выделенного дерева T как подмножество изометрических циклов приходящих по данной хорде.

В качестве примера, рассмотрим следующий граф (см. рис. 2.6).

***Пример 2.9.*** Запишем для данного графа всё множество квазициклов:

$c_0 = \varnothing$, $c_1 = \{e_1,e_2,e_3\}$, $c_2 = \{e_1,e_2,e_4,e_6,e_7\}$, $c_3 = \{e_3,e_4,e_6,e_7\}$,
$c_4 = \{e_5,e_6,e_7\}$, $c_5 = \{e_1,e_2,e_3,e_5,e_6,e_7\}$, $c_6 = \{e_1,e_2,e_4,e_5\}$,
$c_7 = \{e_3,e_4,e_5\}$.

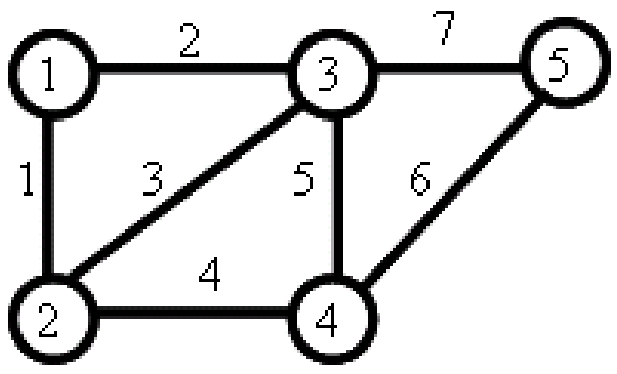

Рис. 2.6. Граф **G**.

Мощность этого подпространства циклов определяется как:

$\nu(G) = k = m - n + 1 = 7 - 5 + 1 = 3$,

$\text{card}\,\mathbf{L}_G^C = 2^{\nu(G)} = 2^3 = 8$.

Для данного графа выделим дерево, состоящее из ребер графа $\mathbf{T} = \{e_1,e_2,e_6,e_7\}$. Тогда хорды образуют множество хорд $\mathbf{H} = \{e_3,e_4,e_5\}$. Запишем в виде однострочных структурных чисел множество номеров квазициклов проходящих по хордам графа.

По хорде $e_3$ проходят циклы $c_1,c_3,c_5,c_7$. Соответственно записывается однострочное структурное число (здесь циклы обозначаются цифрами):

$e_3$: [1 3 5 7].

По хорде $e_4$ проходят циклы $c_2, c_3, c_6, c_7$:

$e_4$: [2 3 6 7].

И наконец, по хорде $e_5$ проходят циклы $c_4, c_5, c_6, c_7$:

$e_5$: [4 5 6 7].

Произведение однострочных структурных чисел характеризует множество следующих базисов из изометрических циклов:

C = [1 3 5 7]×[2 3 6 7]×[4 5 6 7] =

$$= \begin{matrix} 1 & 1 & 1 & 1 & 3 & 3 & 3 & 5 & 5 & 5 & 5 & 7 & 7 & 7 \\ 2 & 3 & 6 & 7 & 2 & 6 & 7 & 2 & 3 & 6 & 7 & 2 & 3 & 6 \end{matrix} \times \begin{matrix} 4 & 5 & 6 & 7 \end{matrix} =$$

=

| | | | | | | | | | | | | | | | | | | | | | | | |
|---|---|---|---|---|---|---|---|---|---|---|---|---|---|---|---|---|---|---|---|---|---|---|---|
| 1 | 1 | 1 | 1 | 1 | 1 | 1 | 1 | 1 | 1 | **1** | 1 | 1 | **1** | 3 | 3 | 3 | 3 | 3 | **3** | **3** | **3** | **3** | **3** |
| 2 | 2 | 2 | 2 | 3 | 3 | 3 | 3 | 6 | 6 | **6** | 7 | 7 | **7** | 2 | 2 | 2 | 2 | 6 | **6** | **6** | **7** | **7** | **7** |
| 4 | 5 | 6 | 7 | 4 | 5 | 6 | 7 | 4 | 5 | **7** | 4 | 5 | **6** | 4 | 5 | 6 | 7 | 4 | **5** | **7** | **4** | **5** | **6** |

| | | | | | | | | | | | | | | | | | |
|---|---|---|---|---|---|---|---|---|---|---|---|---|---|---|---|---|---|
| 5 | 5 | **5** | 5 | **5** | **5** | 5 | **5** | 5 | **5** | 7 | **7** | 7 | **7** | 7 | 7 | 7 | 7 |
| 2 | 2 | **2** | 3 | **3** | **3** | 6 | **6** | 7 | **7** | 2 | **2** | 2 | **3** | 3 | 3 | 6 | 6 |
| 4 | 6 | **7** | 4 | **6** | **7** | 4 | **7** | 4 | **6** | 4 | **5** | 6 | **4** | 5 | 6 | 4 | 5 |

Количество базисов подпространства $L_G^C$ можно определять по формуле:

$$\text{card } L_G^C = (1/k!)(2^k-2^0)(2^k-2^1)...(2^k-2^{k-1}) = 28. \quad (2.28)$$

Если удалить зависимые столбцы (они выделены серым цветом), то получим следующее структурное число:

=

| | | | | | | | | | | | | | | | | | | | | | | | |
|---|---|---|---|---|---|---|---|---|---|---|---|---|---|---|---|---|---|---|---|---|---|---|---|
| 1 | 1 | 1 | 1 | 1 | 1 | 1 | 1 | 1 | 1 | 1 | 1 | 3 | 3 | 3 | 3 | 3 | 5 | 5 | 5 | 5 | 5 | 7 | 7 |
| 2 | 2 | 2 | 2 | 3 | 3 | 3 | 3 | 6 | 6 | 7 | 7 | 2 | 2 | 2 | 2 | 6 | 2 | 2 | 3 | 6 | 7 | 2 | 2 |
| 4 | 5 | 6 | 7 | 4 | 5 | 6 | 7 | 4 | 5 | 4 | 5 | 4 | 5 | 6 | 7 | 4 | 4 | 6 | 4 | 4 | 4 | 4 | 6 |

| | | | |
|---|---|---|---|
| 7 | 7 | 7 | 7 |
| 3 | 3 | 6 | 6 |
| 5 | 6 | 4 | 5 |

Построим структурное число $C_n$ с учетом только простых циклов:

$C_n$ = [1 3 7]×[2 3 6 7]×[4 6 7] =

$$= \begin{matrix} 1 & 1 & 1 & 1 & 3 & 3 & 3 & 7 & 7 & 7 \\ 2 & 3 & 6 & 7 & 2 & 6 & 7 & 2 & 3 & 6 \end{matrix} \times \begin{matrix} 4 & 6 & 7 \end{matrix} =$$

$$= \begin{vmatrix} 1 & 1 & 1 & 1 & 1 & 1 & 1 & 1 & 1 & 1 & 3 & 3 & 3 & 3 & 3 & 3 & 3 & 7 & 7 & 7 & 7 & 7 \\ 2 & 2 & 2 & 3 & 3 & 3 & 6 & 6 & 7 & 7 & 2 & 2 & 2 & 6 & 6 & 7 & 7 & 2 & 2 & 3 & 3 & 6 \\ 4 & 6 & 7 & 4 & 6 & 7 & 4 & 7 & 4 & 6 & 4 & 6 & 7 & 4 & 7 & 4 & 6 & 4 & 6 & 4 & 6 & 4 \end{vmatrix}$$

$$= \begin{vmatrix} 1 & 1 & 1 & 1 & 1 & 1 & 1 & 1 & 3 & 3 & 3 & 3 & 3 & 7 & 7 & 7 \\ 2 & 2 & 2 & 3 & 3 & 3 & 6 & 7 & 2 & 2 & 2 & 6 & 6 & 2 & 2 & 6 \\ 4 & 6 & 7 & 4 & 6 & 7 & 4 & 4 & 4 & 6 & 7 & 4 & 7 & 4 & 6 & 4 \end{vmatrix}$$

Построим структурное число $C_e$ с учетом только изометрических циклов:

$$C_e = [1\ \ 7] \times [7] \times [4\ \ 7] \qquad = \begin{vmatrix} 1 \\ 7 \\ 4 \end{vmatrix}$$

Как видим, количество столбцов структурного числа $C$, состоящего из всех квазициклов графа, соответствует количеству базисов определяемых формулой (2.28). Для данного графа базис, состоящий только из изометрических циклов, единственный (см. рис. 2.7).

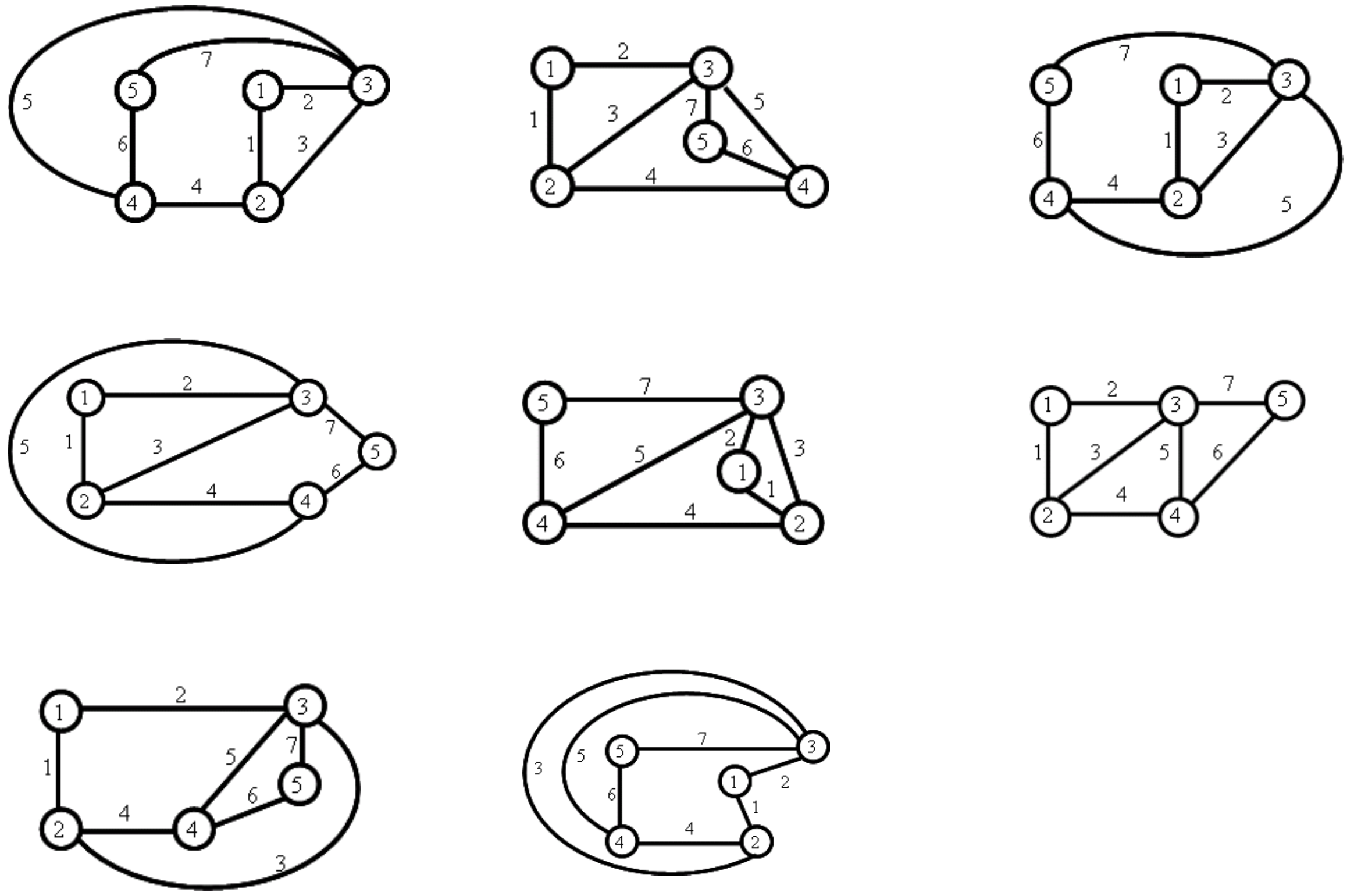

Рис. 2.7. Плоские рисунки графа из простых циклов.

На рисунке 2.7 представлены рисунки графа G. Множество рисунков возникает в результате того, что сумма базисных изометрических циклов определяет обод графа, который можно охарактеризовать как простой цикл. В свою очередь обод также можно рассматривать как границу грани в плоском графе.

С другой стороны, построение структурного числа без цикла $c_5$ можно представить как применение операции взятия обратной производной структурного числа $C$, то есть структурное число $C$ без квазицикла $c_5$.

$$\frac{\delta C}{\delta c_5} = \begin{vmatrix} 1 & 1 & 1 & 1 & 1 & 1 & 1 & 1 & 3 & 3 & 3 & 3 & 3 & 7 & 7 & 7 \\ 2 & 2 & 2 & 3 & 3 & 3 & 6 & 7 & 2 & 2 & 2 & 6 & 6 & 2 & 2 & 6 \\ 4 & 6 & 7 & 4 & 6 & 7 & 4 & 4 & 4 & 6 & 7 & 4 & 7 & 4 & 6 & 4 \end{vmatrix}$$

Структурное числа без циклов $c_2,c_3,c_5,c_6$ можно представить как последовательное применение операции взятия обратной производной структурного числа $C$, то есть структурное число без квазицикла $c_5$ и простых циклов $c_2,c_3,c_6$ запишется в виде:

$$\frac{\delta C}{\delta c_5 \delta c_2 \delta c_3 \delta c_6} = \begin{vmatrix} 1 \\ 7 \\ 4 \end{vmatrix}.$$

Выбранный базис может быть проверен на независимость методами алгебры структурных чисел. Рассмотрим граф $K_5$ представленный на рис. 2.9.

Множество изометрических циклов графа:

$c_1 = \{e_1,e_2,e_5\}$, $c_2 = \{e_1,e_3,e_6\}$, $c_3 = \{e_1,e_4,e_7\}$, $c_4 = \{e_2,e_3,e_8\}$, $c_5 = \{e_2,e_4,e_9\}$,

$c_6 = \{e_3,e_4,e_{10}\}$, $c_7 = \{e_5,e_6,e_8\}$, $c_8 = \{e_5,e_7,e_9\}$, $c_9 = \{e_6,e_7,e_{10}\}$, $c_{10} = \{e_8,e_9,e_{10}\}$.

Выделим дерево графа $T = \{e_1,e_2,e_6,e_9\}$. Тогда, относительно хорд, можно построить однострочные структурные числа, характеризующие циклы, проходящие по ребрам (числа, записанные в однострочном структурном числе, характеризуют номера циклов):

$e_3$: [ $\underline{\mathbf{2}}$ 4 6 ];
$e_4$: [ 3 $\underline{\mathbf{5}}$ 6 ];
$e_5$: [ $\underline{\mathbf{1}}$ 7 8 ];
$e_7$: [ 3 $\underline{\mathbf{8}}$ 9 ];
$e_8$: [ 4 7 $\underline{\mathbf{10}}$ ];
$e_{10}$: [ 6 $\underline{\mathbf{9}}$ 10 ].

Предположим, что выделили следующий базис, состоящий из циклов $b_1 = \{c_1,c_2,c_5,c_8,c_9,c_{10}\}$. Данная система линейно независима, так как имеется нечетное число структурных чисел для выбранной системы циклов.

Заменим цикл $c_1$ проходящий по хорде $e_5$ на цикл $c_7$ проходящий по той же хорде, получим следующую систему циклов $b_2 = \{c_2,c_5,c_7,c_8,c_9,c_{10}\}$.

$e_3$: [ $\underline{\mathbf{2}}$ 4 6 ];
$e_4$: [ 3 $\underline{\mathbf{5}}$ 6 ];
$e_5$: [ 1 $\underline{\mathbf{7}}$ 8 ];
$e_7$: [ 3 $\underline{\mathbf{8}}$ 9 ];

$e_8$: [ 4 7 **10** ];
$e_{10}$: [ 6 **9** 10 ].

Данная система линейно зависима, так как:

$$c_7 \oplus c_8 \oplus c_9 \oplus c_{10} = \{e_5,e_6,e_8\} \oplus \{e_5,e_7,e_9\} \oplus \{e_6,e_7,e_{10}\} \oplus \{e_8,e_9,e_{10}\} = \varnothing.$$

С другой стороны, данная система линейно зависима, так как имеется четное число структурных элементов для выбранной системы циклов [37,77].

| | | |
|---|---|---|
| $e_3$: | 2 | 2 |
| $e_4$: | 5 | 5 |
| $e_5$: | 7 | 8 |
| $e_7$: | 8 | 9 |
| $e_8$: | 10 | 7 |
| $e_{10}$: | 9 | 10 |

Поэтому, при замене циклов в однострочном структурном числе, необходимо применять операцию проверки системы циклов на независимость.

## 2.6 Линейная независимость множества изометрических циклов графа

Прежде чем перейти к рассмотрению алгоритмов проверки графа на планарность и выделения максимально плоского суграфа, мы должны найти возможность определения ранга матрицы, состоящей из множества изометрических циклов. Ранг этой матрицы должен быть равен цикломатическому числу. Для определения ранга обычно используют алгоритм Гаусса [50]. Алгоритм Гаусса служит для определения ранга матрицы, однако небольшая его модификация позволяет определять зависимые системы, состоящие из суграфов, для последующего их исключения из процесса расчета. Приведем соответствующую схему модификации алгоритма Гаусса для определения ранга матрицы и одновременного нахождения подмножества зависимых суграфов.

Вектор $\beta$ называется линейной комбинацией векторов $a_1,a_2, ...,a_s$, если существуют такие числа $l_1,l_2,...,l_s$, что $\beta = l_1a_1+l_2a_2+...+l_sa_s$.

Система векторов $a_1,a_2,...,a_{r-1},a_r$ $(r\geq2)$, (2.29)

называется линейно зависимой, если хотя бы один из этих векторов является линейной комбинацией остальных векторов системы (2.29), и линейно независимой - в противном случае. Такие алгебраические понятия как линейная зависимость векторов, ранг матрицы, обращение матрицы и т.д., справедливы в любом поле.

**Инициализация**. Выбираем дерево в графе. Помечаем все хорды графа.

***Шаг 1.*** Все множество изометрических циклов графа G записываем в виде (0,1)-матрицы. Идем на шаг 2.

***Шаг 2.*** Выделяем первый непомеченный изометрический цикл, который содержит единицу в столбце соответствующем хорде, и объявляем его главным элементом. Помечаем выбранный цикл знаком ⊗. Если количество таких помеченных изометрических циклов равно цикломатическому числу графа, то конец работы алгоритма. Если нет, идем на шаг 3.

***Шаг 3.*** Выбираем непомеченные строки в столбцах, которых имеются единицы, совпадающие с главным элементом выбранной строки. Помечаем эти строки знаком & и приписываем номер цикла содержащего главный элемент. Содержимое таких строк образуется, как результат суммирования по модулю два предыдущего содержания строки и строки содержащей главный элемент. Идем на шаг 2.

Таким образом, применение алгоритма Гаусса для выделения базиса циклов к (0,1)-матрице циклов, имеющий ранг, равный цикломатическому числу графа, позволяет строить диаграмму комбинаций циклов. Ненулевые строки в диаграмме комбинаций циклов соответствуют линейно независимой системе циклов.

***Пример 2.10.*** Сказанное рассмотрим на примере графа $\mathbf{K_5}$, представленного на рис. 2.9.

Множество изометрических циклов состоит из следующих элементов:

$c_1 = \{e_1,e_2,e_5\}$; $c_2 = \{e_1,e_3,e_6\}$; $c_3 = \{e_1,e_4,e_7\}$; $c_4 = \{e_2,e_3,e_8\}$; $c_5 = \{e_2,e_4,e_9\}$;
$c_6 = \{e_3,e_4,e_{10}\}$; $c_7 = \{e_5,e_6,e_8\}$; $c_8 = \{e_5,e_7,e_9\}$; $c_9 = \{e_6,e_7,e_{10}\}$; $c_{10} = \{e_8,e_9,e_{10}\}$.

Выбираем хорды графа $H = \{e_1,e_3,e_4,e_6,e_7,e_8\}$ относительно дерева $T = \{e_2,e_5,e_9,e_{10}\}$. Составляем матрицу изометрических циклов графа и помечаем хорды графа.

| | $e_1$ | $e_2$ | $e_3$ | $e_4$ | $e_5$ | $e_6$ | $e_7$ | $e_8$ | $e_9$ | $e_{10}$ |
|---|---|---|---|---|---|---|---|---|---|---|
| $c_1$ | 1 | 1 | | | 1 | | | | | |
| $c_2$ | 1 | | 1 | | | 1 | | | | |
| $c_3$ | 1 | | | 1 | | | 1 | | | |
| $c_4$ | | 1 | 1 | | | | | 1 | | |
| $c_5$ | | 1 | | 1 | | | | | 1 | |
| $c_6$ | | | 1 | 1 | | | | | | 1 |
| $c_7$ | | | | | 1 | 1 | | 1 | | |
| $c_8$ | | | | | 1 | | 1 | | 1 | |
| $c_9$ | | | | | | 1 | 1 | | | 1 |
| $c_{10}$ | | | | | | | | 1 | 1 | 1 |

Выбираем первую строку с циклом $c_1$, и в качестве главного элемента выбираем ребро $u_1$ как принадлежащее хорде графа. Выбранная строка пересекается со строками $c_2$ и $c_3$ ребром $u_1$. Формируем новое содержание матрицы изометрических циклов записывая в строки содержащие цикл $c_1$ результат сложения по модулю 2 данных строк с 1-ой строкой.

| | $e_1$ | $e_2$ | $e_3$ | $e_4$ | $e_5$ | $e_6$ | $e_7$ | $e_8$ | $e_9$ | $e_{10}$ | |
|---|---|---|---|---|---|---|---|---|---|---|---|
| $c_1$ | 1 | 1 | | | 1 | | | | | | ⊗ |
| $c_2$ | | 1 | 1 | | 1 | 1 | | | | | & $c_1$ |
| $c_3$ | | 1 | | 1 | 1 | | 1 | | | | & $c_1$ |
| $c_4$ | | 1 | 1 | | | | | 1 | | | |
| $c_5$ | | 1 | | 1 | | | | | 1 | | |
| $c_6$ | | | 1 | 1 | | | | | | 1 | |
| $c_7$ | | | | | 1 | 1 | | 1 | | | |

| | | | | | | | | | | |
|---|---|---|---|---|---|---|---|---|---|---|
| $c_8$ | | | | | 1 | | 1 | | 1 | |
| $c_9$ | | | | | | 1 | 1 | | | 1 |
| $c_{10}$ | | | | | | | | 1 | 1 | 1 |

Выбираем следующую непомеченную индексами & и ⊗ строку с циклом $c_4$. В качестве главного элемента выбираем хорду $e_3$. Выбранная строка пересекается со строками $c_2$ и $c_6$. Помечаем строку с циклом $c_6$ индексом &. Формируем новое содержимое строк $c_2$ и $c_6$.

| | $e_1$ | $e_2$ | $e_3$ | $e_4$ | $e_5$ | $e_6$ | $e_7$ | $e_8$ | $e_9$ | $e_{10}$ | |
|---|---|---|---|---|---|---|---|---|---|---|---|
| $c_1$ | **1** | 1 | | | 1 | | | | | | ⊗ |
| $c_2$ | | | | | 1 | 1 | | 1 | | | & $c_1$ $c_4$ |
| $c_3$ | | 1 | | 1 | 1 | | 1 | | | | & $c_1$ |
| $c_4$ | | 1 | **1** | | | | | 1 | | | ⊗ |
| $c_5$ | | 1 | | 1 | | | | | 1 | | |
| $c_6$ | | 1 | | 1 | | | | 1 | | 1 | & $c_4$ |
| $c_7$ | | | | | 1 | 1 | | 1 | | | |
| $c_8$ | | | | | 1 | | 1 | | 1 | | |
| $c_9$ | | | | | | 1 | 1 | | | 1 | |
| $c_{10}$ | | | | | | | | 1 | 1 | 1 | |

Выбираем следующую непомеченную индексами & и ⊗ строку с циклом $c_5$. В качестве главного элемента выбираем хорду $e_4$. Выбранная строка пересекается со строками $c_3$ и $c_6$. Формируем новое содержание строк $c_2$ и $c_6$.

| | $e_1$ | $e_2$ | $e_3$ | $e_4$ | $e_5$ | $e_6$ | $e_7$ | $e_8$ | $e_9$ | $e_{10}$ | |
|---|---|---|---|---|---|---|---|---|---|---|---|
| $c_1$ | **1** | 1 | | | 1 | | | | | | ⊗ |
| $c_2$ | | | | | 1 | 1 | | 1 | | | & $c_1$ $c_4$ |
| $c_3$ | | | | | 1 | | 1 | | 1 | | & $c_1$ $c_5$ |
| $c_4$ | | 1 | **1** | | | | | 1 | | | ⊗ |
| $c_5$ | | 1 | | **1** | | | | | 1 | | ⊗ |
| $c_6$ | | | | | | | | 1 | 1 | 1 | & $c_4$ $c_5$ |
| $c_7$ | | | | | 1 | 1 | | 1 | | | |
| $c_8$ | | | | | 1 | | 1 | | 1 | | |
| $c_9$ | | | | | | 1 | 1 | | | 1 | |
| $c_{10}$ | | | | | | | | 1 | 1 | 1 | |

Выбираем следующую непомеченную индексами & и ⊗ строку с циклом $c_7$. В качестве главного элемента выбираем хорду $e_6$. Выбранная строка пересекается со строками $c_2$ и $c_9$. Помечаем строку с циклом $c_9$ индексом &. Формируем новое содержание строк $c_2$ и $c_9$.

| | $e_1$ | $e_2$ | $e_3$ | $e_4$ | $e_5$ | $e_6$ | $e_7$ | $e_8$ | $e_9$ | $e_{10}$ | |
|---|---|---|---|---|---|---|---|---|---|---|---|
| $c_1$ | **1** | 1 | | | 1 | | | | | | ⊗ |
| $c_2$ | | | | | | | | | | | & $c_1$ $c_4$ $c_7$ ∅ |
| $c_3$ | | | | | 1 | | 1 | | 1 | | & $c_1$ $c_5$ |
| $c_4$ | | 1 | **1** | | | | | 1 | | | ⊗ |
| $c_5$ | | 1 | | **1** | | | | | 1 | | ⊗ |
| $c_6$ | | | | | | | | 1 | 1 | 1 | & $c_4$ $c_5$ |
| $c_7$ | | | | | 1 | **1** | | 1 | | | ⊗ |
| $c_8$ | | | | | 1 | | 1 | | 1 | | |
| $c_9$ | | | | | 1 | | 1 | 1 | | 1 | & $c_7$ |
| $c_{10}$ | | | | | | | | 1 | 1 | 1 | |

Выбираем следующую непомеченную индексами & и ⊗ строку с циклом $c_8$. В качестве главного элемента выбирпем хорду $e_7$. Выбранная строка пересекается со строками $c_3$ и $c_9$. Формируем новое содержание строк $c_3$ и $c_9$.

| | $e_1$ | $e_2$ | $e_3$ | $e_4$ | $e_5$ | $e_6$ | $e_7$ | $e_8$ | $e_9$ | $e_{10}$ | |
|---|---|---|---|---|---|---|---|---|---|---|---|
| $c_1$ | **1** | 1 | | | 1 | | | | | | ⊗ |
| $c_2$ | | | | | | | | | | | & $c_1$ $c_4$ $c_7$ ∅ |
| $c_3$ | | | | | | | | | | | & $c_1$ $c_5$ $c_8$ ∅ |
| $c_4$ | | 1 | **1** | | | | | 1 | | | ⊗ |
| $c_5$ | | 1 | | **1** | | | | | 1 | | ⊗ |
| $c_6$ | | | | | | | | 1 | 1 | 1 | & $c_4$ $c_5$ |
| $c_7$ | | | | | 1 | **1** | | 1 | | | ⊗ |
| $c_8$ | | | | | 1 | | **1** | | 1 | | ⊗ |
| $c_9$ | | | | | | | | 1 | 1 | 1 | & $c_7$ $c_8$ |
| $c_{10}$ | | | | | | | | 1 | 1 | 1 | |

Выбираем следующую непомеченную индексами & и ⊗ строку с циклом $c_{10}$. В качестве главного элемента выбираем хорду $u_8$. Выбранная строка пересекается со строками $c_6$ и $c_9$. Формируем новое содержание строк $c_6$ и $c_9$.

| | $e_1$ | $e_2$ | $e_3$ | $e_4$ | $e_5$ | $e_6$ | $e_7$ | $e_8$ | $e_9$ | $e_{10}$ | |
|---|---|---|---|---|---|---|---|---|---|---|---|
| $c_1$ | **1** | 1 | | | 1 | | | | | | ⊗ |
| $c_2$ | | | | | | | | | | | & $c_1$ $c_4$ $c_7$ ∅ |
| $c_3$ | | | | | | | | | | | & $c_1$ $c_5$ $c_8$ ∅ |
| $c_4$ | | 1 | **1** | | | | | 1 | | | ⊗ |
| $c_5$ | | 1 | | **1** | | | | | 1 | | ⊗ |
| $c_6$ | | | | | | | | | | | & $c_4$ $c_5$ $c_{10}$ ∅ |
| $c_7$ | | | | | 1 | **1** | | 1 | | | ⊗ |
| $c_8$ | | | | | 1 | | **1** | | 1 | | ⊗ |
| $c_9$ | | | | | | | | | | | & $c_7$ $c_8$ $c_{10}$ ∅ |
| $c_{10}$ | | | | | | | | **1** | 1 | 1 | ⊗ |

Таким образом, множество независимых изометрических циклов имеет вид:

| | $e_1$ | $e_2$ | $e_3$ | $e_4$ | $e_5$ | $e_6$ | $e_7$ | $e_8$ | $e_9$ | $e_{10}$ |
|---|---|---|---|---|---|---|---|---|---|---|
| $c_1$ | **1** | 1 | | | 1 | | | | | |
| $c_4$ | | 1 | **1** | | | | | 1 | | |
| $c_5$ | | 1 | | **1** | | | | | 1 | |
| $c_7$ | | | | | 1 | **1** | | 1 | | |
| $c_8$ | | | | | 1 | | **1** | | 1 | |
| $c_{10}$ | | | | | | | | **1** | 1 | 1 |

Естественно, что выбор подмножества независимых циклов из множества изометрических циклов существенно зависит от выбора дерева графа.

Следует заметить, что данный алгоритм можно использовать в задаче проверки базиса на независимость.

## 2.7. Матроиды и структурные числа

Матроиды были введены Уитни [2] в контексте исследований абстрактной теории линейной зависимости. Существует много эквивалентных определений матроида. Для нас наиболее удобным будет следующее определение [54]:

**Определение 2.8.** ***Матроид*** M – это конечное множество S и набор F таких подмножеств множества S, что выполняются следующие условия, называемые условиями независимости:

$$\varnothing \in \mathrm{F} \qquad (2.30)$$

$$\text{Если } \mathrm{Z} \in \mathrm{F} \text{ и } \mathrm{Y} \subseteq \mathrm{Z}, \text{ то } \mathrm{Y} \in \mathrm{F} \qquad (2.31)$$

$$\text{Если Z и Y – члены F□ и } |\mathrm{Z}| = |\mathrm{Y}| + 1, \text{ то}$$

$$\text{существует такое } z \in \mathrm{Z} - \mathrm{Y}, \text{ что } \mathrm{Y} \cup z \in \mathrm{F} \qquad (2.32)$$

В выражении (2.32) вычитание Z – Y понимается как Z /(Z $\cap$ Y). Элементы множества S называются элементами матроида M. Члены набора F, называются множествами матроида M. Максимальное, по включению независимое множество матроида M, называется базой матроида M. Множество баз матроида M, обозначается B(M) или просто B.

Подмножество S, не принадлежащее набору F, называется зависимым. Минимальное по включению зависимое подмножество S, называется циклом матроида M. Набор циклов матроида M, обозначается $\zeta$ (M) или просто $\zeta$ .

Функция ранга $\rho$ матроида M связывает со всяким подмножеством A $\subseteq$ S неотрицательное целое число, определяемое следующим образом:

$\rho$(A) = $max\{ |\mathrm{Z}| : \mathrm{Z} \subseteq \mathrm{A}, \mathrm{Z} \in \mathrm{F}\}$;

где $\rho$(A) – ранг подмножества A. Ранг матроида M, обозначается $\rho$(M) – это ранг множества S.

Пусть S - конечное подмножество векторного пространства. Семейство всех подмножеств линейно независимых векторов в подмножестве S удовлетворяет аксиомам независимости (2.30-2.32). Следовательно, эти подмножества S образуют набор независимых множеств матроида на подмножестве S. Ранг подмножества Z $\subseteq$ S в этом матроиде равен размерности векторного пространства, порожденного Z.

Пусть неориентированный граф с множеством ребер E. Определим на E два матроида. Сначала рассмотрим набор всех F всех подмножеств E, не содержащих циклов. Очевидно, что F удовлетворяет аксиомам (2.30-2.31). Нетрудно показать, что F удовлетворяет и (2.28). Таким образом, F является набором независимых множеств матроида M на E. Каждая база матроида M – это дерево графа G. Ранг любого подмножества Z $\subseteq$ E в этом матроиде равен рангу суграфа графа G, порожденного Z. Кроме того, каждый цикл матроида M является циклом графа G. По этой причине M называется циклическим матроидом.

Рассмотрим теперь семейство $F^*$ всех подмножеств E, не содержащих разрезающих множеств графа G.

Можно показать, что $F^*$ удовлетворяет всем аксиомам независимости (2.26-2.28) и поэтому является семейством независимых множеств матроида $M^*$ на E. Любая база матроида $M^*$ - это кодерево графа G. В этом матроиде ранг подмножества $Z \subseteq E$ равен цикломатическому числу графа G. Кроме того, каждый цикл матроида $M^*$ является разрезающим множеством графа G. Матроид $M^*$ называется матроидом разрезов или матроидом связей графа G.

Определенные таким образом матроиды M и $M^*$, обладают интересным свойством, заключающимся в том, что базы одного из них являются дополнениями в E баз другого.

Здесь следует отметить, что множество баз матроида B(M), точно совпадает со структурным числом T, порожденным умножением однострочных структурных чисел, каждое из которых характеризует независимый единичный разрез. А множество баз матроида $B(M^*)$ точно совпадает со структурным числом $T^d$ порожденного умножением однострочных структурных чисел, каждое из которых характеризует независимый простой цикл (см. пример 2.5).

## 2.8. Циклический матроид и матроид разрезов

***Пример 2.11.*** Рассмотрим следующий граф **G**.

Рассмотрим однострочные структурные числа, характеризующие центральные разрезы графа G, для удобства записи ребра будем обозначать цифрами (см. рис. 2.8).

$s_1$ = [1 2], $s_2$ = [2 3 5], $s_3$ = [4 5], $s_4$ = [1 3 4].

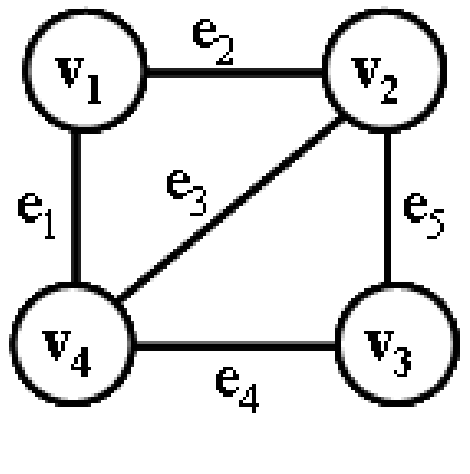

Рис. 2.8. Граф G.

Перемножим однострочные структурные числа $s_2 \times s_3 \times s_4$ по правилам умножения алгебра структурных чисел

| | | | | | | | | | | | | |
|---|---|---|---|---|---|---|---|---|---|---|---|---|
| [2 3 5] | | 2 | 2 | 2 | 2 | 2 | 3 | 3 | 3 | 5 | 5 | |
| × [4 5] | = | 4 | 4 | 5 | 5 | 5 | 4 | 5 | 5 | 4 | 4 | = |
| [1 3 4] | | 1 | 3 | 1 | 3 | 4 | 1 | 1 | 4 | 1 | 3 | |

Здесь 8 и 10 столбцы повторяются четное число раз (отмечено серым цветом) и представляют собой цикл матроида {3,4,5}, и поэтому подлежат удалению.

| | | | | | | | | | |
|---|---|---|---|---|---|---|---|---|---|
| | | 2 | 2 | 2 | 2 | 2 | 3 | 3 | 5 |
| **T** | = | 4 | 4 | 5 | 5 | 5 | 4 | 5 | 4 |
| | | 1 | 3 | 1 | 3 | 4 | 1 | 1 | 1 |

В результате умножения однострочных структурных чисел $x_2 \times x_3 \times x_4$ по правилам умножения алгебра структурных чисел, мы получили структурное число **T** характеризующее базы циклического матроида графа G. Легко показать, что данное структурное число **T** удовлетворяет условиям (2.30 - 2.32) аксиомам независимости. Для проверки условия (2.32) рассмотрим всевозможные сочетания 5-ти элементов по 2.

| | | | | | | | | | |
|---|---|---|---|---|---|---|---|---|---|
| 1 | 1 | 1 | 1 | 2 | 2 | 2 | 3 | 3 | 4 |
| 2 | 3 | 4 | 5 | 3 | 4 | 5 | 4 | 5 | 5 |

Легко убедиться, что все сочетания удовлетворяют условию (2.31). Проверим, удовлетворяется условие (2.32) или нет.

| | | | | |
|---|---|---|---|---|
| 1: | Z = {1,2,4}; | Y = {1,2}; | z = Z – Y = {4}; | z ∪ Y = {1,2,4}∈ B(M); |
| 2: | Z = {1,2,4}; | Y = {1,3}; | z = Z – Y = {2,4}; | z ∪ Y = {1,3,4}∈ B(M); |
| 3: | Z = {1,2,4}; | Y = {1,4}; | z = Z – Y = {2}; | z ∪ Y = {1,2,4}∈ B(M); |
| 4: | Z = {1,2,4}; | Y = {1,5}; | z = Z – Y = {2,4}; | z ∪ Y = {1,4,5}∈ B(M); |
| 5: | Z = {1,2,4}; | Y = {2,3}; | z = Z – Y = {1,4}; | z ∪ Y = {2,3,4}∈ B(M); |
| 6: | Z = {1,2,4}; | Y = {2,4}; | z = Z – Y = {1}; | z ∪ Y = {1,2,4}∈ B(M); |
| 7: | Z = {1,2,4}; | Y = {2,5}; | z = Z – Y = {1,4}; | z ∪ Y = ({1,2,5}∨{2,4,5})∈ B(M); |
| 8: | Z = {1,2,4}; | Y = {3,4}; | z = Z – Y = {1,2}; | z ∪ Y = {2,3,4}∈ B(M); |
| 9: | Z = {1,2,4}; | Y = {3,5}; | z = Z – Y = {1,2,4}; | z ∪ Y = ({1,3,5}∨{2,3,5})∈ B(M); |
| 10: | Z = {1,2,4}; | Y = {4,5}; | z = Z – Y = {1,2}; | z ∪ Y = ({1,4,5}∨{2,4,5})∈ B(M); |
| 11: | Z = {1,2,5}; | Y = {1,2}; | z = Z – Y = {5}; | z ∪ Y = {1,2,5}∈ B(M); |
| 12: | Z = {1,2,5}; | Y = {1,3}; | z = Z – Y = {2,5}; | z ∪ Y = {1,3,5}∈ B(M); |
| 13: | Z = {1,2,5}; | Y = {1,4}; | z = Z – Y = {2,5}; | z ∪ Y = ({1,2,4}∨{1,4,5})∈ B(M); |
| 14: | Z = {1,2,5}; | Y = {1,5}; | z = Z – Y = {2}; | z ∪ Y = {1,2,5}∈ B(M); |
| 15: | Z = {1,2,5}; | Y = {2,3}; | z = Z – Y = {1,5}; | z ∪ Y = {2,3,5}∈ B(M); |
| 16: | Z = {1,2,5}; | Y = {2,4}; | z = Z – Y = {1,5}; | z ∪ Y = {2,4,5}∈ B(M); |
| 17: | Z = {1,2,5}; | Y = {2,5}; | z = Z – Y = {1}; | z ∪ Y = {1,2,5}∈ B(M); |
| 18: | Z = {1,2,5}; | Y = {3,4}; | z = Z – Y = {1,2,5}; | z ∪ Y = ({1,3,4}∨{2,3,4})∈ |
| | | | | $B(M)$; |
| 19: | $Z = \{1,2,5\}$; | $Y = \{3,5\}$; | $z = Z - Y = \{1,2\}$; | $z \cup Y = (\{1,3,5\}\vee\{2,3,5\})\in B(M)$; |
| 20: | $Z = \{1,2,5\}$; | $Y = \{4,5\}$; | $z = Z - Y = \{1,2\}$; | $z \cup Y = (\{1,4,5\}\vee\{2,4,5\})\in B(M)$; |
| 21: | $Z = \{1,3,4\}$; | $Y = \{1,2\}$; | $z = Z - Y = \{3,4\}$; | $z \cup Y = \{1,2,4\}\in B(M)$; |
| 22: | $Z = \{1,3,4\}$; | $Y = \{1,3\}$; | $z = Z - Y = \{4\}$; | $z \cup Y = \{1,3,4\}\in B(M)$; |
| 23: | $Z = \{1,3,4\}$; | $Y = \{1,4\}$; | $z = Z - Y = \{3\}$; | $z \cup Y = \{1,3,4\}\in B(M)$; |
| 24: | $Z = \{1,3,4\}$; | $Y = \{1,5\}$; | $z = Z - Y = \{3,4\}$; | $z \cup Y = (\{1,3,5\}\vee\{1,4,5\})\in B(M)$; |
| 25: | $Z = \{1,3,4\}$; | $Y = \{2,3\}$; | $z = Z - Y = \{1,4\}$; | $z \cup Y = \{2,3,4\}\in B(M)$; |
| 26: | $Z = \{1,3,4\}$; | $Y = \{2,4\}$; | $z = Z - Y = \{1,3\}$; | $z \cup Y = \{2,3,4\}\vee\{1,2,4\})\in B(M)$; |
| 27: | $Z = \{1,3,4\}$; | $Y = \{2,5\}$; | $z = Z - Y = \{1,3,4\}$; | $z \cup Y = (\{1,2,5\}\vee\{2,3,5\})\in B(M)$; |
| 28: | $Z = \{1,3,4\}$; | $Y = \{3,4\}$; | $z = Z - Y = \{1\}$; | $z \cup Y = \{1,3,4\}\in B(M)$; |
| 29: | $Z = \{1,3,4\}$; | $Y = \{3,5\}$; | $z = Z - Y = \{1,4\}$; | $z \cup Y = \{1,3,5\}\in B(M)$; |
| 30: | $Z = \{1,3,4\}$; | $Y = \{4,5\}$; | $z = Z - Y = \{1,3\}$; | $z \cup Y = \{1,4,5\}\in B(M)$; |
| 31: | $Z = \{1,3,5\}$; | $Y = \{1,2\}$; | $z = Z - Y = \{3,5\}$; | $z \cup Y = \{1,2,5\}\in B(M)$; |
| 32: | $Z = \{1,3,5\}$; | $Y = \{1,3\}$; | $z = Z - Y = \{5\}$; | $z \cup Y = \{1,3,5\}\in B(M)$; |
| 33: | $Z = \{1,3,5\}$; | $Y = \{1,4\}$; | $z = Z - Y = \{3,5\}$; | $z \cup Y = (\{1,3,4\}\vee\{1,4,5\})\in B(M)$; |
| 34: | $Z = \{1,3,5\}$; | $Y = \{1,5\}$; | $z = Z - Y = \{3\}$; | $z \cup Y = \{1,3,5\}\in B(M)$; |
| 35: | $Z = \{1,3,5\}$; | $Y = \{2,3\}$; | $z = Z - Y = \{1,5\}$; | $z \cup Y = \{2,3,5\}\in B(M)$; |
| 36: | $Z = \{1,3,5\}$; | $Y = \{2,4\}$; | $z = Z - Y = \{1,3,5\}$; | $z \cup Y = (\{1,2,4\}\vee\{2,3,4\}\vee \vee\{2,4,5\}) \in B(M)$; |
| 37: | $Z = \{1,3,5\}$; | $Y = \{2,5\}$; | $z = Z - Y = \{1,3\}$; | $z \cup Y = \{2,3,5\}\in B(M)$; |
| 38: | $Z = \{1,3,5\}$; | $Y = \{3,4\}$; | $z = Z - Y = \{1,5\}$; | $z \cup Y = \{1,3,4\}\in B(M)$; |
| 39: | $Z = \{1,3,5\}$; | $Y = \{3,5\}$; | $z = Z - Y = \{1\}$; | $z \cup Y = \{1,3,5\}\in B(M)$; |
| 40: | $Z = \{1,3,5\}$; | $Y = \{4,5\}$; | $z = Z - Y = \{1,3\}$; | $z \cup Y = \{1,4,5\}\in B(M)$; |
| 41: | $Z = \{1,4,5\}$; | $Y = \{1,2\}$; | $z = Z - Y = \{4,5\}$; | $z \cup Y = (\{1,2,4\}\vee\{1,2,5\})\in B(M)$; |
| 42: | $Z = \{1,4,5\}$; | $Y = \{1,3\}$; | $z = Z - Y = \{4,5\}$; | $z \cup Y = (\{1,3,4\}\vee\{1,3,5\})\in B(M)$; |
| 43: | $Z = \{1,4,5\}$; | $Y = \{1,4\}$; | $z = Z - Y = \{5\}$; | $z \cup Y = \{1,4,5\}\in B(M)$; |
| 44: | $Z = \{1,4,5\}$; | $Y = \{1,5\}$; | $z = Z - Y = \{4\}$; | $z \cup Y = \{1,4,5\}\in B(M)$; |
| 45: | $Z = \{1,4,5\}$; | $Y = \{2,3\}$; | $z = Z - Y = \{1,4,5\}$; | $z \cup Y = (\{2,3,4\}\vee\{2,3,5\})\in B(M)$; |
| 46: | $Z = \{1,4,5\}$; | $Y = \{2,4\}$; | $z = Z - Y = \{1,5\}$; | $z \cup Y = (\{1,2,4\}\vee\{2,4,5\})\in B(M)$; |
| 47: | $Z = \{1,4,5\}$; | $Y = \{2,5\}$; | $z = Z - Y = \{1,4\}$; | $z \cup Y = (\{1,2,5\}\vee\{2,4,5\})\in B(M)$; |
| 48: | $Z = \{1,4,5\}$; | $Y = \{3,4\}$; | $z = Z - Y = \{1,5\}$; | $z \cup Y = \{1,3,4\}\in B(M)$; |
| 49: | $Z = \{1,4,5\}$; | $Y = \{3,5\}$; | $z = Z - Y = \{1,4\}$; | $z \cup Y = \{1,3,5\}\in B(M)$; |
| 50: | $Z = \{1,4,5\}$; | $Y = \{4,5\}$; | $z = Z - Y = \{1\}$; | $z \cup Y = \{1,4,5\}\in B(M)$; |
| 51: | $Z = \{2,3,4\}$; | $Y = \{1,2\}$; | $z = Z - Y = \{3,4\}$; | $z \cup Y = \{1,2,4\}\in B(M)$; |
| 52: | $Z = \{2,3,4\}$; | $Y = \{1,3\}$; | $z = Z - Y = \{2,4\}$; | $z \cup Y = \{1,3,4\}\in B(M)$; |
| 53: | $Z = \{2,3,4\}$; | $Y = \{1,4\}$; | $z = Z - Y = \{2,3\}$; | $z \cup Y = (\{1,2,4\}\vee\{1,3,4\})\in B(M)$; |
| 54: | $Z = \{2,3,4\}$; | $Y = \{1,5\}$; | $z = Z - Y = \{2,3,4\}$; | $z \cup Y = (\{1,2,5\}\vee\{1,3,5\}\vee$ |

$\vee \{1,4,5\}) \in B(M)$;

55: $Z = \{2,3,4\}$; $Y = \{2,3\}$; $z = Z - Y = \{4\}$; $z \cup Y = \{2,3,4\} \in B(M)$;
56: $Z = \{2,3,4\}$; $Y = \{2,4\}$; $z = Z - Y = \{3\}$; $z \cup Y = \{2,3,4\} \in B(M)$;
57: $Z = \{2,3,4\}$; $Y = \{2,5\}$; $z = Z - Y = \{3,4\}$; $z \cup Y = (\{2,3,5\} \vee \{2,4,5\}) \in B(M)$;
58: $Z = \{2,3,4\}$; $Y = \{3,4\}$; $z = Z - Y = \{2\}$; $z \cup Y = \{2,3,4\} \in B(M)$;
59: $Z = \{2,3,4\}$; $Y = \{3,5\}$; $z = Z - Y = \{2,4\}$; $z \cup Y = \{2,3,5\} \in B(M)$;
60: $Z = \{2,3,4\}$; $Y = \{4,5\}$; $z = Z - Y = \{2,3\}$; $z \cup Y = \{2,4,5\} \in B(M)$;
61: $Z = \{2,3,5\}$; $Y = \{1,2\}$; $z = Z - Y = \{3,5\}$; $z \cup Y = \{1,2,5\} \in B(M)$;
62: $Z = \{2,3,5\}$; $Y = \{1,3\}$; $z = Z - Y = \{2,5\}$; $z \cup Y = \{1,3,5\} \in B(M)$;
63: $Z = \{2,3,5\}$; $Y = \{1,4\}$; $z = Z - Y = \{2,3,5\}$; $z \cup Y = \{1,2,4\} \vee \{1,3,4\}) \in B(M)$;
64: $Z = \{2,3,5\}$; $Y = \{1,5\}$; $z = Z - Y = \{2,3\}$; $z \cup Y = (\{1,2,5\} \vee \{1,3,5\}) \in B(M)$;
65: $Z = \{2,3,5\}$; $Y = \{2,3\}$; $z = Z - Y = \{5\}$; $z \cup Y = (\{2,3,5\} \in B(M)$;
66: $Z = \{2,3,5\}$; $Y = \{2,4\}$; $z = Z - Y = \{3,5\}$; $z \cup Y = (\{2,3,4\} \vee \{2,4,5\}) \in B(M)$;
67: $Z = \{2,3,5\}$; $Y = \{2,5\}$; $z = Z - Y = \{3\}$; $z \cup Y = \{2,3,5\} \in B(M)$;
68: $Z = \{2,3,5\}$; $Y = \{3,4\}$; $z = Z - Y = \{2,5\}$; $z \cup Y = \{2,3,4\} \in B(M)$;
69: $Z = \{2,3,5\}$; $Y = \{3,5\}$; $z = Z - Y = \{2\}$; $z \cup Y = \{2,3,5\} \in B(M)$;
70: $Z = \{2,3,5\}$; $Y = \{4,5\}$; $z = Z - Y = \{2,3\}$; $z \cup Y = \{2,4,5\} \in B(M)$;
71: $Z = \{2,4,5\}$; $Y = \{1,2\}$; $z = Z - Y = \{4,5\}$; $z \cup Y = (\{1,2,4\} \vee \{1,2,5\}) \in B(M)$;
72: $Z = \{2,4,5\}$; $Y = \{1,3\}$; $z = Z - Y = \{2,4,5\}$; $z \cup Y = (\{1,3,4\} \vee \{1,3,5\}) \in B(M)$;
73: $Z = \{2,4,5\}$; $Y = \{1,4\}$; $z = Z - Y = \{2,5\}$; $z \cup Y = (\{1,2,4\} \vee \{1,4,5\}) \in B(M)$;
74: $Z = \{2,4,5\}$; $Y = \{1,5\}$; $z = Z - Y = \{2,4\}$; $z \cup Y = (\{1,2,5\} \vee \{1,4,5\}) \in B(M)$;
75: $Z = \{2,4,5\}$; $Y = \{2,3\}$; $z = Z - Y = \{4,5\}$; $z \cup Y = (\{2,3,4\} \vee \{2,3,5\}) \in B(M)$;
76: $Z = \{2,4,5\}$; $Y = \{2,4\}$; $z = Z - Y = \{5\}$; $z \cup Y = \{2,4,5\} \in B(M)$;
77: $Z = \{2,4,5\}$; $Y = \{2,5\}$; $z = Z - Y = \{4\}$; $z \cup Y = \{2,4,5\} \in B(M)$;
78: $Z = \{2,4,5\}$; $Y = \{3,4\}$; $z = Z - Y = \{2,5\}$; $z \cup Y = \{2,3,4\} \in B(M)$;
79: $Z = \{2,4,5\}$; $Y = \{3,5\}$; $z = Z - Y = \{2,4\}$; $z \cup Y = \{2,3,5\} \in B(M)$;
80: $Z = \{2,4,5\}$; $Y = \{4,5\}$; $z = Z - Y = \{2\}$; $z \cup Y = \{2,4,5\} \in B(M)$.

Как мы видим, все 80 элементов $z \cup Y$ принадлежат множеству B(M) - баз матроида и удовлетворяют условию (2.32).

Проверим удовлетворение условия (2.32) в случае множества $\zeta$ (M) - циклов матроида или нет.

1: $Z = \{1,2,3\}$; $Y = \{1,2\}$; $z = Z - Y = \{3\}$; $z \cup Y = \{1,2,3\} \in \zeta(M)$;
2: $Z = \{1,2,3\}$; $Y = \{1,3\}$; $z = Z - Y = \{2\}$; $z \cup Y = \{1,2,3\} \in \zeta(M)$;
3: $Z = \{1,2,3\}$; $Y = \{1,4\}$; $z = Z - Y = \{2,3\}$; $z \cup Y = (\{1,2,4\} \vee \{1,3,4\}) \in B(M)$;
4: $Z = \{1,2,3\}$; $Y = \{1,5\}$; $z = Z - Y = \{2,3\}$; $z \cup Y = (\{1,2,5\} \vee \{1,3,5\}) \in B(M)$;
5: $Z = \{1,2,3\}$; $Y = \{2,3\}$; $z = Z - Y = \{1\}$; $z \cup Y = \{1,2,3\} \in \zeta(M)$;
6: $Z = \{1,2,3\}$; $Y = \{2,4\}$; $z = Z - Y = \{1,3\}$; $z \cup Y = (\{1,2,3\} \in \zeta(M)) \vee$

| | | | | |
|---|---|---|---|---|
| | | | | $\vee(\{2,3,4\}) \in B(M))$; |
| 7: | $Z = \{1,2,3\}$; | $Y = \{2,5\}$; | $z = Z - Y = \{1,3\}$; | $z \cup Y = (\{1,2,3\} \in \zeta(M)) \vee$ $\vee(\{2,3,5\}) \in B(M))$; |
| 8: | $Z = \{1,2,3\}$; | $Y = \{3,4\}$; | $z = Z - Y = \{1,2\}$; | $z \cup Y = (\{1,3,4\} \vee \{2,3,4\}) \in B(M)$; |
| 9: | $Z = \{1,2,3\}$; | $Y = \{3,5\}$; | $z = Z - Y = \{1,2\}$; | $z \cup Y = (\{1,3,5\} \vee \{1,2,5\}) \in B(M)$; |
| 10: | $Z = \{1,2,3\}$; | $Y = \{4,5\}$; | $z = Z - Y = \{1,2,3\}$; | $z \cup Y = ((\{1,4,5\} \vee \{2,4,5\}) \in$ $\in B(M)) \vee (\{3,4,5\} \in \zeta(M))$; |
| 11: | $Z = \{3,4,5\}$; | $Y = \{1,2\}$; | $z = Z - Y = \{3,4,5\}$; | $z \cup Y = ((\{1,2,4\} \vee \{1,2,5\}) \in$ $\in B(M)) \vee (\{1,2,3\} \in \zeta(M))$; |
| 12: | $Z = \{3,4,5\}$; | $Y = \{1,3\}$; | $z = Z - Y = \{4,5\}$; | $z \cup Y = (\{1,3,4\} \vee \{1,3,5\}) \in B(M)$; |
| 13: | $Z = \{3,4,5\}$; | $Y = \{1,4\}$; | $z = Z - Y = \{3,5\}$; | $z \cup Y = (\{1,3,4\} \vee \{1,4,5\}) \in B(M)$; |
| 14: | $Z = \{3,4,5\}$; | $Y = \{1,5\}$; | $z = Z - Y = \{3,4\}$; | $z \cup Y = (\{1,3,5\} \vee \{1,4,5\}) \in B(M)$; |
| 15: | $Z = \{3,4,5\}$; | $Y = \{2,3\}$; | $z = Z - Y = \{4,5\}$; | $z \cup Y = (\{2,3,4\} \vee \{2,3,5\}) \in B(M)$; |
| 16: | $Z = \{3,4,5\}$; | $Y = \{2,4\}$; | $z = Z - Y = \{3,5\}$; | $z \cup Y = (\{2,3,4\} \vee \{2,4,5\}) \in B(M)$; |
| 17: | $Z = \{3,4,5\}$; | $Y = \{2,5\}$; | $z = Z - Y = \{3,4\}$; | $z \cup Y = (\{2,3,4\} \vee \{2,4,5\}) \in B(M)$; |
| 18: | $Z = \{3,4,5\}$; | $Y = \{3,4\}$; | $z = Z - Y = \{5\}$; | $z \cup Y = \{3,4,5\} \in \zeta(M)$; |
| 19: | $Z = \{3,4,5\}$; | $Y = \{3,5\}$; | $z = Z - Y = \{4\}$; | $z \cup Y = \{3,4,5\} \in \zeta(M)$; |
| 20: | $Z = \{3,4,5\}$; | $Y = \{4,5\}$; | $z = Z - Y = \{3\}$; | $z \cup Y = \{3,4,5\} \in \zeta(M)$; |

Как мы видим, здесь условие (2.32) не удовлетворяется.

***Пример 2.12.*** Для более подробного изучения свойств циклического матроида рассмотрим следующий граф, представленный на рис. 2.9 (здесь цифры характеризуют ребра графа).

Выделим множество изометрических циклов для данного графа:

$c_1 = \{1,2,5\}$, $c_2 = \{1,4,6\}$, $c_3 = \{2,3,7\}$, $c_4 = \{2,4,8\}$,

$c_5 = \{3,4,9\}$, $c_6 = \{7,8,9\}$, $c_7 = \{5,6,8\}$.

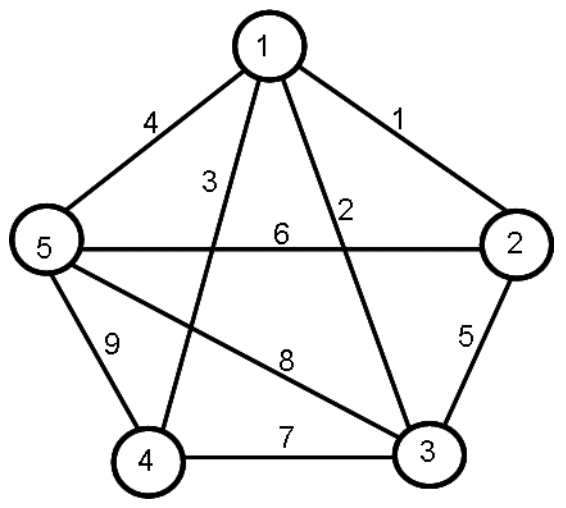

Рис. 2.9. Граф G.

Кроме того, имеется следующая система циклов длиной четыре:

$c_9 = c_1 \oplus c_3 = \{1,2,5\} \oplus \{2,3,7\} = \{1,3,5,7\}$;
$c_{10} = c_3 \oplus c_5 = \{2,3,7\} \oplus \{3,4,9\} = \{2,4,7,9\}$;
$c_{11} = c_1 \oplus c_7 = \{1,2,5\} \oplus \{5,6,8\} = \{1,2,6,8\}$;
$c_{12} = c_2 \oplus c_5 = \{1,4,6\} \oplus \{3,4,9\} = \{1,3,6,9\}$;
$c_{13} = c_2 \oplus c_7 = \{1,4,6\} \oplus \{5,6,8\} = \{1,4,5,8\}$;
$c_{14} = c_3 \oplus c_6 = \{2,3,7\} \oplus \{7,8,9\} = \{2,3,8,9\}$;
$c_{15} = c_4 \oplus c_7 = \{2,4,8\} \oplus \{3,6,8\} = \{2,4,5,6\}$;
$c_{16} = c_6 \oplus c_7 = \{7,8,9\} \oplus \{5,6,8\} = \{5,6,7,9\}$;
$c_{17} = c_5 \oplus c_6 = \{3,4,9\} \oplus \{7,8,9\} = \{3,4,7,8\}$.

Запишем центральные разрезы как однострочные структурные числа:

$s_1 = [1\ \ 2\ \ 3\ \ 4]$;
$s_2 = [1\ \ 5\ \ 6]$;
$s_3 = [2\ \ 5\ \ 7\ \ 8]$;
$s_4 = [3\ \ 7\ \ 9]$;
$s_5 = [4\ \ 6\ \ 8\ \ 9]$.

Произведем умножение первых четырех по правилам умножения структурных чисел

$[1\ \ 2\ \ 3\ \ 4] \times [1\ \ 5\ \ 6] =$

| | | | | | | | | | | |
|---|---|---|---|---|---|---|---|---|---|---|
| 1 | 1 | 2 | 2 | 2 | 3 | 3 | 3 | 4 | 4 | 4 |
| 5 | 6 | 1 | 5 | 6 | 1 | 5 | 6 | 1 | 5 | 6 |

$[1\ \ 2\ \ 3\ \ 4] \times [1\ \ 5\ \ 6] \times [2\ \ 5\ \ 7\ \ 8] =$

| | | | | | | | | | | | | | | | | | | | | | | | | |
|---|---|---|---|---|---|---|---|---|---|---|---|---|---|---|---|---|---|---|---|---|---|---|---|---|
| 1 | 1 | 1 | 1 | 1 | 1 | 1 | 2 | 2 | 2 | 2 | 2 | 2 | 2 | 2 | 3 | 3 | 3 | 3 | 3 | 3 | 3 | 3 | 3 | 3 |
| 5 | 5 | 5 | 6 | 6 | 6 | 6 | 1 | 1 | 1 | 5 | 5 | 6 | 6 | 6 | 1 | 1 | 1 | 1 | 5 | 5 | 6 | 6 | 6 | 6 |
| 2 | 7 | 8 | 2 | 5 | 7 | 8 | 5 | 7 | 8 | 7 | 8 | 5 | 7 | 8 | 2 | 5 | 7 | 8 | 7 | 8 | 2 | 5 | 7 | 8 |

| | | | | | | | | | | |
|---|---|---|---|---|---|---|---|---|---|---|
| 4 | 4 | 4 | 4 | 4 | 4 | 4 | 4 | 4 | 4 | 4 |
| 1 | 1 | 1 | 1 | 5 | 5 | 5 | 6 | 6 | 6 | 6 |
| 2 | 5 | 7 | 8 | 2 | 7 | 8 | 2 | 5 | 7 | 8 |

$[1\ \ 2\ \ 3\ \ 4] \times [1\ \ 5\ \ 6] \times [2\ \ 5\ \ 7\ \ 8] \times [3\ \ 7\ \ 9] =$

| | | | | | | | | | | | | | | | | | | | | | | | | |
|---|---|---|---|---|---|---|---|---|---|---|---|---|---|---|---|---|---|---|---|---|---|---|---|---|
| 1 | 1 | 1 | 1 | 1 | 1 | 1 | 1 | 1 | 1 | 1 | 1 | 1 | 1 | 1 | 1 | 2 | 2 | 2 | 2 | 2 | 2 | 2 | 2 | 2 |
| 5 | 5 | 5 | 5 | 5 | 6 | 6 | 6 | 6 | 6 | 6 | 6 | 6 | 6 | 6 | 6 | 1 | 1 | 1 | 1 | 1 | 5 | 5 | 5 | 5 |
| 7 | 7 | 8 | 8 | 8 | 2 | 2 | 2 | 5 | 5 | 5 | 7 | 7 | 8 | 8 | 8 | 7 | 7 | 8 | 8 | 8 | 7 | 7 | 8 | 8 |
| 3 | 7 | 3 | 7 | 9 | 3 | 7 | 9 | 3 | 7 | 9 | 3 | 9 | 3 | 7 | 9 | 3 | 9 | 3 | 7 | 9 | 3 | 9 | 3 | 7 |
| 2 | 2 | 2 | 2 | 2 | 2 | 2 | 2 | 2 | 3 | 3 | 3 | 3 | 3 | 3 | 3 | 3 | 3 | 3 | 3 | 3 | 3 | 3 | 3 | 3 |
| 5 | 6 | 6 | 6 | 6 | 6 | 6 | 6 | 6 | 1 | 1 | 1 | 1 | 1 | 1 | 1 | 5 | 5 | 5 | 5 | 5 | 6 | 6 | 6 | 6 |
| 8 | 5 | 5 | 5 | 7 | 7 | 8 | 8 | 8 | 2 | 2 | 5 | 5 | 7 | 8 | 8 | 2 | 2 | 7 | 8 | 8 | 2 | 2 | 5 | 5 |
| 9 | 3 | 7 | 9 | 3 | 9 | 3 | 7 | 9 | 7 | 9 | 7 | 9 | 9 | 7 | 9 | 7 | 9 | 9 | 7 | 9 | 7 | 9 | 7 | 9 |
| 3 | 3 | 3 | 4 | 4 | 4 | 4 | 4 | 4 | 4 | 4 | 4 | 4 | 4 | 4 | 4 | 4 | 4 | 4 | 4 | 4 | 4 | 4 | 4 | 4 |
| 6 | 6 | 6 | 1 | 1 | 1 | 1 | 1 | 1 | 1 | 1 | 1 | 1 | 1 | 5 | 5 | 5 | 5 | 5 | 5 | 5 | 5 | 6 | 6 | 6 |
| 7 | 8 | 8 | 2 | 2 | 2 | 5 | 5 | 5 | 7 | 7 | 8 | 8 | 8 | 2 | 2 | 2 | 7 | 7 | 8 | 8 | 8 | 2 | 2 | 2 |
| 9 | 7 | 9 | 3 | 7 | 9 | 3 | 7 | 9 | 3 | 9 | 3 | 7 | 9 | 3 | 7 | 9 | 3 | 9 | 3 | 7 | 9 | 3 | 7 | 9 |

| | | | | | | | |
|---|---|---|---|---|---|---|---|
| 4 | 4 | 4 | 4 | 4 | 4 | 4 | 4 |
| 6 | 6 | 6 | 6 | 6 | 6 | 6 | 6 |
| 5 | 5 | 5 | 7 | 7 | 8 | 8 | 8 |
| 3 | 7 | 9 | 3 | 9 | 3 | 7 | 9 |

Таким образом, каждая база циклического матроида из множества баз циклического матроида, может быть представлена как элемент структурного числа T.

Очевидно также, что множество баз матроида разрезов $B(M^*)$ может быть описано структурным числом $T^{d}$.

## 2.9. Матроид изометрических циклов

Теорема о существовании базиса изометрических циклов [37], говорит о том, что для любого связного графа без петель и кратных ребер линейное подпространство квазициклов имеет независимое подмножество состоящее из изометрических циклов с мощностью равной цикломатическому числу графа, то можно ввести понятие матроида изометрических циклов. Действительно, так как количество изометрических циклов в графе G всегда больше или равно

$$C_\tau \geq m - n + 1 \qquad (2.29)$$

цикломатическому числу графа. И если выбрать базис подпространства циклов состоящий из изомерических циклов, то остальные изометрических циклы не вошедшие в базис линейного подпространства циклов можно представить в виде линейной комбинации базисных изометрических циклов и следовательно менять набор базисных изометрических циклов. Причем один базисный набор будет отличаться от другого одним ребром, а это как раз и есть аксиома независимости (2.32).

Здесь, как раз и проявляется связь теории матроидов с теорией линейной зависимости. Аксиомы теории матроидов выбраны таким образом, чтобы они отражали наиболее характерные свойства независимых множеств линейного пространства. Точнее говоря, отражает свойства независимых подмножеств, содержащихся в некотором конечном подмножестве этого пространства (вспомним, что подмножество $\{e_1,\ldots,e_m\}$ линейного пространства называется независимым, если не существует набора скаляров $\lambda_1,\ldots,\lambda_m$, не всех равных нулю, такого что $\lambda_1 e_1 + \ldots + \lambda_m e_m = 0$). В самом деле, очевидно, что произвольное подмножество независимого множества линейного пространства независимо. Если $|B| = |A| + 1$ для линейно независимых подмножеств A, B линейного пространства, то A порождает пространство размерности $|\mathbf{A}|$, которое может содержать не более чем $|A|$ элементов множества B. Следовательно, существует элемент e $\in$ $B\backslash A$, не принадлежащий этому подпространству. Множество **A** $\cup$ {e} порождает подпространство размерности |**A**| + 1 и, следовательно, является линейно независимым.

***Пример 2.13***. Покажем, как образуется матроид изометрических циклов из множества изометрических циклов графа. С этой целью рассмотрим граф, представленный на рис. 2.7.

Множество изометрических циклов для данного графа имеет вид:

$c_1 = \{e_1,e_2,e_5\}$, $c_2 = \{e_1,e_4,e_6\}$, $c_3 = \{e_2,e_3,e_7\}$, $c_4 = \{e_2,e_4,e_8\}$,
$c_5 = \{e_3,e_4,e_9\}$, $c_6 = \{e_7,e_8,e_9\}$, $c_7 = \{e_5,e_6,e_8\}$.

Выделим в данном графе любое дерево предположим $T = \{e_1,e_2,e_3,e_4\}$, тогда множество хорд $H = \{e_5,e_6,e_7,e_8,e_9\}$. Как известно, цикломатическое число графа равно количеству хорд графа.

Для каждой хорды построим однострочное структурное число характеризующие циклы, проходящие по данной хорде:

$e_5$: $[c_1\ c_7]$;
$e_6$: $[c_2\ c_7]$;
$e_7$: $[c_3\ c_6]$;
$e_8$: $[c_4\ c_6\ c_7]$;
$e_9$: $[c_5\ c_6]$.

Для удобства записи циклы будем обозначать их номерами. Конечно, такое обозначение может ввести путаницу, поэтому каждый раз мы будем оговаривать, что мы подразумеваем под номерами - вершины, ребра или циклы.

Перемножим однострочные структурные числа $e_5{\times}e_6{\times}e_7{\times}e_8{\times}e_9$ по правилам умножения алгебра структурных чисел

$[c_1\ c_7] \times [c_2\ c_7] \times [c_3\ c_6] \times [c_4\ c_6\ c_7] \times [c_5\ c_6] =$

| | | | | | | | | | | | | | | | | |
|---|---|---|---|---|---|---|---|---|---|---|---|---|---|---|---|---|
| | $e_5$: | $c_1$ | $c_1$ | $c_1$ | $c_1$ | $c_1$ | $c_1$ | $c_1$ | $c_1$ | $c_1$ | $c_1$ | $c_1$ | $c_7$ | $c_7$ | $c_7$ | $c_7$ |
| | $e_6$: | $c_2$ | $c_2$ | $c_2$ | $c_2$ | $c_2$ | $c_2$ | $c_2$ | $c_7$ | $c_7$ | $c_7$ | $c_7$ | $c_2$ | $c_2$ | $c_2$ | $c_2$ |
| $B(C_\tau) =$ | $e_7$: | $c_3$ | $c_3$ | $c_3$ | $c_3$ | $c_3$ | $c_6$ | $c_6$ | $c_3$ | $c_3$ | $c_3$ | $c_6$ | $c_3$ | $c_3$ | $c_3$ | $c_6$ |
| | $e_8$: | $c_4$ | $c_4$ | $c_6$ | $c_7$ | $c_7$ | $c_4$ | $c_7$ | $c_4$ | $c_4$ | $c_6$ | $c_4$ | $c_4$ | $c_4$ | $c_6$ | $c_4$ |
| | $e_9$ | $c_5$ | $c_6$ | $c_5$ | $c_5$ | $c_6$ | $c_5$ | $c_5$ | $c_5$ | $c_6$ | $c_5$ | $c_5$ | $c_5$ | $c_6$ | $c_5$ | $c_5$ |

Или в численном виде:

| | | | | | | | | | | | | | | | | |
|---|---|---|---|---|---|---|---|---|---|---|---|---|---|---|---|---|
| | $e_5$: | 1 | 1 | 1 | 1 | 1 | 1 | 1 | 1 | 1 | 1 | 1 | 7 | 7 | 7 | 7 |
| | $e_6$: | 2 | 2 | 2 | 2 | 2 | 2 | 2 | 7 | 7 | 7 | 7 | 2 | 2 | 2 | 2 |
| = | $e_7$: | 3 | 3 | 3 | 3 | 3 | 6 | 6 | 3 | 3 | 3 | 6 | 3 | 3 | 3 | 6 |
| | $e_8$: | 4 | 4 | 6 | 7 | 7 | 4 | 7 | 4 | 4 | 6 | 4 | 4 | 4 | 6 | 4 |
| | $e_9$ | 5 | 6 | 5 | 5 | 6 | 5 | 5 | 5 | 6 | 5 | 5 | 5 | 6 | 5 | 5 |

Что касается циклов матроида, то для нашего графа - это набор линейно зависимых изометрических циклов

$c_1 \oplus c_2 \oplus c_4 \oplus c_7 = \{e_1,e_2,e_5\} \oplus \{e_1,e_4,e_6\} \oplus \{e_2,e_4,e_8\} \oplus \{e_5,e_6,e_8\} = \varnothing$;

$c_3 \oplus c_4 \oplus c_5 \oplus c_6 = \{e_2,e_3,e_7\} \oplus \{e_2,e_4,e_8\} \oplus \{e_3,e_4,e_9\} \oplus \{e_7,e_8,e_9\} = \varnothing$.

Кроме того, существует еще одна линейно зависимая система, состоящая из шести изометрических циклов, характерная тем, что для данного элемента структурного числа изометрических циклов, обод графа есть изометрический цикл

$c_1 \oplus c_2 \oplus c_3 \oplus c_5 \oplus c_6 \oplus c_7 =$

$= \{e_1,e_2,e_5\} \oplus \{e_1,e_4,e_6\} \oplus \{e_2,e_3,e_7\} \oplus \{e_3,e_4,e_9\} \oplus \{e_7,e_8,e_9\} \oplus \{e_5,e_6,e_8\} = \varnothing$.

Например, это столбцы структурного числа $B(C^{\tau})$ с номерами 3,4,5,7,10,14 не содержащие изометрический цикл $c_4$.

## 2.10. Матроид центральных разрезов

Кроме матроида изометрических циклов существует матроид центральных разрезов графа. Построение матроида центральных разрезов графа осуществляется простым способом. Здесь каждый элемент структурного числа S характеризует n-1 независимый набор центральных разрезов графа.

***Пример 2.14***. В качестве примера рассмотрим граф, представленный на рис. 6.9. Множество центральных разрезов графа

$s_1 = \{e_1,e_2,e_3,e_4\}$, $s_2 = \{e_1,e_5,e_6\}$, $s_3 = \{e_2,e_5,e_7,e_8\}$, $s_4 = \{e_3,e_7,e_9\}$, $s_5 = \{e_3,e_4,e_8,e_9\}$.

Выделим в данном графе любое дерево предположим $T = \{e_1,e_2,e_3,e_4\}$.

Для каждой ветви дерева построим однострочное структурное число характеризующие разрезы, включающие данную ветвь дерева:

$e_1$: $[s_1\ s_2]$;
$e_2$: $[s_1\ s_3]$;
$e_3$: $[s_1\ s_4]$;
$e_4$: $[s_1\ s_5]$.

Для удобства записи разрезы будем обозначать их номерами. Перемножим однострочные структурные числа $e_1 \times e_2 \times e_3 \times e_4$ по правилам умножения алгебра структурных чисел

$[s_1\ s_2] \times [s_1\ s_3] \times [s_1\ s_4] \times [s_1\ s_5] =$

| | | | | | | | | | | | | | |
|---|---|---|---|---|---|---|---|---|---|---|---|---|---|
| B(S) = | $e_1$: | $s_1$ | $s_2$ | $s_2$ | $s_2$ | $s_2$ | Или в виде номеров центральных разрезов: | 1 | 2 | 2 | 2 | 2 |
| | $e_2$: | $s_3$ | $s_1$ | $s_3$ | $s_3$ | $s_3$ | | 3 | 1 | 3 | 3 | 3 |
| | $e_3$: | $s_4$ | $s_4$ | $s_1$ | $s_4$ | $s_4$ | | 4 | 4 | 1 | 4 | 4 |
| | $e_4$: | $s_5$ | $s_5$ | $s_5$ | $s_1$ | $s_5$ | | 5 | 5 | 5 | 1 | 5 |

## 2.11. Жадный алгоритм решения оптимизационных задач

Рассмотрим следующую матрицу с действительными неотрицательными коэффициентами:

A=

| | 1 | 2 | 3 |
|---|---|---|---|
| 1 | 7 | 5 | 1 |
| 2 | 3 | 4 | 3 |
| 3 | 2 | 3 | 1 |

Займемся решением следующей оптимизационной задачи.

***Пример 2.15.*** Найти такое подмножество элементов матрицы, что

(а) в каждом столбце и в каждой строке находятся не более одного выбранного элемента

(б) сумма выбранных элементов является наибольшей из всех возможных.

Попробуем решить эту задачу следующим образом: будем выбирать элементы последовательно, причем каждый раз будем выбирать наибольший из элементов, которые мы можем добавить, не нарушая условия (а). Будем действовать так вплоть до момента, пока добавление произвольного элемента не нарушит условия (а). Алгоритм такого типа мы будем называть *жадным*.

Для нашей матрицы A жадный алгоритм находит подмножество

A=

| | | |
|---|---|---|
| **7** | 5 | 1 |
| 3 | **4** | 3 |
| 2 | 3 | **1** |

с суммой 12, что не является правильным решением, поскольку подмножество

A=

| | | |
|---|---|---|
| **7** | 5 | 1 |
| 3 | 4 | **3** |
| 2 | **3** | 1 |

имеет сумму 13. Следовательно, на втором шаге не следовало бы быть жадным, мы выигрываем в конечном результате, выбирая несколько меньший элемент (3 вместо 4).

Возникает вопрос: когда выгодно быть жадным, а когда нет?

Мы будем рассматривать оптимизационные задачи следующего типа.

Даны конечное множество $E$, семейство его подмножеств $J \subseteq \mathcal{P}(E)$ и функция $w : E \to R+$, где R+ обозначает множество вещественных неотрицательных чисел. Найти подмножество $S \in J$ с наибольшей суммой $\sum_{e \in S} w(e)$.

В нашем случае $E$ есть множество позиций матрицы, а $w$ ставит в соответствие позиции

(i,j) матрицы $[a_{ij}]$ число $a_{ij}$, $S \in J \Leftrightarrow$ каждый столбец и каждая строка содержит не более одной позиции из множества S.

Теперь мы можем сформулировать наш вопрос следующим образом, при каких условиях относительно семейства J жадный алгоритм правильно решает задачу для произвольной функции *w*?

Оказывается, что на этот вопрос можно найти простой ответ. А именно, достаточно, чтобы пара (*E,* J) образовывала так матроид.

**Теорема 2.3.** (Радо, Эдмондс, см. [54]). Если M = <E, J> есть матроид, то множество S, найденное жадным алгоритмом, является независимым множеством с наибольшим весом. Напротив, если M = <E, J> не является матроидом, то существует такая функция $w : E \rightarrow R+$, что S не будет независимым множеством с наибольшим весом.

Отметим также, что при обращении упорядочения элементов жадный алгоритм выберет множество *S*, которое не только имеет наименьший вес, но и *i*-й по величине элемент (считая от наименьшего) будет не больше *i*-го по величине элемента произвольной базы.

Теперь мы можем посмотреть на работу жадного алгоритма на примере 2.13 с позиций алгебры структурных чисел.

Запишем индексы каждой строки матрицы **A** как однострочное структурное число и произведем их умножение по правилам алгебры структурных чисел, получим следующее структурное число 3-го порядка, где в нижней строке проставлена сумма соответствующих элементов матрицы А.

$$A = \begin{vmatrix} 1 & 1 & 2 & 2 & 3 & 3 \\ 2 & 3 & 1 & 3 & 1 & 2 \\ 3 & 2 & 3 & 1 & 2 & 1 \end{vmatrix}$$
$$\Sigma \quad 12 \quad 13 \quad 9 \quad 10 \quad 7 \quad 7$$

Первый элемент структурного числа имеет сумму 12. Поменяем местами 2-ую и 3-ью строки матрицы А,

$$A = \begin{vmatrix} 7 & 5 & 1 \\ 2 & 3 & 1 \\ 3 & 4 & 3 \end{vmatrix}$$

структурное число не изменится. Но здесь первый элемент структурного числа уже имеет сумму 13.

$$\mathbf{A} = \begin{vmatrix} 1 & 1 & 2 & 2 & 3 & 3 \\ 2 & 3 & 1 & 3 & 1 & 2 \\ 3 & 2 & 3 & 1 & 2 & 1 \end{vmatrix}$$
$$\Sigma \quad 13 \quad 12 \quad 10 \quad 9 \quad 7 \quad 7$$

Во всех других случаях оптимизации мы получим оценку сложности алгоритма в виде многочлена относительно размерности задачи, что отнюдь не является типичной ситуацией для оптимизационных задач. В большой степени эффективность жадных алгоритмов вызывается тем, что элемент, один раз, включенный в оптимальное решение, остается в нем до конца. Здесь не бывает «проверки всех возможностей» и выбора последовательности включения элементов, характерной для алгоритмов с возвратом, что обычно приводит к экспоненциальному росту числа шагов при росте размерности задачи.

В качестве примера рассмотрим алгоритм Краскала [2], как пример работы жадного алгоритма с точки зрения методов алгебры структурных чисел. Следует отметить, что алгоритм Краскала сыграл большую роль в развитии теории оптимизации.

***Пример 2.14.*** Пусть задан граф G(V,E) с заданными весами ребер, записанными в круглых скобках (см. рис. 2.10). Будем искать дерево с минимальной суммой его весов.

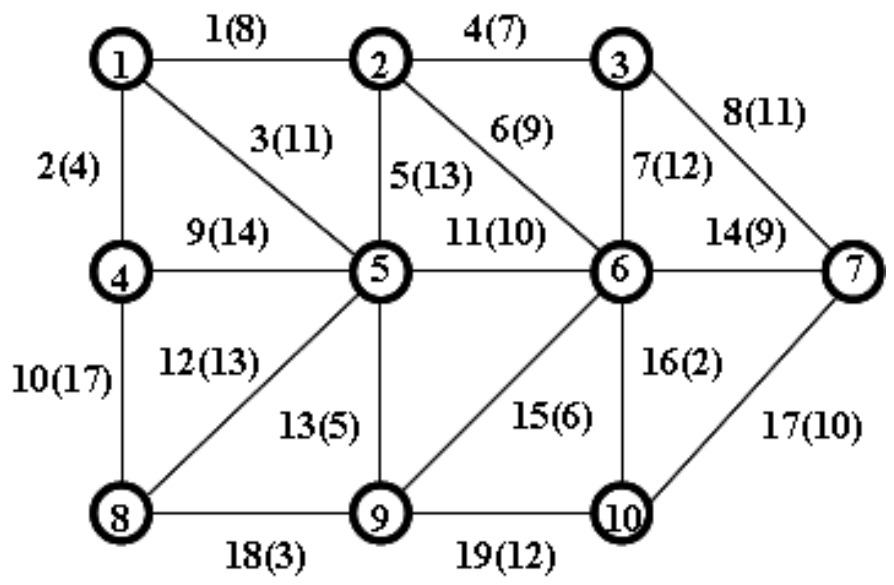

Рис. 2.10. Граф G с заданными весами ребер (ребра обозначены цифрами).

Запишем однострочные структурные числа, характеризующие центральные разрезы в виде списка ребер упорядоченного по не убыванию весов.

$s_1$ = [2(4) 1(8) 3(11)];

$s_2$ = [4(7) 1(8) 6(9) 5(13)];

$s_3$ = [4(7) 8(11) 7(12)];

$s_4$ = [2(4) 9(14) 10(17)];

$s_5$ = [13(5) 11(10) 3(11) 5(13) 12(13) 9(14)];

$s_6$ = [16(2) 15(6) 6(9) 14(9) 11(10) 7(12)];

$s_7$ = [14(9) 17(10) 8(11) ];

$s_8$ = [18(3) 12(13) 10(17)];

$s_9$ = [18(3) 13(5) 15(6) 19(12)];

$s_{10}$ = [16(2) 17(10) 19(12)].

**шаг 1**: На первых местах стоят четыре пары ребер, это ребро $e_2$ в разрезах $s_1$ и $s_4$, это ребро $e_4$ в разрезах $s_2$ и $s_3$, это ребро $e_{16}$ в разрезах $s_6$ и $s_{10}$, это ребро $e_{18}$ в разрезах $s_8$ и $s_9$. Ребро $e_2$ и вершины $v_1$ и $v_4$ образуют блок $Б_1$. Ребро $e_4$ и вершины $v_2$ и $v_3$ образуют блок $Б_2$. Ребро $e_{16}$ и вершины $v_6$ и $v_{10}$ образуют блок $Б_3$. Ребро $e_{18}$ и вершины $v_8$ и $v_9$ образуют блок $Б_4$. Теперь однострочные структурные числа характеризуют не только центральные разрезы, но и реберные разрезы (см. рис. 2.11).

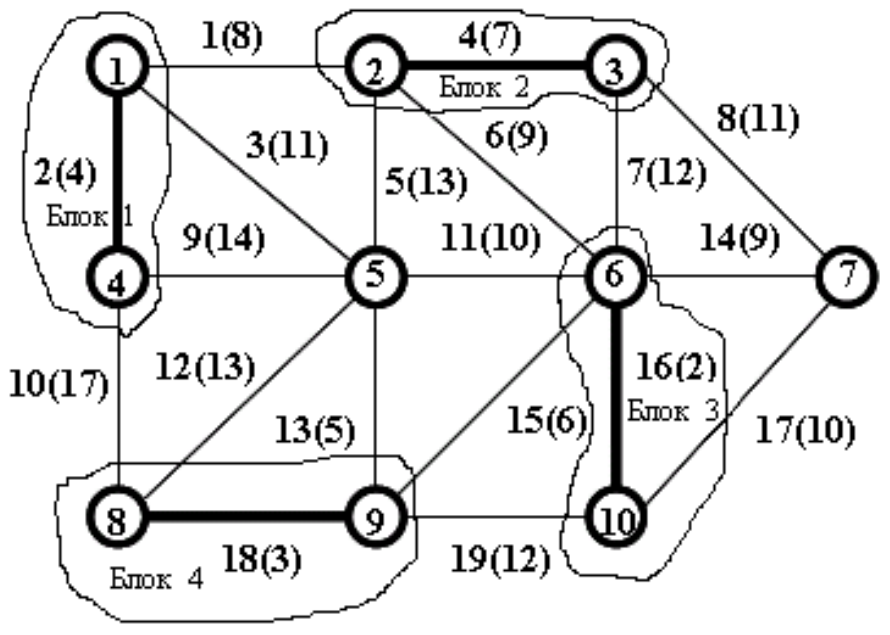

Рис. 2.11. Блоки графа G с минимальными весами ребер.

$Б_1$ = [1(8) 3(11) 9(14) 10(17)];

$Б_2$ = [1(8) 6(9) 8(11) 7(12) 5(13)];

$Б_3$ = [15(6) 6(9) 14(9) 11(10) 17(10) 7(12) 19(12];

$Б_4$ = [13(5) 15(6) 19(12) 12(13) 10(17)];

$s_5$ = [13(5) 11(10) 3(11) 5(13) 12(13) 9(14)];

$s_7$ = [14(9) 17(10) 8(11) ].

**шаг 2**: Вершины $v_5$ и $v_7$ оказались не включенные в блоки. Ребро $e_{13}$ с минимальным весом 5 принадлежащее разрезу $s_5$, соединено с блоком $Б_4$, поэтому оно включается в новый блок $Б_4$. Ребро $e_{14}$ с минимальным весом 9 принадлежащее разрезу $s_7$,соединено с блоком $Б_3$, поэтому оно включается в новый блок $Б_3$. Формируется новая запись однострочных структурных чисел

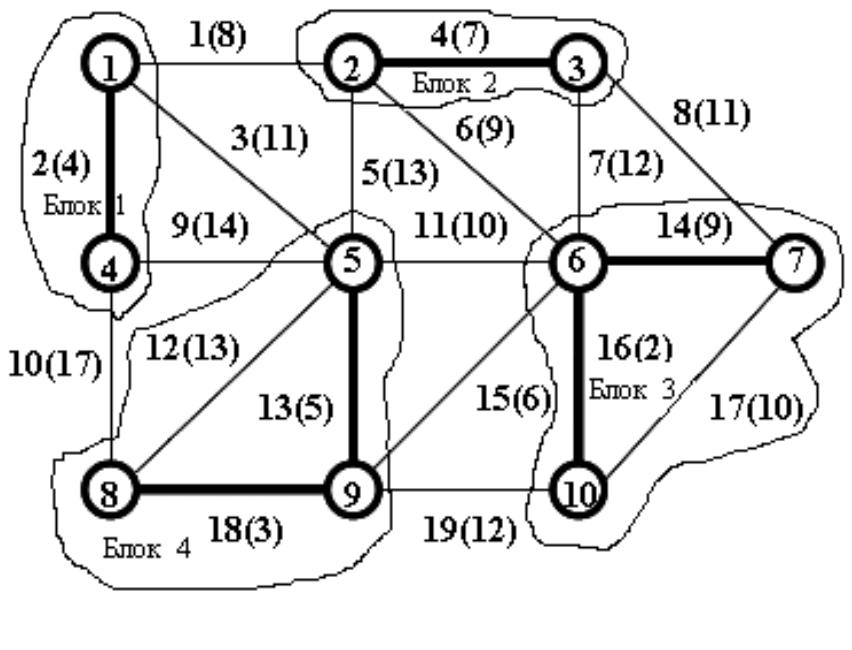

Рис. 2.12. Новые блоки графа G с минимальными весами ребер

$Б_1$ = [1(8) 3(11) 9(14) 10(17)];

$Б_2$ = [1(8) 6(9) 8(11) 7(12) 5(13)];

$Б_3$ = [15(6) 11(10) 3(11) 19(12) 5(13) 9(14) 10(17)];

$Б_4$ = [15(6) 6(9) 11(10) 8(11) 7(12) 19(12].

**шаг 3**: Ребро $e_1$ с минимальным весом 8 принадлежит блоку $Б_1$ и блоку $Б_2$, соединяем их в один новый блок $Б_{1,2}$. Ребро $e_{15}$ с минимальным весом 6 принадлежит блоку $Б_3$ и блоку $Б_4$, соединяем их в один новый блок $Б_{3,4}$. Формируется новая запись однострочных структурных чисел

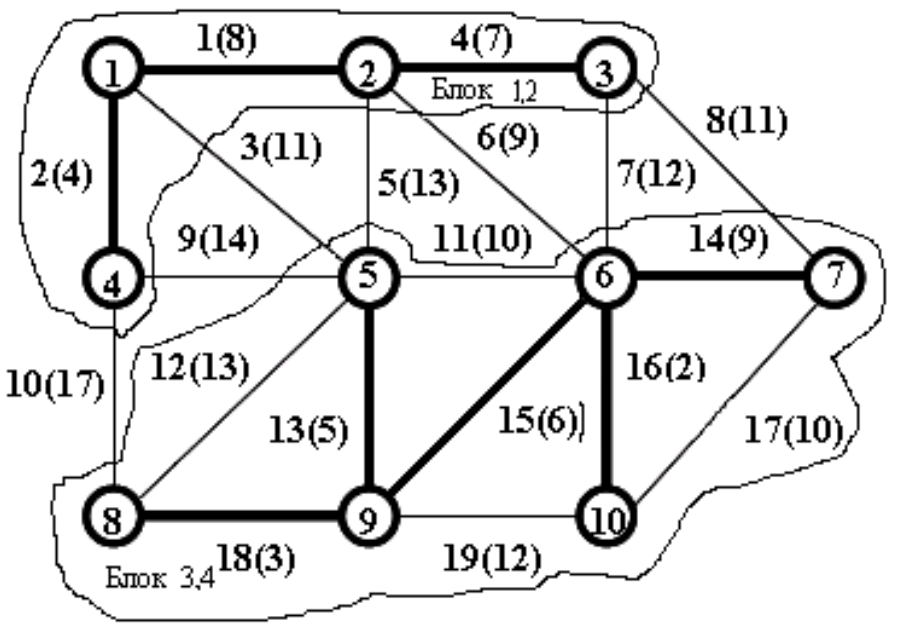

Рис. 2.13. Блоки графа G с минимальными весами ребер.

$Б_{1,2}$ = [6(9) 3(11) 8(11) 7(12) 5(13) 9(14) 10(17)];

$Б_{3,4}$ = [6(9) 3(11) 8(11) 7(12) 5(13) 9(14) 10(17)].

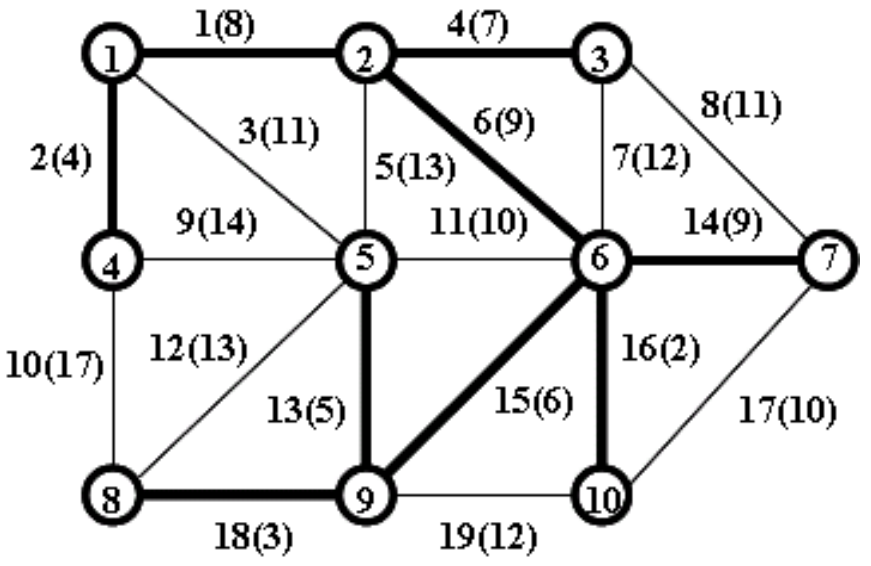

Рис. 2.14. Дерево графа G с минимальным суммарным весом ребер.

**шаг 4**: И наконец происходит объединение блоков $Б_{1,2}$ и $Б_{3,4}$ по ребру $e_6$. В результате получено дерево T = {1,2,4,6,13,14,15,16,18] с минимальным суммарным весом

= 8 + 4 + 7 + 9 + 5 + 9 + 6 + 2 + 3 = 53

Анализируя сказанное, можно сделать следующий вывод:

а) жадный алгоритм рассматривает только первые элементы однострочных структурных чисел, не производя их умножения;

б) неугодные элементы (то есть элементы стоящие не на первых местах) автоматически и последовательно исключаются из рассмотрения.

Например, в алгоритме Краскала или Прима для примера 2.14, ребра $e_{12}(13)$ и $e_{17}(10)$ на шаге 2, ребро $e_{11}(10)$ и $e_{19}(12)$ на шаге 3, и ребра $e_3(11)$, $e_5(13)$, $e_7(12)$, $e_8(11)$, $e_9(14)$, $e_{10}(17)$ на шаге 4.

## 2.12. Алгоритм «прыгающая строка»

Поставим следующий вопрос: возможно ли, получить подмножество баз матроида, не производя умножение структурных чисел, а рассматривая только первые элементы в однострочных структурных числах. Рассмотрим следующий алгоритм типа «прыгающая строка».

***Пример 2.15.*** Предположим, что нам нужно выделить все базы матроида изометрических циклов для графа, представленного на рис. 2.9 (см. пример 2.10).

Множество единичных циклов для данного графа:

$c_1 = \{e_1,e_2,e_5\}$, $c_2 = \{e_1,e_4,e_6\}$, $c_3 = \{e_2,e_3,e_7\}$, $c_4 = \{e_2,e_4,e_8\}$,

$c_5 = \{e_3,e_4,e_9\}$, $c_6 = \{e_7,e_8,e_9\}$, $c_7 = \{e_5,e_6,e_8\}$.

В примере 2.10 мы рассмотрели умножение однострочных структурных чисел $u_5 \times u_6 \times u_7 \times u_8 \times u_9$ по правилам умножения алгебра структурных чисел и нашли структурное число $W = [c_1\ c_7] \times [c_2\ c_7] \times [c_3\ c_6] \times [c_4\ c_6\ c_7] \times [c_5\ c_6]$. Рассмотрим несколько иной алгоритм вычисления структурного числа W под названием «прыгающая строка».

В данном алгоритме последовательно формируется строка-счетчик указывающая на номер последовательного расположения элемента в однострочном структурном числе. Например, запись строки-счетчика (1,1,1,2,1) указывает на 1-ый элемент в 1-ом однострочном структурном числе, затем на 1-ый элемент во 2-ом однострочном структурном числе, затем на 1-ый элемент во 3-ем однострочном структурном числе, затем на 2-ой элемент в 4-ом однострочном структурном числе, затем на 1-ый элемент в 5-ом однострочном структурном числе. Таким образом, выделено следующее подмножество изометрических циклов, записанное в порядке выбора $<c_1,c_2,c_3,c_6,c_5>$. Данный алгоритм последовательно перебирает всевозможные комбинации, не исключая случай повторяющихся циклов.

| | | | | | | | | | | | | |
|---|---|---|---|---|---|---|---|---|---|---|---|---|
| **Номера столбцов** | 1 | 2 | 3 | 4 | 5 | | **Номера циклов в элементе струк-турного числа** | | | | | |
| **Количество номеров циклов в однострочном структурном числ** | 2 | 2 | 2 | 3 | 2 | | | | | | | |
| **Первоначальное положение счетчика, характерно тем, что присутствуют только первые номера в однострочных структурных числах** | 1 | 1 | 1 | 1 | 1 | → | 1 | 2 | 3 | 4 | 5 | **Запись** |
| **В записи нет повторяющихся номеров +1 в последнюю ячейку** | | | | | +1 | | | | | | | |
| **Новая запись счетчика** | 1 | 1 | 1 | 1 | 2 | → | 1 | 2 | 3 | 4 | 6 | **Запись** |
| **В записи нет повторяющихся номеров +1 в последнюю ячейку** | | | | | +1 | | | | | | | |
| **Новая запись счетчика** | 1 | 1 | 1 | 1 | 3 | | | | | | | |
| **В 5-ой ячейки счетчика число 3 больше чем количество элементов в 5-ом однострочном структурном числе, +1 в 4-ую ячейку и все последующие = 1** | | | | +1 | | | | | | | | |
| **Новая запись счетчика** | 1 | 1 | 1 | 2 | 1 | → | 1 | 2 | 3 | 6 | 5 | Запись |
| **В записи нет повторяющихся номеров +1 в последнюю ячейку** | | | | | +1 | | | | | | | |
| **Новая запись счетчика** | 1 | 1 | 1 | 2 | 2 | → | 1 | 2 | 3 | 6 | 6 | |
| **В записи есть повторяющиеся элементы, в ячейку счетчика соответствующую последнему** | | | | | | | | | | | | |

| | | | | | | |
|---|---|---|---|---|---|---|
| **повторяющемуся элементу +1** | | | | | +1 | |
| **Новая запись счетчика** | 1 | 1 | 1 | 2 | 3 | |
| **В 5-ой ячейки счетчика число 3 больше чем количество элементов в 5-ом однострочном структурном числе, +1 в 4-ую ячейку и все последующие = 1** | | | | +1 | | |
| **Новая запись счетчика** | 1 | 1 | 1 | 3 | 1 | → 1 2 3 7 5 Запись |
| **В записи нет повторяющихся номеров +1 в последнюю ячейку** | | | | | +1 | |
| **Новая запись счетчика** | 1 | 1 | 1 | 3 | 2 | → 1 2 3 7 6 Запись |
| **В записи нет повторяющихся номеров +1 в последнюю ячейку** | | | | | +1 | |
| **Новая запись счетчика** | 1 | 1 | 1 | 3 | 3 | |
| **В 5-ой ячейки счетчика число 3 больше чем количество элементов в 5-ом однострочном структурном числе, +1 в 4-ую ячейку и все последующие = 1** | | | | +1 | | |
| **Новая запись счетчика** | 1 | 1 | 1 | 4 | 1 | |
| **В 4-ой ячейки счетчика число 4 больше чем количество элементов в 4-ом однострочном структурном числе, +1 в 3-ью ячейку и все последующие = 1** | | | +1 | | | |
| **Новая запись счетчика** | 1 | 1 | 2 | 1 | 1 | → 1 2 6 4 5 Запись |
| **В записи нет повторяющихся номеров +1 в последнюю ячейку** | | | | | +1 | |
| **Новая запись счетчика** | 1 | 1 | 2 | 1 | 2 | → 1 2 6 4 6 |
| **В записи есть повторяющиеся элементы, в ячейку счетчика соответствующую последнему повторяющемуся элементу +1** | | | | | +1 | |
| **Новая запись счетчика** | 1 | 1 | 2 | 1 | 3 | |
| **В 5-ой ячейки счетчика число 3 больше чем количество элементов в 5-ом однострочном структурном числе, +1 в 4-ую ячейку и все последующие = 1** | | | | +1 | | |
| **Новая запись счетчика** | 1 | 1 | 2 | 2 | 1 | → 1 2 6 6 |
| **В записи есть повторяющиеся элементы, в ячейку счетчика соответствующую последнему повторяющемуся элементу +1** | | | | +1 | | |
| **Новая запись счетчика** | 1 | 1 | 2 | 3 | 1 | → 1 2 6 7 5 Запись |
| **В записи нет повторяющихся номеров +1 в последнюю ячейку** | | | | | +1 | |
| **Новая запись счетчика** | 1 | 1 | 2 | 3 | 2 | → 1 2 6 7 6 |
| **В записи есть повторяющиеся элементы, в ячейку счетчика** | | | | | | |
| **соответствующую последнему повторяющемуся элементу +1** | | | | | +1 | |
| **Новая запись счетчика** | 1 | 1 | 2 | 3 | 3 | |
| **В 5-ой ячейки счетчика число 3 больше чем количество элементов в 5-ом однострочном структурном числе, +1 в 4-ую ячейку и все последующие = 1** | | | | +1 | | |
| **Новая запись счетчика** | 1 | 1 | 2 | 4 | 1 | |
| **В 4-ой ячейки счетчика число 4 больше чем количество элементов в 4-ом однострочном структурном числе, +1 в 3-ью ячейку и все последующие = 1** | | | +1 | | | |
| **Новая запись счетчика** | 1 | 1 | 3 | 1 | 1 | |
| **В 3-ей ячейки счетчика число 3 больше чем количество элементов в 3-ем однострочном структурном числе, +1 в 2-ую ячейку и все последующие = 1** | | +1 | | | | |
| **Новая запись счетчика** | 1 | 2 | 1 | 1 | 1 | → 1 7 3 4 5 Запись |
| **В записи нет повторяющихся номеров +1 в последнюю ячейку** | | | | | +1 | |
| **Новая запись счетчика** | 1 | 2 | 1 | 1 | 2 | → 1 7 3 4 6 Запись |
| **В записи нет повторяющихся номеров +1 в последнюю ячейку** | | | | | +1 | |
| **Новая запись счетчика** | 1 | 2 | 1 | 1 | 3 | |
| **В 5-ой ячейки счетчика число 3 больше чем количество элементов в 5-ом однострочном структурном числе, +1 в 4-ую ячейку и все последующие = 1** | | | | +1 | | |
| **Новая запись счетчика** | 1 | 2 | 1 | 2 | 1 | → 1 7 3 6 5 Запись |
| **В записи нет повторяющихся номеров +1 в последнюю ячейку** | | | | | +1 | |
| **Новая запись счетчика** | 1 | 2 | 1 | 2 | 2 | → 1 7 3 6 6 |
| **В записи есть повторяющиеся элементы, в ячейку счетчика соответствующую последнему повторяющемуся элементу +1** | | | | | +1 | |
| **Новая запись счетчика** | 1 | 2 | 1 | 2 | 3 | |
| **В 5-ой ячейки счетчика число 3 больше чем количество элементов в 5-ом однострочном структурном числе, +1 в 4-ую ячейку и все последующие = 1** | | | | +1 | | |
| **Новая запись счетчика** | 1 | 2 | 1 | 3 | 1 | → 1 7 3 7 |
| **В записи есть повторяющиеся элементы, в ячейку счетчика соответствующую последнему повторяющемуся элементу +1** | | | | +1 | | |
| **Новая запись счетчика** | 1 | 2 | 1 | 4 | 1 | |
| **В 4-ой ячейки счетчика число 4 больше чем количество элементов в 4-ом однострочном структурном числе, +1 в 3-ью ячейку и все последующие = 1** | | | +1 | | | |
| **Новая запись счетчика** | 1 | 2 | 2 | 1 | 1 | → 1 7 6 4 5 Запись |
| **В записи нет повторяющихся номеров +1 в последнюю ячейку** | | | | | +1 | |
| **Новая запись счетчика** | 1 | 2 | 2 | 1 | 2 | → 1 7 6 4 6 |
| **В записи есть повторяющиеся элементы, в ячейку счетчика соответствующую последнему повторяющемуся элементу +1** | | | | | +1 | |
| **Новая запись счетчика** | 1 | 2 | 2 | 1 | 3 | |
| **В 5-ой ячейки счетчика число 3 больше чем количество элементов в 5-ом однострочном структурном числе, +1 в 4-ую ячейку и все последующие = 1** | | | | +1 | | |
| **Новая запись счетчика** | 1 | 2 | 2 | 2 | 1 | → 1 7 6 6 |
| **В записи есть повторяющиеся элементы, в ячейку счетчика соответствующую последнему повторяющемуся элементу +1 и последующие = 1** | | | | +1 | | |
| **Новая запись счетчика** | 1 | 2 | 2 | 3 | 1 | → 1 7 6 7 |
| **В записи есть повторяющиеся элементы, в ячейку счетчика соответствующую последнему повторяющемуся элементу +1 и последующие = 1** | | | | +1 | | |
| **Новая запись счетчика** | 1 | 2 | 2 | 4 | 1 | |
| **В 4-ой ячейки счетчика число 4 больше чем количество элементов в 4-ом однострочном структурном числе, +1 в 3-ью ячейку и все последующие = 1** | | | +1 | | | |
| **Новая запись счетчика** | 1 | 2 | 3 | 1 | 1 | |
| **В 3-ей ячейки счетчика число 3 больше чем количество элементов в 3-ом однострочном структурном числе, +1 в 2-ую ячейку и все последующие = 1** | | +1 | | | | |
| **Новая запись счетчика** | 1 | 3 | 1 | 1 | 1 | |
| **Во 2-ой ячейки счетчика число 3 больше чем количество элементов в 2-ом однострочном структурном числе, +1 в 1-ую ячейку и все последующие = 1** | | | | | | |
| **Новая запись счетчика** | 2 | 1 | 1 | 1 | 1 | → 7 2 3 4 5 Запись |
| **В записи нет повторяющихся номеров +1 в последнюю ячейку** | | | | | +1 | |
| **Новая запись счетчика** | 2 | 1 | 1 | 1 | 2 | → 7 2 3 4 6 Запись |
| **В записи нет повторяющихся номеров +1 в последнюю ячейку** | | | | | +1 | |
| **Новая запись счетчика** | 2 | 1 | 1 | 1 | 3 | |
| **В 5-ой ячейки счетчика число 3 больше чем количество элементов в 5-ом однострочном структурном числе, +1 в 4-ую ячейку и все последующие = 1** | | | | +1 | | |
| **Новая запись счетчика** | 2 | 1 | 1 | 2 | 1 | → 7 2 3 6 5 Запись |
| **В записи нет повторяющихся номеров +1 в последнюю ячейку** | | | | | +1 | |
| **Новая запись счетчика** | 2 | 1 | 1 | 2 | 2 | → 7 2 3 6 6 |
| **В записи есть повторяющиеся элементы, в ячейку счетчика соответствующую последнему повторяющемуся элементу +1** | | | | | +1 | |
| **Новая запись счетчика** | 2 | 1 | 1 | 2 | 3 | |
| **В 5-ой ячейки счетчика число 3 больше чем количество элементов в 5-ом однострочном структурном числе, +1 в 4-ую ячейку и все последующие = 1** | | | | +1 | | |
| **Новая запись счетчика** | 2 | 1 | 1 | 3 | 1 | → 7 2 3 7 |
| **В записи есть повторяющиеся элементы, в ячейку счетчика соответствующую последнему повторяющемуся элементу +1 и последующие = 1** | | | | +1 | | |
| **Новая запись счетчика** | 2 | 1 | 1 | 4 | 1 | |
| **В 4-ой ячейки счетчика число 4 больше чем количество элементов в 4-ом однострочном структурном числе, +1 в 3-ью ячейку и все последующие = 1** | | | +1 | | | |
| **Новая запись счетчика** | 2 | 1 | 2 | 1 | 1 | → 7 2 6 4 5 Запись |
| **В записи нет повторяющихся номеров +1 в последнюю ячейку** | | | | | +1 | |
| **Новая запись счетчика** | 2 | 1 | 2 | 1 | 2 | → 7 2 6 4 6 Запись |
| **В записи нет повторяющихся номеров +1 в последнюю ячейку** | | | | | +1 | |
| **Новая запись счетчика** | 2 | 1 | 2 | 1 | 3 | |
| **В 5-ой ячейки счетчика число 3 больше чем количество элементов в 5-ом однострочном структурном числе, +1 в 4-ую ячейку и все последующие = 1** | | | | +1 | | |
| **Новая запись счетчика** | 2 | 1 | 2 | 2 | 1 | → 7 2 6 6 |
| **В записи есть повторяющиеся элементы, в ячейку счетчика** | | | | | | |
| **соответствующую последнему повторяющемуся элементу +1 и последующие = 1** | | | | +1 | | |
| **Новая запись счетчика** | 2 | 1 | 2 | 3 | 1 | → 7 2 6 7 |
| **В записи есть повторяющиеся элементы, в ячейку счетчика соответствующую последнему повторяющемуся элементу +1 и последующие = 1** | | | | +1 | | |
| **Новая запись счетчика** | 2 | 1 | 2 | 4 | 1 | |
| **В 4-ой ячейки счетчика число 4 больше чем количество элементов в 4-ом однострочном структурном числе, +1 в 3-ью ячейку и все последующие = 1** | | | +1 | | | |
| **Новая запись счетчика** | 2 | 1 | 3 | 1 | 1 | |
| **В 3-ей ячейки счетчика число 3 больше чем количество элементов в 3-ом однострочном структурном числе, +1 в 2-ую ячейку и все последующие = 1** | | +1 | | | | |
| **Новая запись счетчика** | 2 | 2 | 1 | 1 | 1 | → 7 7 |
| **В записи есть повторяющиеся элементы, в ячейку счетчика соответствующую последнему повторяющемуся элементу +1 и последующие = 1** | | +1 | | | | |
| **Новая запись счетчика** | 2 | 3 | 1 | 1 | 1 | |
| **В 2-ой ячейки счетчика число 3 больше чем количество элементов в 2-ом однострочном структурном числе, +1 в 1-ую ячейку и все последующие = 1** | +1 | | | | | |
| **Новая запись счетчика** | 3 | 1 | 1 | 1 | 1 | |
| **В 1-ой ячейки счетчика число 3 больше чем количество элементов в 1-ом однострочном структурном числе. Конец работы алгоритма.** | | | | | | |

Записи в этом методе представляют собой элементы структурного числа W. Уникальным свойством этого алгоритма, является возможность нахождения только одного структурного элемента из множества, а именно первого столбца.

Будем называть элемент структурного числа *первым*, если он характеризуется минимальным набором номеров в соответствующих однострочных структурных числах расположенных в определенном порядке. Например, в нашем примере первый структурный элемент – это элемент структурного числа:

$u_5$: 1-ый номер записи в однострочном структурном числе $c_1$

$u_6$: 1-ый номер записи в однострочном структурном числе $c_2$

$u_7$: 1-ый номер записи в однострочном структурном числе $c_3$

$u_8$: 1-ый номер записи в однострочном структурном числе $c_4$

$u_9$: 1-ый номер записи в однострочном структурном числе $c_5$

## 2.13. Задача о назначении. Венгерский алгоритм

Рассмотрим связь методов алгебры структурных чисел и решение задачи о назначении на следующем примере [63].

***Пример 2.16***. Рассмотрим задачу о назначении, представленную следующей матрицей стоимости для ребер двудольного графа []:

| | $e_1$ | $e_2$ | $e_3$ | $e_4$ | $e_5$ |
|---|---|---|---|---|---|
| $v_1$ | 7 | 2 | 1 | 9 | 4 |
| $v_2$ | 9 | 6 | 9 | 5 | 5 |
| $v_3$ | 3 | 8 | 3 | 1 | 8 |
| $v_4$ | 7 | 9 | 4 | 2 | 2 |
| $v_5$ | 8 | 4 | 7 | 4 | 8 |

Необходимо выбрать элементы матрицы так чтобы их сумма была минимальной, и при этом из каждой строки и столбца был выбран только один элемент.

Запишем однострочные структурные числа для каждой строки матрицы, расположив номера столбцов с их весами в порядке не убывания весов:

$v_1$: [$e_3$(1) $e_2$(2) $e_5$(4) $e_1$(7) $e_4$(9)];

$v_2$: [$e_4$(5) $e_5$(5) $e_2$(6) $e_1$(9) $e_3$(9)];

$v_3$: [$e_4$(1) $e_1$(3) $e_3$(3) $e_2$(8) $e_5$(8)];

$v_4$: [$e_4$(2) $e_5$(2) $e_3$(4) $e_1$(7) $e_2$(9)];

$v_5$: [$e_2$(4) $e_4$(4) $e_3$(7) $e_1$(8) $e_5$(8)].

Расставим однострочные структурные числа в последовательности не убывания весов в каждом однострочном структурном числе:

$v_1$: [$e_3$(1) $e_2$(2) $e_5$(4) $e_1$(7) $e_4$(9)];

$v_3$: [$e_4$(1) $e_1$(3) $e_3$(3) $e_2$(8) $e_5$(8)];

$v_4$: [$e_4$(2) $e_5$(2) $e_3$(4) $e_1$(7) $e_2$(9)];

$v_5$: [$e_2$(4) $e_4$(4) $e_3$(7) $e_1$(8) $e_5$(8)];

$v_2$: [$e_4$(5) $e_5$(5) $e_2$(6) $e_1$(9) $e_3$(9)].

Алгоритмом «прыгающая строка» выделим первый элемент структурного числа.

$v_1$: [**$\mathbf{e_3}$(1)** $e_2$(2) $e_5$(4) $e_1$(7) $e_4$(9)];

$v_3$: [**$\mathbf{e_4}$(1)** $e_1$(3) $e_3$(3) $e_2$(8) $e_5$(8)];

$v_4$: [$e_4(2)$ $\mathbf{e_5(2)}$ $e_3(4)$ $e_1(7)$ $e_2(9)$];

$v_5$: [$\mathbf{e_2(4)}$ $e_4(4)$ $e_3(7)$ $e_1(8)$ $e_5(8)$];

$v_2$: [$e_4(5)$ $e_5(5)$ $e_2(6)$ $\mathbf{e_1(9)}$ $e_3(9)$].

Сумма весов выбранных элементов равна 17. Рассмотрим возможность получения меньшей суммы весов, если удалить максимальный элемент из выбранного подмножества. Выбираем в данном элементе структурного числа элемент матрицы с максимальным весом, это $e_1(9)$ в строке $v_2$ и удаляем его из рассмотрения. Алгоритмом «прыгающая строка» выделим следующий первый элемент структурного числа.

$v_1$: [$\mathbf{e_3(1)}$ $e_2(2)$ $e_5(4)$ $e_1(7)$ $e_4(9)$];

$v_3$: [$\mathbf{e_4(1)}$ $e_1(3)$ $e_3(3)$ $e_2(8)$ $e_5(8)$];

$v_4$: [$e_4(2)$ $\mathbf{e_5(2)}$ $e_3(4)$ $e_1(7)$ $e_2(9)$];

$v_5$: [$e_2(4)$ $e_4(4)$ $e_3(7)$ $\mathbf{e_1(8)}$ $e_5(8)$];

$v_2$: [$e_4(5)$ $e_5(5)$ $\mathbf{e_2(6)}$ $e_3(9)$ ].

Сумма весов выбранных элементов равна 18. Будем осуществлять дальнейшее построение процесса. Выбираем в данном элементе структурного числа элемент матрицы с максимальным весом, это $e_1(8)$ в строке $v_5$ и удаляем его из рассмотрения. Алгоритмом «прыгающая строка» выделим следующий первый элемент структурного числа.

$v_1$: [$\mathbf{e_3(1)}$ $e_2(2)$ $e_5(4)$ $e_1(7)$ $e_4(9)$];

$v_3$: [$\mathbf{e_4(1)}$ $e_1(3)$ $e_3(3)$ $e_2(8)$ $e_5(8)$];

$v_4$: [$e_4(2)$ $e_5(2)$ $e_3(4)$ $\mathbf{e_1(7)}$ $e_2(9)$];

$v_5$: [$\mathbf{e_2(4)}$ $e_4(4)$ $e_3(7)$ $e_5(8)$ ];

$v_2$: [$e_4(5)$ $\mathbf{e_5(5)}$ $e_2(6)$ $e_3(9)$ ].

Сумма весов выбранных элементов снова равна 18. Выбираем в данном элементе структурного числа элемент матрицы с максимальным весом, это $e_1(7)$ в строке $v_4$ и удаляем его из рассмотрения. Алгоритмом «прыгающая строка» выделим следующий первый элемент структурного числа.

$v_1$: [$\mathbf{e_3(1)}$ $e_2(2)$ $e_5(4)$ $e_1(7)$ $e_4(9)$];

$v_3$: [$e_4(1)$ $\mathbf{e_1(3)}$ $e_3(3)$ $e_2(8)$ $e_5(8)$];

$v_4$: [$\mathbf{e_4(2)}$ $e_5(2)$ $e_3(4)$ $e_2(9)$ ];

$v_5$: [$\mathbf{e_2(4)}$ $e_4(4)$ $e_3(7)$ $e_5(8)$ ];

$v_2$: [$e_4(5)$ $\mathbf{e_5(5)}$ $e_2(6)$ $e_3(9)$ ].

Сумма весов выбранных элементов равна 15. Рассмотрим возможность построения другого решения. Выбираем в данном элементе структурного числа элемент матрицы с максимальным весом, это $e_5$(5) в строке $v_2$ и удаляем его из рассмотрения. Алгоритмом «прыгающая строка» выделим первый элемент структурного числа:

$v_1$: [**$e_3$(1)** $e_2$(2) $e_5$(4) $e_1$(7) $e_4$(9)];

$v_3$: [$e_4$(1) **$e_1$(3)** $e_3$(3) $e_2$(8) $e_5$(8)];

$v_4$: [$e_4$(2) **$e_5$(2)** $e_3$(4) $e_2$(9) ];

$v_5$: [**$e_2$(4)** $e_4$(4) $e_3$(7) $e_5$(8) ];

$v_2$: [**$e_4$(5)** $e_2$(6) $e_3$(9) ].

Сумма весов выбранных элементов равна 15. Продолжаем поиск, выбираем в данном элементе структурного числа, элемент матрицы с максимальным весом, это $e_4$(5) в строке $v_2$ и удаляем его из рассмотрения. Алгоритмом «прыгающая строка» выделим первый элемент структурного числа:

$v_1$: [**$e_3$(1)** $e_2$(2) $e_5$(4) $e_1$(7) $e_4$(9)];

$v_3$: [$e_4$(1) **$e_1$(3)** $e_3$(3) $e_2$(8) $e_5$(8)];

$v_4$: [$e_4$(2) **$e_5$(2)** $e_3$(4) $e_2$(9) ];

$v_5$: [e2(4) **$e_4$(4)** $e_3$(7) $e_5$(8) ];

$v_2$: [**$e_2$(6)** $e_3$(9) ].

Сумма весов выбранных элементов равна 16. Приостанавливаем работу алгоритма, так как дальнейший поиск приведет только к увеличению суммы весов. А последующее удаление $e_3$(9) к появлению пустой строки $v_2$.

Модифицируем данный подход, после выбора первого элемента структурного числа в первый раз и определения элемент матрицы с максимальным весом, это $e_1$(9) во строке $v_2$:

$v_1$: [**$e_3$(1)** $e_2$(2) $e_5$(4) $e_1$(7) $e_4$(9)];

$v_3$: [**$e_4$(1)** $e_1$(3) $e_3$(3) $e_2$(8) $e_5$(8)];

$v_4$: [$e_4$(2) **$e_5$(2)** $e_3$(4) $e_1$(7) $e_2$(9)];

$v_5$: [**$e_2$(4)** $e_4$(4) $e_3$(7) $e_1$(8) $e_5$(8)];

$v_2$: [$e_4$(5) $e_5$(5) $e_2$(6) **$e_1$(9)** $e_3$(9)].

Удалим все элементы матрицы, имеющие вес 9 из рассмотрения. Алгоритмом «бегущая строка» выделим первый элемент структурного числа:

$v_1$: [**$e_3$(1)** $e_2$(2) $e_5$(4) $e_1$(7) ];

$v_3$: [**$e_4$(1)** $e_1$(3) $e_3$(3) $e_2$(8) $e_5$(8)];

$v_4$: [$e_4$(2) **$e_5$(2)** $e_3$(4) $e_1$(7) ];

$v_5$: [$e_2$(4) $e_4$(4) $e_3$(7) **$e_1$(8)** $e_5$(8)];

$v_2$: [$e_4$(5) $e_5$(5) **$e_2$(6)** ].

Сумма весов выбранных элементов равна 18. Выбираем в данном элементе структурного числа, элемент матрицы с максимальным весом, это $e_1(8)$ в строке $v_5$ и удаляем все элементы матрицы, имеющие вес равный 8 из рассмотрения. Алгоритмом «прыгающая строка» выделим первый элемент структурного числа.

$v_1$: [**$e_3(1)$** $e_2(2)$ $e_5(4)$ $e_1(7)$ ];

$v_3$: [**$e_4(1)$** $e_1(3)$ $e_3(3)$ ];

$v_4$: [$e_4(2)$ $e_5(2)$ $e_3(4)$ **$e_1(7)$** ];

$v_5$: [**$e_2(4)$** $e_4(4)$ $e_3(7)$ ];

$v_2$: [$e_4(5)$ **$e_5(5)$** $e_2(6)$ ].

Сумма весов выбранных элементов равна 18. Выбираем в данном элементе структурного числа элемент матрицы с максимальным весом, это $e_1(7)$ в строке $v_4$ и удаляем все элементы матрицы, имеющие вес равный 7 из рассмотрения. Алгоритмом «прыгающая строка» выделим первый элемент структурного числа.

$v_1$: [**$e_3(1)$** $e_2(2)$ $e_5(4)$ ];

$v_3$: [$e_4(1)$ **$e_1(3)$** $e_3(3)$ ];

$v_4$: [**$e_4(2)$** $e_5(2)$ $e_3(4)$ ];

$v_5$: [**$e_2(4)$** $e_4(4)$ ];

$v_2$: [$e_4(5)$ **$e_5(5)$** $e_2(6)$ ].

Сумма весов выбранных элементов равна 15. Если теперь удалить из рассмотрения все элементы матрицы с весом 5 и выше, то строка $v_2$ будет пустой.

Можно произвести и следующую модификацию метода. В исходной матрице выделить все минимальные элементы по столбцам, а затем по строкам:

| | $e_1$ | $e_2$ | $e_3$ | $e_4$ | $e_5$ |
|---|---|---|---|---|---|
| $v_1$ | 7 | **2** | **1** | 9 | 4 |
| $v_2$ | 9 | 6 | 9 | **5** | **5** |
| $v_3$ | **3** | 8 | 3 | **1** | 8 |
| $v_4$ | 7 | 9 | 4 | **2** | **2** |
| $v_5$ | 8 | **4** | 7 | **4** | 8 |

Построим однострочные структурные числа из выбранных элементов матрицы, и выберем элементы структурного числа с минимальным весом алгоритмом «прыгающая строка». Если паросочетания не произойдет, то будем добавлять другие минимальные элементы до образования паросочетания.

$v_1$: [$e_3(1)$ $e_2(2)$ ];

$v_3$: [$e_4(1)$ $e_1(3)$ ];

$v_4$: [$e_4(2)$ $e_5(2)$ ];

$v_5$: [$e_2(4)$ $e_4(4)$ ];

$v_2$: [$e_4(5)$ $e_5(5)$ ].

В данном случае это всего два элемента структурного числа с суммарным весом равным 15 (см. рис. 2.15):

| | | |
|---|---|---|
| $v_1$: | $e_3(1)$ | $e_3(1)$ |
| $v_3$: | $e_1(3)$ | $e_1(3)$ |
| $v_4$: | $e_4(2)$ | $e_5(2)$ |
| $v_5$: | $e_2(4)$ | $e_2(4)$ |
| $v_2$: | $e_5(5)$ | $e_4(5)$ |

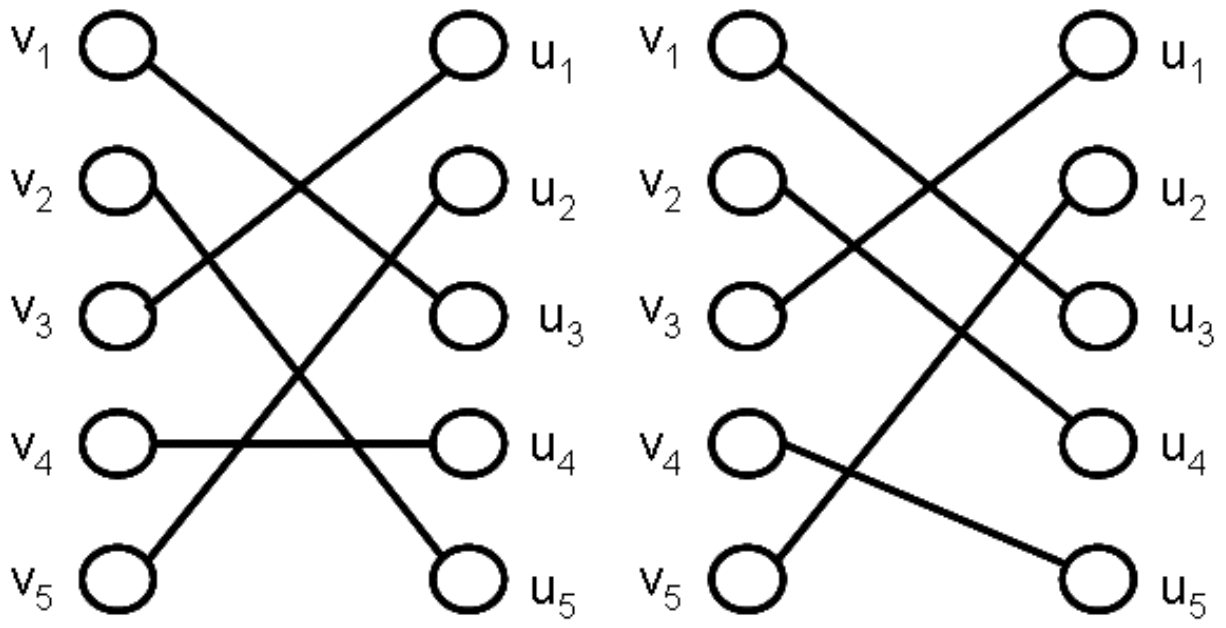

Рис. 2.15. Трасверсали для двух элементов структурного числа с минимальным весом.

Мы рассмотрели с вами подход, основанный на попытке просмотра всего множества трансверсалей как бы сверху, отсекая элементы, приводящие к увеличению суммарного веса.

Рассмотрим венгерский алгоритм, который подробно описан в литературе [], оттуда взят и пример. Основная идея этого метода заключена в том, что паросочетания строятся как бы снизу, рассматривая только перспективные проекты.

Выберем в каждом столбце нашей матрицы минимальный элемент и пометим его. Это соответствует следующему варианту паросочетаний (см. рис. 2.16)

| | | | | |
|---|---|---|---|---|
| 7 | **2** | **1** | 9 | 4 |
| 9 | 6 | 9 | 5 | 5 |
| **3** | 8 | 3 | **1** | 8 |
| 7 | 9 | 4 | 2 | **2** |
| 8 | 4 | 7 | 4 | 8 |

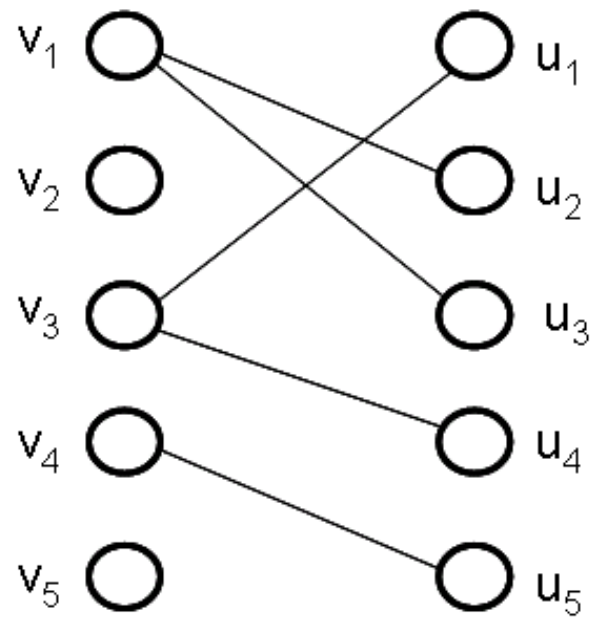

Рис. 2.16. Вариант паросочетания.

Произведем выборку выбранных элементов матрицы так чтобы, начиная с первого столбца выбранные элементы, находились в столбцах и в строках только по одному разу. Это соответствует следующему варианту паросочетаний и выбору трасверсалей (см. рис. 2.17).

| | | | | |
|---|---|---|---|---|
| 7 | **2** | **1** | 9 | 4 |
| 9 | 6 | 9 | 5 | 5 |
| **3** | 8 | 3 | **1** | 8 |
| 7 | 9 | 4 | 2 | **2** |
| 8 | 4 | 7 | 4 | 8 |

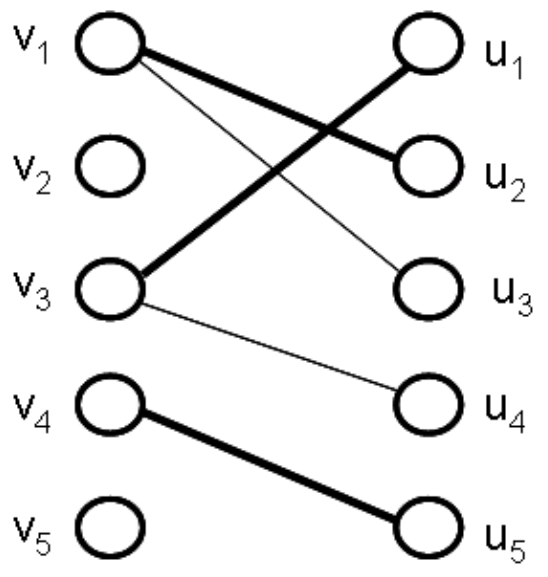

Рис. 2.17. Выбор трансверсалей.

Строим вектор невязки и вектор соседней вершины, рассчитывая параметры α и β для строк и столбцов.

| | $\alpha$ \ $\beta$ | 3 | 2 | 1 | 1 | 2 |
|---|---|---|---|---|---|---|
| | 0 | 7 | **2** | **1** | 9 | 4 |
| * | 0 | 9 | 6 | 9 | 5 | 5 |
| | 0 | **3** | 8 | 3 | **1** | 8 |
| | 0 | 7 | 9 | 4 | 2 | **2** |
| * | 0 | 8 | 4 | 7 | 4 | 8 |

| | | | | | |
|---|---|---|---|---|---|
| вектор невязки | 5 | 2 | 6 | 3 | 3 |
| вектор соседних вершин | 5 | 5 | 5 | 5 | 2 |

В результате расчета в паросочетание, добавляется новое ребро $v_5e_2$, и появляется новый чередующийся маршрут $v_5e_2v_1e_3$. Затем определяем новый набор трансверсалей (см. рис. 2.18).

| | | | | | |
|---|---|---|---|---|---|
| | 7 | **2** | 1 | 9 | 4 |
| | 9 | 6 | 9 | 5 | 5 |
| | **3** | 8 | 3 | 1 | 8 |
| | 7 | 9 | 4 | 2 | **2** |
| | 8 | **4** | 7 | 4 | 8 |

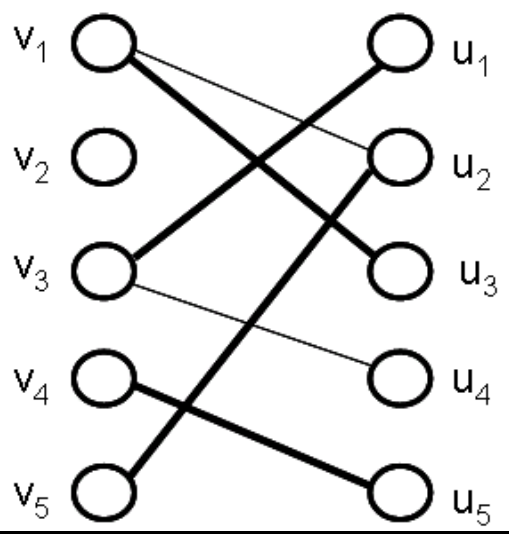

Рис. 2.18. Новый набор трансверсалей

Строим вектор невязки и вектор соседней вершины, рассчитывая параметры α и β для строк и столбцов

| | $\beta$ | 4 | 3 | 2 | 2 | 3 |
|---|---|---|---|---|---|---|
| | -1 | 7 | **2** | **1** | 9 | 4 |
| * | +1 | 9 | 6 | 9 | 5 | 5 |
| | -1 | **3** | 8 | 3 | **1** | 8 |

| | | | | | |
|---|---|---|---|---|---|
| -1 | 7 | 9 | 4 | 2 | **2** |
| +1 | 8 | **4** | 7 | 4 | 8 |
| $\alpha$ | | | | | |
| вектор невязки | 4 | 2 | 6 | 2 | 1 |
| вектор соседних вершин | 2 | 2 | 2 | 2 | 2 |

В результате расчета в паросочетание добавляется новое ребро ($v_2$,$e_5$). В данном паросочетании, появляется новый маршрут ($v_2$,$e_5$,$v_4$) который является перспективным для построения новых чередующихся маршрутов. Поэтому вновь строим вектор невязки и вектор соседней вершины, рассчитывая параметры $\alpha$ и $\beta$ для строк и столбцов.

| | | | | | | | |
|---|---|---|---|---|---|---|---|
| | | 4,5 | 3,5 | 2,5 | 2,5 | 3,5 | $\beta$ |
| | -1,5 | 7 | **2** | **1** | 9 | 4 | |
| * | +1,5 | 9 | 6 | 9 | 5 | **5** | |
| | -1,5 | **3** | 8 | 3 | **1** | 8 | |
| * | -1,5 | 7 | 9 | 4 | 2 | **2** | |
| | +0,5 | 8 | **4** | 7 | 4 | 8 | |
| | $\alpha$ | | | | | | |
| | вектор невязки | 3 | 1 | 3 | 1 | 0 | |
| | вектор соседних вершин | 2 | 2 | 4 | 2,4 | 2,4 | |

В результате расчета в паросочетание добавляется два новых ребра ($v_2$,$e_4$) и ($v_4$,$e_4$). В данном паросочетании появляются два новых чередующихся маршрута ($v_2$,$e_4$,$v_4$,$e_5$) и ($v_2$,$u_5$,$v_4$,$u_4$), строим новый набор трансверсалей, которые и являются решениями для нашего примера (см. рис. 2.19).

| | | | | |
|---|---|---|---|---|
| 7 | **2** | 1 | 9 | 4 |
| 9 | 6 | 9 | **5** | **5** |
| **3** | 8 | 3 | 1 | 8 |
| 7 | 9 | 4 | 2 | **2** |
| 8 | **4** | 7 | **4** | 8 |

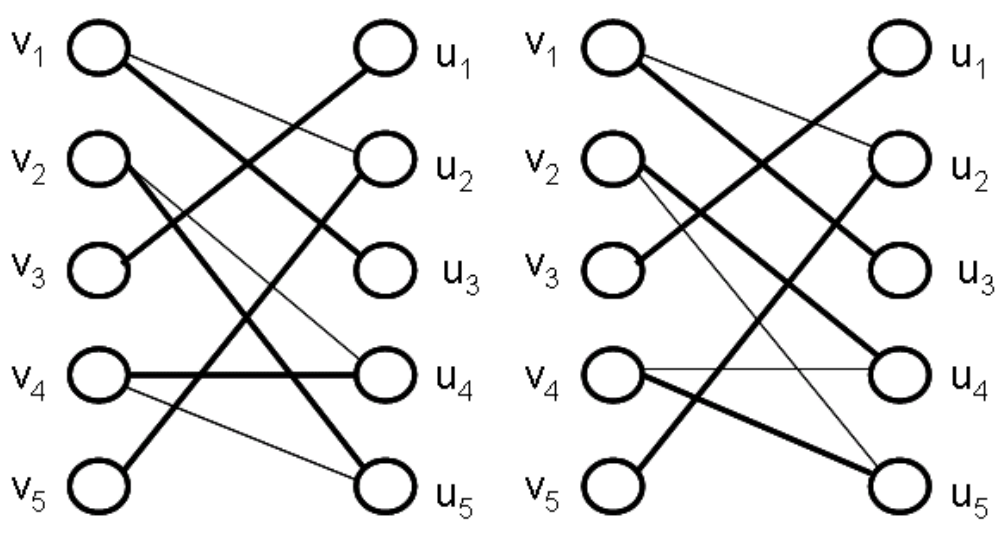

Рис. 2.19. Построение трансверсалей для нашего примера.

Таким образом, принцип отсечения, применяемый к вычислению структурного числа, приводит к построению венгерского алгоритма для решения задачи о назначении.

## 2.14. Программа «прыгающая строка»

```
program Gall19;

type
        TMasy = array[1..1000] of integer;
        TMass = array[1..4000] of integer;
var
        F1,F2,F3,F4 : text;
        i,ii,j,jj,K,K1,Np,Nv,Kzikl,M,MakLin : integer;
        Ziklo,Nr,Nh,KKK,K9,Priz,L,KK77,KM : integer;
        Masy: TMasy;
        Mass: TMass;
        Masi: TMass;
        MasyT: TMasy;
        MassT : TMass;
        MasMdop : TMasy;
        MasMy1 : TMasy;
        MasMs1 : TMass;
        MasMy2 : TMasy;
        MasMs2 : TMass;
        MasMy3 : TMasy;
        MasMs3 : TMass;
        MasMy4 : TMasy;
        MasMs4 : TMass;
        MasMy5 : TMasy;
        MasMs5 : TMass;
        MasMcg : TMasy;
        MasMcg1 : TMasy;
        Mass1 : TMass;
{*************************************************************}
{ процедура  предназначена для построения  структурного числа                }
{*************************************************************}
{*************************************************************}
procedure FormStryk(var Nh,KKK : integer;
                    var MasMcg : TMasy;
                    var MasMcg1 : TMasy;
                    var MasMy1 : TMasy;
                    var MasMs1 : TMass;
                    var MasMy2 : TMasy;
                    var MasMs2 : TMass;
                    var MasMy3 : TMasy;
                    var MasMs3 : TMass;
                    var MasMy4 : TMasy;
                    var MasMs4 : TMass;
                    var MasMy5 : TMasy;
                    var MasMs5 : TMass;
                    var Mass1 : TMass;
                    var MasyT : TMasy;
                    var MassT : TMass);
{Процедура формирования элементов структурного числа              }
{                                                   }
{ Nv - количество элементов структурного числа;                  }
{                                                   }
var i,j,max,k,kp,kip,kuku : integer;
label 1,2,3;
begin
{формирование массива количества элементов}
  for i:= 1 to Nh do MasMcg1[i]:=MasMy3[i+1]-MasMy3[i];
{определение максимального количества элементов }
  max:=0;
```

```
  for i:=1 to Nh do if MasMcg1[i]>max then max:=MasMcg1[i];
{ перестановка количества циклов по не убыванию в MasMy4,
соответсвующие индексы хранятся в массиве MasMs4}
  k:=0;
  for i:=1 to max do
  begin {1}
   for j:=1 to Nh do if MasMcg1[j]=i then
   begin {2}
    k:=k+1;
    MasMy4[k]:=i;
    MasMs4[k]:=j;
    Mass1[k]:=MasMcg[MasMs4[k]];
   end; {2}
  end; {1}
 {перезапись расположения 'элементов с учетом их количества}
  kp:=0;
  MasMy1[1]:=1;
  for i:=1 to Nh do
  begin {3}
   kip:=MasMs4[i];
   for j:=MasMy3[kip] to MasMy3[kip+1]-1 do
   begin {4}
    kp:=kp+1;
    MasMs1[kp]:= MasMs3[j];
    {writeln(F2,'i = ',i,' ','kp = ',kp,' ','MasMs1[kp] =',MasMs1[kp]); }
   end; {4}
  MasMy1[i+1]:=kp+1;
  end; {3}
  for i:=1 to Nh+1 do MasMy5[i]:=MasMy1[i];
  for i:=1 to MasMy1[Nh+1]-1 do MasMs5[i]:=MasMs1[i];
  {for i:=1 to Nh+1 do
         begin
            if i<>Nh+1 then write(F2,MasMy1[i],' ');
            if i=Nh+1 then writeln(F2,MasMy1[i]);
         end;
  for I:=1 to Nh do
        begin
         for j:=MasMy1[i] to MasMy1[i+1]-1 do
         begin
            if j<>MasMy1[i+1]-1 then write(F2,MasMs1[j],' ');
            if j=MasMy1[i+1]-1 then writeln(F2,MasMs1[j]);
         end;
        end; }
{алгоритм прыгающая строка}
{формирование массива счетчика MasMy2}
  for i:=1 to Nh do MasMy2[i]:=1;
  {for i:=1 to Nh do
         begin
            if i<>Nh then write(F2,MasMy2[i],' ');
            if i=Nh then writeln(F2,MasMy2[i]);
         end; }
{формирование структурного числа алгоритмом прыгающая строка}
  KKK:=0;
  MasyT[1]:=1;
2: for i:=1 to Nh do
   MasMs2[i]:=MasMs1[MasMy1[i]-1+MasMy2[i]];
  for i:=2 to Nh do
  begin  {5}
   for j:=1 to i-1 do
   begin {6}
    if MasMs2[i]=MasMs2[j] then
    begin {7}
     kuku:=i;
```

```
        MasMy2[kuku]:=MasMy2[kuku]+1;
        goto 1;
       end; {7}
      end; {6};
     end; {5}
   {запись элемента структурного числа}
     KKK:=KKK+1;
     for i:=1 to Nh do
        MassT[Nh*(KKK-1)+i]:=MasMs2[i];
     MasyT[KKK+1]:=MasyT[KKK]+Nh;
     MasMy2[Nh]:=MasMy2[Nh]+1;
   1: if MasMy2[1]>MasMy4[1] then goto 3;
     for i:=1 to Nh do
     begin {8}
      if MasMy2[i]>MasMy4[i] then
      begin {9}
       MasMy2[i-1]:=MasMy2[i-1]+1;
       for j:=i to Nh do MasMy2[j]:=1;
       goto 1;
      end; {9}
     end; {8}
     goto 2;
   3:  {for I:=1 to Nh do
       begin
        if i<>Nh then write(F2,MasMy1[i],' ');
        if i=Nh  then writeln(F2,MasMy1[i]);
       end; }
   end;{FormStryk}
   {*********************************************************************}
   label 1,2,3;
   begin
          assign(F1,'D:\Delphi1\DYD\10a.dyd');
          reset(F1);
          assign(F2,'D:\Delphi1\KUK\10a.kuk');
          Rewrite(F2);
          readln(F1,Nh);
          for I:=1 to Nh do
          begin  {1}
            if I<>Nh then read(F1,MasMcg[I]);
            if I=Nh then readln(F1,MasMcg[I]);
          end;{1}
          {for  i:=1 to Nh do
          begin
            if I<>Nh then write(F2,MasMcg[I],' ');
            if I=Nh then writeln(F2,MasMcg[I]);
          end;  }
          for I:=1 to Nh+1 do
          begin  {2}
            if I<>Nh+1 then read(F1,MasMy3[I]);
            if I=Nh+1 then readln(F1,MasMy3[I]);
          end;{2}
          for I:=1 to Nh do
          begin {3}
            for j:= MasMy3[i] to MasMy3[i+1]-1 do
            begin {4}
              if j<>MasMy3[i+1]-1 then read(F1,MasMs3[j]);
              if j=MasMy3[i+1]-1 then readln(F1,MasMs3[j]);
            end;  {4}
          end; {3}
          close (F1);
          { Создаём новый файл и открываем его в режиме "для чтения и записи"}
          FormStryk(Nh,KKK,MasMcg,MasMcg1,MasMy1,MasMs1,MasMy2,MasMs2,
            MasMy3,MasMs3,MasMy4,MasMs4,MasMy5,MasMs5,Mass1,MasyT,MassT);
```

```
        writeln(F2,'Длина элемента структурного числа = ',Nh,'.');
        writeln(F2,'Количество элементов структурного числа = ',KKK,'.');
        writeln(F2,'Однострочные структурные числа:');
        {for I:=1 to Nh+1 do
        begin
          if i<>Nh+1 then write(F2,MasMy5[i],' ');
          if i=Nh+1  then writeln(F2,MasMy5[i]);
        end; }
        for I:=1 to Nh do
        begin {17}
          write(F2,'элементы однострочного структурного числа ',Mass1[i],': [');
          for j:=MasMy5[i] to MasMy5[i+1]-1 do
          begin {18}
            if j<>MasMy5[i+1]-1 then write(F2,MasMs5[j],',');
            if j=MasMy5[i+1]-1  then writeln(F2,MasMs5[j],'];');
          end; {18}
        end; {17}
        for I:=1 to KKK do
        begin {17}
          write(F2,i,' - элемент структурного числа = {');
          for j:=MasyT[i] to MasyT[i+1]-1 do
          begin {18}
            if j<>MasyT[i+1]-1 then write(F2,MassT[j],',');
            if (j=MasyT[i+1]-1) and (i=KKK)
              then writeln(F2,MassT[j],'}.');
            if (j=MasyT[i+1]-1) and (i<>KKK)
              then writeln(F2,MassT[j],'};');
          end; {18}
        end; {17}
        close (F2);
        writeln('Конец расчета');
end.
```

## 2.15. Файлы программы Gall19

### Входной файл 10a.dyd

```
10                                  {количество столбцов в структурном числе}
1 2 3 4 5 6 7 8 9 10                {порядок расположения строк}
1 3 5 7 10 12 14 16 18 21 23        {массив указателей}
5 7                                 {1 строка элементов}
6 10                                {2 строка элементов}
8 9                                 {3 строка элементов}
1 6 7                               {4 строка элементов}
9 10                                {5 строка элементов}
4 6                                 {6 строка элементов}
2 4                                 {7 строка элементов}
6 9                                 {8 строка элементов}
2 3 6                               {9 строка элементов}
1 3                                 {10 строка элементов}
```

### Выходной файл 10a.kuk

```
Длина элемента структурного числа = 10.
Количество элементов структурного числа = 2.
Однострочные структурные числа:
элементы однострочного структурного числа 1: [5,7];
```

элементы однострочного структурного числа 2: [6,10];
элементы однострочного структурного числа 3: [8,9];
элементы однострочного структурного числа 5: [9,10];
элементы однострочного структурного числа 6: [4,6];
элементы однострочного структурного числа 7: [2,4];
элементы однострочного структурного числа 8: [6,9];
элементы однострочного структурного числа 10: [1,3];
элементы однострочного структурного числа 4: [1,6,7];
элементы однострочного структурного числа 9: [2,3,6];
1 - элемент структурного числа = {5,6,8,10,4,2,9,1,7,3};
2 - элемент структурного числа = {5,10,8,9,4,2,6,1,7,3}.

## Комментарии

Математический аппарат алгебры структурных чисел позволяет наглядно и просто описывать матроидные структуры. Он незаменим при решении многих задач дискретной оптимизации. С его помощью, осуществляется выбор линейно независимых элементов системы. Для выделения независимой системы элементов представлен алгоритм «прыгающая строка» и программа на языке Паскаль.

# Глава 3. ТОПОЛОГИЧЕСКИЙ РИСУНОК ПЛОСКОГО ГРАФА

## 3.1. Идентификация топологического рисунка графа на плоскости

Пусть G = (V,E) – несеперабельный граф с пронумерованным множеством ребер $E = \{e_1, e_2,...,e_m\}$ и $V = \{v_1,v_2,...,v_n\}$ вершин, причем card V = n и card E = m. Обычно граф G представляется матрицей инциденций или матрицей смежностей. Графически граф может быть представлен диаграммой, в которой вершина изображена точкой или кружком, а ребро – отрезком линии, соединяющим вершины [84].

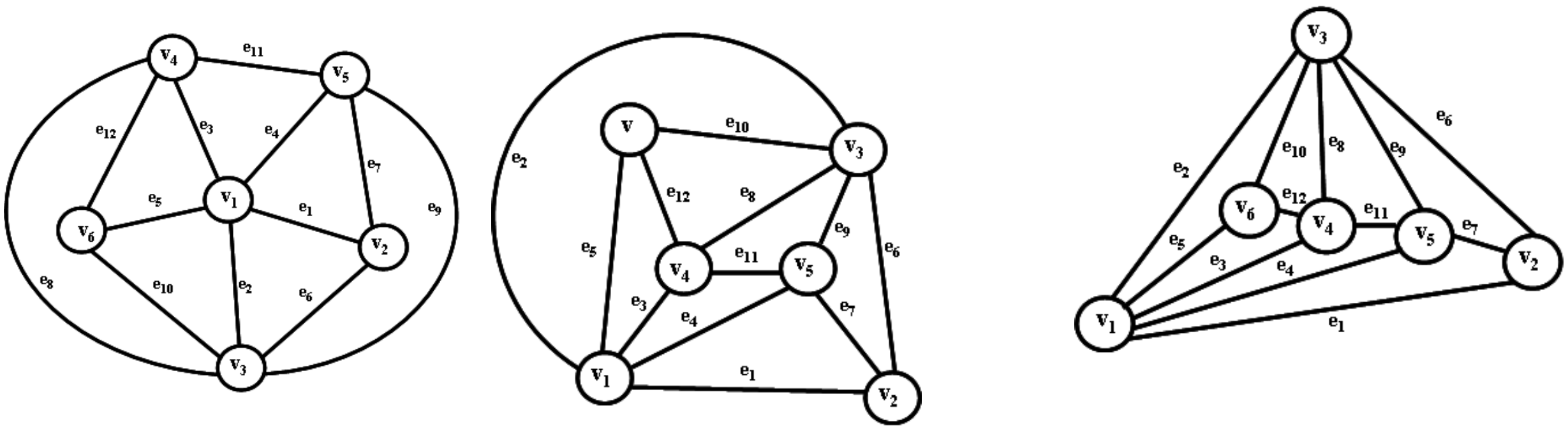

Рис. 3.1. Различные диаграммы графа G.

В случае планарного графа всегда имеется возможность проведения соединений (ребер графа) без пересечений. Причем такое представление не зависит от расположения вершин и характера проведения соединений (ребер). Такое представление планарного графа называется *плоским изображением графа* (см. рис. 3.6).

Но существуют структуры, которые являются общими для любого плоского изображения графа. Рассмотрим множество простых циклов являющихся границами граней плоского изображения. В качестве примера, рассмотрим граф G представленный на рис. 3.6. Запишем множество граничных циклов в виде элементов пространства суграфов:

$c_1 = \{e_2,e_5,e_{10}\}$; $c_2 = \{e_3,e_5,e_{12}\}$; $c_3 = \{e_8,e_{10},e_{12}\}$; $c_4 = \{e_8,e_9,e_{11}\}$;
$c_5 = \{e_6,e_7,e_9\}$; $c_6 = \{e_3,e_4,e_{11}\}$; $c_7 = \{e_1,e_4,e_7\}$; $c_8 = \{e_1,e_2,e_6\}$.

Цикломатическое число определяет количество независимых циклов графа $\nu(G) = m - n + 1$. Кольцевая сумма независимых циклов определяет обод.

Если задать направление обхода ребер в циклах, с соблюдением условия планарности Маклейна [], то можно записать циклы как кортежи вершин (см. рис. 3.7):

$c_1 = \{e_2,e_5,e_{10}\} \to \langle v_1,v_3,v_6 \rangle$; $c_2 = \{e_3,e_5,e_{12}\} \to \langle v_1,v_6,v_4 \rangle$; $c_3 = \{e_8,e_{10},e_{12}\} \to \langle v_4,v_6,v_3 \rangle$;
$c_4 = \{e_8,e_9,e_{11}\} \to \langle v_4,v_3,v_5 \rangle$; $c_5 = \{e_6,e_7,e_9\} \to \langle v_3,v_2,v_5 \rangle$; $c_6 = \{e_3,e_4,e_{11}\} \to \langle v_4,v_5,v_1 \rangle$;
$c_7 = \{e_1,e_4,e_7\} \to \langle v_1,v_5,v_2 \rangle$; $c_8 = \{e_1,e_2,e_6\} \to \langle v_1,v_2,v_3 \rangle$.

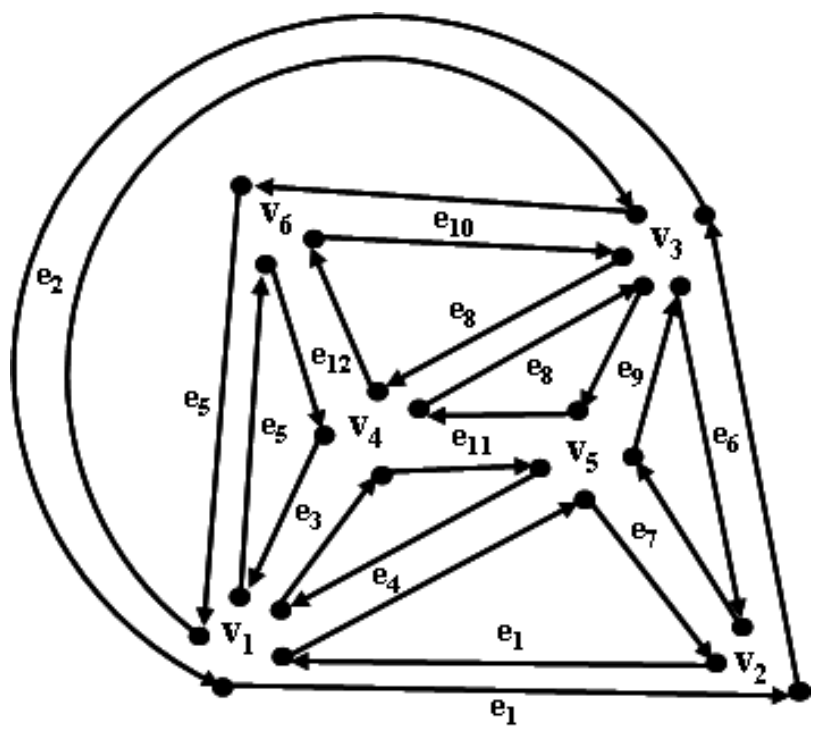

Рис. 3.2. Задание направления обхода ребер в циклах.

С другой стороны заданное подмножество циклов с направлением обхода ребер, порождает (индуцирует) определенный циклический порядок расположения смежных вершин для каждой вершины (см. рис. 3.3).

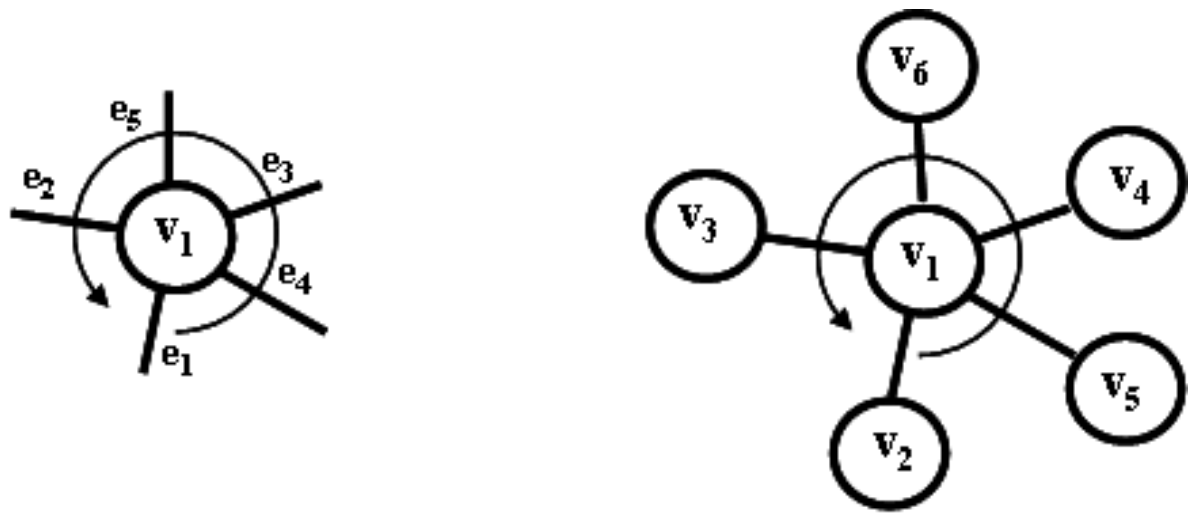

Рис. 3.3. Вращение вершины $v_1$.

Необходимым понятием для описания плоского рисунка графа **G**, является понятие о вращении вершин графа, введенное Г. Рингелем [75].

**Определение 3.1.** Для данного графа G вращение вершины A графа **G** - это ориентированный циклический порядок (или циклическая перестановка) всех ребер, инцидентных вершине A.

Пометим вершины графа числами 1,2,...,n. Если вершина **A** имеет степень три, а 1,2,3 - три вершины, смежные с вершиной A, то имеются две различные возможности для задания вращения вершины A. Будем описывать вращение вершины, указывая циклический порядок вершин смежных с вершиной A, вместо ребер, инцидентных этой вершине. Тогда указанные две возможности это:

( 1 2 3 ) = ( 2 3 1 ) = ( 3 1 2 )

и ( 3 2 1 ) = ( 2 1 3 ) = ( 1 3 2 ).

Вообще, число возможных вращений вершины степени n равно:

$$(n-1)! = 1\ 2\ ...\ n-1. \qquad (3.1)$$

**Определение 3.2.** Вращение $\sigma$ графа G - это вращение всех вершин графа G.

Запись $(G, \sigma)$ будет обозначать граф G с некоторым вращением $\sigma$.

Граф G с вращением часто бывает удобно изображать на плоскости таким образом,

чтобы, читая ребра, инцидентные некоторой вершине, или вершины, смежные с этой вершиной, по часовой стрелке (или против часовой стрелки), мы получили вращение в этой вершине.

Вращение графа можно описывать и представлять следующим образом. Обозначим вершины числами 1,2,...,n. Затем выпишем циклическую перестановку соседей для каждой вершины *i*. Эта перестановка порождается вращением вершины *i*, которое является циклической перестановкой ребер, инцидентных вершине *i*. Например, графу с вращением, показанному на рис. 3.9, соответствует схема:

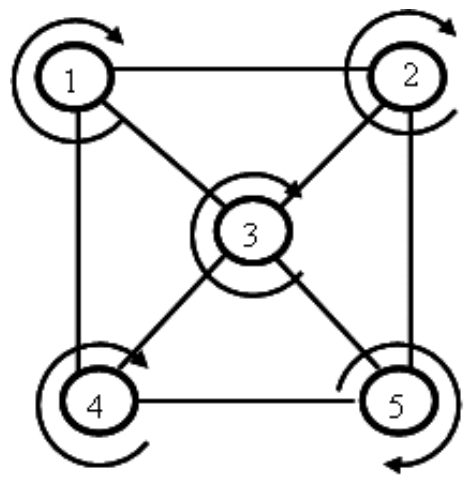

Рис. 3.4. Граф и его вращение вершин.

Диаграмма 3.1 (для рис. 3.4), $\sigma_i$ – вращение вершины $v_i$:

| | | | | |
|---|---|---|---|---|
| $\sigma_1$: | 2 | 3 | 4 | |
| $\sigma_2$: | 5 | 3 | 1 | |
| $\sigma_3$: | 1 | 2 | 5 | 4 |
| $\sigma_4$: | 1 | 3 | 5 | |
| $\sigma_5$: | 4 | 3 | 2 | |

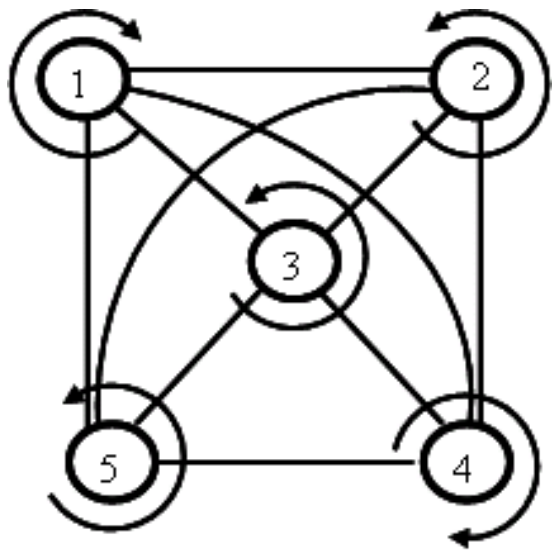

Рис. 3.5. Граф $K_5$ и его вращение вершин.

Диаграмма 3.2 (для рис. 3.5):

| | | | | |
|---|---|---|---|---|
| $\sigma_1$: | 2 | 4 | 3 | 5 |
| $\sigma_2$: | 4 | 5 | 3 | 1 |
| $\sigma_3$: | 2 | 1 | 5 | 4 |
| $\sigma_4$: | 3 | 1 | 2 | 5 |
| $\sigma_5$: | 2 | 3 | 1 | 4 |

представляет вращение полного графа $K_5$, показанное на рис. 3.5.

Пусть $V_1$ - вершина, инцидентная ребру $e_1$ в графе **G** с вращением $(G, \sigma)$. Мы построим в графе **G** замкнутый маршрут:

$$V_1, e_1, V_2, e_2, V_3, e_3, ..., \qquad (3.2)$$

где вершина $V_2$ - второй конец ребра $e_1$, а ребро $e_2$ следует за ребром $e_1$ во вращении вершины $V_2$, определяемом вращением $\sigma$. Затем определяется $V_3$, как вершина, инцидентная ребру $e_2$ и не равная $V_2$. После этого в качестве $e_3$ выбирается ребро, следующее за ребром $e_2$ во вращении вершины $V_3$ и т.д. Закончим процесс в точности перед тем моментом, когда должна повториться пара $V_1,e_1$. Она должна повториться, ибо граф **G** конечный, а наш процесс однозначно определен и в обратном направлении, а именно, если часть $V_{t-1},e_t,V_t,\ldots$ известна, то ребро $e_{t-1}$ определяется вращением вокруг вершины $V_{t-1}$. Мы назовем такой замкнутый маршрут циклом, порожденным вершиной $V_1$ и ребром $e_1$ и индуцированным вращением $\sigma$.

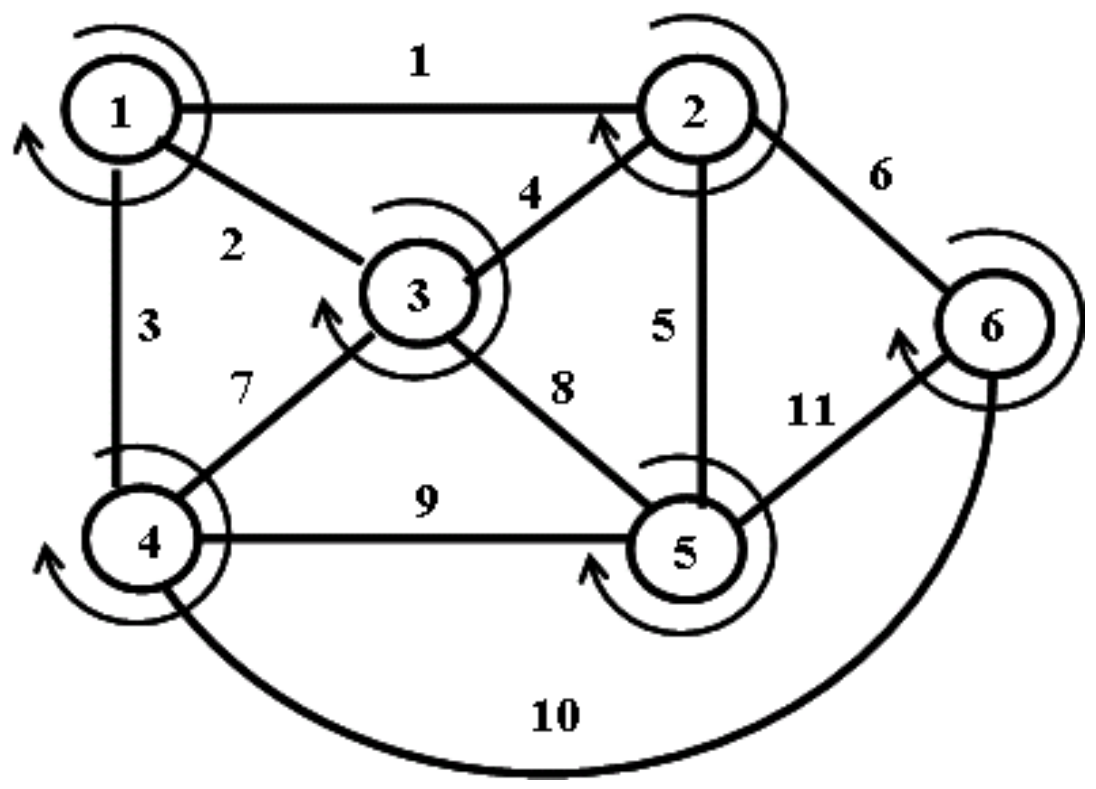

Рис. 3.6. Граф G и его вращение вершин $\sigma$.

Например, на рис. 3.5 показан цикл $\langle x_2,x_3,x_1,x_5,x_4,x_3,x_2,x_5,x_3,x_4,x_1,x_3,x_5,x_1,x_2\rangle$ индуцированный графом G с вращением $\sigma$. Еще два цикла индуцируются этим вращением: $v_4,v_5,v_2,v_4$ и $v_4,v_2,v_1,v_4$. Заметим, что во вращении каждое ребро появляется в точности дважды, второй раз - всегда в противоположном направлении.

Если граф планарен и имеется вращение, описывающее плоский рисунок, то циклы, индуцированные вращением, суть простые циклы. Например, для плоского графа **G** с вращением, представленным на рис. 3.4, имеем следующую систему индуцированных циклов.

Диаграмма вращений 3.3 (для рис. 3.4, но запись произведена в реберном представлении):

| | | | | |
|---|---|---|---|---|
| $\sigma_1$: | $e_1$ | $e_2$ | $e_3$ | |
| $\sigma_2$: | $e_6$ | $e_5$ | $e_4$ | $e_1$ |
| $\sigma_3$: | $e_2$ | $e_4$ | $e_8$ | $e_7$ |
| $\sigma_4$: | $e_3$ | $e_7$ | $e_9$ | $e_{10}$ |
| $\sigma_5$: | $e_9$ | $e_8$ | $e_5$ | $e_{11}$ |
| $\sigma_6$: | $e_6$ | $e_{10}$ | $e_{11}$ | |

Индуцированные циклы: $\langle v_1,v_3,v_2\rangle,\langle v_1,v_4,v_3\rangle,\langle v_2,v_3,v_5\rangle,\langle v_5,v_3,v_4\rangle,\langle v_2,v_5,v_6\rangle,\langle v_6,v_5,v_4\rangle$, $\langle v_2,v_6,v_4,v_1\rangle$ для вершин.

Или в виде множества ребер: $\{e_1,e_2,e_4\},\{e_2,e_3,e_7\},\{e_7,e_8,e_9\},\{e_4,e_5,e_8\},\{e_5,e_6,e_{11}\},\{e_9,e_{10},e_{11}\}$, $\{e_1,e_3,e_6,e_{10}\}$.

В случае непланарных графов, понятия только вращения вершин уже не достаточно, так как вращение вершин $\sigma$ не позволяет описывать пересечение ребер. В качестве примера, рассмотрим рисунки непланарного графа G (см. рис. 3.12), заданного одним и тем же вращением.

На рис. 3.7. вращение вершин в обоих случаях одинаково и производится по часовой стрелке, однако пересечение ребер различно. Например, на рисунке слева ребро $(v_7,v_{10})$ пересекается с ребрами $(v_2,v_6)$ и $(v_8,v_9)$, а на рисунке справа ребро $(v_7,v_{10})$ пересекается только с ребром $(v_2,v_6)$.

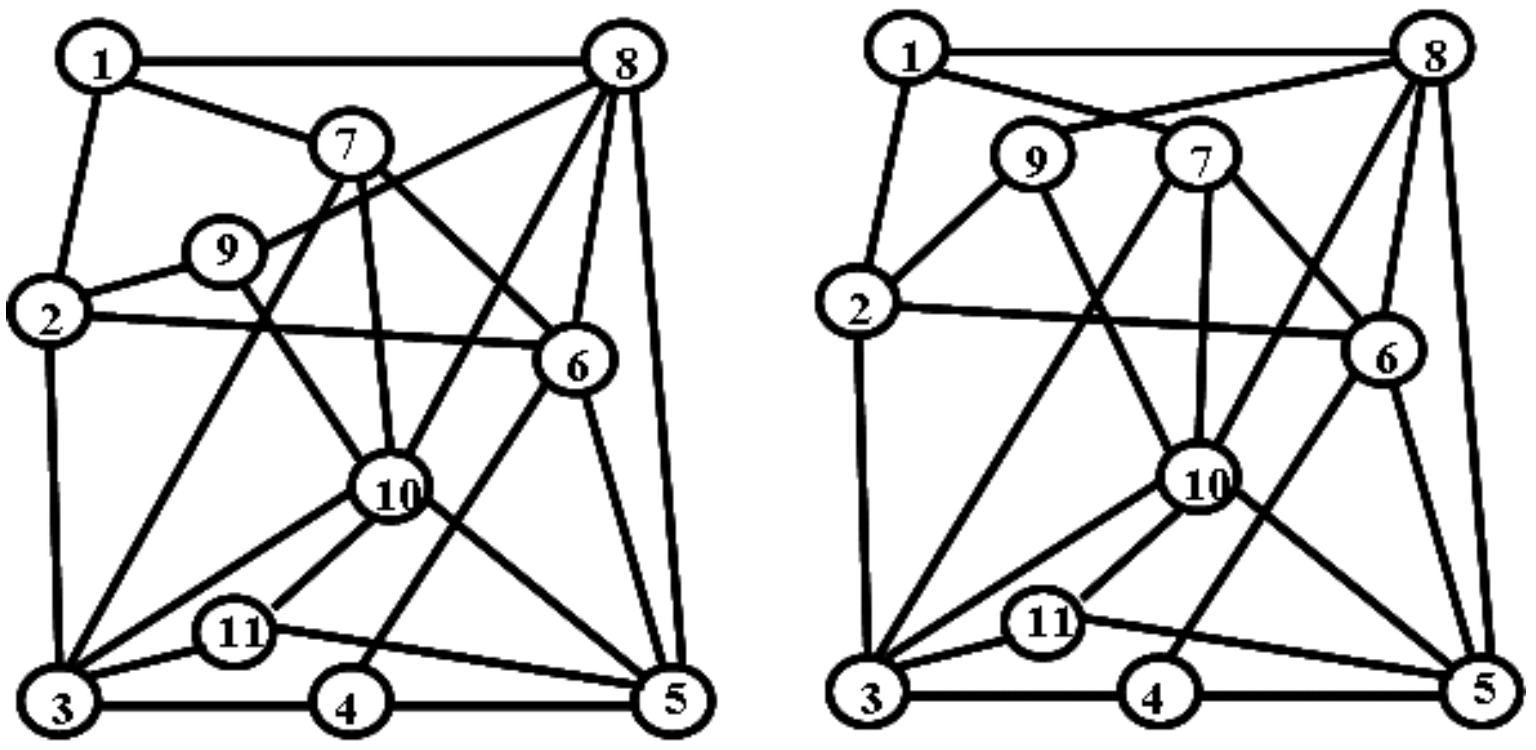

Рис. 3.7. Рисунок графа G с пересекающимися ребрами и его вращение вершин.

Диаграмма вращения 3.4:

| | | | | | | |
|---|---|---|---|---|---|---|
| $\sigma_1$: | 8 | 7 | 2 | | | |
| $\sigma_2$: | 1 | 9 | 6 | 3 | | |
| $\sigma_3$: | 2 | 7 | 10 | 11 | 4 | |
| $\sigma_4$: | 3 | 6 | 5 | | | |
| $\sigma_5$: | 4 | 11 | 10 | 6 | 8 | |
| $\sigma_6$: | 7 | 8 | 5 | 4 | 2 | |
| $\sigma_7$: | 1 | 6 | 10 | 3 | | |
| $\sigma_8$: | 5 | 6 | 10 | 9 | 1 | |
| $\sigma_9$: | 8 | 10 | 2 | | | |
| $\sigma_{10}$: | 9 | 7 | 8 | 5 | 11 | 3 |
| $\sigma_{11}$: | 3 | 10 | 5 | | | |

Теперь предположим, что планарный граф G вложен в евклидову плоскость, причем выполняется такое свойство:

(а) *Граница любой внутренней грани — простой цикл графа G.*

Дадим определение топологического рисунка графа.

**Определение 3.3.** *Топологическим рисунком планарного графа* будем называть одинаково направленное вращение всех его вершин $\sigma$, индуцирующее множество простых циклов графа, удовлетворяющее критерию планарности Маклейна [] и записывать его как $(G,\sigma)$.

Топологический рисунок графа позволяет осуществлять операции с рисунком графа, не производя никаких геометрических построений в пространстве.

Так как, вращение вершин графа индуцирует замкнутые маршруты, то топологический рисунок может быть описан системой изометрических циклов.

Как видно из приведенного примера, в случае описания рисунка непланарного графа одно и то же вращение вершин может характеризовать различное пересечение ребер. На рис. 3.2 вращение вершин производится по часовой стрелке и может быть записано в виде множества циклических подмножеств:

$\sigma_1 = \langle x_8,x_7,x_2\rangle$, $\sigma_2 = \langle x_1,x_9,x_6,x_3\rangle$, $\sigma_3 = \langle x_2,x_7,x_{10},x_{11},x_4\rangle$, $\sigma_4 = \langle x_3,x_6,x_5\rangle$,

$\sigma_5 = \langle x_4,x_{11},x_{10},x_6,x_8\rangle$, $\sigma_6 = \langle x_7,x_8,x_5,x_4,x_2\rangle$, $\sigma_7 = \langle x_1,x_6,x_{10},x_3\rangle$,

$\sigma_8 = \langle x_5,x_6,x_{10},x_9,x_1\rangle$, $\sigma_9 = \langle x_8,x_{10},x_2\rangle$, $\sigma_{10} = \langle x_9,x_7,x_8,x_5,x_{11},x_3\rangle$, $\sigma_{11} = \langle x_3x_{10},x_5\rangle$.

Для описания рисунка графа с пересекающимися ребрами можно ввести понятие мнимой вершины - т.е. вершины, которая характеризует топологическое пересечение ребер. То есть, введение мнимой вершины определяет пересечение пары ребер. И тогда удается описывать топологический рисунок графа с помощью понятия вращения вершин с его свойством индуцировать простые циклы.

**Определение 3.4.** *Мнимая вершина – это топологическое местоположение пересечения двух ребер.*

На рис. 3.7 представлен рисунок графа с мнимыми вершинами. Для идентификации мнимых вершин применяется обозначение отличное от обозначения истинных вершин графа. Таким образом, запись $(G,\sigma)$ будет обозначать граф **G** с некоторым вращением $\sigma$ и одновременно будет характеризовать топологический рисунок графа на плоскости.

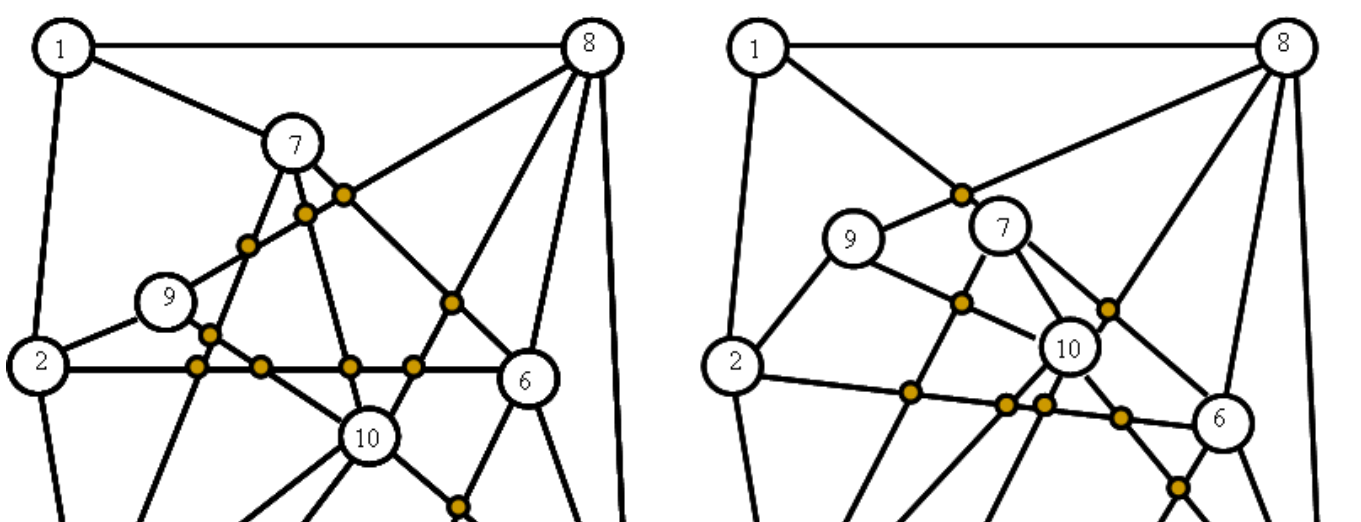

Рис. 3.8. Рисунок графа **G** с введенными мнимыми вершинами.

Таким образом, рисунок графа с $n_e$ исходными и $n_f$ мнимыми вершинами $n = n_e + n_f$ можно описывать и представлять следующим образом:

- трехместным предикатом, устанавливающим соответствие между ребром и его концевой парой вершин в виде матрицы смежностей, матрицы инциденций или логической записью соответствия между ребром и вершинами;

- вращением вершин $\sigma$ записанное в виде множества циклических подмножеств, которые в свою очередь индуцируют систему простых циклов

$$\sigma = \{ \sigma_1, \sigma_2, \ , \sigma_n \}, \text{ где } \sigma_i = <v_1, v_2, \ldots, v_{i-1}, v_{i+1}, \ldots, v_k>, \ v_i \in V. \qquad (3.3)$$

- количеством мнимых вершин $n_f$.

**Определение 3.5.** *Будем говорить, что топологический рисунок непланарного графа – граф с заданным вращением вершин $\sigma$ и заданным количеством мнимых вершин и обозначать как $(G, \sigma, n_f)$.*

Для описания заданного рисунка графа на плоскости при заданном трехмерном предикате **P** и вращении вершин $\sigma$ этого вполне достаточно. Однако, существует и обратная задача, когда вращения вершин исходного графа и мнимых вершин нужно определить при заданном трехместном предикате **P**.

## 3.2. Функционал Маклейна

Пусть G(V,E) – несепарабельный граф с множеством вершин $V=\{v_1,v_2,...,v_n\}$ и ребер $E=\{e_1,e_2,...,e_m\}$, где n-количество вершин графа и m-количество ребер графа G.

Пусть $L_G$ - множество всех суграфов этого графа []. Относительно операции сложения:

$$(V,E_1;P) \oplus (V,E_2;P) \overset{\text{def}}{=} (V,(E_1 \cup E_2) \backslash (E_1 \cap E_2);P). \qquad (3.4)$$

Это множество, как известно, образует абелеву 2-группу, которую можно рассматривать как векторное пространство над полем из двух элементов GF(2). Размерность этого пространства, называемого пространством суграфов графа G, конечно и равно m ($\dim L_G = m$). В качестве базиса этого пространства выберем множество однореберных суграфов $(u_1,u_2,...,u_m)$. Тогда в этом базисе каждому элементу Y пространства $L_G$ однозначно сопоставляется последовательность координат $(a_1,a_2,...,a_m)$, где $a_i \in \{0,1\}$. При этом оказывается, что ребро $u_i$ входит в суграф Y, если $a_i = 1$, и не входит в данный суграф - в противном случае. В дальнейшем для удобства будем отождествлять пространство суграфов $L_G$ и его координатное пространство.

Напомним, что суграф называется квазициклом, если все его вершины имеют четную валентность (в данном случае валентность совпадает с локальной степенью вершины). Множество графа G образует, как легко можно видеть, подпространство квазициклов C пространства $L_G$ [21]. Известно, что размерность подпространства $L_G^C$ совпадает с цикломатическим числом $\nu(G) = m - n + 1$ графа G, а порядок группы $L_G^C$ равен $2^{\nu(G)}$.

Теперь мы можем сформулировать критерий планарности Маклейна [55].

**Теорема 3.1** [55]. Граф G планарен тогда и только тогда, когда существует такой базис подпространства квазициклов, где каждое ребро принадлежит не более, чем двум циклам.

Выберем некоторый базис $(c_1,c_2,...,c_k)$, где $k = \nu(G)$ - размерность подпространства квазициклов. Рассмотрим матрицу циклов C, строки которой соответствуют элементам указанного базиса:

$$C = \begin{Vmatrix} a_{11} & a_{12} & \dots & a_{1m} \\ a_{21} & a_{22} & \dots & a_{2m} \\ \dots & \dots & \dots & \dots \\ a_{k1} & a_{k2} & \dots & a_{km} \end{Vmatrix}$$

Элементы $a_{ij}$ этой матрицы принадлежат полю GF(2)={0,1}. Очевидно, что указанный базис удовлетворяет условию Маклейна тогда и только тогда, когда в каждом столбце матрицы C содержится не более двух единиц. Рассмотрим другой базис этого подпространства $(c'_1,c'_2,...,c'_k)$, которому соответствует матрица циклов $C'$. Тогда эти матрицы связаны отношением:

$$C' = TC, \qquad (3.5)$$

где T - невырожденная матрица (матрица перехода от базиса C к базису C'). Поскольку каждая невырожденная матрица разлагается в произведение элементарных матриц $P_1,P_2,...,P_s$, а умножение слева на элементарную матрицу равносильно выполнению одной элементарной операции над строками, то из (3.5) следует, что:

$$C' = P_s\, P_{s-1} \dots P_1 C. \qquad (3.6)$$

Таким образом, каждый базис в этом пространстве получается из данного базиса при помощи цепочки элементарных преобразований. А на матричном языке проблема распознавания планарности сводится к нахождению такой матрицы в классе эквивалентных матриц (т.е. матриц, которые получаются друг из друга при помощи элементарных преобразований над строками), у которой в каждом столбце содержится не более двух единиц [].

*Вектором циклов по ребру* – будем называть вектор $P_e$, где каждая составляющая определяет количество циклов проходящим по соответствующему ребру графа.

Указанный критерий позволяет разработать методику определения планарности графа, сводя проблему планарности к отысканию минимума некоторого функционала на множестве базисов подпространства квазициклов. Что позволяет применить к решению задачи методы дискретной оптимизации. Определим следующий функционал на матрице C, соответствующий базису подпространства квазициклов (и будем его впредь называть функционалом Маклейна):

$$F(C)=\sum_{i=1}^{m}(s_i-1)(s_i-2)=\sum_{i=1}^{m}s_i^2-3\sum_{i=1}^{m}s_i+2m\,, \qquad (3.7)$$

где $s_i$ – количество циклов проходящих по i-му ребру графа.

Очевидно, что матрица циклов C соответствует базису Маклейна (т.е. базису, удовлетворяющему условию Маклейна) тогда и только тогда, когда F(C) = 0.

Функционал F(C) принимает целое неотрицательное значение и проблема отыскания базиса Маклейна, таким образом, является частным случаем следующей задачи дискретной оптимизации: найти минимум F(C) на множестве матриц C размера $k \times m$ и ранга $k = \nu(G)$, покрывающим все множество вершин графа.

Очевидно, что не для любого графа G минимум F(C) будет равен нулю, согласно критерию Маклейна. Отметим, что нулевое значение данный функционал принимает только в случае планарного графа. В общем случае, решение указанной задачи минимизации позволяет построить алгоритм выделения максимально плоского суграфа. Поэтому, вначале рассмотрим свойства функционала F(C), где C пробегает множество указанных матриц.

В произвольном линейном пространстве отсутствуют метрические характеристики, такие как длина, расстояние и т.д. Однако, их использование становится возможным, если в линейное пространство удается ввести следующую операцию.

**Определение 3.6.** Пусть в линейном пространстве каждой паре элементов $x$ и $y$ поставлено в соответствие вещественное число $(x,y)$, называемое скалярным произведением так, что выполнены следующие аксиомы:

$$(x,y)=(y,x); \qquad (3.8)$$

$$(\lambda x,y)=\lambda(x,y) \qquad (3.9)$$

$$(x+z,y)=(x,y)+(z,y); \qquad (3.10)$$

$$(x,x)\ge 0 \text{ причем } (x,x)=0 \leftrightarrow x=0. \qquad (3.11)$$

Тогда говорят, что задано *евклидово пространство E* [14].

Введем понятие модуля пересечения на множестве квазициклов. Для двух квазициклов $c_i$ и $c_j$ модуль пересечения, будем обозначать его $(c_i,c_j)$, определяется из следующей формулы:

$$(c_i,c_j) = |c_i \cap c_j|, \qquad (3.12)$$

где: $|c_i|$ - длина квазицикла i , $|c_i \cap c_j|$ - мощность множества пересечения квазициклов $c_i$ и $c_j$.

Как видно, модуль пересечения квазициклов удовлетворяет всем аксиомам скалярного произведения евклидового пространства и поэтому пространство суграфов является *евклидовым пространством*.

Полезным инструментом при исследовании свойств набора базисных циклов для решения задачи построения плоского графа и укладки его на плоскость является матрица Грамма.

**Определение 3.7.** Будем называть матрицей Грамма следующую матрицу, образованную набором модулей пересечения всевозможных пар простых циклов:

**Г**(C) =

| | $c_1$ | $c_2$ | | $c_{\nu-1}$ | $c_\nu$ |
|---|---|---|---|---|---|
| $c_1$ | $(c_1,c_1)$ | $(c_1,c_2)$ | … | $(c_1,c_{\nu-1})$ | $(c_1,c_\nu)$ |
| $c_2$ | $(c_2,c_1)$ | $(c_2,c_2)$ | … | $(c_2,c_{\nu-1})$ | $(c_2,c_\nu)$ |
| $c_3$ | $(c_3,c_1)$ | $(c_3,c_2)$ | … | $(c_3,c_{\nu-1})$ | $(c_3,c_\nu)$ |
| … | … | … | … | … | … |
| $c_{\nu-1}$ | $(c_{\nu-1},c_1)$ | $(c_{\nu-1},c_2)$ | … | $(c_{\nu-1},c_{\nu-1})$ | $(c_{\nu-1},c_\nu)$ |
| $c_\nu$ | $(c_\nu,c_1)$ | $(c_\nu,c_2)$ | … | $(c_\nu,c_{\nu-1})$ | $(c_\nu,c_\nu)$ |

Введем для произвольно заданного базиса Q(G), $Q(G) = \{c_i\ ;\ i = [1,\nu]\}$ следующую оценку:

$$\Gamma(C) = \sum_{i=1}^{k}\sum_{j=1}^{k}(c_i,c_j) = (c_1,c_1)\oplus(c_1,c_2)\oplus...\oplus(c_k,c_k),$$

где $c_i$ = m - мерный вектор:

$c_i =(a_{i1},a_{i2},\ldots,a_{im})$, $a_{ip} \in \{0,1\}$.

Учитывая также, что

$$(ci,cj) = \sum_{p=1}^{m} a_{ip}a_{jp}\ . \tag{3.13}$$

Записываем Г(C) в следующем виде:

$$\Gamma(C) = \sum_{i=1}^{k}\sum_{j=1}^{k}(\sum_{p=1}^{m} a_{ip}aj_p) = \sum_{p=1}^{m}(\sum_{i=1}^{k}\sum_{j=1}^{k} a_{ip}a_{jp}) = \tag{3.14}$$

Легко видеть, что:

$$\sum_{i=1}^{k}\sum_{j=1}^{k} a_{ip}a_{jp} = s_p^2\ , \tag{3.15}$$

где $S_p$ - количество циклов, содержащих p-ое ребро.

Подставляя выражение (3.15) в (3.14), получим:

$$\Gamma(C) = \sum_{p=1}^{m} s_p^2\ . \tag{3.16}$$

Между функционалом Маклейна (3.6) и (3.15) имеется следующая связь:

$$F(C) = \Gamma(C) - 3\sum_{i=1}^{k} l_i + 2m. \tag{3.17}$$

Таким образом, функционал Маклейна можно записать через модуль пересечения квазициклов ($l_i$ - длина i-го квазицикла):

$$F(C)=\sum_{i=1}^{k}\sum_{j=1}^{k}(c_i,c_j)-3\sum_{i=1}^{k}l_i+2m. \qquad (3.18)$$

Следующее выражение, справедливое для плоских графов, связывает количество ребер плоского графа и суммарную длину квазициклов, принадлежащих выбранному базису. Для базиса Маклейна (τ – система циклов [43,45]):

$$\sum_{i=1}^{k}l_i+l_0-2m=0, \qquad (3.19)$$

где $l_0$ - длина обода графа.

Введенное выражение (3.19) позволяет построить целый спектр функционалов, эквивалентных функционалу Маклейна и характеризующих любой базис с учетом обода графа. Например:

$$F_1(C)=\sum_{j=1}^{m}s_j^2-2\sum_{i=1}^{k}l_i+l_0; \qquad (3.20)$$

$$F_2(C)=\sum_{j=1}^{m}s_j^2-\sum_{i=1}^{k}l_i+2l_0-2m; \qquad (3.21)$$

$$F_3(C)=\sum_{j=1}^{m}s_j^2+3l_0-4m. \qquad (3.22)$$

***Пример 3.1.*** В качестве иллюстрации, рассмотрим граф G, представленный на рис. 3.14. Для данного графа G существует базис Маклейна, функционал которого равен нулю.

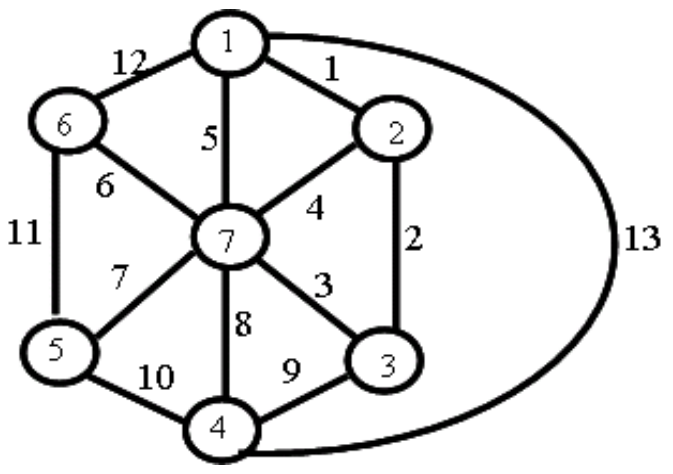

Рис. 3.9. Граф G.

Матрица независимых циклов имеет вид:

| | | $e_1$ | $e_2$ | $e_3$ | $e_4$ | $e_5$ | $e_6$ | $e_7$ | $e_8$ | $e_9$ | $e_{10}$ | $e_{11}$ | $e_{12}$ | $e_{13}$ | |
|---|---|---|---|---|---|---|---|---|---|---|---|---|---|---|---|
| | $c_1$ | 1 | | | 1 | 1 | | | | | | | | | **3** |
| | $c_2$ | | 1 | 1 | 1 | | | | | | | | | | **3** |
| | $c_3$ | | | 1 | | | | | 1 | 1 | | | | | **3** |
| C= | $c_4$ | | | | | | | 1 | 1 | | 1 | | | | **3** |
| | $c_5$ | | | | | | 1 | 1 | | | | 1 | | | **3** |
| | $c_6$ | | | | | 1 | 1 | | | | | | 1 | | **3** |
| | $c_7$ | | | | | | | | | | 1 | 1 | 1 | 1 | **4** |
| | $\mathbf{V_u}$ | **1** | **1** | **2** | **2** | **2** | **2** | **2** | **2** | **1** | **1** | **2** | **2** | **1** | |

↑

количество циклов, проходящих по данному ребру $e_3$.

Базис Маклейна, представленный в виде циклов, проходящих по соответствующим ребрам, имеет вид:

$c_1 = \{e_1, e_4, e_5\}$; $c_2 = \{e_2, e_3, e_4\}$; $c_3 = \{e_3, e_8, e_9\}$; $c_4 = \{e_7, e_8, e_{10}\}$;

$c_5 = \{e_6, e_7, e_{11}\}$; $c_6 = \{e_5, e_6, e_{12}\}$; $c_7 = \{e_{10}, e_{11}, e_{12}, e_{13}\}$; $c_0 = \{e_1, e_2, e_9, e_{13}\}$.

Обод плоского графа $c_0$ можно представить как кольцевую сумму элементов базиса.

Вращение графа, индуцируемое данным базисом циклов и ободом графа, представлено диаграммой, полученной по правилу формирования вращения вершин:

$\sigma_1$: 2 7 6 4
$\sigma_2$: 3 7 1
$\sigma_3$: 2 4 7
$\sigma_4$: 5 7 3 1
$\sigma_5$: 6 7 4
$\sigma_6$: 7 5 1
$\sigma_7$: 1 2 3 4 5 6

Таким образом, задача проверки графа на планарность может быть сведена к задаче поиска базиса подпространства квазициклов, у которого функционал Маклейна равен нулю. Рисунок такого графа на плоскости определяется вращением $(G, \sigma)$, индуцированным таким базисом циклов и ободом графа.

***Пример 3.2.*** Вычислим функционал Маклейна по формуле (3.18) для нашего графа G (см. рис. 3.9).

Матрица Грама для графа представленного на рис. 3.8 и циклов $c_1, c_2, c_3, c_4, c_5, c_6, c_7$ имеет вид:

Г(G) =

| | $c_1$ | $c_2$ | $c_3$ | $c_4$ | $c_5$ | $c_6$ | $c_7$ |
|---|---|---|---|---|---|---|---|
| $c_1$ | 3 | 1 | 0 | 0 | 0 | 1 | 0 |
| $c_2$ | 1 | 3 | 1 | 0 | 0 | 0 | 0 |
| $c_3$ | 0 | 1 | 3 | 1 | 0 | 0 | 0 |
| $c_4$ | 0 | 0 | 1 | 3 | 1 | 0 | 1 |
| $c_5$ | 0 | 0 | 0 | 1 | 3 | 1 | 1 |
| $c_6$ | 1 | 0 | 0 | 0 | 1 | 3 | 1 |
| $c_7$ | 0 | 0 | 0 | 1 | 1 | 1 | 4 |

$$\Gamma(C) = \sum_{i=1}^{7}\sum_{j=1}^{7} a_{ij} = 5+5+5+6+6+6+7 = 40,$$

$$F(C) = \Gamma(C) - 3\sum_{i=1}^{k} l_i + 2m = 40 - 3\times 22 + 2\times 13 = 40 - 66 + 26 = 0.$$

Ранее было показано, что структурное число C, полученное как произведение однострочных структурных чисел, характеризующих простые циклы, проходящие по хордам графа на выделенном дереве, определяет все множество базисов, состоящих из простых циклов, и его построение не зависит от выбора дерева графа.

Так как любой столбец структурного числа С характеризует систему независимых циклов, то соответствующему столбцу можно поставить в соответствие значение функционала Маклейна. Однако, ввиду того, что величина функционала Маклейна является интегральной характеристикой всего базиса, не удается для каждого цикла графа поставить в соответствие число $\lambda_i$ так, чтобы сумма всех $\lambda_i$ (i = 1,2,...,$\nu$(G)) определяла какой-то функционал. Однако можно поставить множеству изометрических циклов в соответствие значение функционала Маклейна.

Напомним, что алгебраической обратной производной структурного числа называется структурное число $\delta C/\delta\alpha$, равное:

$$\frac{\delta C}{\delta \alpha} = C \text{ - где столбцы, содержащие элемент } \alpha \text{, опущены.} \qquad (3.23)$$

Воспользовавшись способом записи структурного числа в виде семейства множеств, можно записать обратную производную как:

$$\frac{\delta C}{\delta \alpha} = \{c_i \mid \alpha \in c_i, c_i \in C\}. \qquad (3.24)$$

Алгебраическая обратная производная структурного числа С эквивалентна рассмотрению множеств базисов, не содержащих базисы с номером $\alpha$ .

Таким образом, множество базисов, состоящих из простых циклов (или суграфов) можно записывать и хранить в компактном виде как результат произведения однострочных структурных чисел, состоящих из номеров суграфов, включающих хорды графа (т.е. в виде структурного числа С).

## 3.3. Свойство скалярного произведения суграфов

Рассмотрим свойства скалярного произведения суграфов (3.12). Предположим, что заданы вектора x,y,z характеризующие простые циклы графа. Естественно, что простые циклы x,y,z могут иметь общие ребра, а могут и не иметь. Нас интересует проверка выполнения аксиомы (3.11) и аксиомы (3.3). Обозначим высказывание (x,y) = 0 если циклы x и y не имеют общих ребер и (x,y) = 1, если циклы x и y имеют общие ребра. Соответственно это относится и к скалярным произведениям (y,z) и (x,z)

Построим следующую таблицу истинности для высказываний (3.4) и (3.9).

| $(x,y)$ | $(y,z)$ | $(x,z)$ | | $(x,y) + (y,z) \geq (x,z)$ | $(x \oplus z, y) = (x, y) +$ |
|---|---|---|---|---|---|
| 0 | 0 | 0 | 1 случай | выполнено | выполнено |
| 0 | 0 | 1 | 2 случай | не выполнено | выполнено |
| 0 | 1 | 0 | 3 случай | выполнено | выполнено |
| 0 | 1 | 1 | 4 случай | выполнено | выполнено |
| 1 | 0 | 0 | 5 случай | выполнено | выполнено |
| 1 | 0 | 1 | 6 случай | выполнено | выполнено |
| 1 | 1 | 0 | 7 случай | выполнено | выполнено |
| 1 | 1 | 1 | 8 случай | выполнено | выполнено |

Рассмотрим граф $G_1$, представленный на рис. 3.15. Данный граф определяется следующим базисом простых циклов:

$c_1 = \{e_1,e_2,e_3\}$; $c_2 = \{e_3,e_6,e_7\}$; $c_3 = \{e_4,e_5,e_6\}$;
$c_4 = \{e_7,e_8,e_9\}$; $c_5 = \{e_5,e_{10},e_{11}\}$.

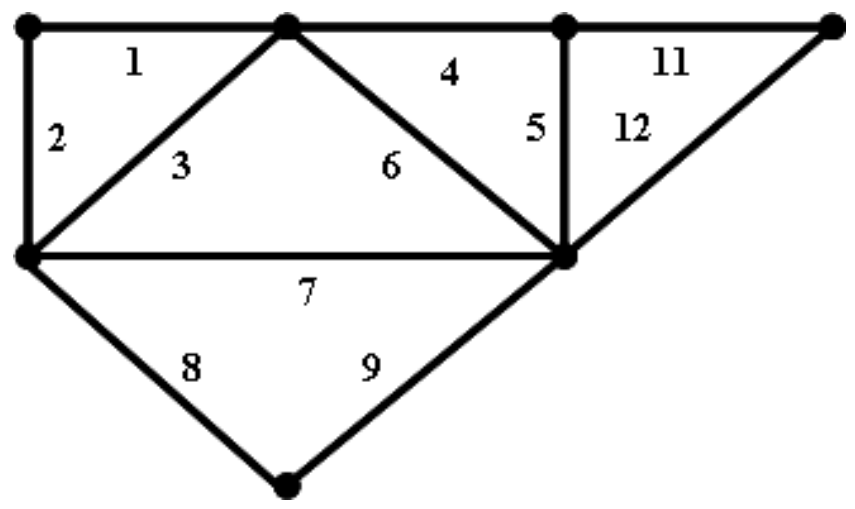

Рис. 3.10. Граф $G_1$.

**Случай 1.**

$x = c_1 = \{e_1,e_2,e_3\}$;
$y = c_3 = \{e_4,e_5,e_6\}$;
$z = c_4 = \{e_7,e_8,e_9\}$.

$x \oplus z == \{e_1,e_2,e_3\} \oplus \{e_7,e_8,e_9\} = \{e_1,e_2,e_3,e_7,e_8,e_9\}$;

$|(x \oplus z, y)| = |\{e_1,e_2,e_3,e_7,e_8,e_9\} \cap \{e_4,e_5,e_6\}| = |0|$;

$|(x \oplus z, y)| = |(x, y)| + |(z, y)|$; $|0| = |0| + |0|$;

$|(x, y)| + |(z, y)| \geq |(x, z)|$; $|0| + |0| = |0|$.

**Случай 2.**

$x = c_1 = \{e_1,e_2,e_3\}$;
$y = c_5 = \{e_5,e_{10},e_{11}\}$;
$z = c_2 = \{e_3,e_6,e_7\}$.

$x \oplus z == \{e_1,e_2,e_3\} \oplus \{e_3,e_6,e_7\} = \{e_1,e_2,e_6,e_7\}$;

$|(x \oplus z, y)| = |\{e_1,e_2,e_6,e_7\} \cap \{e_5,e_{10},e_{11}\}| = |0|$;

$|(x \oplus z, y)| = |(x, y)| + |(z, y)|$; $|0| = |0| + |0|$;

$|(x, y)| + |(z, y)| \geq |(x, z)|$; $|0| + |0| < |1|$.

**Случай 3.**

$x = c_1 = \{e_1,e_2,e_3\}$;
$y = c_3 = \{e_4,e_5,e_6\}$;
$z = c_5 = \{e_5,e_{10},e_{11}\}$.

$x \oplus z = \{e_1,e_2,e_3\} \oplus \{e_5,e_{10},e_{11}\} = \{e_1,e_2,e_3,e_5,e_{10},e_{11}\}$;

$|(x \oplus z, y)| = |\{e_1,e_2,e_3,e_5,e_{10},e_{11}\} \cap \{e_4,e_5,e_6\}| = |1|$;

$|(x \oplus z, y)| = |(x, y)| + |(z, y)|$; $|1| = |0| + |1|$;

$|(x, y)| + |(z, y)| \geq |(x, z)|$; $|0| + |1| = |1|$.

**Случай 4.**

$x$ = $c_1$ = {$e_1$,$e_2$,$e_3$};
$y$ = $c_3$ = {$e_4$,$e_5$,$e_6$};
$z$ = $c_2$ = {$e_3$,$e_6$,$e_7$}.

$$x \oplus z == \{e_1, e_2, e_3\} \oplus \{e_3, e_6, e_7\} = \{e_1, e_2, e_6, e_7\};$$

$$|(x \oplus z, y)| = |\{e_1, e_2, e_6, e_7\} \cap \{e_4, e_5, e_6\}| = |1|;$$

$$|(x \oplus z, y)| = |(x, y)| + |(z, y)|; \qquad |1| = |0| + |1|;$$

$$|(x, y)| + |(z, y)| \geq |(x, z)|; \qquad |0| + |1| = |1|.$$

**Случай 5.**

$x$ = $c_1$ = {$e_1$,$e_2$,$e_3$};
$y$ = $c_2$ = {$e_3$,$e_6$,$e_7$};
$z$ = $c_5$ = {$e_5$,$e_{10}$,$e_{11}$}.

$$x \oplus z == \{e_1, e_2, e_3\} \oplus \{e_5, e_{10}, e_{11}\} = \{e_1, e_2, e_3, e_5, e_{10}, e_{11}\};$$

$$|(x \oplus z, y)| = |\{e_1, e_2, e_3, e_5, e_{10}, e_{11}\} \cap \{e_3, e_6, e_7\}| = |1|;$$

$$|(x \oplus z, y)| = |(x, y)| + |(z, y)|; \qquad |1| = |1| + |0|;$$

$$|(x, y)| + |(z, y)| \geq |(x, z)|; \qquad |1| + |0| > |0|.$$

**Случай 6.**

$x$ = $c_2$ = {$e_3$,$e_6$,$e_7$};
$y$ = $c_1$ = {$e_1$,$e_2$,$e_3$};
$z$ = $c_4$ = {$e_7$,$e_8$,$e_9$}.

$$x \oplus z == \{e_3, e_6, e_7\} \oplus \{e_7, e_8, e_9\} = \{e_3, e_6, e_8, e_9\};$$

$$|(x \oplus z, y)| = |\{e_3, e_6, e_8, e_9\} \cap \{e_1, e_2, e_3\}| = |1|;$$

$$|(x \oplus z, y)| = |(x, y)| + |(z, y)|; \qquad |1| = |1| + |0|;$$

$$|(x, y)| + |(z, y)| \geq |(x, z)|; \qquad |1| + |0| > |0|.$$

**Случай 7.**

$x$ = $c_1$ = {$e_1$,$e_2$,$e_3$};
$y$ = $c_2$ = {$e_3$,$e_6$,$e_7$};
$z$ = $c_4$ = {$e_7$,$e_8$,$e_9$}.

$$x \oplus z == \{e_1, e_2, e_3\} \oplus \{e_7, e_8, e_9\} = \{e_1, e_2, e_3, e_7, e_8, e_9\};$$

$$|(x \oplus z, y)| = |\{e_1, e_2, e_3, e_7, e_8, e_9\} \cap \{e_3, e_6, e_7\}| = |2|;$$

$$|(x \oplus z, y)| = |(x, y)| + |(z, y)|; \qquad |2| = |1| + |1|;$$

$$|(x, y)| + |(z, y)| \geq |(x, z)|; \qquad |1| + |1| = |2|.$$

**Случай 8.**

$x$ = $c_1$ = {$e_1$,$e_2$,$e_3$};
$y$ = $c_2$ = {$e_1$,$e_3$,$e_5$};
$z$ = $c_3$ = {$e_4$,$e_5$,$e_6$}.

$$x \oplus z == \{e_1, e_2, e_3\} \oplus \{e_4, e_5, e_6\} = \{e_1, e_2, e_3, e_4, e_5, e_6\};$$

$$|(x \oplus z, y)| = |\{e_1, e_2, e_3, e_4, e_5, e_6\} \cap \{e_1, e_3, e_5\}| = |3|;$$

$$|(x \oplus z, y)| = |(x, y)| + |(z, y)|; \qquad |3| = |2| + |1|;$$

$$|(x, y)| + |(z, y)| \geq |(x, z)|; \qquad |2| + |1| = |3|.$$

Мы рассмотрели все случаи скалярного произведения циклов. Следовательно, пространство евклидово, так как выполняются аксиомы (3.8 – 3.11). Но не выполняется метрическая аксиома Фреше (3.3), следовательно, пространство не метрическое.

## 3.4. Алгоритм построения вращений вершин графа для заданной системы циклов

Для создания топологического рисунка графа нам необходимо иметь алгоритм построения вращений вершин для заданной системы циклов, удовлетворяющих нулевому значению функционала Маклейна. Предположим, что такая система циклов найдена и существует обод, тогда:

**Шаг 1**. [**Циклическая вершинная запись подмножества циклов**].

Представим запись циклов, в виде ориентированных ребер используя вектор $V_u$ для выделенного базиса с нулевым значением функционала Маклейна. Запись будем производить с учетом того, что по каждому ребру мы проходим дважды, как в прямом, так и в обратном направлении. Идем на шаг 2.

**Шаг 2. [Просмотр списка циклов]**. Последовательно просматриваем изометрические циклы базиса и обод данного базиса. Выбираем непомеченный цикл и идеи на шаг 3. Если список всех циклов окончен то получаем таблицу вращения вершин для плоского рисунка графа.

**Шаг 3**. [**Окаймляющая запись вершин**].

Просматриваем ориентированную запись вершин выбранного цикла. В таблице вращений вершин для каждой вершины записываем окаймляющие ее вершины из данной ориентированной записи вершин. В случае необходимости объединяем записи для каждой вершины в таблице вращений. Идем на шаг 3.

***Пример 3.3***. Построить таблицу вращений вершин графа для максимально плоского рисунка суграфа представленного на рис. 3.10.

Согласно шагу 1 алгоритма, запишем систему независимых циклов и обод для плоского графа в виде ориентированных ребер через вершины (вращение по часовой стрелке):

$c_1 = \{e_1,e_2,e_6\} = \langle v_1,v_2,v_3\rangle$; $c_2 = \{e_1,e_5,e_8\} = \langle v_1,v_6,v_2\rangle$;
$c_3 = \{e_2,e_3,e_9\} = \langle v_1,v_3,v_4\rangle$; $c_4 = \{e_3,e_4,e_{11}\} = \langle v_1,v_4,v_5\rangle$;
$c_5 = \{e_4,e_5,e_{12}\} = \langle v_1,v_5,v_6\rangle$; $c_6 = \{e_6,e_7,e_{10}\} = \langle v_2,v_5,v_3\rangle$;
$c_7 = \{e_9,e_{10},e_{11}\} = \langle v_3,v_5,v_4\rangle$; $c_0 = \{e_7,e_8,e_{12}\} = \langle v_2,v_6,v_5\rangle$.

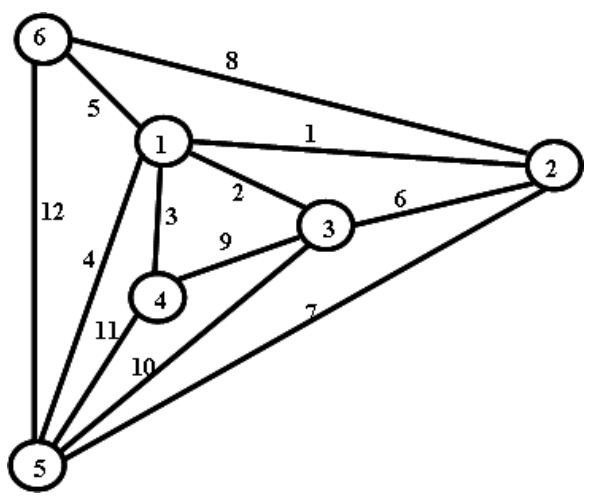

Рис. 3.11. Граф G.

Согласно шагу 2 алгоритма, просматриваем все циклы и выбираем первый цикл $c_1$ = $<e_1,e_2,e_3>$ = $<v_1,v_2,v_3>$. Строим таблицу вращений вершин по принципу помещения в строку данной вершины смежных вершин находящихся слева и справа в вершинной записи цикла:

| | | |
|---|---|---|
| вершина 1: | 3 | 2 |
| вершина 2: | 1 | 3 |
| вершина 3: | 2 | 1 |
| вершина 4: | | |
| вершина 5: | | |
| вершина 6: | | |

Выбираем второй цикл $c_2 = \{e_1,e_5,e_8\}$ = $<v_1,v_6,v_2>$ и достраиваем таблицу вращений вершин согласно ориентированной записи цикла:

| | | | |
|---|---|---|---|
| вершина 1: | 3 | 2 | 6 |
| вершина 2: | 6 | 1 | 3 |
| вершина 3: | 2 | 1 | |
| вершина 4: | | | |
| вершина 5: | | | |
| вершина 6: | 1 | 2 | |

Выбираем третий цикл $c_3 = \{e_2,e_3,e_9\}$ = $<v_1,v_3,v_4>$ и достраиваем таблицу вращений вершин:

| | | | | |
|---|---|---|---|---|
| вершина 1: | 4 | 3 | 2 | 6 |
| вершина 2: | 6 | 1 | 3 | |
| вершина 3: | 2 | 1 | 4 | |
| вершина 4: | 3 | 1 | | |
| вершина 5: | | | | |
| вершина 6: | 1 | 2 | | |

Выбираем четвертый цикл $c_4 = \{e_3,e_4,e_{11}\}$ = $<v_1,v_4,v_5>$ и достраиваем таблицу вращений вершин:

| | | | | | |
|---|---|---|---|---|---|
| вершина 1: | 5 | 4 | 3 | 2 | 6 |
| вершина 2: | 6 | 1 | 3 | | |
| вершина 3: | 2 | 1 | 4 | | |
| вершина 4: | 3 | 1 | 5 | | |
| вершина 5: | 4 | 1 | | | |
| вершина 6: | 1 | 2 | | | |

Выбираем пятый цикл $c_5 = \{e_4,e_5,e_{12}\}$ = $<v_1,v_5,v_6>$ и достраиваем таблицу вращений вершин:

| | | | | | | |
|---|---|---|---|---|---|---|
| вершина 1: | 6 | 5 | 4 | 3 | 2 | 6 |
| вершина 2: | 6 | 1 | 3 | | | |
| вершина 3: | 2 | 1 | 4 | | | |
| вершина 4: | 3 | 1 | 5 | | | |
| вершина 5: | 4 | 1 | 6 | | | |
| вершина 6: | 5 | 1 | 2 | | | |

Выбираем шестой цикл $c_6 = \{e_6,e_7,e_{10}\} = <v_2,v_5,v_3>$ и достраиваем таблицу вращений вершин:

| | | | | | | |
|---|---|---|---|---|---|---|
| вершина 1: | 4 | 3 | 2 | 6 | 5 | 4 |
| вершина 2: | 6 | 1 | 3 | 5 | | |
| вершина 3: | 5 | 2 | 1 | 4 | | |
| вершина 4: | 3 | 1 | 5 | | | |
| вершина 5: | 4 | 1 | 6 | | 2 | 3 |
| вершина 6: | 5 | 1 | 2 | | | |

Выбираем седьмой цикл $c_7 = \{e_9,e_{10},e_{11}\} = <v_3,v_5,v_4>$ и достраиваем таблицу вращений вершин, объединяя элементы для пятой вершины:

| | | | | | | |
|---|---|---|---|---|---|---|
| вершина 1: | 4 | 3 | 2 | 6 | 5 | 4 |
| вершина 2: | 6 | 1 | 3 | 5 | | |
| вершина 3: | 4 | 5 | 2 | 1 | 4 | |
| вершина 4: | 5 | 3 | 1 | 5 | | |
| вершина 5: | 2 | 3 | 4 | 1 | 6 | |
| вершина 6: | 5 | 1 | 2 | | | |

Выбираем обод $c_0 = \{e_7,e_8,e_{12}\} = <v_2,v_6,v_5>$ и окончательно достраиваем таблицу вращений вершин:

| | | | | | | |
|---|---|---|---|---|---|---|
| вершина 1: | 4 | 3 | 2 | 6 | 5 | 4 |
| вершина 2: | 6 | 1 | 3 | 5 | 6 | |
| вершина 3: | 5 | 2 | 1 | 4 | 5 | |
| вершина 4: | 5 | 3 | 1 | 5 | | |
| вершина 5: | 2 | 3 | 4 | 1 | 6 | 2 |
| вершина 6: | 5 | 1 | 2 | 5 | | |

По исчерпанию циклов, алгоритм заканчивает свою работу.

### 3.5. Алгоритм построения системы изометрических циклов для заданной таблицы вращения вершин

Если задан топологический рисунок плоского графа, то задана и таблица вращений вершин, тогда необходимо иметь алгоритм построения изометрических циклов для заданного вращения вершин. Предположим, что задано вращение вершин.

**Шаг 1**. [**Поиск непомеченной пары вершин**].

Последовательно просматриваем вершины с вращением. В таблице вращения вершин находим первую встречную непомеченную пару вершин $v_j,v_k$ для выбранной вершины $x_i$. Специальным способом помечаем выбранные вершины. Если таких вершин больше нет, то конец работы алгоритма. Иначе, образуем тройку выбранных вершин $<v_j,v_i,v_k>$. Идем на шаг

2.

**Шаг 2**. [**Формирование изометрических циклов**].

В созданном кортеже рассматриваем только три последних элемента, обозначаемые как $<\ldots,v_j,v_i,v_k>$ Ищем в списке вершины $v_k$ вершину $v_i$. Выбираем следом за $v_i$ стоящую вершину $v_l$ в циклической последовательности вершины $v_k$, помечаем вершины $v_k$ и $v_l$. Если номера первой вершины в кортеже и номер вершины $v_l$ совпадают, то образовался изометрический цикл, начинаем строить другой изометрический цикл, идя на шаг 1. Иначе, идем на шаг 3.

**Шаг 3.** [**Образование кортежа для выбранных вершин**].

Образуем кортеж выбранных вершин $<\ldots,v_j,v_i,v_k,v_l>$. Идем на шаг 2.

По окончанию работы алгоритма получаем систему изометрических циклов и обод плоского графа.

***Пример 3.4***. В качестве примера рассмотрим пошаговую работу алгоритма для топологического рисунка графа, представленного на рис. 3.13. Топологический рисунок описан следующим вращением вершин:

$v_1$: $v_4$ $v_3$ $v_2$ $v_6$ $v_5$;
$v_2$: $v_6$ $v_1$ $v_3$ $v_5$;
$v_3$: $v_5$ $v_2$ $v_1$ $v_4$;
$v_4$: $v_5$ $v_3$ $v_1$;
$v_5$: $v_2$ $v_3$ $v_4$ $v_1$ $v_6$;
$v_6$: $v_5$ $v_1$ $v_2$.

Согласно шагу 1 алгоритма, выделяем первую встречную непомеченную пару вершин $v_j = v_4$ и $v_k = v_3$ в строке для 1-ой вершины. Образуем циклический кортеж $<v_4,v_1,v_3>$, помечаем специальным образом вершины $v_j = v_4$ и $v_k = v_3$ расположив между ними вершину $v_1$

1: 4 $_1$ 3 2 6 5. Идем на шаг 2.

Здесь удобно представлять вершины только натуральными числами.

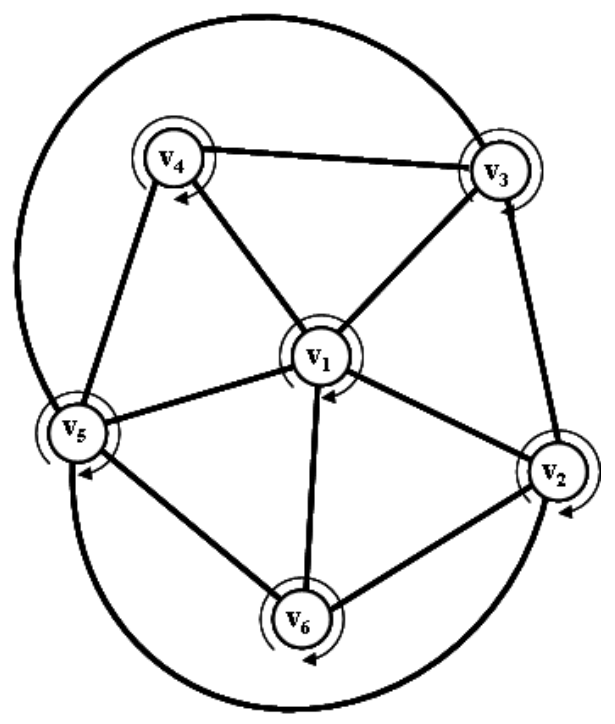

Рис. 3.17. Граф.

Согласно шагу 2 алгоритма, идем на строку для 3-ей вершины как стоящей справа после разметки. Выбираем $e_k = e_1$ и следующую за ней в циклической последовательности вершину

$e_l = e_4$. Помечаем выбранные вершины, введя вершину $e_3$ между ними

3: 5 2 1 $_{\mathbf{3}}$ 4. Идем на шаг 3.

Согласно шагу 3 алгоритма, образуем циклический кортеж $\langle v_4,v_1,v_3,v_4\rangle$. Идем на шаг 2.

Согласно шагу 2 алгоритма, идем на строку для 4-ой вершины. Выбираем $v_k = v_3$ и следующую за ней в циклической последовательности вершину $v_l = v_1$. Помечаем выбранные вершины, введя вершину $v_4$ между ними

4: 5 3 $_{\mathbf{4}}$ 1.

Здесь номера первой и последней вершин в кортеже $\langle v_4,v_1,v_3,v_4\rangle$ совпали. Тем самым образовался изометрический цикл $\langle v_4,v_1,v_3,v_4\rangle$ или $\{e_2,e_3,e_9\}$ в реберной записи. Идем на шаг 1. Теперь таблица вращений вершин имеет вид:

1: 4 $_{\mathbf{1}}$ 3 2 6 5;<br>
2: 6 1 3 5;<br>
3: 5 2 1 $_{\mathbf{3}}$ 4;<br>
4: 5 3 $_{\mathbf{4}}$ 1;<br>
5: 2 3 4 1 6;<br>
6: 5 1 2.

Согласно шагу 1 алгоритма, выделяем первую встречную непомеченную пару вершин $v_j = v_3$ и $v_k = v_2$ в строке для 1-ой вершины. Образуем кортеж $\langle v_3,v_1,v_2\rangle$, помечаем вершины $v_j = v_3$ и $v_k = v_2$ расположив между ними вершину $v_1$

1: 4 $_{\mathbf{1}}$ 3 $_{\mathbf{1}}$ 2 6 5. Идем на шаг 2.

Согласно шагу 2 алгоритма, идем на строку для 2-ой вершине. Выбираем $v_k = v_1$ и следующую за ней в циклической последовательности вершину $v_l = v_3$. Помечаем выбранные вершины, введя вершину $v_2$ между ними

2: 6 1 $_{\mathbf{2}}$ 3 5. Идем на шаг 3.

Согласно шагу 3 алгоритма, образуем кортеж $\langle v_3,v_1,v_2,v_3\rangle$. Идем на шаг 2.

Согласно шагу 2 алгоритма, идем на строку для 3-ей вершины. Выбираем $v_k = v_2$ и следующую за ней в циклической последовательности вершину $v_l = v_1$. Помечаем выбранные вершины, введя вершину $v_3$ между ними

3: 5 2 $_{\mathbf{3}}$ 1 $_{\mathbf{3}}$ 4

Здесь номера первой и последней вершин в кортеже совпали. Тем самым образовался изометрический цикл $\langle v_3,v_1,v_2,v_3\rangle$ или $\{e_1,e_2,e_6\}$. Идем на шаг 1. После проведенных операций таблица вращений вершин имеет вид:

1: 4 $_{\mathbf{1}}$ 3 $_{\mathbf{1}}$ 2 6 5;<br>
2: 6 1 $_{\mathbf{2}}$ 3 5;<br>
3: 5 2 $_{\mathbf{3}}$ 1 $_{\mathbf{3}}$ 4;<br>
4: 5 3 $_{\mathbf{4}}$ 1<br>
5: 2 3 4 1 6<br>
6: 5 1 2

Согласно шагу 1 алгоритма, выделяем первую непомеченную пару вершин $v_j = v_2$ и $v_k =$

$v_6$ в строке для 1-ой вершины. Образуем кортеж $<v_2,v_1,v_6>$, помечаем вершины $v_j = v_2$ и $v_k = v_6$ расположив между ними вершину $x_1$

1: 4 $_{\mathbf{1}}$ 3 $_{\mathbf{1}}$ 2 $_{\mathbf{1}}$ 6  5. Идем на шаг 2.

Согласно шагу 2 алгоритма, идем на строку для 6-ой вершины. Выбираем $v_k = v_l$ и следующую за ней в циклической последовательности вершину $v_l = v_2$. Помечаем выбранные вершины, введя вершину $v_6$ между ними

6: 5  1 $_{\mathbf{6}}$ 2. Идем на шаг 3.

Согласно шагу 3 алгоритма, образуем кортеж $<v_2,v_1,v_6,v_2>$. Идем на шаг 2.

Согласно шагу 2 алгоритма, идем на строку для 2-ой вершины. Выбираем $v_k = v_6$ и следующую за ней в циклической последовательности вершину $v_l = v_1$. Помечаем выбранные вершины, введя вершину $v_2$ между ними

2: 6 $_{\mathbf{2}}$ 1 $_{\mathbf{2}}$ 3  5

Здесь номера первой и последней вершин в кортеже совпали. Тем самым образовался изометрический цикл $<v_2,v_1,v_6,v_2>$ или $\{e_1,e_5,e_8\}$. Идем на шаг 1. Таблица вращений вершин имеет вид:

1: 4 $_{\mathbf{1}}$ 3 $_{\mathbf{1}}$ 2 $_{\mathbf{1}}$ 6  5;
2: 6 $_{\mathbf{2}}$ 1 $_{\mathbf{2}}$ 3  5;
3: 5  2 $_{\mathbf{3}}$ 1 $_{\mathbf{3}}$ 4;
4: 5  3 $_{\mathbf{4}}$ 1;
5: 2  3  4  1  6;
6: 5  1 $_{\mathbf{6}}$ 2.

Согласно шагу 1 алгоритма, выделяем первую встречную непомеченную пару вершин $v_j = v_6$ и $v_k = v_5$ в строке для 1-ой вершины. Образуем кортеж $<v_6,v_1,v_5>$, помечаем вершины $v_j = v_6$ и $v_k = v_5$ расположив между ними вершину $v_1$

1: 4 $_{\mathbf{1}}$ 3 $_{\mathbf{1}}$ 2 $_{\mathbf{1}}$ 6 $_{\mathbf{1}}$ 5. Идем на шаг 2.

Согласно шагу 2 алгоритма, идем на строку для 5-ой вершины. Выбираем $v_k = v_l$ и следующую за ней в циклической последовательности вершину $v_l = v_6$. Помечаем выбранные вершины, введя вершину $v_5$ между ними

5: 2  3  4  1 $_{\mathbf{5}}$ 6. Идем на шаг 3.

Согласно шагу 3 алгоритма, образуем кортеж $<v_6,v_1,v_5,v_6>$. Идем на шаг 2.

Согласно шагу 2 алгоритма, идем на 6-ую вершину. Выбираем $v_k = v_5$ и следующую за ней в циклической последовательности вершину $v_l = v_1$. Помечаем выбранные вершины, введя вершину $x_6$ между ними

6: 5 $_{\mathbf{6}}$ 1 $_{\mathbf{6}}$ 2

Здесь номера первой и последней вершин в кортеже совпали. Тем самым образовался изометрический цикл $<v_6,v_1,v_5,v_6>$ или $\{e_4,e_5,e_{12}\}$. Идем на шаг 1. Таблица вращений вершин имеет вид:

1: 4 $_{\mathbf{1}}$ 3 $_{\mathbf{1}}$ 2 $_{\mathbf{1}}$ 6 $_{\mathbf{1}}$ 5;
2: 6 $_{\mathbf{2}}$ 1 $_{\mathbf{2}}$ 3   5;
3: 5   2 $_{\mathbf{3}}$ 1 $_{\mathbf{3}}$ 4;
4: 5   3 $_{\mathbf{4}}$ 1;
5: 2   3   4   1 $_{\mathbf{5}}$ 6;
6: 5 $_{\mathbf{6}}$ 1 $_{\mathbf{6}}$ 2.

Согласно шагу 1 алгоритма, выделяем первую встречную непомеченную пару вершин $v_j = v_5$ и $v_k = v_4$ в строке для 1-ой вершины. Образуем кортеж $\langle v_5,v_1,v_4\rangle$, помечаем вершины $v_j = v_5$ и $v_k = v_4$ расположив между ними вершину $v_1$

1: 4 $_{\mathbf{1}}$ 3 $_{\mathbf{1}}$ 2 $_{\mathbf{1}}$ 6 $_{\mathbf{1}}$ 5 $_{\mathbf{1}}$. Идем на шаг 2.

Согласно шагу 2 алгоритма, идем на 4-ую вершину. Выбираем $v_k = v_1$ и следующую за ней в циклической последовательности вершину $v_l = v_5$. Помечаем выбранные вершины, введя вершину $v_4$ между ними

4: 5   3 $_{\mathbf{4}}$ 1 $_{\mathbf{4}}$. Идем на шаг 3.

Согласно шагу 3 алгоритма, образуем кортеж $\langle v_5,v_1,v_4,v_5\rangle$. Идем на шаг 2.

Согласно шагу 2 алгоритма, идем на строку для 5-ой вершины. Выбираем $v_k = v_4$ и следующую за ней в циклической последовательности вершину $v_l = v_1$. Помечаем выбранные вершины, введя вершину $v_5$ между ними

5: 2   3   4 $_{\mathbf{5}}$ 1 $_{\mathbf{5}}$ 6

Здесь номера первой и последней вершин в кортеже совпали. Тем самым образовался изометрический цикл $\langle v_5,v_1,v_4,v_5\rangle$ или $\{e_3,e_4,e_{11}\}$. Идем на шаг 1. Таблица вращений вершин имеет вид:

1: 4 $_{\mathbf{1}}$ 3 $_{\mathbf{1}}$ 2 $_{\mathbf{1}}$ 6 $_{\mathbf{1}}$ 5 $_{\mathbf{1}}$;
2: 6 $_{\mathbf{2}}$ 1 $_{\mathbf{2}}$ 3   5;
3: 5   2 $_{\mathbf{3}}$ 1 $_{\mathbf{3}}$ 4;
4: 5   3 $_{\mathbf{4}}$ 1 $_{\mathbf{4}}$;
5: 2   3   4 $_{\mathbf{5}}$ 1 $_{\mathbf{5}}$ 6;
6: 5 $_{\mathbf{6}}$ 1 $_{\mathbf{6}}$ 2.

Согласно шагу 1 алгоритма, выделяем первую непомеченную пару вершин $v_j = v_3$ и $v_k = v_5$ в строке для 2-ой вершины. Образуем кортеж $\langle v_3,v_2,v_5\rangle$, помечаем вершины $v_j = v_5$ и $v_k = v_4$ расположив между ними вершину $v_2$

2: 6 $_{\mathbf{2}}$ 1 $_{\mathbf{2}}$ 3 $_{\mathbf{2}}$ 5. Идем на шаг 2.

Согласно шагу 2 алгоритма, идем на строку для 5-ой вершины. Выбираем $v_k = v_2$ и следующую за ней в циклической последовательности вершину $v_l = v_3$. Помечаем выбранные вершины, введя вершину $v_5$ между ними

5: 2 $_{\mathbf{5}}$ 3   4 $_{\mathbf{5}}$ 1 $_{\mathbf{5}}$ 6. Идем на шаг 3.

Согласно шагу 3 алгоритма, образуем кортеж $\langle v_3,v_2,v_5,v_3\rangle$. Идем на шаг 2.

Согласно шагу 2 алгоритма, идем на строку для 3-ей вершины. Выбираем $v_k = v_5$ и следующую за ней в циклической последовательности вершину $v_l = v_2$. Помечаем выбранные вершины, введя вершину $v_3$ между ними

3: $5_{3}\ 2_{3}\ 1_{3}\ 4$

Здесь номера первой и последней вершин в кортеже совпали. Тем самым образовался изометрический цикл $<v_3,v_2,v_5,v_3>$ или $\{u_6,u_7,u_{10}\}$. Идем на шаг 1. Таблица вращений вершин имеет вид:

1: $4_{1}\ 3_{1}\ 2_{1}\ 6_{1}\ 5_{1}$;<br>
2: $6_{2}\ 1_{2}\ 3_{2}\ 5$;<br>
3: $5_{3}\ 2_{3}\ 1_{3}\ 4$;<br>
4: $5\ \ 3_{4}\ 1_{4}$;<br>
5: $2_{5}\ 3\ \ 4_{5}\ 1_{5}\ 6$;<br>
6: $5_{6}\ 1_{6}\ 2$.

Согласно шагу 1 алгоритма, выделяем первую встречную непомеченную пару вершин $v_j = v_5$ и $v_k = v_6$ в строке для 2-ой вершины. Образуем кортеж $<v_5,v_2,v_6>$, помечаем вершины $v_j = v_5$ и $v_k = v_6$ расположив между ними вершину $x_2$

2: $6_{2}\ 1_{2}\ 3_{2}\ 5_{2}$. Идем на шаг 2.

Согласно шагу 2 алгоритма, идем на строку для 6-ой вершины. Выбираем $v_k = v_2$ и следующую за ней в циклической последовательности вершину $v_l = v_5$. Помечаем выбранные вершины, введя вершину $v_6$ между ними

6: $5_{6}\ 1_{6}\ 2_{6}$. Идем на шаг 3.

Согласно шагу 3 алгоритма, образуем кортеж $<v_5,v_2,v_6,v_5>$. Идем на шаг 2.

Согласно шагу 2 алгоритма, идем на строку для 5-ой вершины. Выбираем $v_k = v_6$ и следующую за ней в циклической последовательности вершину $v_l = v_2$. Помечаем выбранные вершины, введя вершину $v_5$ между ними

5: $2_{5}\ 3\ \ 4_{5}\ 1_{5}\ 6_{5}$.

Здесь номера первой и последней вершин в кортеже совпали. Тем самым образовался изометрический цикл $<v_5,v_2,v_6,v_5>$ или $\{e_7,e_8,e_{12}\}$. Идем на шаг 1. Таблица вращений вершин имеет вид:

1: $4_{1}\ 3_{1}\ 2_{1}\ 6_{1}\ 5_{1}$;<br>
2: $6_{2}\ 1_{2}\ 3_{2}\ 5_{2}$;<br>
3: $5_{3}\ 2_{3}\ 1_{3}\ 4$<br>
4: $5\ \ 3_{4}\ 1_{4}$;<br>
5: $2_{5}\ 3\ \ 4_{5}\ 1_{5}\ 6_{5}$;<br>
6: $5_{6}\ 1_{6}\ 2_{6}$.

Согласно шагу 1 алгоритма, выделяем первую встречную непомеченную пару вершин $v_j = v_4$ и $v_k = v_5$ в строке для 3-ей вершины. Образуем кортеж $<v_4,v_3,v_5>$, помечаем вершины $v_j = v_4$ и $v_k = v_5$ расположив между ними вершину $v_3$

3: $5_{3}\ 2_{3}\ 1_{3}\ 4_{3}$. Идем на шаг 2.

Согласно шагу 2 алгоритма, идем на строку для 5-ой вершины. Выбираем $v_k = v_3$ и следующую за ней в циклической последовательности вершину $v_l = v_4$. Помечаем выбранные вершины, введя вершину $v_5$ между ними

$5:\ 2_{\mathbf{5}}\ 3_{\mathbf{5}}\ 4_{\mathbf{5}}\ 1_{\mathbf{5}}\ 6_{\mathbf{5}}$. Идем на шаг 3.

Согласно шагу 3 алгоритма, образуем кортеж $\langle v_4,v_3,v_5,v_4 \rangle$. Идем на шаг 2.

Согласно шагу 2 алгоритма, идем на строку для 4-ой вершины. Выбираем $v_k = v_5$ и следующую за ней в циклической последовательности вершину $v_l = v_3$. Помечаем выбранные вершины, введя вершину $v_4$ между ними

$4:\ 5_{\mathbf{4}}\ 3_{\mathbf{4}}\ 1_{\mathbf{4}}$.

Здесь номера первой и последней вершин в кортеже совпали. Тем самым образовался изометрический цикл $\langle v_4,v_3,v_5,v_4 \rangle$ или $\{e_9,e_{10},e_{11}\}$. Идем на шаг 1. Таблица вращений вершин имеет вид:

$1:\ 4_{\mathbf{1}}\ 3_{\mathbf{1}}\ 2_{\mathbf{1}}\ 6_{\mathbf{1}}\ 5_{\mathbf{1}}$;
$2:\ 6_{\mathbf{2}}\ 1_{\mathbf{2}}\ 3_{\mathbf{2}}\ 5_{\mathbf{2}}$;
$3:\ 5_{\mathbf{3}}\ 2_{\mathbf{3}}\ 1_{\mathbf{3}}\ 4_{\mathbf{3}}$;
$4:\ 5_{\mathbf{4}}\ 3_{\mathbf{4}}\ 1_{\mathbf{4}}$;
$5:\ 2_{\mathbf{5}}\ 3_{\mathbf{5}}\ 4_{\mathbf{5}}\ 1_{\mathbf{5}}\ 6_{\mathbf{5}}$;
$6:\ 5_{\mathbf{6}}\ 1_{\mathbf{6}}\ 2_{\mathbf{6}}$.

Согласно шагу 1 алгоритма, непомеченных вершин в таблице вращения вершин нет. Конец работы алгоритма

Итак, выделены следующие изометрические циклы:

$c_1 = \{e_1,e_2,e_6\} \rightarrow \langle v_3,v_1,v_2 \rangle$; $c_2 = \{e_1,e_5,e_8\} \rightarrow \langle v_2,v_1,v_6 \rangle$;
$c_3 = \{e_2,e_3,e_9\} \rightarrow \langle v_4,v_1,v_3 \rangle$; $c_4 = \{e_3,e_4,e_{11}\} \rightarrow \langle v_5,v_1,v_4 \rangle$;
$c_5 = \{e_4,e_5,e_{12}\} \rightarrow \langle v_6,v_1,v_5 \rangle$; $c_6 = \{e_6,e_7,e_{10}\} \rightarrow \langle v_3,v_2,v_5 \rangle$;
$c_7 = \{e_9,e_{10},e_{11}\} \rightarrow \langle v_5,v_4,v_3 \rangle$; $c_0 = \{e_7,e_8,e_{12}\} \rightarrow \langle v_5,v_2,v_6 \rangle$.

Из теоремы Маклейна [] следует, что любое неориентированное ребро в плоском рисунке графа можно представить в виде двух ориентированных ребер, например ребро $e_4$ для графа на рис. 3.16 можно записать $e_4 = (v_1,v_5)$ или $(v_5,v_1)$. Тогда сложение двух ориентированных циклов можно представит в виде сложения ориентированных ребер, например:

$c_4 \oplus c_5 = \{e_3,e_4,e_{11}\} \oplus \{e_4,e_5,e_{12}\} \rightarrow \langle v_1,v_4,v_5 \rangle + \langle v_1,v_5,v_6 \rangle \rightarrow (v_1,v_4) +$
$+ (v_4,v_5) + (\underline{v}_5,\underline{v}_1) + (\underline{v}_1,\underline{v}_5) + (v_5,v_6) + (v_6,v_1) \rightarrow (v_1,v_4) + (v_4,v_5) +$
$+ (v_5,v_6) + (v_6,v_1) = \langle v_1,v_4,x_5,v_6 \rangle$.

В этой операции разнонаправленные ориентированные ребра уничтожаются.

## Комментарии

Современные методы теории графов описывают рисунок графа при помощи диаграмм с точностью до изоморфизма. При этом изображение графа строится геометрическими методами используя метрические свойства представления. Однако современные средства визуализации должны предоставлять наряду с геометрическими методами топологические способы описания рисунка графа.

Теория вращения вершин, предложенная американскими математиками Дж. Янгсом и Герхардом Рингелем позволяет создать математическую модель для топологического описания рисунка графа. Причем данная модель позволяет связать центральные разрезы и простые циклы в единое целое. Так как вращение вершин в графе индуцирует (порождает) простые циклы и обратно индуцированная система циклов порождает вращение вершин. Если индуцированная система циклов удовлетворяет комбинаторному условию теоремы Маклейна, то граф планарен, то есть он может быть представлен на плоскости без пересечения ребер.

Приводятся формулы для вычисления значения функционала Маклейна, представляющего собой интегральную характеристику рассматриваемой системы циклов. Введение дополнительных (мнимых) вершин характеризующих местоположение двух пересекающихся ребер позволяет описывать топологические рисунки непланарных графов, используя методы теории вращения вершин. Ненулевое значение функционала Маклейна говорит о том, что рассматриваемая система циклов описывает непланарный суграф графа G, а нулевое значение функционала описывает планарный суграф. Установлено влияние элементов матрицы Грамма составленной из скалярных произведений суграфов на свойства интегральной характеристики системы циклов.

Представлены алгоритмы построения вращений вершин графа для заданной системы циклов и построения системы изометрических циклов для заданной таблицы вращения вершин.

Следует заметить, что вращение вершин индуцирует простые циклы. В свою очередь, система независимых простых циклов с нулевым значением функционала Маклейна индуцирует вращение вершин графа. Таким образом, устанавливается природная связь между подпространством разрезов S(G) и подпространством циклов C(G) пространства суграфов.

## Глава 4 . Векторная алгебра пересечений

### 4.1. Пересечение оьтрезков

Рассмотрим следующую задачу. Пусть заданы два отрезка характеризующие связи в пространстве $R^2$ Первый отрезок описывается концевыми точками *a* и *b*, второй отрезок описывается концевыми точками *c* и *d*. Пусть заданы координаты точек. Каждая пара точек характеризует открытое множество точек, где точки являются границами интервала. Такое открытое множество точек будем называть отрезком или соединением. На рис. 4.1 отрезки пересекаются, а на рис. 4.2 они не пересекаются. Нас интересует ответ на следующий вопрос: существует ли условие определяющее пересечения или не пересечения отрезков.

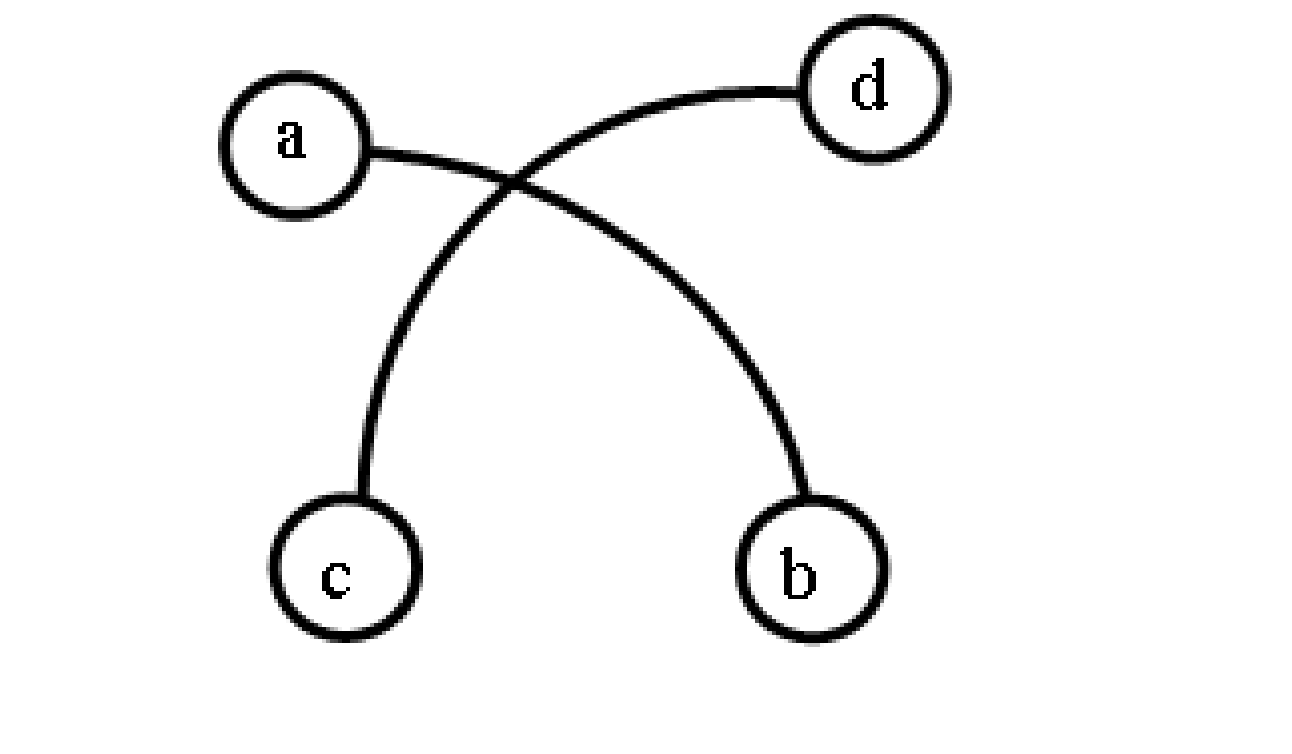

Рис. 4.1. Пересечение  соединений.

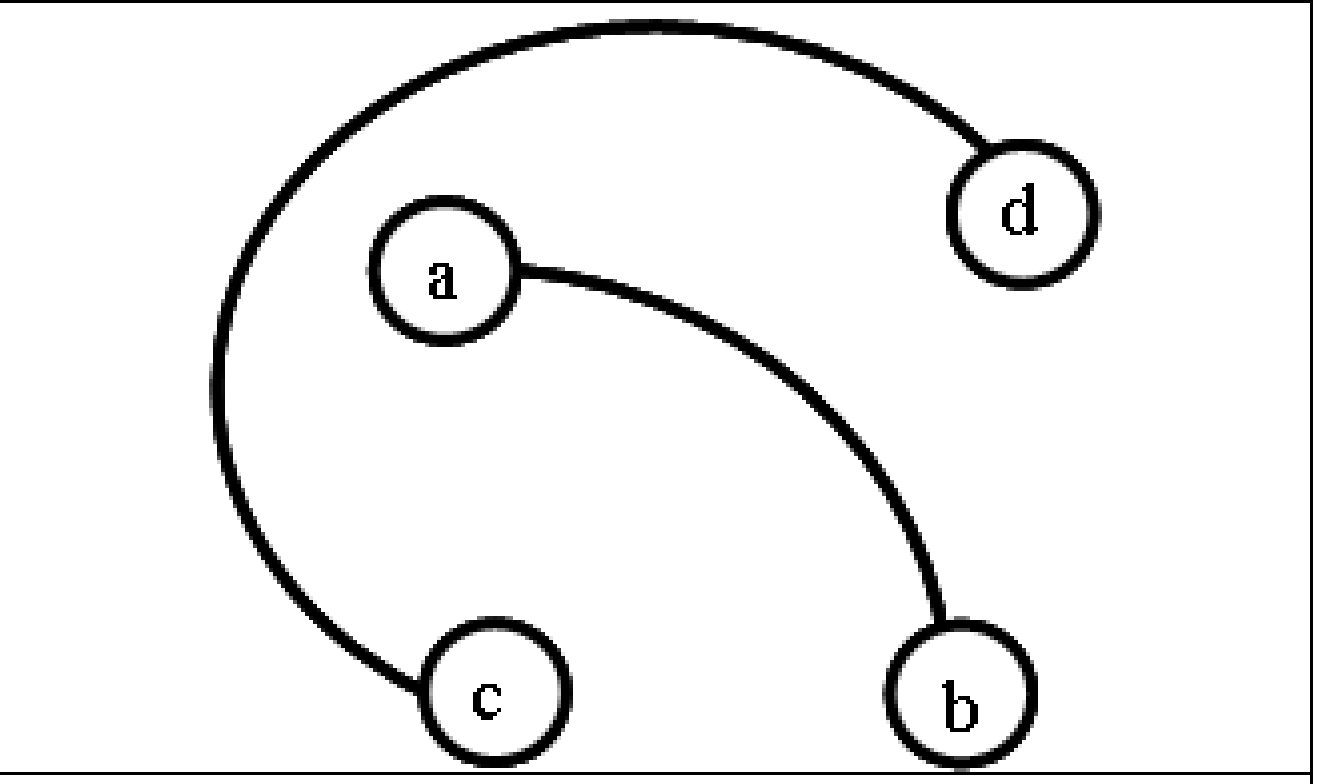

Рис. 4.2. Соединения не пнрнснкается

Очевидно, что для определения пересечения следует решить два уравнения описывающих отрезки. Если решение определяет координаты точки принадлежащгй отрезкам, то отрезки пересекаются. Иначе, отрезки не пересекаются. Но для этого решения нужно перебрать множество функций описывающих отрезки. Это очень трудоемкая задача.

Представим отрезки в виде векторов с заданными координанами. И рассмотрим пересечение двух отрезков как результат векторного произведения.

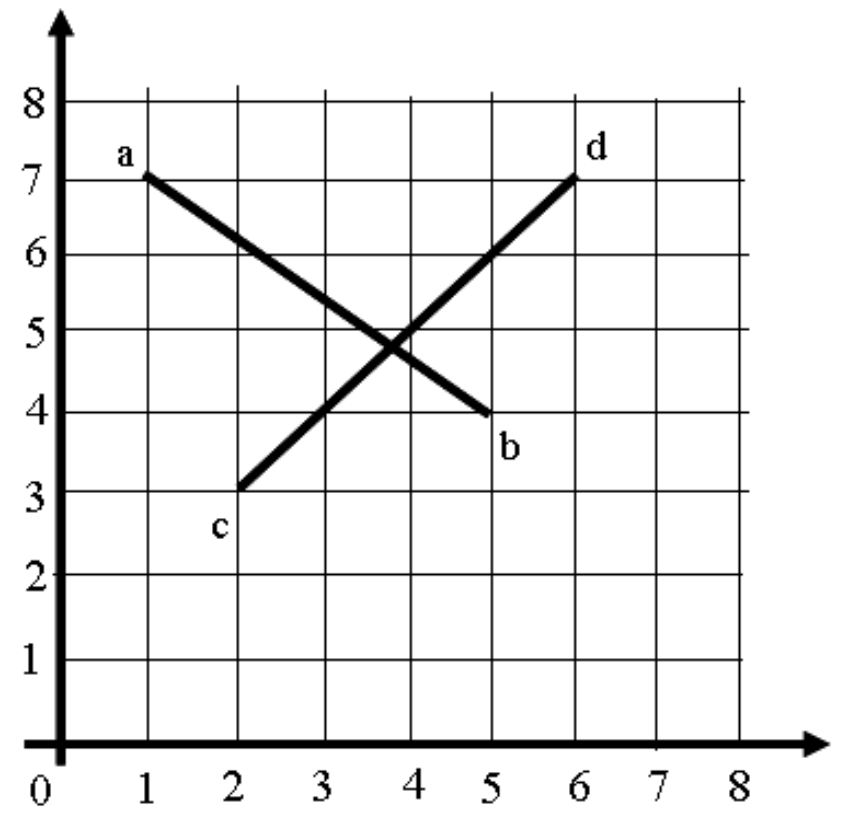

Воспользуемся свойствами векторного произведения векторов. Будем перемножать следующие векторы:

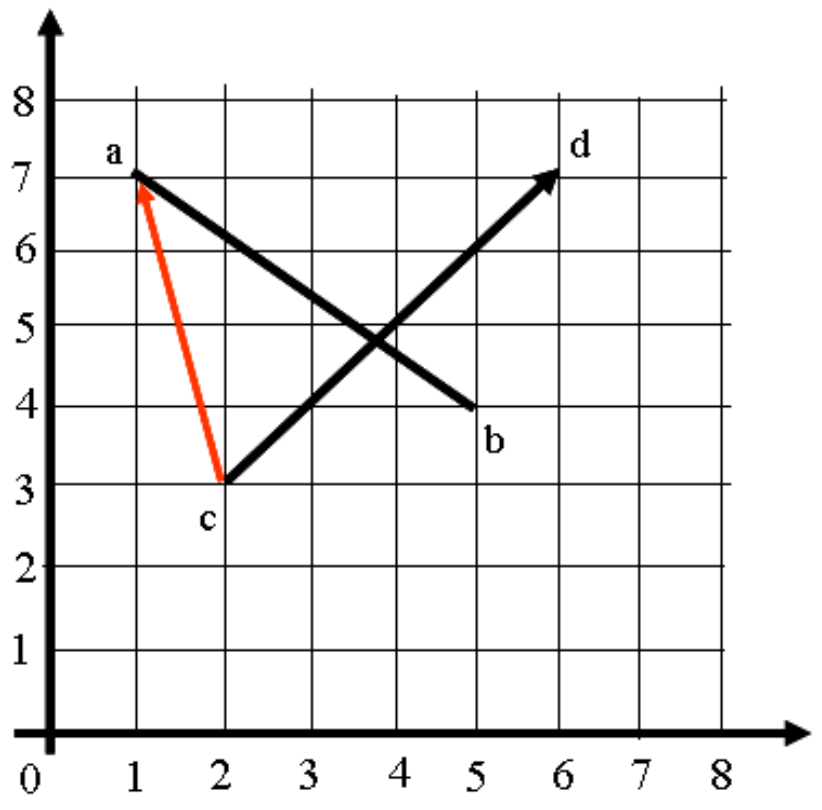

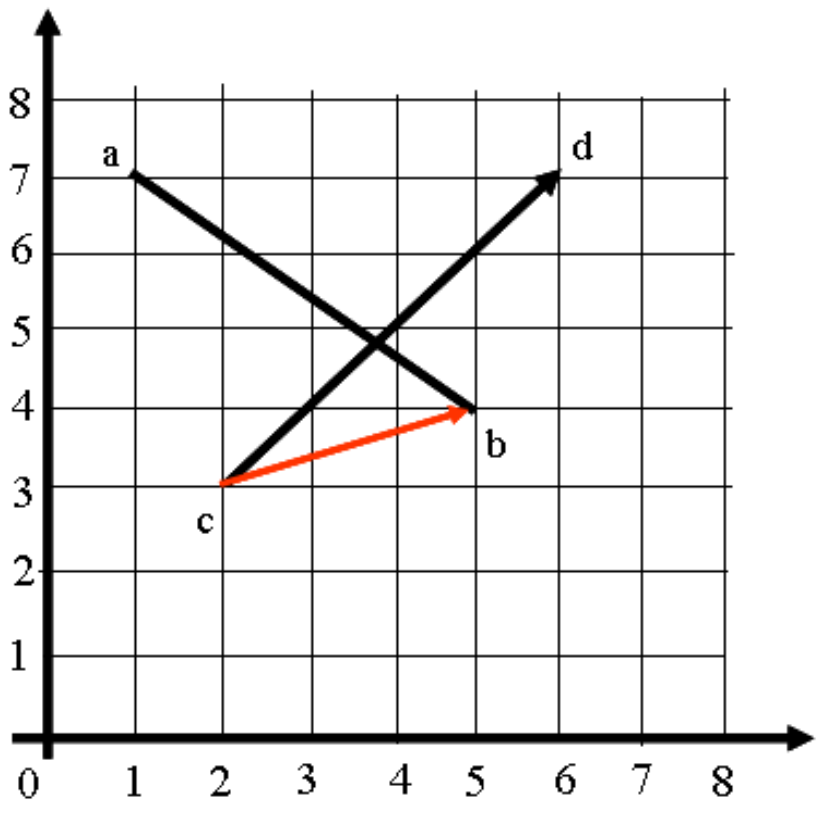

$$\vec{cd}\times\vec{ca}=\begin{Vmatrix} x_d-x_c & x_a-x_c \\ y_d-y_c & y_a-y_c \end{Vmatrix}=\begin{Vmatrix} 6-2 & 1-2 \\ 7-3 & 7-3 \end{Vmatrix}=16+4=+20\,;$$

$$\vec{cd}\times\vec{cb}=\begin{Vmatrix} x_d-x_c & x_b-x_c \\ y_d-y_c & y_b-y_c \end{Vmatrix}=\begin{Vmatrix} 6-2 & 4-2 \\ 7-3 & 4-3 \end{Vmatrix}=4-8=-4\,;$$

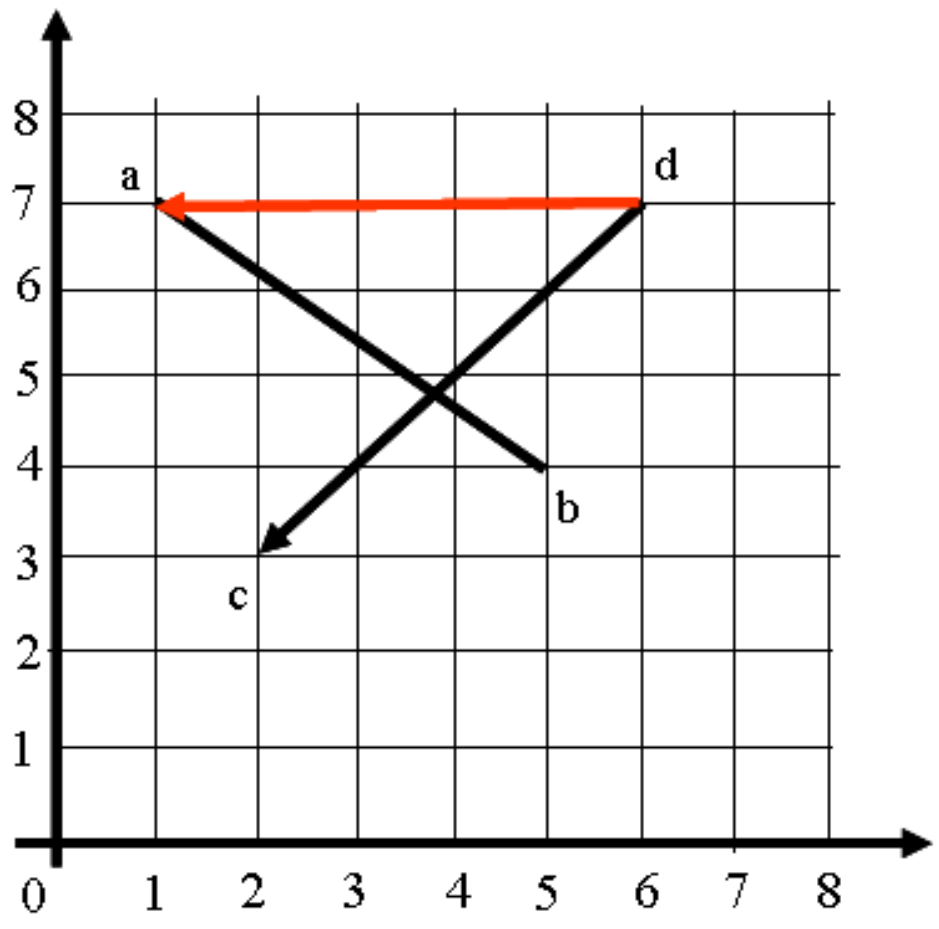

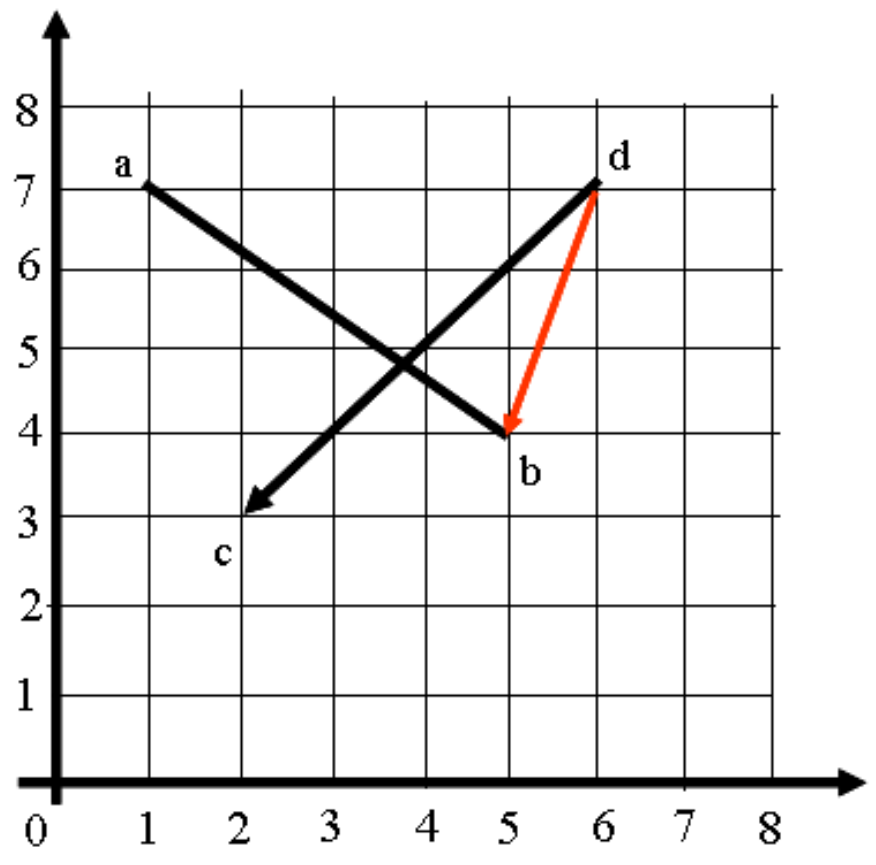

$$\vec{dc}\times\vec{da}=\begin{Vmatrix} x_c-x_d & x_a-x_d \\ y_c-y_d & y_a-y_d \end{Vmatrix}=\begin{Vmatrix} 2-6 & 1-6 \\ 3-7 & 7-7 \end{Vmatrix}=0-20=-20\,;$$

$$\vec{dc}\times\vec{db}=\begin{Vmatrix} x_c-x_d & x_b-x_d \\ y_c-y_d & y_b-y_d \end{Vmatrix}=\begin{Vmatrix} 2-6 & 5-6 \\ 3-7 & 4-7 \end{Vmatrix}=12-4=+8\,.$$

Векторное произведение первых двух векторов имеет различные знаки. Векторное произведение следующих двух векторов также имеет различные знаки. Следовательно, отрезки пересекаются. Если хотя бы одна пара векторов имеет одинаковые знаки, то отрезки не пересекаются.

Очевидно, что векторное произведение отрезков заданных векторами, также не определяет свойсва пересечения отрезков. Очевидно, нежно иметь другую математическую модель для описания этого явления.

**Определение 4.1**. Будем считать отрезки не пересекаемыми, если они одновременно присутствуют в топологическом рисунке максимально плоского суграфа.

Например, если представить все множество отрезков для точек *a*,*b*,*c*,*d* в виде полного графа $K_4$ (см. рис. 4.3). то в топологическом рисунке максимально плоского суграфа графа $K_4$ (см. рис. 4.4) одновременно присутствуют ребра (*a*,*b*) и (*c*,*d*).

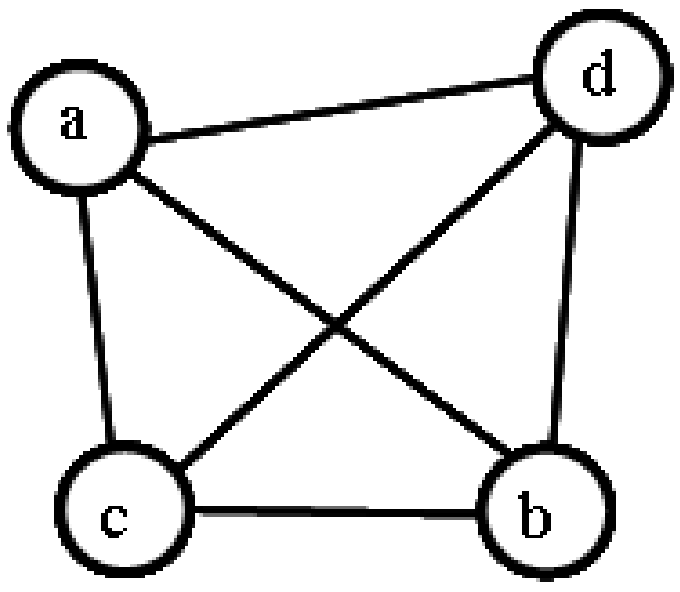

Рис. 4.3. Полный граф $K_4$.

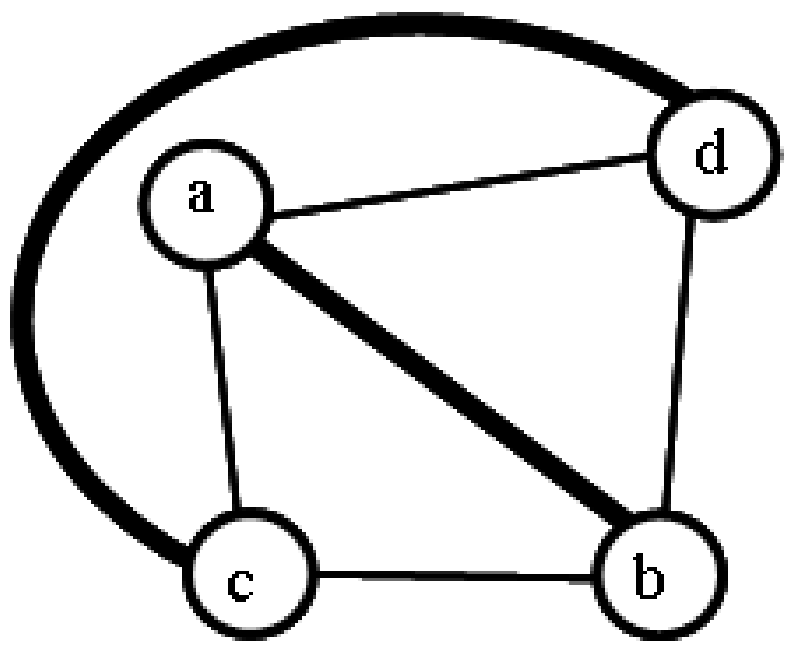

Рис. 4.4. Максимально плоский суграф графа $K_4$.

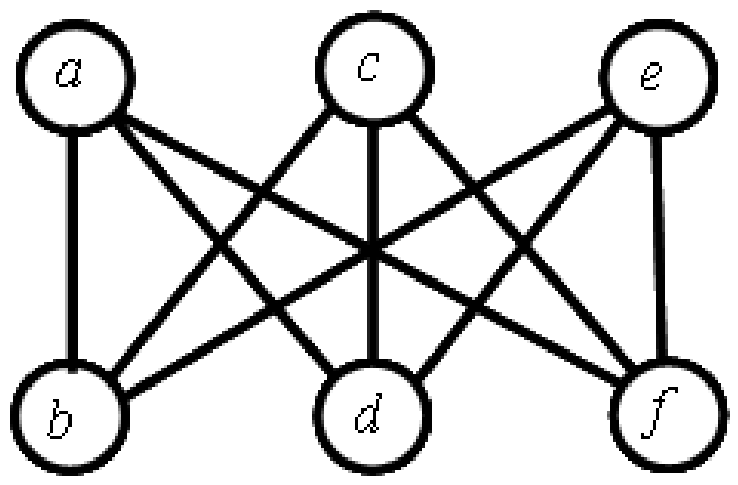

Рис. 4.5. Полный граф $K_{3,3}$

Рассмотрим задачу «три дома т три колодца», в которой требуется провести дорожки от каждого домика к каждому колодцу без перечечения.

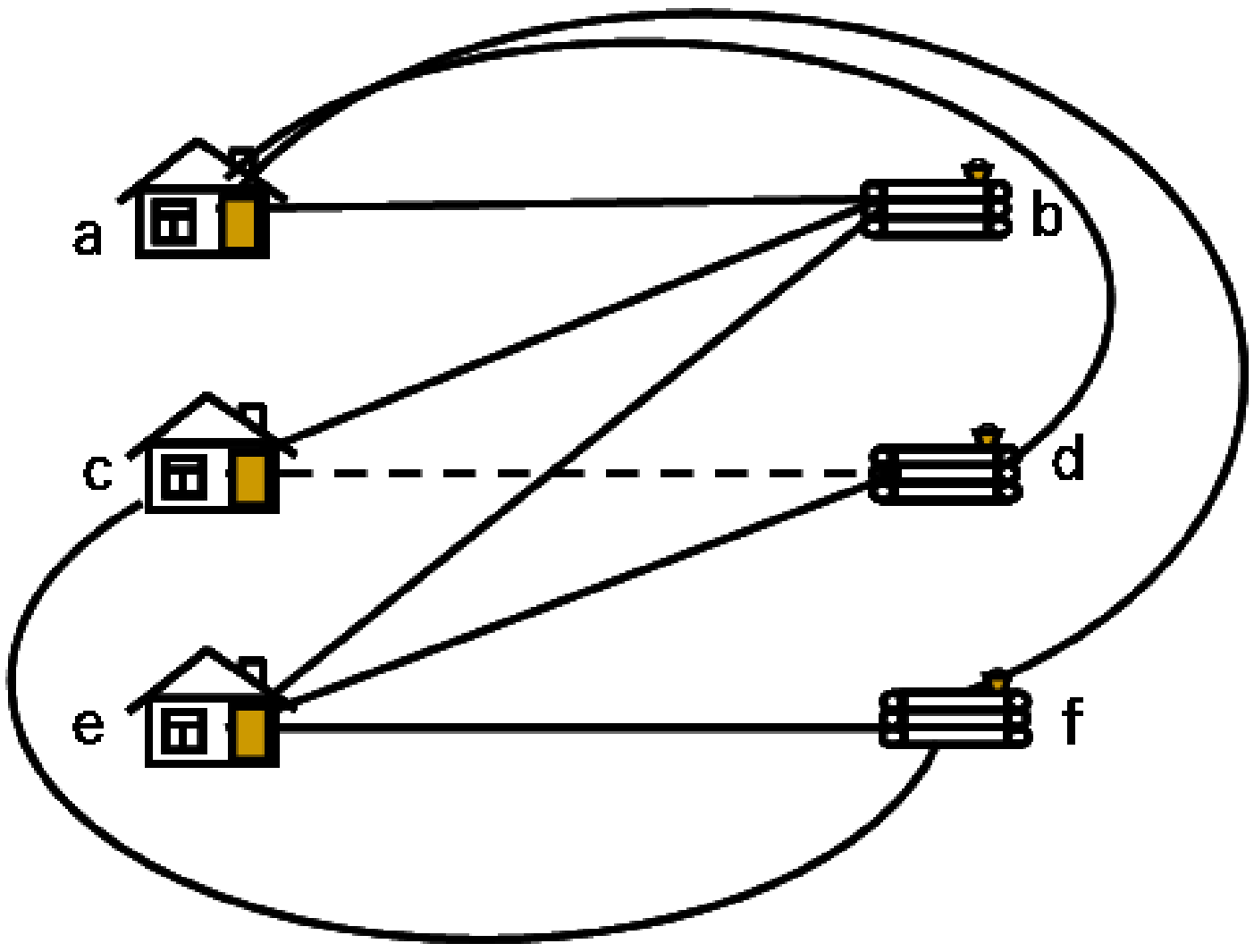

Рис. 4.6. Задача «три дома и три колодца».

Задачу «три дома и три колодца» можно представить в виде полного двудольного графа $K_{3,3}$ (см. рис. 2.5 и рис. 2.6). Известно, что для полного двудольного графа $K_{3,3}$ не существует топологического рисунка максимально плоского суграфа с участием всех ребер графа (см. рис. 4.7).

**4.2. Координатно-базисныя система векторов**

Существует и другая математическая модель для определения пересечения соединений.

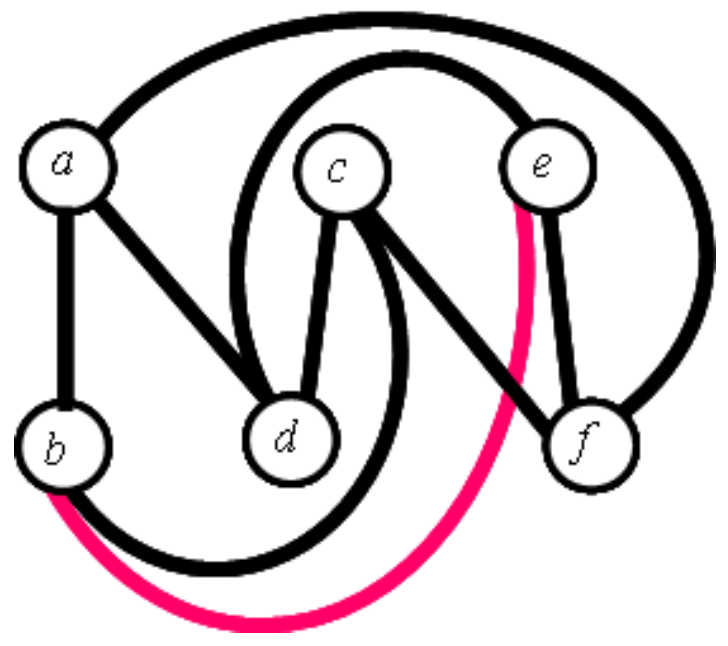

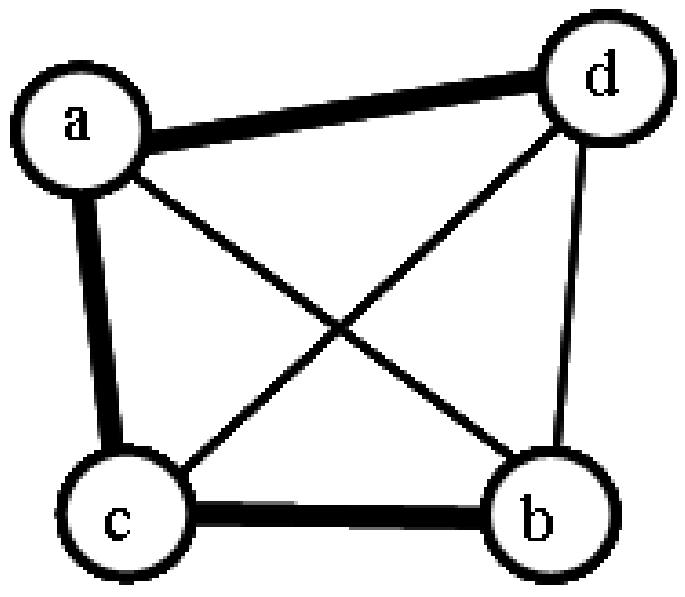

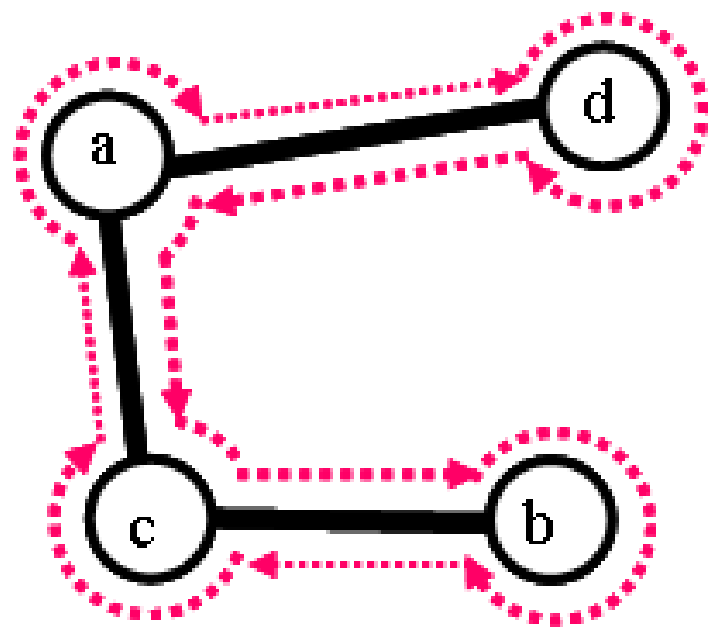

Рис. 4.7. Максимально плоский суграфа графа $K_{3,3}$.

Рис. 4.8. Дерево графа $K_4$.

Рис. 4.9. Обход дерева графа.

Будем рассматривать способ проведения соединений на основе векторной алгебры пересечений [34-36]. С этой целью, построим полный граф $K_4$ (см. рис. 4.8). Выделим произвольное дерево графа. Произведем циклический обход дерева графа (см. рис. 4.9). Циклический обход дерева индуцирует (порождает) замкнутую систему векторов (см. рис. 4.10 и рис. 4.11). Будем называть такую замкнутую систему ***координатно-базисной системой векторов*** (сокращенно КБС).

Существует два вида изображения координатно-базисной системы векторов. Изображение координатно-базисной системы в виде ориентированных петель не предполагает разделения вершин (точек) на части (см. рис. 4.10), такое представление КБС будем называть *веерным*. Изображение замкнутой координатно-базисной системы векторов (см. рис. 4.11), такое представление предполагает разделение вершин (точек), такое представление будем называть *циклическим*.

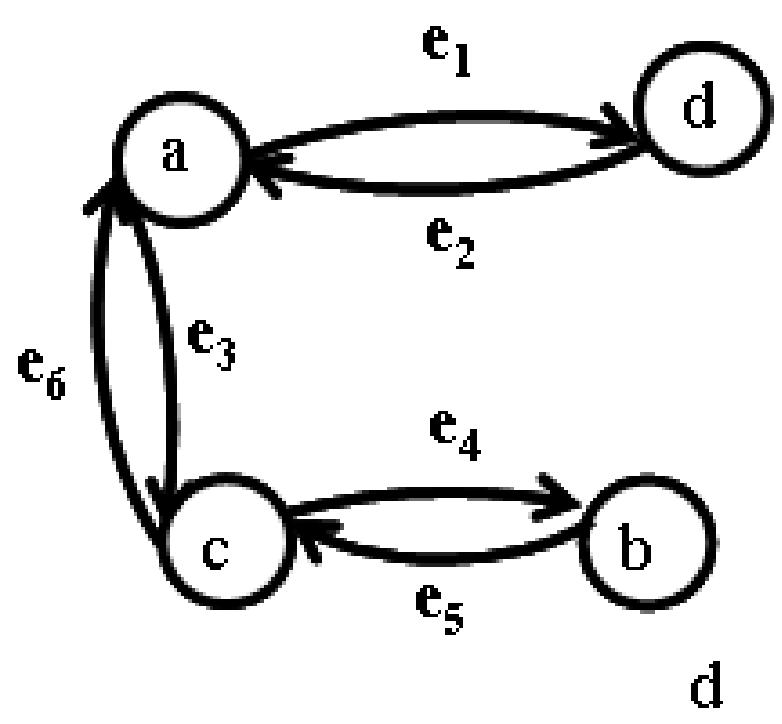

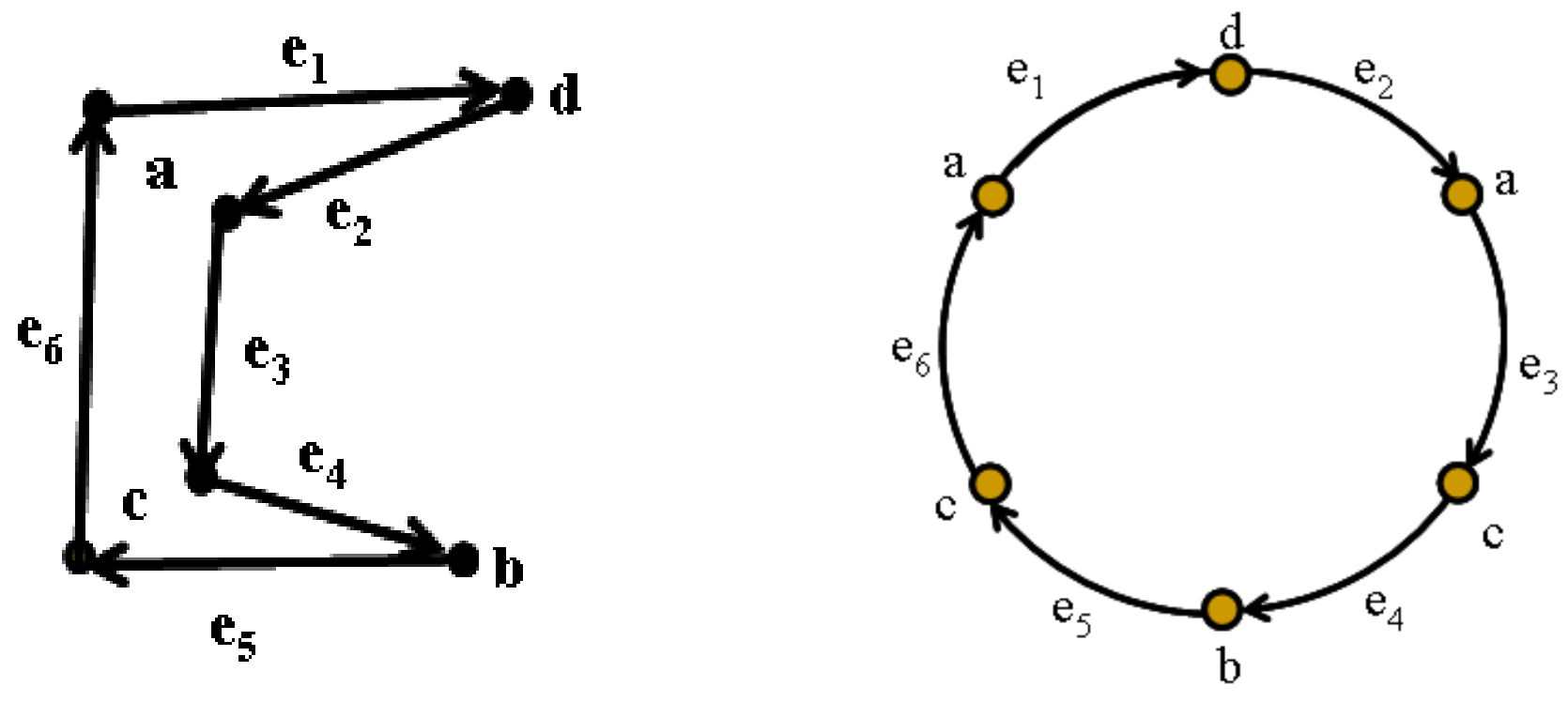

Рис. 4.10. Веерное представление КБС.

Рис. 4.11. Циклическое представление координатно-базисной системы векторов.

При построении рисунка графа анализ отношения пересечения ребер производится в топологическом пространстве, в котором метрические свойства не определены. Поэтому для полного и непротиворечивого описания рисунка нужно определить координатно-базисную систему и установить проекции всех соединений на координатно-базисную систему с целью определения пересечения ребер по их проекциям.

**Определение 4.2.** *Координатно-базисной системой* (*КБС*) будем называть систему базисных ориентированных векторов (орт), полученную в результате циклического обхода ребер произвольно выбранного дерева графа G.

**Определение 4.3.** Прообраз вершин в координатно-базисной системе будем называть *полюсом*.

Будем учитывать тот факт, что неориентированное ребро исходного графа G можно представлять двумя разнонаправленными ориентированными дугами $q_1$ и $q_2$.

**Определение 4.4.** Под *ориентированной петлей* будем понимать два разнонаправленных ориентированных вектора, соединяющие пару вершин.

Например, для нашего примера (см. рис. 4.10 и рис. 4.11), последовательность векторов $<e_1,e_2,e_3,e_4,e_5,e_6>$ представляет собой координатно-базисную систему (КБС). Ребра $\{e_1,e_2\},\{e_3,e_4\},\{e_5,e_6\}$ – ориентированные петли.

Соединение (a,b) представляется двумя дугами Соединение (*a,b*) представляется двумя дугами $(\overline{a,b})$ и $(\overline{b,a})$ имеющими соответствующие проекции на КБС в виде подмножеств $\{e_1,e_2,e_3,e_4\}$ и $\{e_5,e_6\}$. $(\overline{a,b})$ и $(\overline{b,a})$ имеющими соответствующие проекции на КБС в виде подмножеств $\{e_1,e_2,e_3,e_4\}$ и $\{e_5,e_6\}$. Соединение (*c,d*) представляется двумя дугами $(\overline{c,d})$ и $(\overline{d,c})$ имеющими соответствующие проекции на КБС в виде подмножеств $\{e_6,e_1\}$ и $\{e_2,e_3,e_4,e_5\}$ (см. рис. 4.13).

**Определение 4.5**. Объединение проекций соединения равно координатно-базисной системе векторов.

$$np(\overline{a,b}) \cup np(\overline{b,a}) = \{e_1,e_2,e_3,e_4\} \cup \{e_5,e_6\} = \{e_1,e_2,e_3,e_4,e_5,e_6\}.$$

$$np(\overline{c,d}) \cup np(\overline{d,c}) = \{e_6,e_1\} \cup \{e_2,e_3,e_4,e_5\} = \{e_1,e_2,e_3,e_4,e_5,e_6\}.$$

**Определение 4.6.** Будем говорить, что существует *пересечение двух топологических векторов* в рисунке неориентированного графа, если существует непустое пересечение (в множественном смысле) проекций их двух топологических векторов на топологическую координатно-базисную систему и не одно из этих двух проекций не является подмножеством другого

$$U_i \cap U_j \neq \varnothing \text{ ; } U_i \not\subset U_j \vee U_j \not\subset U_i. \quad (4.1)$$

Здесь $U_i$ - проекция дуги одного соединения $u_i$, $U_j$ - проекция дуги другого соединения $u_j$.

Рассмотрим проекции соединений (*a,b*) и (*c,d*) на координатно-базисную систему векторов (см. рис. 4.12). Удалим ориентированные петли из проекций соединений и результат запишем в множественном виде:

пр($a$,$b$) = {$e_3$,$e_4$};
пр($b$,$a$) = {$e_5$,$e_6$};
пр($c$,$d$) = {$e_2$,$e_3$};
пр($d$,$c$) = {$e_6$,$e_1$}.

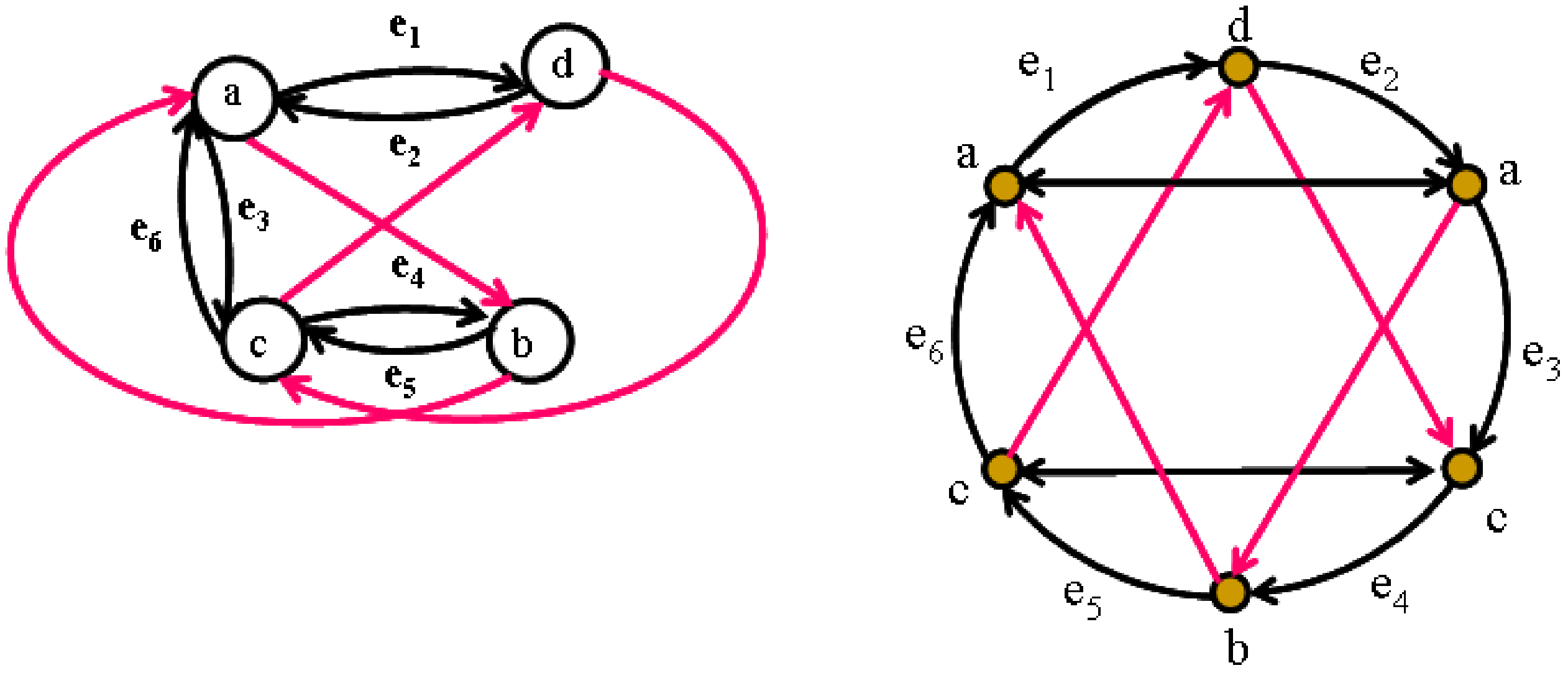

Рис. 4.12. Проекция соединений на координатно-базисную систему.

Как видно связи $\{a,b\}$ и $\{c,d\}$ пересекаются, так как пересекаются их проекции

$\{a,b\} \cap \{c,d\}$ = пр($a$,$b$) $\cap$ пр($c$,$d$) = $\{e_3,e_4\}\{e_2,e_3\}=\{e_3\}$;

а вот соединения $\{a,b\}$ и $\{d,c\}$ не пересекаются, так как пересечение проекций пусто

$\{a,b\} \cap \{d,c\}$ = пр($a$,$b$) $\cap$ пр($d$,$c$) = $\{e_3,e_4\}\{e_6,e_1\}=\varnothing$.

## 4.3. Преобразование векторов при смене координатно-базисной системы

При смене координатно-базисной системы преобразование координат осуществляется как преобразование тензоров 2-рода.

Построим координатно-базисную систему и проекцию соединений для графа $K_4$ при смене топологической координатно-базисной системы. В качестве основы новой топологической координатно-базисной системы положим дерево (см. рис. 4.14).

Матрица проекций старой топологической координатно-базисной системы на новую топологическую координатно-базисную систему имеет вид:

Матрица проекций векторов старой КБС на новую координатно-базисную систему.

| | $x_1$ | $x_2$ | $x_3$ | $x_4$ | $x_5$ | $x_6$ |
|---|---|---|---|---|---|---|
| $e_1$ | 1 | 1 | | | | |
| $e_2$ | | | 1 | 1 | 1 | 1 |
| $e_3$ | 1 | | | | | |
| $e_4$ | | | | 1 | | |
| $e_5$ | | | | | 1 | |
| $e_6$ | | | | | | 1 |

Результат получения проекций на новую КБС представлен на рис. 4.15.

Таким образом, координатно-базисная система (КБС) служит визуальной основой для построения рисунка графа на плоскости. Преобразование координат векторов при смене координатно-базисной системы производится по законам преобразования тензоров 2-го рода, используя сложение по модулю 2. Отсюда можно определить систему простых циклов для построения диаграмму вращения вершин описывающую топологический рисункок графа.

Матрица проекций соединений на КБС

| | $e_1$ | $e_2$ | $e_3$ | $e_4$ | $e_5$ | $e_6$ |
|---|---|---|---|---|---|---|
| **(a,b)** | | | **1** | **1** | | |
| **(b,a)** | | | | | **1** | **1** |
| **(c,d)** | **1** | | | | | **1** |
| **(d,c)** | | **1** | **1** | | | |

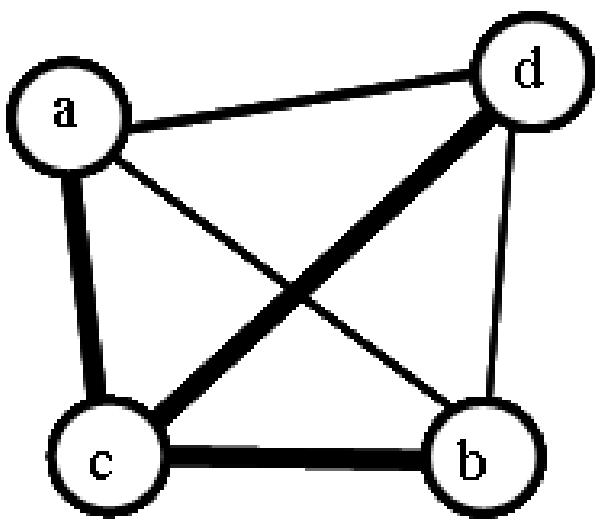

Рис. 4.13. Дерево графа.

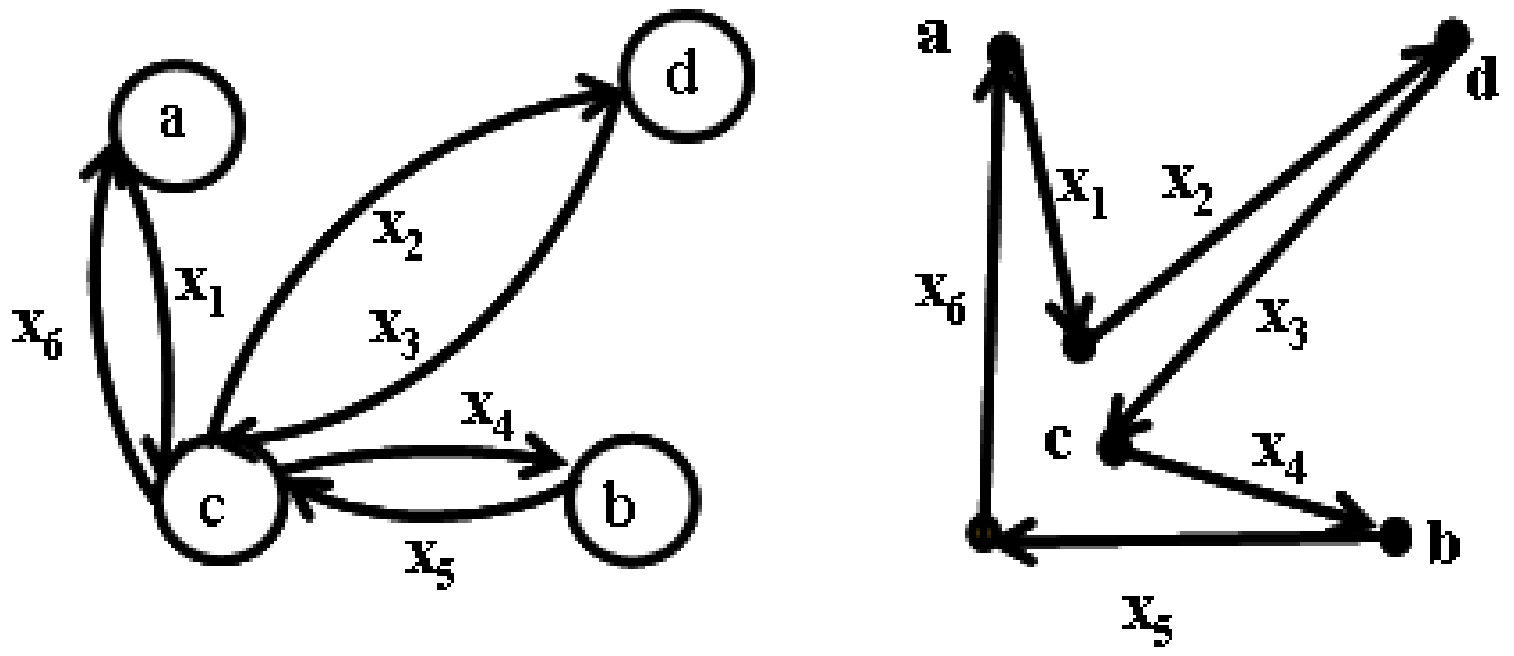

Рис. 4.14. Новое координатно-базисная система.

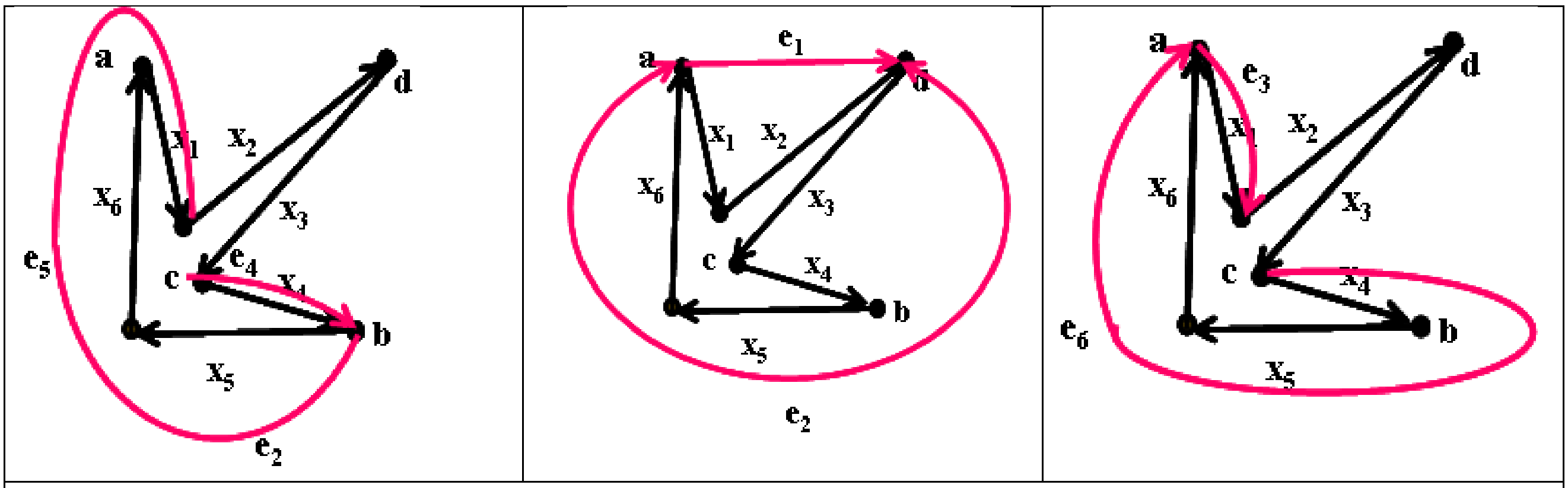

Рис. 4.15. Проекции векторов старой КБС на новую КБС.

Умножение матриц

| | $e_1$ | $e_2$ | $e_3$ | $e_4$ | $e_5$ | $e_6$ |
|---|---|---|---|---|---|---|
| **(a,b)** | | | **1** | **1** | | |
| **(b,a)** | | | | | **1** | **1** |
| **(c,d)** | **1** | | | | | **1** |
| **(d,c)** | | **1** | **1** | | | |

| | $x_1$ | $x_2$ | $x_3$ | $x_4$ | $x_5$ | $x_6$ |
|---|---|---|---|---|---|---|
| **$e_1$** | **1** | **1** | | | | |
| **$e_2$** | | | **1** | **1** | **1** | **1** |
| **$e_3$** | **1** | | | | | |
| **$e_4$** | | | | **1** | | |
| **$e_5$** | | | | | **1** | |
| **$e_6$** | | | | | | **1** |

Результат умножения

| $x_1$ | $x_2$ | $x_3$ | $x_4$ | $x_5$ | $x_6$ |
|---|---|---|---|---|---|

| | | | | | | |
|---|---|---|---|---|---|---|
| **(a,b)** | **1** | | | **1** | | |
| **(b,a)** | | | | | **1** | **1** |
| **(c,d)** | **1** | **1** | | | | **1** |
| **(d,c)** | **1** | | **1** | **1** | **1** | **1** |

Удаляем разнонаправленные вектора ($x_4$,$x_5$) и ($x_1$,$x_6$) и получаем проекции связей на новую КБС.

| | **$x_1$** | **$x_2$** | **$x_3$** | **$x_4$** | **$x_5$** | **$x_6$** |
|---|---|---|---|---|---|---|
| **(a,b)** | **1** | | | **1** | | |
| **(b,a)** | | | | | **1** | **1** |
| **(c,d)** | | **1** | | | | |
| **(d,c)** | | | **1** | | | |

## Комментарии

В данной главе рассмотрены элементы векторной алгебры пересечения отрекков.

Для решения задачи определения пересечения отрезклв с заданными координатами концов, в клеточном пространстве, существует волновой алгоритм Ли. Данный алгоритм зависит от порядка проведения соединений и не всегда позволяет достигать желаемых результатов. Алгоритм подробно описан в классических работах посвещенных автоматизированным методам проведения соединений [5,28]. Однако, существует ряд задач, где требуется определить пересечение отрезков в топологическом простраестве, без привязки к метрики.

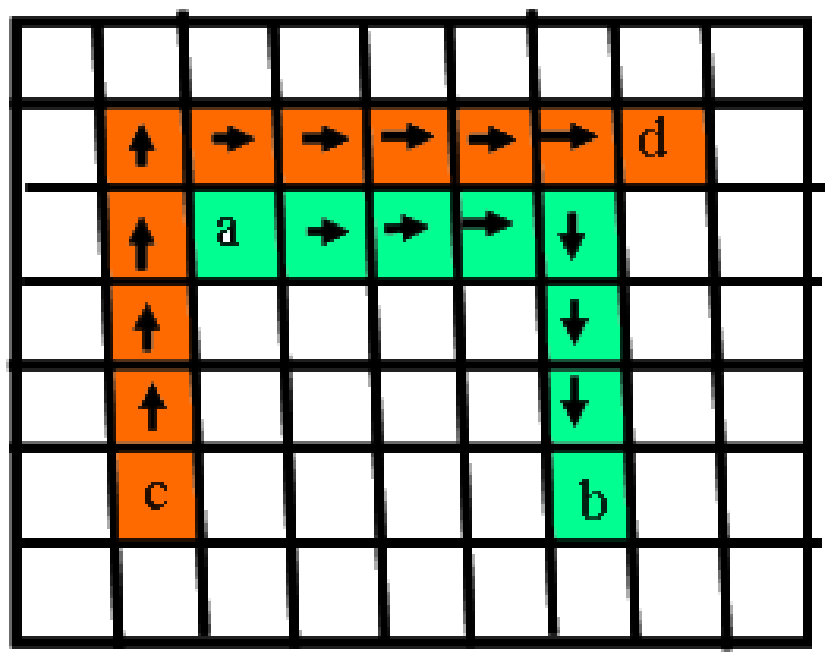

# Глава 5. ПРОВЕРКА ГРАФА НА ПЛАНАРНОСТЬ

## 5.1. Действия над суграфами

Будем рассматривать модифицированный метод определения планарности и построения вращения вершин с вычислительной сложностью соизмеримой с алгоритмом Хопкрофта-Тарьяна [7,42].

Перед тем как перейти к описанию алгоритма, рассмотрим вопросы необходимых действий над суграфами.

Любой суграф может быть записан как множество его рёбер. Что касается ориентированных циклов, то здесь возможно три вида записи цикла, например (см. рис. 5.1):

$c_4 = \{e_4,e_5,e_6,e_9,e_{27}\} \to \langle v_2,v_3,v_4,v_{15},v_{14},v_2\rangle \to \langle v_2,v_3\rangle + \langle v_3,v_4\rangle + \langle v_4,v_{15}\rangle +$

$+ \langle v_{15},v_{14}\rangle + \langle v_{14},v_2\rangle$.

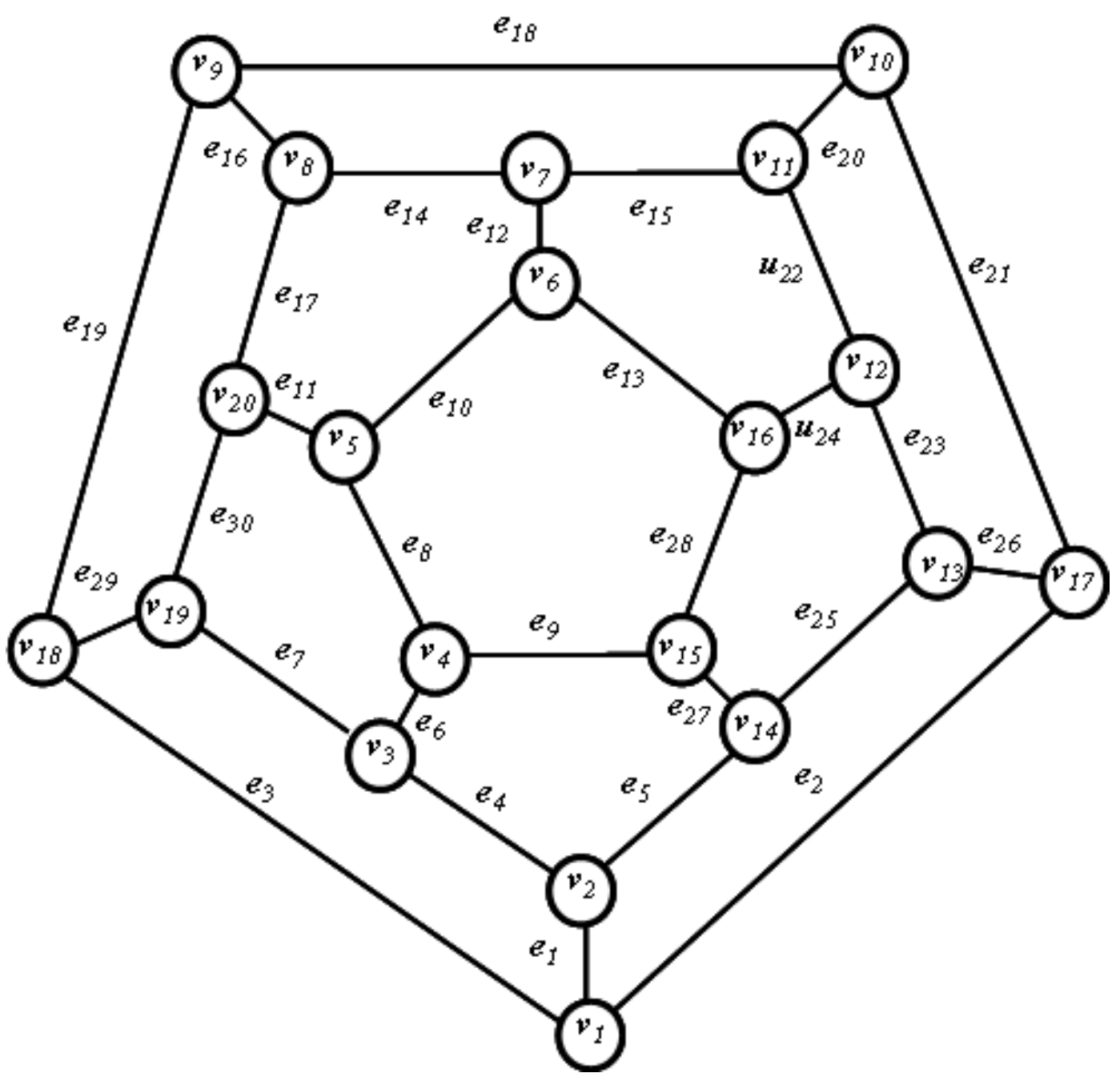

Рис. 5.1. Граф $G_1$.

Первая запись $\{e_4,e_5,e_6,e_9,e_{27}\}$ характеризуется множеством рёбер. Вторая запись характеризует кортеж, состоящий из последовательности вершин. Третья запись представляет цикл в виде последовательного сложения ориентированных дуг.

Известно также, что любое неориентированное ребро в графе можно представить в виде двух разнонаправленных ориентированных рёбер [41]. Например, ребро $e_4$ можно записать как $e_4 = \{v_2,v_3\} \to \langle v_2,v_3\rangle + \langle v_3,v_2\rangle$. Тогда сложение двух ориентированных циклов можно представить в виде последовательного сложения ориентированных рёбер, например (см. рис. 5.1):

$c_4 \oplus c_5 = \{e_4,e_5,e_6,e_9,e_{27}\} \oplus \{e_8,e_9,e_{10},e_{13},e_{28}\} \to \langle v_2,v_3,v_4,v_{15},v_{14},v_2\rangle +$
$\langle v_4,v_5,v_6,v_{16},v_{15},v_4\rangle \to \langle v_2,v_3\rangle + \langle v_3,v_4\rangle + \langle \mathbf{v_4,v_{15}}\rangle + \langle v_{15},v_{14}\rangle +$

$+ \langle v_{14}, v_2 \rangle + \langle v_4, v_5 \rangle + \langle v_5, v_6 \rangle + \langle v_6, v_{16} \rangle + \langle v_{16}, v_{15} \rangle + \langle \mathbf{v_{15}, v_4} \rangle =$
$= \langle v_2, v_3, v_4, v_5, v_6, v_{16}, v_{15}, v_{14}, v_2 \rangle$.

В этой операции разнонаправленные ориентированные рёбра попарно сокращаются.

Другой важной операцией является включение ребра в цикл.

Включение ребра в цикл разбивает цикл на две части. Например, в произвольный цикл $c_i = \langle v_{18}, v_9, v_{10}, v_{11}, v_{12}, v_{13}, v_{17}, v_1, v_{18} \rangle$ нужно включить ребро $\{v_{10}, v_{17}\}$. Для выполнения этой операции необходимо чтобы концевые вершины вновь введённого ребра (или совокупности рёбер если подключается маршрут) принадлежали циклу. Тогда к содержимому цикла можно добавить ориентированные дуги данного ребра (или совокупности рёбер):

$c_i = \langle v_{18}, v_9 \rangle + \langle v_9, v_{10} \rangle + \langle v_{10}, v_{11} \rangle + \langle v_{11}, v_{12} \rangle + \langle v_{12}, v_{13} \rangle + \langle v_{13}, v_{17} \rangle + \langle v_{17}, v_1 \rangle +$
$+ \langle v_1, v_{18} \rangle + \langle \mathbf{v_{17}, v_{10}} \rangle + \langle \mathbf{v_{10}, v_{17}} \rangle$

При этом содержимое цикла разбивается на части, где концевые вершины ребра участвуют в формировании маршрутов:

$\langle \mathbf{v_{17}}, v_1 \rangle + \langle v_1, v_{18} \rangle + \langle v_{18}, v_9 \rangle + \langle v_9, \mathbf{v_{10}} \rangle$;
$\langle \mathbf{v_{10}}, v_{11} \rangle + \langle v_{11}, v_{12} \rangle + \langle v_{12}, v_{13} \rangle + \langle v_{13}, \mathbf{v_{17}} \rangle$.

Затем с соблюдением последовательности к выделенным частям присоединяются дуги. Таким образом, образуются два новых обруча:

$\langle \mathbf{v_{17}}, v_1 \rangle + \langle v_1, v_{18} \rangle + \langle v_{18}, v_9 \rangle + \langle v_9, \mathbf{v_{10}} \rangle + \langle \mathbf{v_{10}, v_{17}} \rangle$;
$\langle \mathbf{v_{10}}, v_{11} \rangle + \langle v_{11}, v_{12} \rangle + \langle v_{12}, v_{13} \rangle + \langle v_{13}, \mathbf{v_{17}} \rangle + \langle \mathbf{v_{17}, v_{10}} \rangle$.

Введя необходимые обозначения и операции, продолжим рассмотрение вопросов дальнейшего развития методов и алгоритмов проверки графа на планарность с одновременным получением топологического рисунка графа.

## 5.2. Распределение обратных путей по блокам

Рассмотрим следующий пример проверки на планарность графа $G_1$ представленного на рис. 5.1. Выделим DFS-дерево методом поиска в глубину (см. рис. 5.2).

Построим фундаментальную матрицу циклов для выбранного дерева.

Используя фундаментальную матрицу циклов, выберем самый длинный цикл, образованный ветвями дерева и одной хордой (см. рис. 5.2). Пусть это будет цикл $\langle v_4, v_5, v_6, v_7, v_8, v_9, v_{10}, v_{11}, v_{12}, v_{13}, v_{14}, v_4 \rangle$. Данный цикл будем называть *опорным*.

**Определение 5.1.** *Обратным ребром* – называется ориентированный маршрут, состоящий из одной хорды [14].

**Определение 5.2.** *Обратным путём* – называется ориентированный маршрут, состоящий из последовательно расположенных ветвей дерева и одной, и только одной, хорды. Причем, концевые вершины такого обратного пути не всегда принадлежат опорному циклу[14].

Очевидно, что обратные ребра являются частным случаем обратных путей.

Как известно, матрица фундаментальных циклов состоит из двух подматриц: единичной подматрицы и подматрицы $\pi$ состоящей только из ветвей дерева.

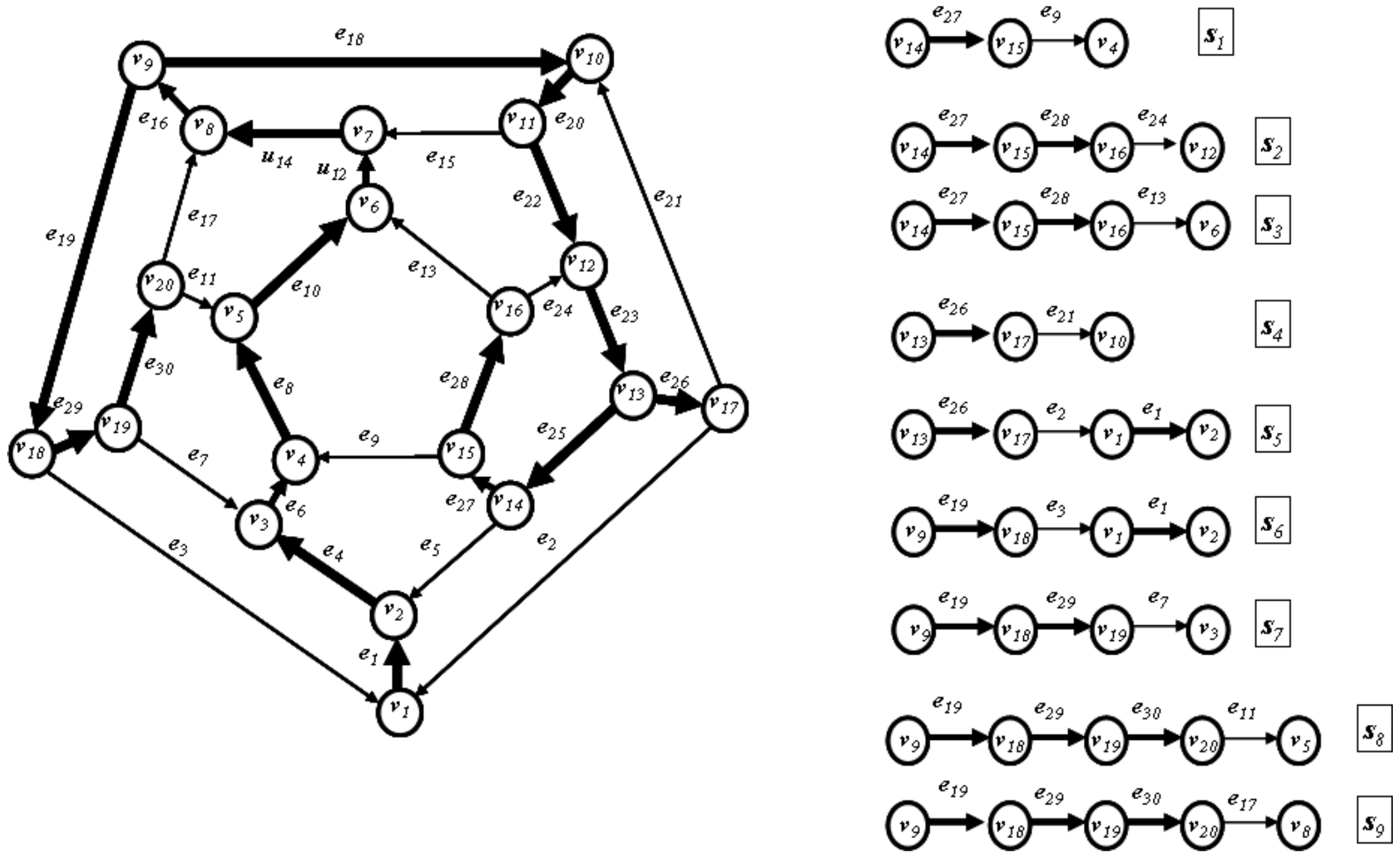

Рис. 5.2. DFS-дерево графа $G_1$, опорный цикл и обратные пути.

Фундаментальная подматрица $\pi$ имеет вид:

| | | $e_4$ | $e_6$ | $e_8$ | $e_{10}$ | $e_{12}$ | $e_{14}$ | $e_{16}$ | $e_{18}$ | $e_{20}$ | $e_{22}$ | $e_{23}$ | $e_{25}$ | $e_{27}$ | $e_{28}$ | $e_1$ | $e_{26}$ | $e_{19}$ | $e_{29}$ | $e_{30}$ |
|---|---|---|---|---|---|---|---|---|---|---|---|---|---|---|---|---|---|---|---|---|
| $s_5$ | $e_2$ | 1 | 1 | 1 | 1 | 1 | 1 | 1 | 1 | 1 | 1 | 1 | | | | **1** | **1** | | | |
| $s_6$ | $e_3$ | 1 | 1 | 1 | 1 | 1 | 1 | 1 | | | | | | | | **1** | | **1** | | |
| | $e_5$ | 1 | 1 | 1 | 1 | 1 | 1 | 1 | 1 | 1 | 1 | 1 | 1 | | | | | | | |
| $s_7$ | $e_7$ | | 1 | 1 | 1 | 1 | 1 | 1 | | | | | | | | | | **1** | **1** | |
| $s_1$ | $e_9$ | | | 1 | 1 | 1 | 1 | 1 | 1 | 1 | 1 | 1 | 1 | **1** | | | | | | |
| $s_8$ | $e_{11}$ | | | | 1 | 1 | 1 | 1 | | | | | | | | | | **1** | **1** | **1** |
| $s_3$ | $e_{13}$ | | | | | 1 | 1 | 1 | 1 | 1 | 1 | 1 | 1 | **1** | **1** | | | | | |
| | $e_{15}$ | | | | | | 1 | 1 | 1 | 1 | | | | | | | | | | |
| $s_9$ | $e_{17}$ | | | | | | | 1 | | | | | | | | | | **1** | **1** | **1** |
| $s_4$ | $e_{21}$ | | | | | | | | | 1 | 1 | 1 | | | | | **1** | | | |
| $s_2$ | $e_{24}$ | | | | | | | | | | | 1 | 1 | **1** | **1** | | | | | |

Разбиваем подматрицу $\pi$ на две части: часть, состоящую из ветвей дерева принадлежащих опорному циклу, и часть, состоящую из ветвей дерева не принадлежащих опорному циклу (выделены серым цветом). Строим обратные пути. Для этого объединяем ветви дерева, не включённые в опорный цикл и соответствующую хорду:

$s_1 = \{e_{27}, e_9\} \to \langle v_{14}, v_{15}, v_4 \rangle \to \langle v_{14}, v_{15} \rangle + \langle v_{15}, v_4 \rangle$;

$s_2 = \{e_{27}, e_{28}, e_{24}\} \to \langle v_{14}, v_{15}, v_{16}, v_{12} \rangle \to \langle v_{14}, v_{15} \rangle + \langle v_{15}, v_{16} \rangle + \langle v_{16}, v_{12} \rangle$;

$s_3 = \{e_{27}, e_{28}, e_{13}\} \to \langle v_{14}, v_{15}, v_{16}, v_6 \rangle \to \langle v_{14}, v_{15} \rangle + \langle v_{15}, v_{16} \rangle + \langle v_{16}, v_6 \rangle$;

$s_4 = \{e_{26}, e_{21}\} \to \langle v_{13}, v_{17}, v_{10} \rangle \to \langle v_{13}, v_{17} \rangle + \langle v_{17}, v_{10} \rangle$;

$s_5 = \{e_{26}, e_2, e_1\} \to \langle v_{13}, v_{17}, v_1, v_2 \rangle \to \langle v_{13}, v_{17} \rangle + \langle v_{17}, v_1 \rangle + \langle v_1, v_3 \rangle$;

$s_6 = \{e_{19}, e_3, e_1\} \to \langle v_9, v_{18}, v_1, v_2 \rangle \to \langle v_9, v_{18} \rangle + \langle v_{18}, v_1 \rangle + \langle v_1, v_2 \rangle$;

$s_7 = \{e_{19}, e_{29}, e_7\} \to \langle v_9, v_{18}, v_{19}, v_3 \rangle \to \langle v_9, v_{18} \rangle + \langle v_{18}, v_{19} \rangle + \langle v_{19}, v_3 \rangle$;

$s_8 = \{e_{19}, e_{29}, e_{30}, e_{11}\} \rightarrow \langle v_9, v_{18}, v_{19}, v_{20}, v_5\rangle \rightarrow \langle v_9, v_{18}\rangle + \langle v_{18}, v_{19}\rangle + \langle v_{19}, v_{20}\rangle + \langle v_{20}, v_5\rangle$;

$s_9 = \{e_{19}, e_{29}, e_{30}, e_{17}\} \rightarrow \langle v_9, v_{18}, v_{19}, v_{20}, v_8\rangle \rightarrow \langle v_9, v_{18}\rangle + \langle v_{18}, v_{19}\rangle + \langle v_{19}, v_{20}\rangle + \langle v_{20}, v_8\rangle$.

Объединяем обратные пути в блоки. Объединение в блоки осуществляется по первой вершине:

Блок $B_1$ включает следующие обратные пути:

$\{s_1 = \langle v_{14}, v_{15}, v_4\rangle, s_2 = \langle v_{14}, v_{15}, v_{16}, v_{12}\rangle, s_3 = \langle v_{14}, v_{15}, v_{16}, v_6\rangle\}$.

Блок $B_2$ включает следующие обратные пути:

$\{s_4 = \langle v_{13}, v_{17}, v_{10}\rangle, s_5 = \langle v_{13}, v_{17}, v_1, v_2\rangle\}$.

Блок $B_3$ включает следующие обратные пути:

$\{s_6 = \langle v_9, v_{18}, v_1, v_2\rangle, s_7 = \langle v_9, v_{18}, v_{19}, v_3\rangle, s_8 = \langle v_9, v_{18}, v_{19}, v_{20}, v_5\rangle, s_9 = \langle v_9, v_{18}, v_{19}, v_{20}, v_8\rangle\}$.

В блоке обратных путей $B_1$ выбираем самый длинный путь, и называем его *главным обратным путём для блока $B_1$*. В данном случае это будет путь $s_3 = \langle v_{14}, v_{15}, v_{16}, v_6\rangle$. Следующий обратный путь $s_2 = \langle v_{14}, v_{15}, v_{16}, v_{12}\rangle$ имеет последовательность вершин $\langle v_{14}, v_{15}, v_{16}\rangle$ уже включённую в главный обратный путь. Остаётся только блочное обратное ребро, состоящее из вершин $\langle v_{16}, v_{12}\rangle$. Следующий обратный путь $s_1 = \langle v_{14}, v_{15}, v_4\rangle$ имеет последовательность вершин $\langle v_{14}, v_{15}\rangle$ уже включённую в главный обратный путь. Остаётся только блочное обратное ребро, состоящее из вершин $\langle v_{15}, v_4\rangle$.

В блоке обратных путей $B_2$ выбираем самый длинный путь, и называем его *главным обратным путём для блока $B_2$*. В данном случае это будет путь $s_5 = \langle v_{13}, v_{17}, v_1, v_2\rangle$. Следующий обратный путь $s_4 = \langle v_{13}, v_{17}, v_{10}\rangle$ имеет последовательность вершин $\langle v_{13}, v_{17}\rangle$ уже включённую в главный обратный путь. Остаётся только блочное обратное ребро, состоящее из вершин $\langle v_{17}, v_{10}\rangle$.

В блоке обратных путей $B_3$ выбираем самый длинный путь, и называем его *главным обратным путем для блока $B_3$*. В данном случае выбираем путь $s_9 = \langle v_9, v_{18}, v_{19}, v_{20}, v_8\rangle$. Следующий обратный путь $s_8 = \langle v_9, v_{18}, v_{19}, v_{20}, v_5\rangle$ имеет последовательность вершин $\langle v_9, v_{18}, v_{19}, v_{20}\rangle$ уже включённую в главный обратный путь. Остаётся только блочное обратное ребро, состоящее из вершин $\langle v_{20}, v_5\rangle$. Следующий обратный путь $s_7 = \langle v_9, v_{18}, v_{19}, v_3\rangle$ имеет последовательность вершин $\langle v_9, v_{18}, v_{19}\rangle$ уже включённую в главный обратный путь. Остаётся только блочное обратное ребро состоящее из вершин $\langle v_{19}, v_3\rangle$. Наконец, обратный путь $s_6 = \langle v_9, v_{18}, v_1, v_2\rangle$ имеет последовательность вершин $\langle v_9, v_{18}\rangle$ уже включённую в главный обратный путь, но вершины $v_1$ и $v_2$ не включены в главный обратный путь. Поэтому формируем частичный обратный путь $\langle v_{18}, v_1, v_2\rangle$ из концевой вершины $v_{18}$ и остатка $v_1, v_2$.

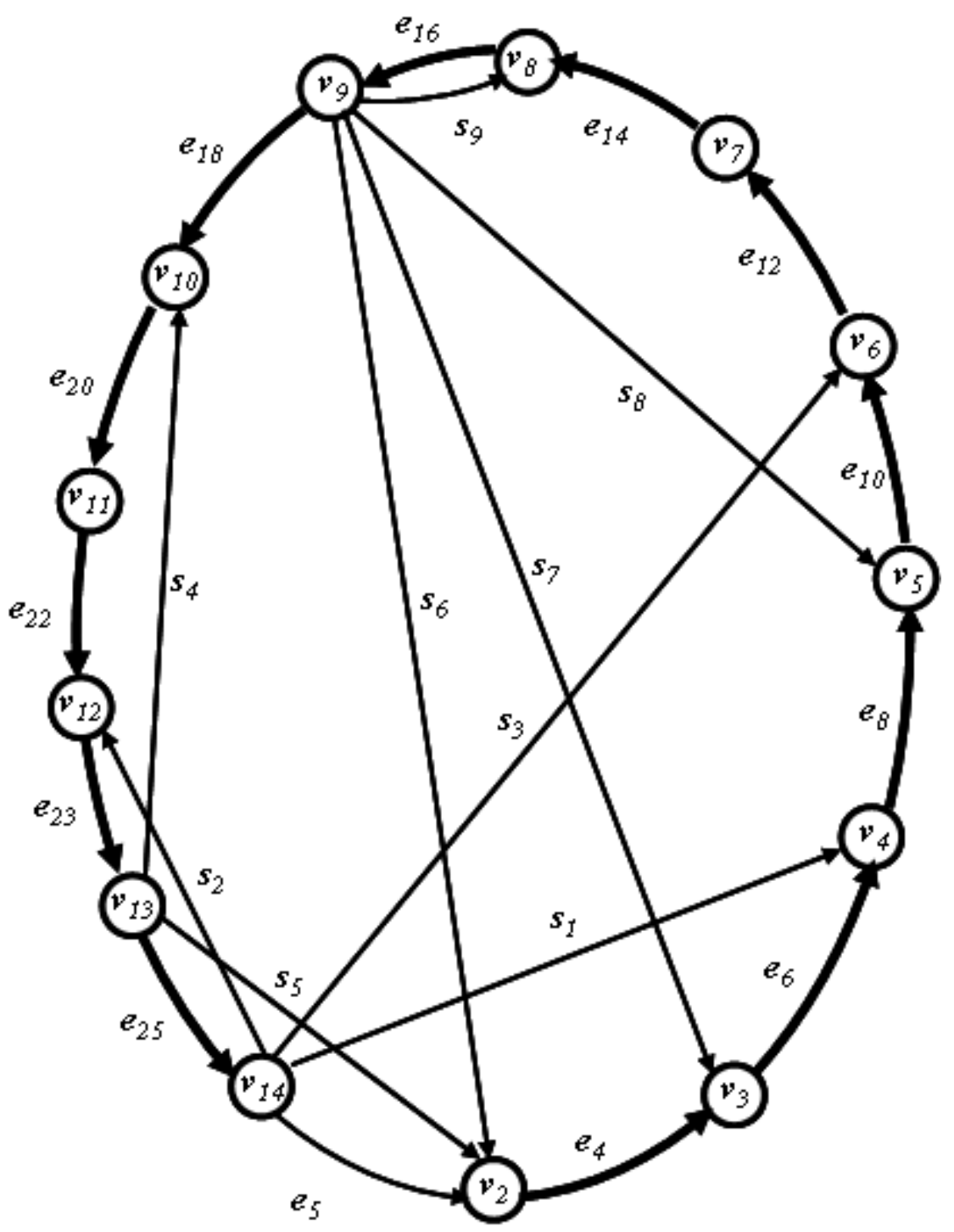

Рис. 5.3. Порождение обратных путей

Таким образом, можно сказать, что разбиение на блоки породило систему обратных маршрутов:

Блок $B_1$ {$s_3 = <v_{14},v_{15},v_{16},v_6>$ – главный обратный путь;
$e_{24} = <v_{16},v_{12}>$ – блочное обратное ребро;
$e_9 = <v_{15},v_4>$} – блочное обратное ребро.

Блок $B_2$ {$s_5 = <v_{13},v_{17},v_1,v_2>$ – главный обратный путь;
$e_{21} = <v_{17},v_{10}>$} – блочное обратное ребро.

Блок $B_3$ {$s_9 = <v_9,v_{18},v_{19},v_{20},v_8>$ – главный обратный путь;
$e_{11} = <v_{20},v_5>$ – блочное обратное ребро;
$e_7 = <v_{19},v_3>$ – блочное обратное ребро;
$s_6 = <v_{18},v_1,v_2>$} – частичный обратный путь.

И обратные рёбра: $e_5 = <v_{14},v_2>$; $e_{15} = <v_{11},v_7>$.

Рисунок 5.3 показывает, что если перенести блок $B_1$ вне опорного цикла то граф можно уложить на плоскости без пересечения связей.

Произведем расположение обратных маршрутов в следующем порядке: вначале разместим главные обратные пути, затем частичные обратные пути, затем блочные обратные рёбра и, наконец, обратные рёбра. Последовательно просматривая порядок расположения обратных маршрутов, производим их включение в обручи с дальнейшим разбиением. В нашем случае порядок расположения обратных маршрутов имеет вид: M = $<s_9,s_3,s_5,s_6,e_{11},e_7,e_{24},e_9,e_{21},e_5,e_{13},e_{15}>$. Кортеж M будем называ*ть списком очередности обратных маршрутов.*

## 5.3. Построение простых циклов

Опорный цикл разбивает поверхность $R^2$ на две области: внутреннюю и внешнюю. Границы этих областей могут быть описаны двумя симметричными циклами (обручами):

$c_1 = \langle v_2,v_3,v_4,v_5,v_6,v_7,v_8,v_9,v_{10},v_{11},v_{12},v_{13},v_{14},v_2\rangle$,
и
$c_2 = \langle v_2,v_{14},v_{13},v_{12},v_{11},v_{10},v_9,v_8,v_7,v_6,v_5,v_4,v_3,v_2\rangle$.

**Определение 5.3.** *Обруч* – это простой цикл, образованный кольцевой суммой простых циклов, каждый из которых имеет, по крайней мере, одно общее ребро с другим простым циклом.

Заметим, что простой цикл можно рассматривать как частный случай обруча. В свою очередь, обратные ребра и обратные пути разбивают поверхности на грани, где границами граней служат обручи.

Применим операцию включения обратного ребра в обруч графа. Рассмотрим главный обратный путь $s_9 = \langle v_9,v_{18},v_{19},v_{20},v_8\rangle$ в порядке расположения обратных маршрутов. Концевые вершины пути $v_9$ и $v_8$ включаются и в цикл $c_1$ и в цикл $c_2$. Учитывая то, что первый обратный маршрут может быть включен в любой цикл, выбираем цикл $c_1$. Включаем обратный маршрут $s_9$ в состав цикла $c_1$. Производим необходимое преобразование. В результате получим два новых цикла:

$c_1$: $\langle v_2,v_3\rangle + \langle v_3,v_4\rangle + \langle v_4,v_5\rangle + \langle v_5,v_6\rangle + \langle v_6,v_7\rangle + \langle v_7,v_8\rangle + \langle v_8,v_9\rangle + \langle v_9,v_{10}\rangle +$
$+ \langle v_{10},v_{11}\rangle + \langle v_{11},v_{12}\rangle + \langle v_{12},v_{13}\rangle + \langle v_{13},v_{14}\rangle + \langle v_{14},v_2\rangle +$
$+ \langle v_9,v_{18}\rangle + \langle v_{18},v_{19}\rangle + \langle v_{19},v_{20}\rangle + \langle v_{20},v_8\rangle +$
$+ \langle v_8,v_{20}\rangle + \langle v_{20},v_{19}\rangle + \langle v_{19},v_{18}\rangle + \langle v_{18},v_9\rangle =$
$= \langle v_2,v_3\rangle + \langle v_3,v_4\rangle + \langle v_4,v_5\rangle + \langle v_5,v_6\rangle + \langle v_6,v_7\rangle + \langle v_7,v_8\rangle + \langle v_8,v_{20}\rangle + \langle v_{20},v_{19}\rangle +$
$+ \langle v_{19},v_{18}\rangle + \langle v_{18},v_9\rangle + \langle v_9,v_{10}\rangle + \langle v_{10},v_{11}\rangle + \langle v_{11},v_{12}\rangle + \langle v_{12},v_{13}\rangle + \langle v_{13},v_{14}\rangle +$
$+ \langle v_{14},v_2\rangle + \langle v_9,v_{18}\rangle + \langle v_{18},v_{19}\rangle + \langle v_{19},v_{20}\rangle + \langle v_{20},v_8\rangle + \langle v_8,v_9\rangle$;

$c_1' = \langle v_2,v_3,v_4,v_5,v_6,v_7,v_8,v_{20},v_{19},v_{18},v_9,v_{10},v_{11},v_{12},v_{13},v_{14},v_2\rangle$;
$c_3 = \langle v_9,v_{18},v_{19},v_{20},v_8,v_9\rangle$.

Таким образом, система сформированных обручей имеет вид:

$c_1' = \langle v_2,v_3,v_4,v_5,v_6,v_7,v_8,v_{20},v_{19},v_{18},v_9,v_{10},v_{11},v_{12},v_{13},v_{14},v_2\rangle$;
$c_3 = \langle v_9,v_{18},v_{19},v_{20},v_8,v_9\rangle$;
$c_2 = \langle v_2,v_{14},v_{13},v_{12},v_{11},v_{10},v_9,v_8,v_7,v_6,v_5,v_4,v_3,v_2\rangle$.

Будем размещать главный обратный цикл $s_5 = \langle v_{13},v_{17},v_1,v_2\rangle$. Концевые вершины $v_2,v_{13}$ включаются и в цикл $c_1'$ и в цикл $c_2$. Возникла неопределенность. Тогда пробуем разместить главный обратный путь $s_3 = \langle v_{14},v_{15},v_{16},v_6\rangle$, обе концевые вершины $v_{14},v_6$ обратного пути $s_3 = \langle v_{14},v_{15},v_{16},v_6\rangle$ включаются и в цикл $c_1'$ и в цикл $c_2$. Снова возникла неопределенность. Пробуем поместить путь $s_6 = \langle v_{18},v_1,v_2\rangle$. Он размещается только в $c_1$ , так как только там имеются совместно две концевые вершины $\{v_{18},v_2\}$. В результате получим два новых цикла:

$c_1'$: $\langle v_2,v_3\rangle + \langle v_3,v_4\rangle + \langle v_4,v_5\rangle + \langle v_5,v_6\rangle + \langle v_6,v_7\rangle + \langle v_7,v_8\rangle + \langle v_8,v_{20}\rangle + \langle v_{20},v_{19}\rangle$
$+ \langle v_{19},v_{18}\rangle + \langle v_{18},v_9\rangle + \langle v_9,v_{10}\rangle + \langle v_{10},v_{11}\rangle + \langle v_{11},v_{12}\rangle + \langle v_{12},v_{13}\rangle + \langle v_{13},v_{14}\rangle + \langle v_{14},v_2\rangle +$
$+ \langle v_{18},v_1\rangle + \langle v_1,v_2\rangle +$

$+ \langle v_2,v_1\rangle + \langle v_1,v_{18}\rangle =$
$= \langle v_2,v_3\rangle + \langle v_3,v_4\rangle + \langle v_4,v_5\rangle + \langle v_5,v_6\rangle + \langle v_6,v_7\rangle + \langle v_7,v_8\rangle + \langle v_8,v_{20}\rangle + \langle v_{20},v_{19}\rangle +$
$+\langle v_{19},v_{18}\rangle + \langle v_{18},v_1\rangle + \langle v_1,v_2\rangle +$
$+ \langle v_{18},v_9\rangle + \langle v_9,v_{10}\rangle + \langle v_{10},v_{11}\rangle + \langle v_{11},v_{12}\rangle + \langle v_{12},v_{13}\rangle + \langle v_{13},v_{14}\rangle + \langle v_{14},v_2\rangle +$
$+ \langle v_2,v_1\rangle + \langle v_1,v_{18}\rangle$;
$c_1'' = \langle v_2,v_3,v_4,v_5,v_6,v_7,v_8,v_{20},v_{19},v_{18},v_1,v_2\rangle$;
$c_4 = \langle v_{18},v_9,v_{10},v_{11},v_{12},v_{13},v_{14},v_2,v_1,v_{18}\rangle$.

Система обручей имеет вид:

$c_1'' = \langle v_2,v_3,v_4,v_5,v_6,v_7,v_8,v_{20},v_{19},v_{18},v_1,v_2\rangle$;
$c_4 = \langle v_{18},v_9,v_{10},v_{11},v_{12},v_{13},v_{14},v_2,v_1,v_{18}\rangle$;
$c_3 = \langle v_9,v_{18},v_{19},v_{20},v_8,v_9\rangle$;
$c_2 = \langle v_2,v_{14},v_{13},v_{12},v_{11},v_{10},v_9,v_8,v_7,v_6,v_5,v_4,v_3,v_2\rangle$.

Размещаем обратный маршрут $s_3 = \langle v_{14},v_{15},v_{16},v_6\rangle$ в $c_2$. В результате получим два новых цикла:

$c_2' = \langle v_{14},v_{13}v_{12},v_{11},v_{10},v_9,v_8,v_7,v_6,v_{16},v_{15},v_{14}\rangle$;
$c_5 = \langle v_2,v_{14},v_{15},v_{16},v_6,v_5,v_4,v_3,v_2\rangle$.

Система обручей имеет вид:

$c_2' = \langle v_{14},v_{13}v_{12},v_{11},v_{10},v_9,v_8,v_7,v_6,v_{16},v_{15},v_{14}\rangle$;
$c_5 = \langle v_2,v_{14},v_{15},v_{16},v_6,v_5,v_4,v_3,v_2\rangle$;
$c_1'' = \langle v_2,v_3,v_4,v_5,v_6,v_7,v_8,v_{20},v_{19},v_{18},v_1,v_2\rangle$;
$c_4 = \langle v_{18},v_9,v_{10},v_{11},v_{12},v_{13},v_{14},v_2,v_1,v_{18}\rangle$;
$c_3 = \langle v_9,v_{18},v_{19},v_{20},v_8,v_9\rangle$.

Поочерёдно вставляем обратные маршруты из списка очерёдности в образованные обручи до полного исчерпания. Если список очередности исчерпан и все обратные маршруты уложены в обручах, то граф планарен. В противном случае, если список очередности не исчерпан, граф непланарен.

Окончательная система простых циклов имеет вид:

$c_2''' = \langle v_{12},v_{11},v_7,v_6,v_{16},v_{12}\rangle \rightarrow \{e_{22},e_{15},e_{12},e_{13},e_{24}\}$;
$c_{12} = \langle v_{11},v_{10},v_9,v_8,v_7,v_{11}\rangle \rightarrow \{e_{20},e_{18},e_{16},e_{14},e_{15}\}$;
$c_1'''' = \langle v_2,v_3,v_{19},v_{18},v_1,v_2\rangle \rightarrow \{e_4,e_7,e_{29},e_3,e_1\}$;
$c_{11} = \langle v_3,v_4,v_5,v_{20},v_{19},v_3\rangle \rightarrow \{e_6,e_8,e_{11},e_{30},e_7\}$;
$c_{10} = \langle v_5,v_6,v_7,v_8,v_{20},v_5\rangle \rightarrow \{e_{10},e_{12},e_{14},e_{17},e_{11}\}$;
$c_4'' = \langle v_{18},v_9,v_{10},v_{17},v_1,v_{18}\rangle \rightarrow \{e_{19},e_{18},e_{21},e_2,e_3\}$;
$c_9 = \langle v_{10},v_{11},v_{12},v_{13},v_{17},v_{10}\rangle \rightarrow \{e_{20},e_{22},e_{23},e_{26},e_{21}\}$;
$c_5' = \langle v_2,v_{14},v_{15},v_4,v_3,v_2\rangle \rightarrow \{e_5,e_{27},e_9,e_6,e_4\}$;
$c_8 = \langle v_{15},v_{16},v_6,v_5,v_4,v_{15}\rangle \rightarrow \{e_{28},e_{13},e_{10},e_8,e_9\}$;
$c_7 = \langle v_{14},v_{13},v_{12},v_{16},v_{15},v_{14}\rangle \rightarrow \{e_{25},e_{23},e_{24},e_{28},e_{27}\}$;
$c_6 = \langle v_2,v_1,v_{17},v_{13},v_{14},v_2\rangle \rightarrow \{e_1,e_2,e_{26},e_{25},e_5\}$;
$c_3 = \langle v_9,v_{18},v_{19},v_{20},v_8,v_9\rangle \rightarrow \{e_{19},e_{29},e_{30},e_{17},e_{16}\}$.

Так как кольцевая сумма циклов есть пустое множество, то данная система циклов характеризует топологический рисунок плоского графа [15]. Диаграмма вращения вершин представлена на рис. 5.4, а топологический рисунок представлен на рис. 5.5.

$\hbar_1 : v_2\ \ v_{17}\ \ v_{18}$
$\hbar_2 : v_3\ \ v_{14}\ \ v_1$
$\hbar_3 : v_{19}\ \ v_4\ \ v_2$
$\hbar_4 : v_5\ \ v_{15}\ \ v_3$
$\hbar_5 : v_{20}\ \ v_6\ \ v_4$
$\hbar_6 : v_7\ \ v_{16}\ \ v_5$
$\hbar_7 : v_8\ \ v_{11}\ \ v_6$
$\hbar_8 : v_9\ \ v_7\ \ v_{20}$
$\hbar_9 : v_{10}\ \ v_8\ \ v_{18}$
$\hbar_{10} : v_9\ \ v_{17}\ \ v_{11}$
$\hbar_{11} : v_{10}\ \ v_{12}\ \ v_7$
$\hbar_{12} : v_{11}\ \ v_{13}\ \ v_{16}$
$\hbar_{13} : v_{12}\ \ v_{17}\ \ v_{14}$
$\hbar_{14} : v_{15}\ \ v_{13}\ \ v_2$
$\hbar_{15} : v_{16}\ \ v_{14}\ \ v_4$
$\hbar_{16} : v_6\ \ v_{12}\ \ v_{15}$
$\hbar_{17} : v_{13}\ \ v_{10}\ \ v_1$
$\hbar_{18} : v_9\ \ v_{19}\ \ v_1$
$\hbar_{19} : v_{18}\ \ v_{20}\ \ v_3$
$\hbar_{20} : v_8\ \ v_5\ \ v_{19}$

Рис. 5.4. Диаграмма вращения вершин

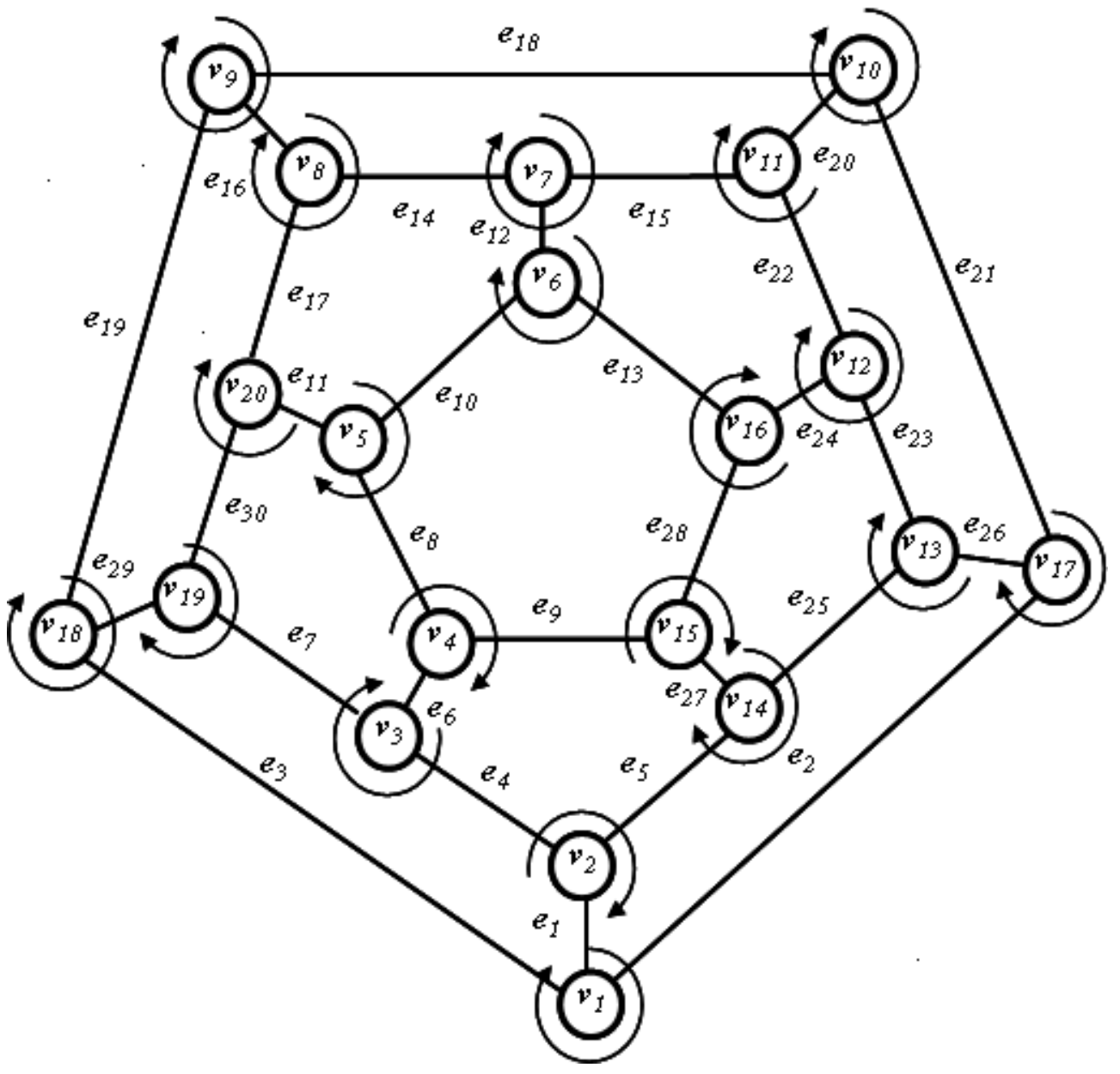

Рис. 5.5. Топологический рисунок графа

## 5.4. Задача выбора опорного цикла

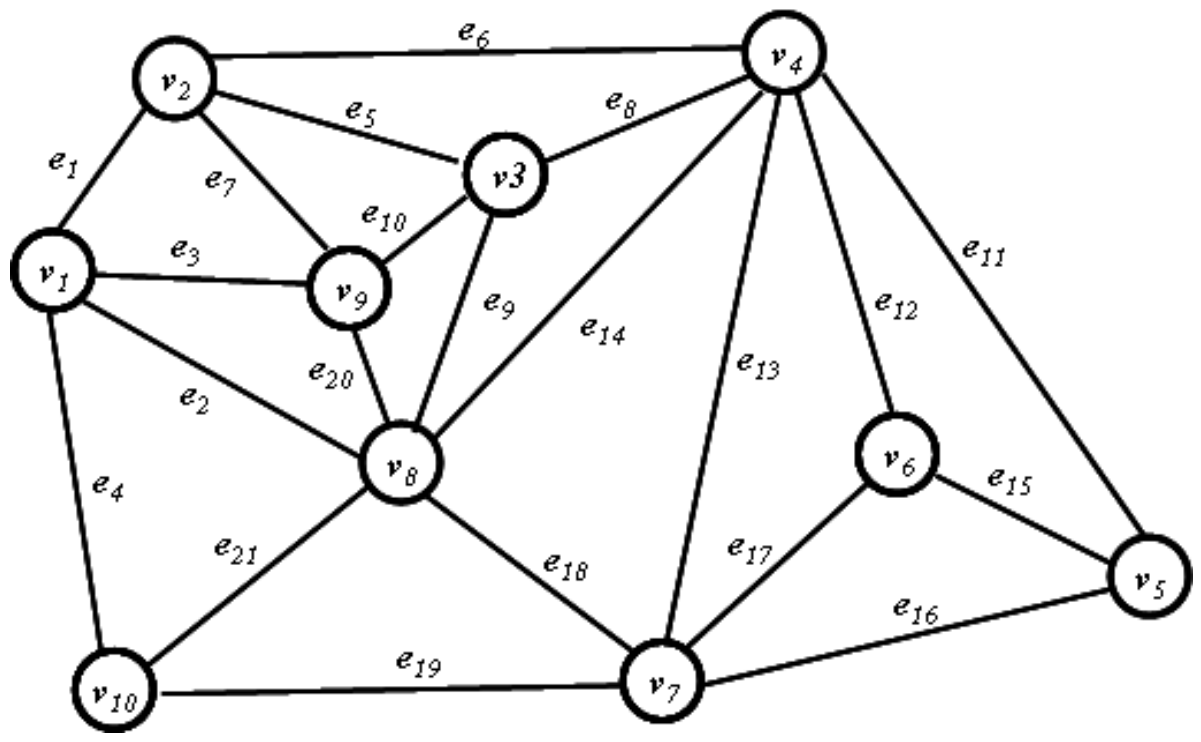

Рис. 5.6. Граф $G_2$

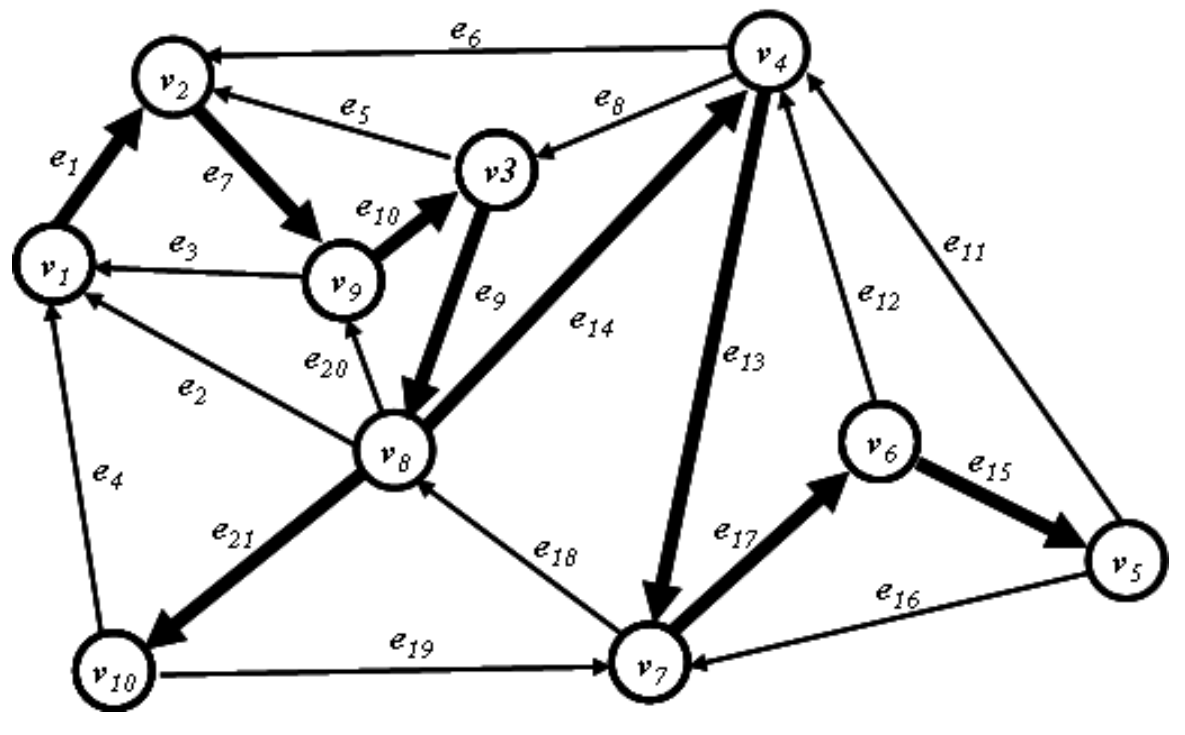

Рис. 5.7. DFS-дерево графа $G_2$

В предыдущем примере рассматривается случай удачного выбора опорного цикла. В этом случае все концевые вершины обратных маршрутов и обратных рёбер принадлежат опорному циклу. Рассмотрим случай, когда концевые вершины обратных маршрутов не принадлежат опорному циклу или обратные маршруты представляют собой замкнутые маршруты.

Рассмотрим следующий граф, представленный на рис. 5.6. Выделим DFS-дерево (см. рис.

5.7). Выберем самый длинный цикл $<v_1,v_2,v_9,v_3,v_8,v_{10},v_1>$ и назовем его опорным циклом. Будем рассматривать только подматрицу $\pi$. Выделенный опорный цикл разбивает подматрицу $\pi$ на две части. Ветви дерева, вошедшие в опорный цикл, образуют подмножество $R_1$. Ветви дерева, не вошедшие в опорный цикл – подмножество $R_2$ (серый цвет). Так как построение дерева определяет ориентированный граф, то в качестве элементов подматрицы будут выступать и ориентированные рёбра.

| | **ветви дерева (подматрица $\pi$)** | | | | | | | | |
|---|---|---|---|---|---|---|---|---|---|
| | $<v_1,v_2>$ | $<v_2,v_9>$ | $<v_9,v_3>$ | $<v_3,v_8>$ | $<v_8,v_{10}>$ | $<v_8,v_4>$ | $<v_4,v_7>$ | $<v_7,v_6>$ | $<v_6,v_5>$ |
| | $e_1$ | $e_7$ | $e_{10}$ | $e_9$ | $e_{21}$ | $e_{14}$ | $e_{13}$ | $e_{17}$ | $e_{15}$ |
| $e_2 = <v_8,v_1>$ | **1** | **1** | **1** | **1** | | | | | |
| $e_3 = <v_9,v_1>$ | **1** | **1** | | | | | | | |
| $e_4 = <v_{10},v_1>$ | **1** | **1** | **1** | **1** | **1** | | | | |
| $e_5 = <v_3,v_2>$ | | **1** | **1** | | | | | | |
| $e_6 = <v_4,v_2>$ | | **1** | **1** | **1** | | **1** | | | |
| $e_8 = <v_4,v_3>$ | | | | **1** | | **1** | | | |
| $e_{11} = <v_5,v_4>$ | | | | | | | **1** | **1** | **1** |
| $e_{12} = <v_6,v_4>$ | | | | | | | **1** | **1** | |
| $e_{16} = <v_5,v_7>$ | | | | | | | | **1** | **1** |
| $e_{18} = <v_7,v_8>$ | | | | | | **1** | **1** | | |
| $e_{19} = <v_7,v_{10}>$ | | | | | **1** | **1** | **1** | | |
| $e_{20} = <v_8,v_9>$ | | | **1** | **1** | | | | | |

Будем строить обратные пути. С этой целью выпишем все фундаментальные циклы в виде множества ребер:

$c_1 = \{e_2,\mathbf{e_1},\mathbf{e_7},\mathbf{e_{10}},\mathbf{e_9}\}$; $c_2 = \{e_3,\mathbf{e_1},\mathbf{e_7}\}$; $c_3 = \{e_4,\mathbf{e_1},\mathbf{e_7},\mathbf{e_{10}},\mathbf{e_9},\mathbf{e_{21}}\}$; $c_4 = \{e_5,\mathbf{e_7},\mathbf{e_{10}}\}$;
$c_5 = \{e_6,\mathbf{e_7},\mathbf{e_{10}},\mathbf{e_9},\mathrm{e}_{14}\}$; $c_6 = \{e_8,\mathbf{e_9},\mathrm{e}_{14}\}$; $c_7 = \{e_{11},\mathrm{e}_{13},\mathrm{e}_{17},\mathrm{e}_{15}\}$; $c_8 = \{e_{12},\mathrm{e}_{13},\mathrm{e}_{17}\}$;
$c_9 = \{e_{16},\mathrm{e}_{17},\mathrm{e}_{15}\}$; $c_{10} = \{e_{18},\mathrm{e}_{14},\mathrm{e}_{13}\}$; $c_{11} = \{e_{19},\mathbf{e_{21}},\mathrm{e}_{14},\mathrm{e}_{13}\}$; $c_{12} = \{e_{20},\mathbf{e_{10}},\mathbf{e_9}\}$.

Курсивом выделены хорды, жирным шрифтом выделены рёбра дерева, вошедшие в опорный цикл, простым шрифтом – рёбра дерева, не вошедшие в опорный цикл. Исключим из описания циклов ребра, вошедшие в опорный цикл. Получим следующие части циклов:

$s_1 = \{e_2\}$; $s_2 = \{e_3\}$; $c_3$ – опорный цикл; $s_4 = \{e_5\}$;
$s_5 = \{e_6,\mathrm{e}_{14}\}$; $s_6 = \{e_8,\mathrm{e}_{14}\}$; $s_7 = \{e_{11},\mathrm{e}_{13},\mathrm{e}_{17},\mathrm{e}_{15}\}$; $s_8 = \{e_{12},\mathrm{e}_{13},\mathrm{e}_{17}\}$;
$s_9 = \{e_{16},\mathrm{e}_{17},\mathrm{e}_{15}\}$; $s_{10} = \{e_{18},\mathrm{e}_{14},\mathrm{e}_{13}\}$; $s_{11} = \{e_{19},\mathrm{e}_{14},\mathrm{e}_{13}\}$; $s_{12} = \{e_{20}\}$.

Запишем полученные части циклов как обратные рёбра в виде ориентированных маршрутов:

$s_1 = \{e_2\} \rightarrow <v_8,v_1>$; $s_2 = \{e_3\} \rightarrow <v_9,v_1>$; $s_4 = \{e_5\} \rightarrow <v_3,v_2>$;
$s_5 = \{e_6,\mathrm{e}_{14}\} \rightarrow <v_4,v_2> + <v_8,v_4> = <v_8,v_4,v_2>$;
$s_6 = \{e_8,\mathrm{e}_{14}\} \rightarrow <v_4,v_3> + <v_8,v_4> = <v_8,v_4,v_3>$;
$s_7 = \{e_{11},\mathrm{e}_{13},\mathrm{e}_{17},\mathrm{e}_{15}\} \rightarrow <v_5,v_4> + <v_4,v_7> + <v_7,v_6> + <v_6,v_5> = <v_5,v_4,v_7,v_6,v_5>$;
$s_8 = \{e_{12},\mathrm{e}_{13},\mathrm{e}_{17}\} \rightarrow <v_6,v_4> + <v_4,v_7> + <v_7,v_6> = <v_6,v_4,v_7,v_6>$;
$s_9 = \{e_{16},\mathrm{e}_{17},\mathrm{e}_{15}\} \rightarrow <v_5,v_7> + <v_7,v_6> + <v_6,v_5>= <v_5,v_7,v_6,v_5>$;
$s_{10} = \{e_{18},\mathrm{e}_{14},\mathrm{e}_{13}\} \rightarrow <v_7,v_8> + <v_8,v_4> + <v_4,v_7>= <v_7,v_8,v_4,v_7>$;

$s_{11} = \{e_{19}, e_{14}, e_{13}\} \rightarrow \langle v_8, v_4\rangle + \langle v_4, v_7\rangle + \langle v_7, v_{10}\rangle += \langle v_8, v_4, v_7, v_{10}\rangle$;
$s_{12} = \{e_{20}\} \rightarrow \langle v_8, v_9\rangle$.

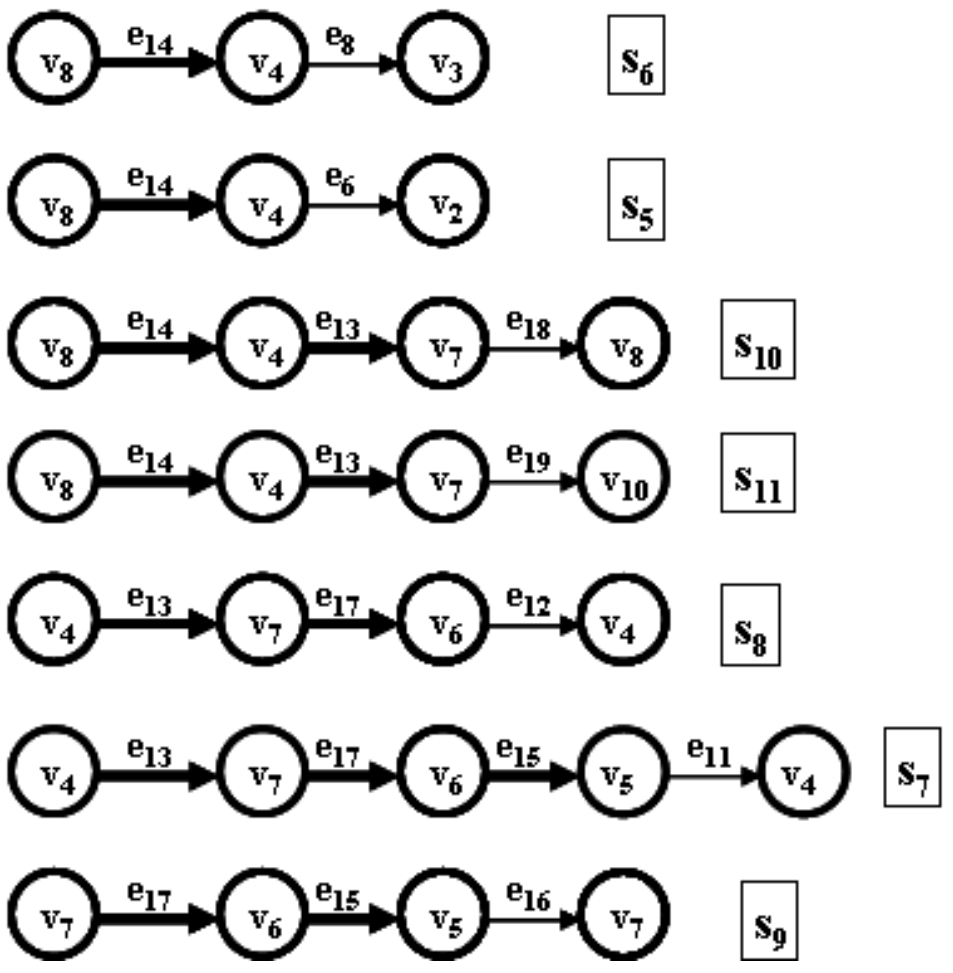

Рис. 5.8. Построение обратных путей для выбранного опорного цикла

Так как при построении топологического рисунка существование обратных путей в виде замкнутых маршрутов невозможно по определению, будем проводить дальнейшие преобразования, увеличивая длину опорного цикла.

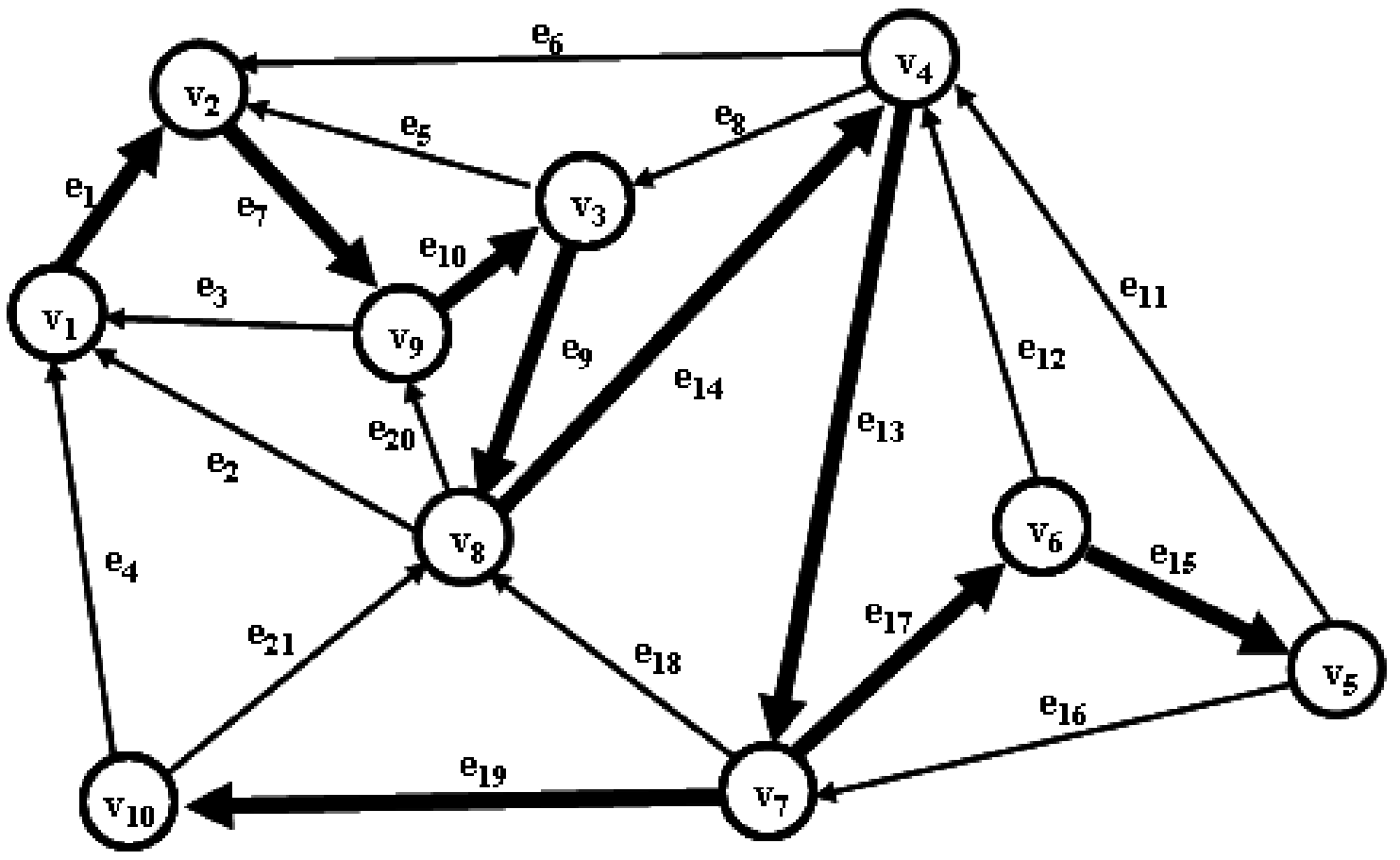

Рис. 5.9. Новое DFS-дерево графа $G_2$ и новый опорный цикл

Для преобразования ищем единственный элемент в строке подматрицы $R_1$. Такая подстрока имеется – это хорда $e_{19}$. Выбираем обратный маршрут $s_{11}$, так как в нем количество ветвей дерева, не вошедших в опорный цикл больше, чем в других случаях. Производим преобразования в подматрице $\pi$. Объявляем хорду $e_{19}$ ребром дерева, принадлежащим опорному циклу, и включаем в подматрицу $\pi$, а ребро дерева $e_{21}$ объявляем хордой и исключаем из подматрицы $\pi$. Ветви дерева $e_{13}$ и $e_{14}$ включаем в опорный цикл. Тем самым преобразуем DFS-дерево графа $G_2$ (см. рис. 5.9).

Для нового дерева построим подматрицу $\pi$:

| | **ветви дерева (подматрица $\pi$)** | | | | | | | | |
|---|---|---|---|---|---|---|---|---|---|
| | <v$_1$,v$_2$> | <v$_2$,v$_9$> | <v$_9$,v$_3$> | <v$_3$,v$_8$> | <v$_8$,v$_4$> | <v$_4$,v$_7$> | <v$_7$v$_{10}$> | <v$_7$,v$_6$> | <v$_6$,v$_5$> |
| | $e_1$ | $e_7$ | $e_{10}$ | $e_9$ | $e_{14}$ | $e_{13}$ | $e_{19}$ | $e_{17}$ | $e_{15}$ |
| $e_2$ = | **1** | **1** | **1** | **1** | | | | | |
| $e_3$ = | **1** | **1** | | | | | | | |
| $e_4$ = | **1** | **1** | **1** | **1** | **1** | **1** | **1** | | |
| $e_5$ = | | **1** | **1** | | | | | | |
| $e_6$ = | | **1** | **1** | **1** | **1** | | | | |
| $e_8$ = | | | | **1** | **1** | | | | |
| $e_{11}$ = | | | | | | **1** | | **1** | **1** |
| $e_{12}$ = | | | | | | **1** | | **1** | |
| $e_{16}$ = | | | | | | | | **1** | **1** |
| $e_{18}$ = | | | | | **1** | **1** | | | |
| $e_{21}$ = | | | | | **1** | **1** | **1** | | |
| $e_{20}$ = | | | **1** | **1** | | | | | |

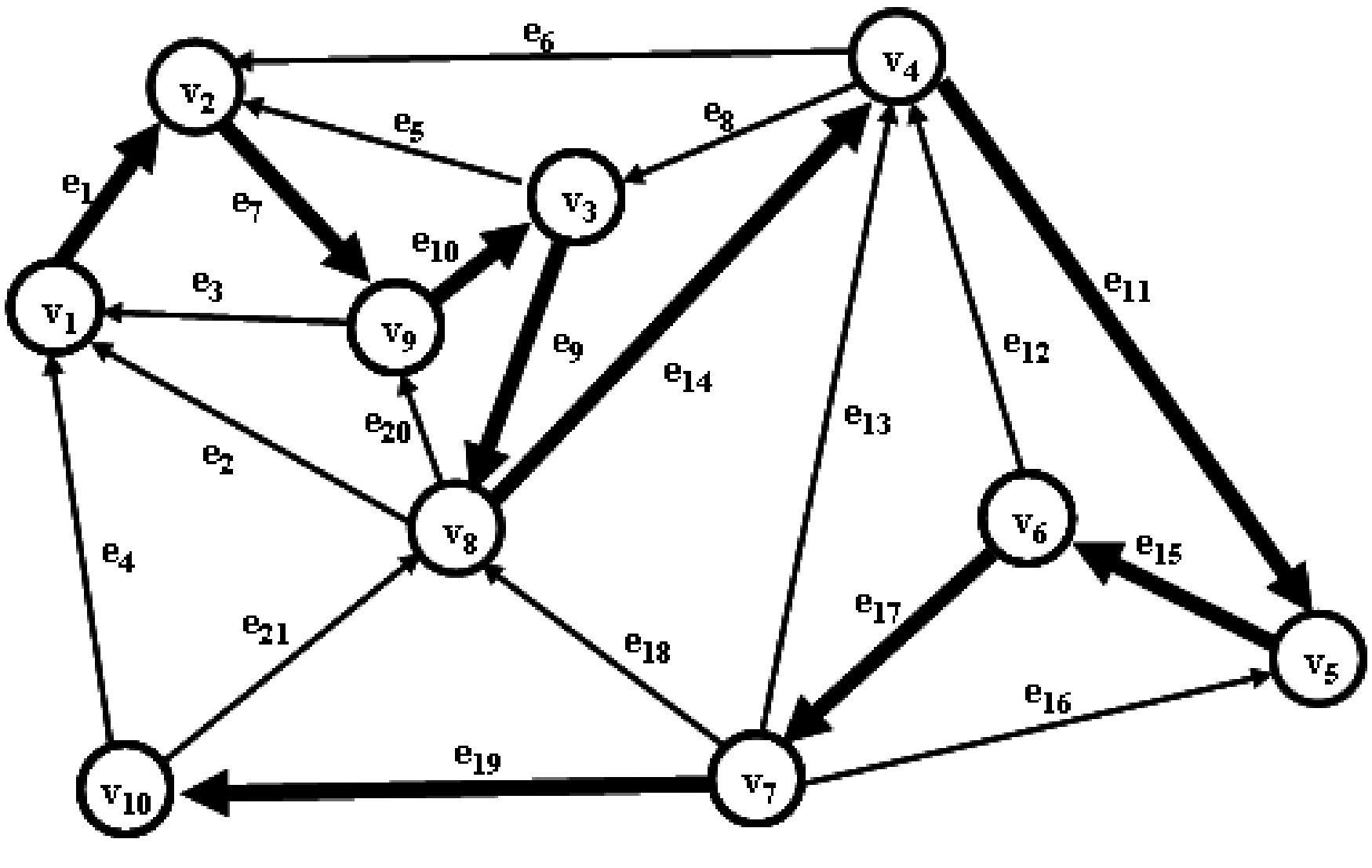

Рис. 5.10. Новое DFS-дерево и новый опорный цикл

Продолжим увеличивать длину опорного цикла. Для преобразования выбираем обратный маршрут характеризующийся строкой матрицы с хордой $e_{11}$, так как в ней имеется единственный элемент в подматрице $R_1$ и количество ветвей дерева не вошедших в опорный цикл больше, чем в других случаях. Производим преобразования в подматрице $\pi$. Объявляем хорду $e_{11}$ ребром дерева, принадлежащим опорному циклу, и включаем в подматрицу $\pi$, а ребро дерева $e_{13}$ объявляем хордой и исключаем из подматрицы $\pi$. Ветви дерева $e_{15}$ и $e_{17}$ включаем в опорный цикл, изменяя ориентацию. Тем самым преобразуем DFS-дерево графа $G_2$ (см. рис. 5.10).

Окончательно подматрица $\pi$ принимает вид:

| | ветви дерева (подматрица $\pi$) | | | | | | | | |
|---|---|---|---|---|---|---|---|---|---|
| | $\langle v_1,v_2\rangle$ | $\langle v_2,v_9\rangle$ | $\langle v_9,v_3\rangle$ | $\langle v_3,v_8\rangle$ | $\langle v_8,v_4\rangle$ | $\langle v_4,v_5\rangle$ | $\langle v_5 v_6\rangle$ | $\langle v_6,v_7\rangle$ | $\langle v_7,v_{10}\rangle$ |
| | $e_1$ | $e_7$ | $e_{10}$ | $e_9$ | $e_{14}$ | $e_{11}$ | $e_{15}$ | $e_{17}$ | $e_{19}$ |
| $e_2$ = | **1** | **1** | **1** | **1** | | | | | |
| $e_3$ = | **1** | **1** | | | | | | | |
| $e_4$ = | **1** | **1** | **1** | **1** | **1** | **1** | **1** | **1** | **1** |
| $e_5$ = | | **1** | **1** | **1** | | | | | |
| $e_6$ = | | **1** | **1** | **1** | **1** | | | | |
| $e_8$ = | | | | **1** | **1** | | | | |
| $e_{13}$ = | | | | | | **1** | **1** | **1** | |
| $e_{12}$ = | | | | | | **1** | **1** | | |
| $e_{16}$ = | | | | | | | **1** | **1** | |
| $e_{18}$ = | | | | | **1** | **1** | **1** | **1** | |
| $e_{21}$ = | | | | | **1** | **1** | **1** | **1** | **1** |
| $e_{20}$ = | | | **1** | **1** | | | | | |

В результате получаем множество, состоящее только из обратных рёбер.

Формируем множество обратных рёбер:

$\{e_2\} \to \langle v_8,v_1\rangle$; $\{e_3\} \to \langle v_9,v_1\rangle$; $\{e_5\} \to \langle v_3,v_2\rangle$;
$\{e_6\} \to \langle v_4,v_2\rangle$ ; $\{e_8\} \to \langle v_4,v_3\rangle$ ;
$\{e_{13}\} \to \langle v_7,v_4\rangle$; $\{e_{12}\} \to \langle v_6,v_4\rangle$;
$\{e_{16}\} \to \langle v_7,v_5\rangle$; $\{e_{18}\} \to \langle v_7,v_8\rangle$;
$\{e_{20}\} \to \langle v_8,v_9\rangle$; $\{e_{21}\} \to \langle v_{10},v_8\rangle$.

Будем формировать простые циклы графа, встраивая обратные рёбра в обручи (циклы).

$c_1 = \langle v_1,v_2,v_9,v_3,v_4,v_5,v_6,v_7,v_8,v_{10},v_1\rangle$;
$c_2 = \langle v_1,v_{10},v_8,v_7,v_6,v_5,v_4,v_3,v_9,v_2,v_1\rangle$

Вставим в цикл $c_1$ обратное ребро $\{e_2\} \to \langle v_8,v_1\rangle$.

$c_1$: $\langle v_1,v_2\rangle + \langle v_2,v_9\rangle + \langle v_9,v_3\rangle + \langle v_3,v_4\rangle + \langle v_4,v_5\rangle + \langle v_5,v_6\rangle + \langle v_6,v_7\rangle + \langle v_7,v_8\rangle +$
$+ \langle v_8,v_{10}\rangle + \langle v_{10},v_1\rangle + \langle v_8,v_1\rangle + \langle v_1,v_8\rangle = (\langle v_1,v_2\rangle + \langle v_2,v_9\rangle + \langle v_9,v_3\rangle + \langle v_3,v_4\rangle + \langle v_4,v_5\rangle + \langle v_5,v_6\rangle + \langle v_6,v_7\rangle + \langle v_7,v_8\rangle + \langle v_8,v_1\rangle) +$
$+ (\langle v_8,v_{10}\rangle + \langle v_{10},v_1\rangle + \langle v_1,v_8\rangle)$.

Образуются два новых обруча:

$c_1' = \langle v_1,v_2,v_9,v_3,v_4,v_5,v_6,v_7,v_8,v_1\rangle$;
$c_3 = \langle v_8,v_{10},v_1,v_8\rangle$.

Система обручей имеет вид:

$c_1' = \langle v_1,v_2,v_9,v_3,v_4,v_5,v_6,v_7,v_8,v_1\rangle$;
$c_2 = \langle v_1,v_{10},v_8,v_7,v_6,v_5,v_4,v_3,v_9,v_2,v_1\rangle$;
$c_3 = \langle v_8,v_{10},v_1,v_8\rangle$.

Продолжаем процесс введения обратных рёбер до полного их исчерпания. Окончательно система обручей имеет вид:

$c_1^{''''} = \langle v_6,v_7,v_8,v_4,v_6\rangle$;
$c_2^{''''} = \langle v_4,v_2,v_1,v_{10},v_7,v_4\rangle$;

$c_3 = \langle v_8, v_{10}, v_1, v_8 \rangle$;
$c_4 = \langle v_8, v_7, v_{10}, v_8 \rangle$;
$c_5 = \langle v_1, v_2, v_9, v_1 \rangle$;
$c_6 = \langle v_3, v_9, v_2, v_3 \rangle$;
$c_7 = \langle v_9, v_8, v_1, v_9 \rangle$;
$c_8 = \langle v_4, v_3, v_2, v_4 \rangle$;
$c_9 = \langle v_9, v_3, v_8, v_9 \rangle$;
$c_{10} = \langle v_8, v_3, v_4, v_8 \rangle$;
$c_{11} = \langle v_7, v_6, v_5, v_7 \rangle$;
$c_{12} = \langle v_4, v_5, v_6, v_4 \rangle$;
$c_{13} = \langle v_5, v_4, v_7, v_5 \rangle$.

В результате получен плоский топологический рисунок графа $G_2$ состоящий из простых циклов (см. рис. 5.11).

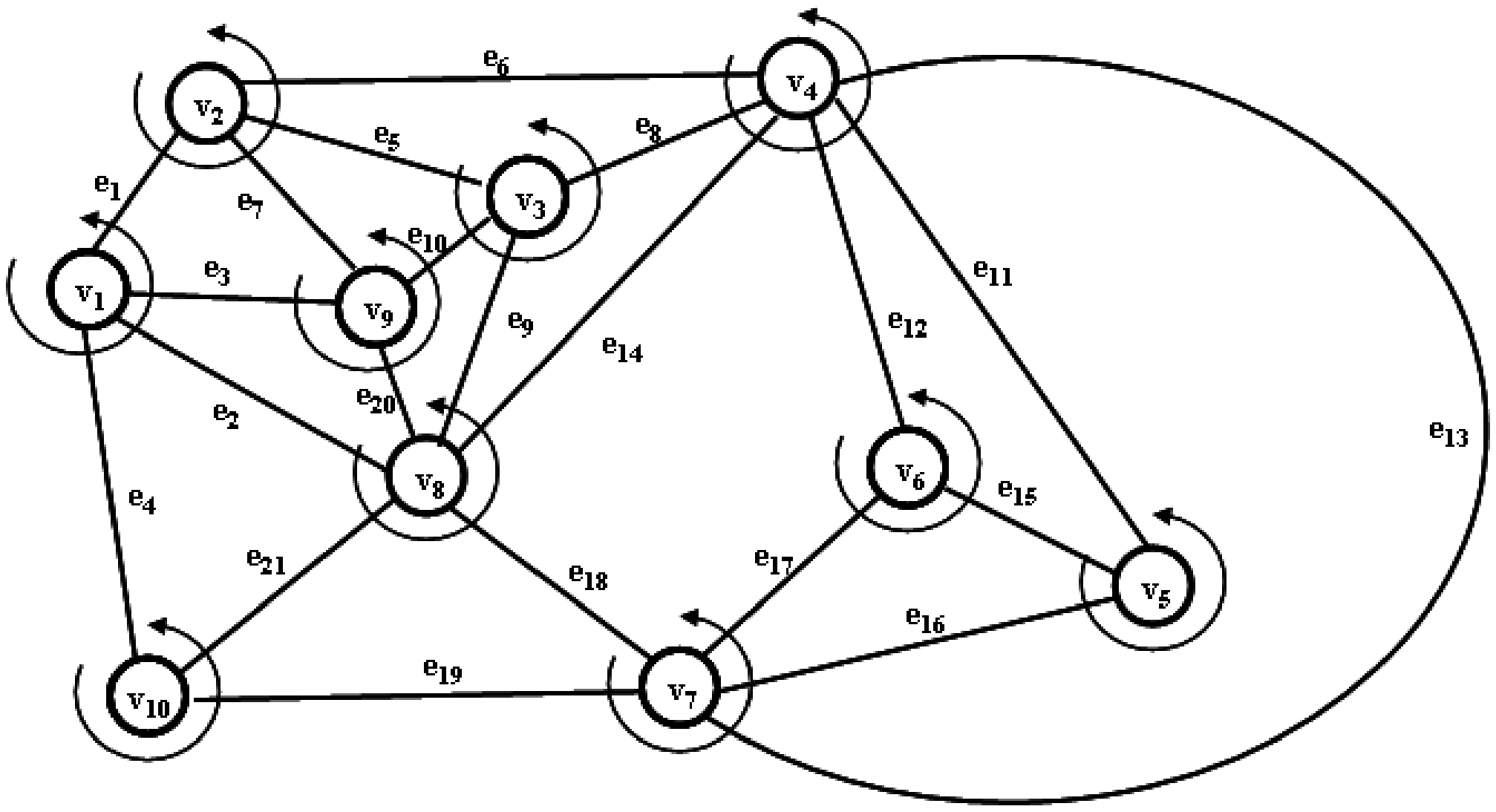

Рис. 5.11. Топологический рисунок графа $G_2$.

## 5.5. Обоснование алгоритма

Обоснование алгоритма основано на следствии теоремы Маклейна вытекающей из построения топологического рисунка графа.

**Следствие 5.1**. В плоском топологическом рисунке несепарабельного графа по ребрам маршрута соединяющего две вершины, принадлежащие простому циклу, проходят два простых цикла, образованные кольцевой суммой циклов принадлежащих инцидентным ребрам внутренних вершин.

Доказательство. Проведение такого маршрута разбивает выделенный цикл на две части. Тогда, исходя из теоремы Маклейна и как следствие построения топологического рисунка кольцевая сумма циклов проходящих по инцидентным ребрам внутренних вершин, для каждой части разбиения, формирует два простых цикла проходящих по ребрам маршрута. В случае отсутствия внутренних вершин, происходит кольцевое суммирование частей выделенного цикла и ребер маршрута.

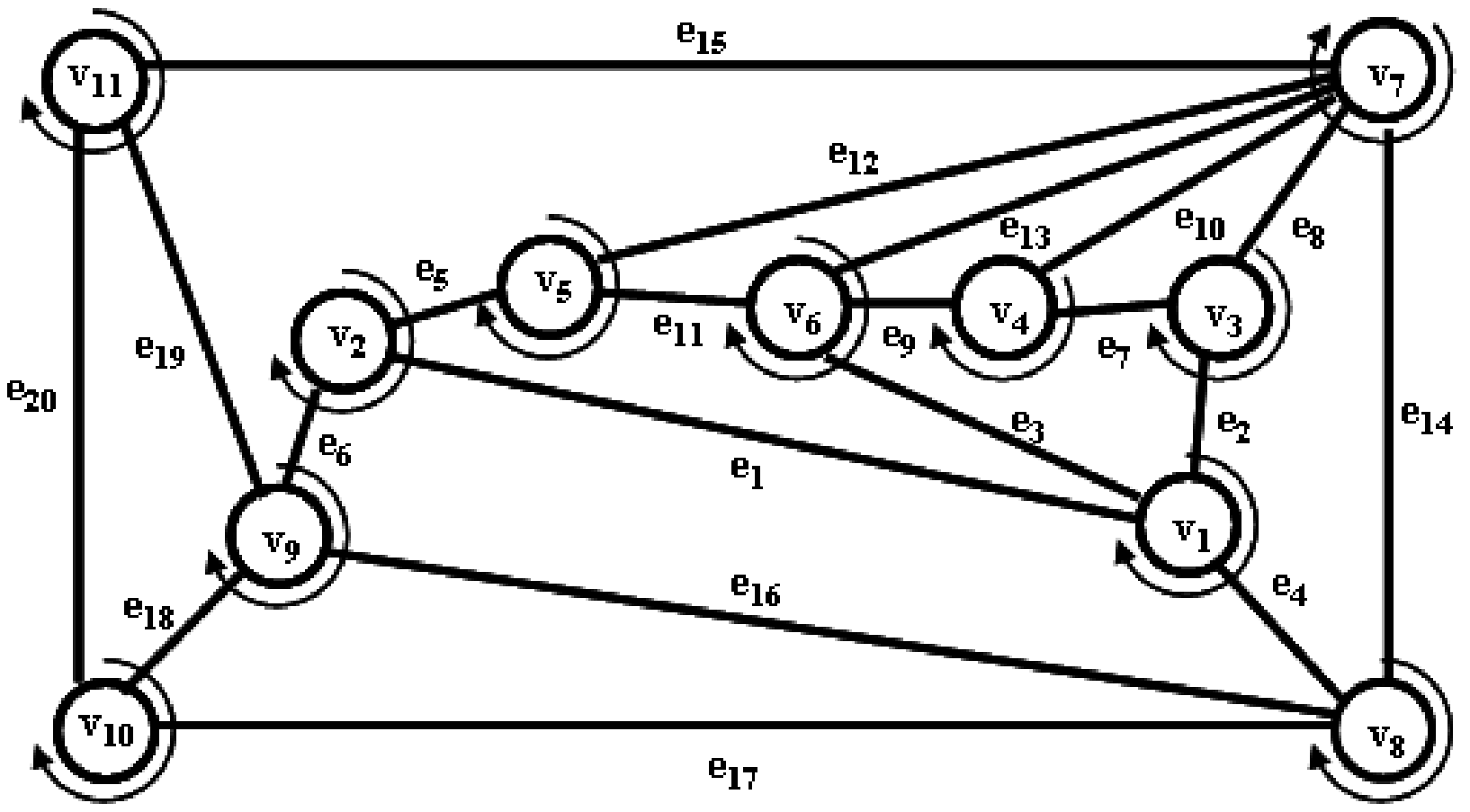

Рис. 5.12. Топологический рисунок плоского графа.

Например, для топологического рисунка графа представленного на рис. 5.12, вращение вершин можно описать диаграммой:

$v_1$: $e_3$ $e_2$ $e_4$ $e_1$
$v_2$: $e_6$ $e_5$ $e_1$
$v_3$: $e_8$ $e_2$ $e_7$
$v_4$: $e_{13}$ $e_7$ $e_9$
$v_5$: $e_{12}$ $e_{11}$ $e_5$
$v_6$: $e_{11}$ $e_{13}$ $e_9$ $e_3$
$v_7$: $e_{15}$ $e_{14}$ $e_8$ $e_{10}$ $e_{13}$ $e_{12}$
$v_8$: $e_4$ $e_{14}$ $e_{17}$ $e_{16}$
$v_9$: $e_{19}$ $e_6$ $e_{16}$ $e_{18}$
$v_{10}$: $e_{20}$ $e_{18}$ $e_{17}$
$v_{11}$: $e_{15}$ $e_{19}$ $e_{20}$

Множество простых циклов и обод:

$c_1 = \{e_1,e_3,e_5,e_{11}\}$; $c_2 = \{e_1,e_4,e_6,e_{16}\}$; $c_3 = \{e_2,e_3,e_7,e_9\}$;
$c_4 = \{e_2,e_4,e_8,e_{14}\}$; $c_5 = \{e_5,e_6,e_{12},e_{15},e_{19}\}$; $c_6 = \{e_7,e_8,e_{10}\}$;
$c_7 = \{e_9,e_{10},e_{13}\}$; $c_8 = \{e_{11},e_{12},e_{13}\}$;
$c_9 = \{e_{16},e_{17},e_{18}\}$; $c_{10} = \{e_{18},e_{19},e_{20}\}$; $c_0 = \{e_{14},e_{15},e_{17},e_{20}\}$.

Выделим цикл $\{e_{14},e_{15},e_{17},e_{20}\}$. Проведем разрезающий маршрут $\langle e_3,e_4,e_{13}\rangle$ вершины которого $x_7$ и $x_8$ принадлежат циклу. Тогда подмножество вершин $\{v_2,v_5,v_9\}$ характеризует внутренность левой части, а подмножество вершин $\{v_3,v_4\}$ характеризует внутренность правой части. В свою очередь инцидентные ребра левой части графа образуют подмножество $\{e_1,e_5,e_6,e_{11},e_{12},e_{16},e_{18},e_{19}\}$, а инцидентные ребра для правой части образуют подмножество $\{e_2,e_7,e_8,e_9,e_{10}\}$. Соответственно подмножество циклов $\{c_1,c_2,c_5,c_8,c_9,c_{10}\}$ характеризует левую часть, а подмножество $\{c_3,c_4,c_6,c_7\}$ характеризует правую часть. Тогда можно сформировать два цикла проходящих по ребрам маршрута и части выделенного цикла:

$c_1 \oplus c_2 \oplus c_5 \oplus c_8 \oplus c_9 \oplus c_{11} = \{e_1,e_3,e_5,e_{11}\} \oplus \{e_1,e_4,e_6,e_{16}\} \oplus \{e_5,e_6,e_{12},e_{15},e_{19}\} \oplus$

$\oplus \{e_{11},e_{12},e_{13}\} \oplus \{e_{16},e_{17},e_{18}\} \oplus \{e_{18},e_{19},e_{20}\} = \{e_3,e_4,e_{15},e_{13},e_{17},e_{20}\}$;

$c_3 \oplus c_4 \oplus c_6 \oplus c_7 = \{e_2,e_3,e_7,e_9\} \oplus \{e_2,e_4,e_8,e_{14}\} \oplus \{e_7,e_8,e_{10}\} \oplus \{e_9,e_{10},e_{13}\} =$

= {$e_3, e_4, e_{13}, e_{14}$}.

Следствие создает основу для создания рекурсивного процесса построения топологического рисунка графа. С этой целью выделяется простой цикл и проводится маршрут между его вершинами, разрезающий выделенный цикл на два цикла меньшей длины принадлежащими данному маршруту. Процесс определяется и управляется выделенным деревом графа и продолжается рекурсивно до появления простых неразделяемых циклов. В свою очередь выделенные циклы индуцируют (порождают) вращение его вершин и определяют плоский топологический рисунок графа.

### 5.6. Удлинения опорного цикла

Опишем алгоритм удлинения опорного цикла.

**шаг 1.** [Выбор максимального по длине цикла]. Строим матрицу фундаментальных циклов графа. Выбираем максимальный по длине цикл.

**шаг 2.** [Построение обратных маршрутов]. Формируем обратные маршруты графа. Для этого объединяем ветви дерева, не включенные в опорный цикл и соответствующую хорду.

**шаг 3.** [Проверка обратных маршрутов на замкнутость]. Проверяем обратные пути на замкнутость и соответствие концевых вершин опорному циклу. Если таких обратных путей нет, конец работы алгоритма.

**шаг 4.** [Построение нового опорного цикла]. Находим обратный маршрут, у которого количество ветвей дерева не вошедших в опорный цикл превышает количество вошедших. Одну из ветвей дерева входящих в опорный цикл объявляем хордой, а хорду данного обратного маршрута объявляем ветвью дерева входящей в опорный цикл. Включаем ветви дерева данного обратного маршрута не входящие в опорный цикл во множество ветвей дерева входящих в опорный цикл с возможной их переориентацией. В случае необходимости переориентируем обратные ветви, иначе – окончательный опорный цикл выбран; конец работы алгоритма.

**шаг 5.** [Построение новой матрицы фундаментальных циклов]. Строим новое DFS-дерево и новую матрицу фундаментальных циклов графа и идем на шаг 4.

Конец алгоритма.

В общем случае, количество операций для выделения дерева графа поиском в глубину определяется величиной равной рангу графа $\rho(\mathbf{G}) = \mathbf{n-1}$ и не может быть больше количества ребер в графе. Поэтому вычислительная сложность данного алгоритма может быть описана как $f_1(m) = O(m)$.

### 5.7. Алгоритм метода нитей

Маршрут начало и конец которого принадлежат различным вершинам, будем называть *нитью*. Тогда проведение маршрута (нити) разбивает цикл, в составе которого находятся вершины характеризующие начало и конец маршрута, на два цикла меньшей длины, при этом удовлетворяется критерий планарности Маклейна. Невозможность проведения нити означает, что вершины характеризующие начало и конец маршрута принадлежат не одному циклу, а разным. Данный процесс и является обоснованием применения метода нитей для проверки графа на планарность.

Опишем алгоритм проверки графа на планарность методом нитей [20].

**шаг 1.** [Построение DFS-дерева]. Методом поиска в глубину строим DFS-дерево графа.

**шаг 2.** [Выбор опорного цикла]. Строим матрицу фундаментальных циклов графа и выбираем опорный цикл.

**шаг 3.** [Формирование обратных путей]. Формируем обратные пути графа относительно выбранного опорного цикла. Для этого объединяем не включённые в опорный цикл ветви дерева и соответствующую хорду. Если концевые вершины сформированных обратных путей не принадлежат опорному циклу и имеются обратные пути в виде замкнутых маршрутов, то производим удлинение опорного цикла и идём на шаг 2.

**шаг 4.** [Построение блоков обратных путей]. Обратные пути, имеющие одинаковые вершины в качестве истока, объединяем в блоки. Ранжируем обратные маршруты в блоках.

**шаг 5.** [Включение обратных маршрутов в обручи графа]. Последовательно подключаем обратные маршруты в обручи графа согласно рангу. При подключении образуем новые обручи графа.

**шаг 6.** [Построение системы простых циклов и топологического рисунка плоского графа]. В результате построения, если подключены все обратные пути и обратные рёбра, получаем систему простых циклов графа. Строим диаграмму вращения вершин. Если имеются неподключенные обратные маршруты, то граф непланарен.

Конец алгоритма.

В свою очередь, количество последовательно подключаемых обратных маршрутов (нитей) не может быть больше, чем количество рёбер графа. Поэтому вычислительную сложность алгоритма включения обратных маршрутов можно определить как $f_2(m) = O(m)$.

Последовательное применение алгоритмов можно оценить как аддитивную сумму их вычислительных сложностей $f_1(m) + f_2(m) = O(m)$.

### 5.8. Оценка сложности метода

Исходной информацией для проведения расчетов служит связанный список вершин графа (аналог матрицы смежностей графа без учета нулевых элементов) и количество ребер *m*.

Процесс построения такого списка состоит из перечисления всех ребер графа представленных в виде $e_k = (v_i,v_j) \vee (v_j,v_i)$, $k \in \{1,\ldots,m\}$ с записью текущего номера ребра в ячейку соответствующую записи вершин $(v_i,v_j)$ для прямого ребра, и в ячейку соответствующую записи вершин $(v_j,v_i)$ для обратного ребра. Вычислительная сложность процесса может быть представлена как функция $f_1 = m$. В крайнем случае, $m = n(n\text{-}1)/2$ и тогда вычислительная сложность данного процесса определяется как $O(n^2)$.

Трудоемкость процесса выделения DFS-дерева можно определить как перечисление всех веощин $v \in$ V для неоткрытых вершин, и первое, что она делает – это помечает переданную в качестве параметра вершину. Но для каждой вершины осуществляется просмотр и проверка на смежность всех остальных вершин. Вычислительная сложность процесса может быть представлена как функция $f_2 = n \times n = n^2$. Итого $O(n^2)$, где $n$ - количество вершин в графе. Параллельно определяются хорды и ветви дерева.

Процесс построения фундаментальной подматрицы циклов $\pi$ относительно множества хорд H = $\{h_1,h_2,\ldots,h_{m-n+1}\}$ осуществляется процедурой поиска в ширину для каждой хорды (число которых равно цикломатическому числу графа) на ациклическом подграфе, состоящем только из ветвей дерева, и определением пути между двумя вершинами. Вычислительная сложность процесса может быть представлена функцией $f_3 = (m\text{-}n\text{+}1)n = mn - n^2 + n$. В наихудшем случае, $m = n(n\text{-}1)/2$, тогда $f_3 = (n^3 - 3n^2)/2 + n$. Итого $O(n^3)$.

Задача выделения опорного цикла состоит из операции определения количества ребер в каждом цикле и выборе максимального по длине цикла. Длина максимально возможного цикла может быть определена как ранг графа (n-1). Вычислительная сложность процесса может быть представлена функцией $f_4 = (m - n + 1)(n - 1) = m(n - 1) - n^2 + 2n - 1$. В крайнем случае, $m = n(n\text{-}1)/2$, тогда $f_4 = (n^3 - 4n^2 + 5\text{n} - 2)/2$. Итого $O(n^3)$.

Процесс разделения элементов подматрицы $\pi$ на ветви дерева, принадлежащие и не принадлежащие опорному циклу, касается всех элементов подматрицы циклов, число которых может быть определено как $(m\text{-}n\text{+}1)(n\text{-}1)$. Параллельно с этим определяется количество ветвей дерева принадлежащих и не принадлежащих опорному циклу, проходящему по данной хорде. Тогда вычислительную сложность процесса можно представить функцией $f_5 = (m\text{-}n\text{+}1)(n\text{-}1) = m(n - 1)\text{-}n^2\text{+}2n\text{-}1 = (\text{n}^3\text{-}4\text{n}^2\text{+}5\text{n-}2)/2$. Итого $O(n^3)$.

Увеличение длины опорного цикла приводит к преобразованию всех элементов подматрицы циклов $\pi$ и добавлению к ней элементов. Общее число преобразуемых элементов может быть определено как $(m\text{-}n\text{+}1)(n\text{-}1)$. Но так как операция увеличения длины опорного цикла производится для всех циклов, число которых равно количеству хорд графа, то выражение $(m\text{-}n\text{+}1)(n\text{-}1)$ нужно умножить на величину $(m\text{-}n\text{+}1)$. Тогда вычислительную сложность процесса можно представить функцией $f_6 = (m\text{-}n\text{+}1)(m\text{-}n\text{+}1)(n\text{-}1) = (n^5\text{-}7n^4\text{+}19n^3\text{-}25n^2\text{+}16n\text{-}4)/2$. Итого $O(n^5)$.

Формирование обратных ветвей и блоков обратных ветвей может быть произведено для множества хорд графа, количество которых равно числу m-n+1. Параллельно определяем общие вершины для обратных маршрутов входящих в блок. Тогда вычислительную сложность процесса можно представить функцией $f_7 = (m\text{-}n+1) = (n^2\text{-}3n+2)/2$. Итого $O(n^2)$.

Вычислительный процесс формирования множества простых циклов начинается с записи опорного цикла в виде последовательности ориентированных ребер для соответствующей хорды с последовательным порядком подключения входящих ветвей дерева. В результате формируется запись цикла в векторном виде. Для формирования другого цикла производится запись векторов цикла в обратном порядке. Очевидно, что функцию вычислительной сложности можно записать как $f_8=2(n\text{-}1)$. Итого $O(n)$.

Дальнейшее формирование простых циклов производится путем выбора единственного из вновь образованных циклов, содержащего концевые вершины обратного маршрута. Количество вновь образованных циклов не может превышать цикломатическое число графа и в наихудшем случае равно ему. Далее выбранный цикл разбивается на две части относительно концевых вершин с подключением соответствующих ориентированных обратных маршрутов. Число просматриваемых при этом ориентированных ребер можно определить как $2(n\text{-}1)$. Параллельно определяется количество обратных маршрутов еще не участвовавших в процессе формирования простых циклов. В конце расчета, если количество оставшихся обратных маршрутов не равно нулю, делается вывод о том, что исходный граф не планарный. Учитывая, что количество простых циклов равно цикломатическому числу графа, то вычислительную сложность процесса можно представить функцией $f_9 = (m\text{-}n+1)(m\text{-}n+1)(n\text{-}1) = (n^5\text{-}7n^4+19n^3\text{-}25n^2+16n\text{-}4)/2$. Итого $O(n^5)$.

Вычислительный процесс построения диаграммы вращения вершин, характеризующей топологический рисунок графа, заключается в последовательном просмотре всех ориентированных ребер в простом цикле или ободе и записью для каждой вершины смежной вершины соответствующего ориентированного ребра. Вычислительную сложность алгоритма можно представить функцией $f_{10} = (m\text{-}n+2)(n\text{-}1)=(n^3\text{-}4n^2+7n\text{-}4)/2$. Итого $O(n^3)$.

Общую вычислительную сложность алгоритма можно представить аддитивным сложением всех вычислительных сложностей процедур $f = \sum_{k=1}^{10} f_k$ . Данная функция позволяет оценить общую вычислительную сложность алгоритма как $O(n^5)$.

## Комментарии

Рассматривается модифицированный алгоритм проверки графа на планарность с одновременным построением математических структур для описания топологического рисунка плоского графа. Основой расчета является выделение DFS-дерева графа. Затем производится выделение опорного цикла, который производит разбиение поверхности $R^2$ на внутреннюю и внешнюю области. Выделенное дерево позволяет построить фундаментальную матрицу циклов. В свою очередь, выделение опорного цикла разбивает подматрицу $\pi$ фундаментальной матрицы циклов на две подматрицы: с ветвями дерева входящими, и с ветвями дерева не входящими в опорный цикл. Подматрица, состоящая только из ветвей дерева не вошедших в опорный цикл, позволяет формировать обратные маршруты. В случае существования обратных замкнутых маршрутов и обратных маршрутов с концевыми вершинами, не принадлежащими опорному циклу, алгоритмом удлинения опорного цикла строится новый опорный цикл, что порождает новое построение DFS-дерева графа и новую фундаментальную матрицу циклов. В случае отсутствия обратных замкнутых маршрутов и обратных маршрутов с концевыми вершинами, не принадлежащими опорному циклу, производится объединение обратных маршрутов в блоки и их ранжирование. Обратные маршруты последовательно просматриваются согласно их рангу, и в случае, если их концевые вершины принадлежат единственному обручу, производится размещение данного обратного маршрута с разбиением данного обруча на два обруча меньшей длины. Процесс продолжается рекурсивно до полного исчерпания списка обратных маршрутов. В этом случае граф планарен. Если список обратных маршрутов не исчерпан, то граф не планарен. Вычислительная сложность алгоритма оценивается как $O(n^5)$. Полученная система простых циклов графа индуцирует вращение вершин для описания топологического рисунка плоского графа. Топологический рисунок плоской части графа позволяет описывать процесс планаризации алгебраическими методами, не производя никаких геометрических построений на плоскости. Получение вращения вершин графа сразу решает две важнейшие задачи теории графов: задачу проверки графа на планарность и задачу построения топологического рисунка плоского графа.

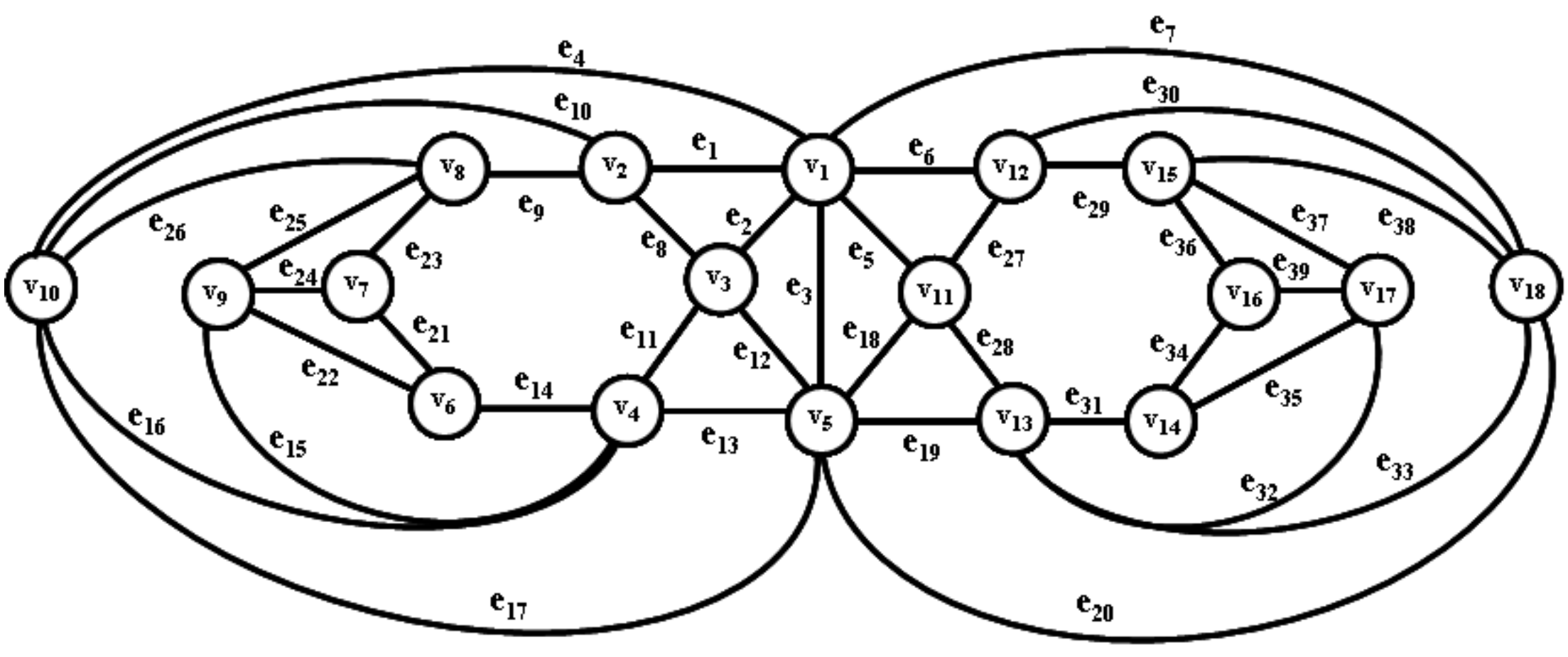

5.13. Топологический рисунок графа $G_3$.

Следует также заметить, что в общем случае, не всякий базис планарного графа состоит из изометрических циклов. Достаточно привести пример графа, у которого в базисе существуют простые, но не изометрические циклы. Рассмотрим следующий топологический рисунок графа $G_3$ (см. рис. 5.13).

В данном плоском рисунке (см. рис. 5.13) присутствуют не изометрические базисные циклы $\{e_8,e_9,e_{11},e_{14},e_{21},e_{23}\}$ и $\{e_{27},e_{28},e_{29},e_{31},e_{34},e_{36}\}$. Поэтому обод топологического рисунка состоит из трех простых циклов объединенных в квазицикл.

## Заключение

Рассматриваются вопросы создания математической модели для описания топологического рисунка графа на основе теории вращения вершин.

В первой главе на основе метрических свойств эвклидова пространства графа, введено понятие изометрического цикла графа. Рассмотрены основные свойства изометрических циклов графа. Указана необходимость построения минимальных маршрутов относительно всех ребер графа, для алгоритма выделения множества изометрических циклов графа. Описан способ получения множества изометрических циклов для подграфов определенных путем исключения ребер из полного графа.

Вторая глава посвящена вопросам применения методов алгебры структурных чисел в теории графов. Рассмотрены вопросы описания матроидов центральных разрезов и изометрических чисел методами алгебры структурных чисел. Представлен алгоритм «бегущая строка» для построения множества структурных чисел. Приведены примеры решения некоторых известных оптимизационных задач с его помощью.

В третьей главе, используя методы теории вращения вершин, построена математическая модель топологического рисунка графа. Показана связь между вращением вершин и построением индуцированной этим вращением системы простых циклов и обратно. Для оценки подмножества изометрических циклов вводится интегральная характеристика в виде функционала Маклейна, построенного на основе комбинаторных свойств системы циклов. Для выделения подмножества с заданным значением целевой функции представлен алгоритм наискорейшего спуска с применением методом алгебры структурных чисел.

Четвертая глава посвящена описанию модифицированного алгоритма Хопкрофта-Тарьяна для проверки графа на планарность, с целью его применимости для построения топологического рисунка графа. Для улучшения вычислительных возможностей представлен способ удлинения базового цикла графа.

В работе показано, что представленную математическую модель можно применить для решения задачи построения рисунка непланарного графа с учетом введения дополнительных вершин характеризующих пересечение рёбер. Естественно, что создание математической модели построения рисунка необходимо для дальнейшего развития и решения таких задач теории графов: проверка планарности графа; выделение максимально плоского суграфа; определение толщины графа; получение графа с минимальным количеством пересечений и т.д.

## СПИСОК ИСПОЛЬЗУЕМОЙ ЛИТЕРАТУРЫ

1. Алексеев В.В. Графы. Модели вычислений. Структуры данных. / В.В. Алексеев, В.А. Таланов - Нижний Новгород, 2005. – 307 с.
2. Асанов М.О. Дискретная математика: графы, матроиды, алгоритмы / М.О. Асанов, В.А. Баранский, В.В. Расин Ижевск: Регулярная и хаотическая динамика, 2001. – 288 с.
3. Апанович З.В. Средства для работы с графами большого объема: построение и оптимизация компоновочных планов / З.В. Апанович // Системная информатика: Сб.науч.тр. – Новосибирск: Изд-во СОРАН, 2006. – Вып. 10. Методы и модели современного программирования. – С. 7–58.
4. Апанович З.В. От рисования графов к визуализации информации. / З.В. Апанович - Новосибирск, РАН, 2007. – 24 с.
5. Базилевич Р.П. Декомпозиционные и топологические методы автоматизированного конструирования электронных устройств. Львов: Вища шк., 1981.
6. Беллерт С. Анализ и синтез электрических цепей методом структурных чисел / С. Беллерт, Г. Возняцки - М.: Мир, 1972. - 332 с.
7. Гэри М. Вычислительные машины и труднорешаемые задачи: Пер. с англ. / М. Гэри., Д. Джонсон - М.: Мир.- 1982. – 416 с.
8. Емеличев В.А. Многогранники, графы и оптимизация./ В.А. Емеличев., М.М. Ковалев, М.К. Кравцов. - М.: Наука, 1981.– 344 с.
9. Емеличев В.А. Лекции по теории графов / В.А. Емеличев, О.И. Мельников, В.И.Сарванов, Р.И. Тышкевич - М.: Наука. ГРФМЛ.-1990. – 384 с.
10. Зыков А.А. Теория конечных графов. / А.А. Зыков - Новосибирск: ГРФМЛ.- 1963.- 542 с.
11. Зыков А.А. Основы теории графов. / А.А. Зыков – М.: Наука, ГРФМЛ, 1987. – 384 с.
12. Касьянов В.Н. Графы в программировании: обработка, визуализация и применение. / В.Н. Касьянов, В.А. Евстигнеев - СПб.: БХВ-Петербург, 2003. – 1104 с.
13. Козин И.В., Фрагментарные модели для некоторых экстремальных задач на графах / Козин И.В., Полюга С.И. // Математические машины и системы. – 2014. – № 1. – С. 143 – 150.
14. Курапов С.В. Методы построения топологического рисунка графа. / С.В. Курапов, А.В.Толок // Автоматика и телемеханика - № 9. - 2013.- С.78-97.
15. Курапов С.В., Давидовский М.В. Проверка планарности и построение топологического рисунка плоского графа (поиском в глубину) //Прикладная дискретная математика. - 2016. № 2(32). С.100 - 114.

16. Курапов С.В. Топологические методы построения рисунка графа / С.В. Курапов, В.А. Чеченя //Науковий журнал Запорізького національного технічного університету: Радіоелектроніка, Інформатика, управління. - Запоріжжя: ЗНТУ, 2013, № 1(18). - С.72-81.
17. Курапов С.В., Козин И.В., Полюга С.И. Эволюционно-фрагментарный алгоритм нахождения максимального планарного суграфа // Прикладная дискретная математика. - 2015. № 3(29). С.74 - 82.
18. Курапов С.В., Толок А.В. “Построение топологического рисунка максимально плоского суграфа не планарного графа”, *Автомат. и телемех.*, 2018, № 5, 24–45с; *Autom. Remote Control*, **79**:5 (2018), 793–810
19. Курапов С.В., Давидовский М.В., Сгадов С.А. Выделение плоской части графа. Матроиды и алгебра структурных чисел // Компоненты и технологии. - 2018. № 9. с. 56 - 59.
20. Курапов С.В., Давидовский М.В., Толок А.В. Модифицированный алгоритм проверки планарности графа и построение топологического рисунка. Метод нитей. / Научная визуализация, 2018, том 10, номер 4, страницы 53 – 74.
21. Курапов С.В., Давидовский М.В., Толок А.В. Метод визуализации рисунка непланарного графа. / Научная визуализация, 2019, том 11, номер 2, страницы 126 - 142, DOI: 10.26583/sv.10.4.05
22. Курапов С.В., Давидовский М.В., Толок А.В. Генерация топологического рисунка плоской части непланарного графа. / Научная визуализация, 2020, том 12, номер 1, страницы 90 - 102, DOI: 10.26583/sv.12.1.08
23. Курапов С.В., Давидовский М.В., Клиценко А.А. Изометрические циклы графа // Вісник Запорізького національного університету: Збірник наук. статей. Фіз.-мат. науки. - № 1. - Запоріжжя: ЗНУ. - 2016.- С.121-137.
24. Курапов С.В., Давидовский М.В. Единичные и изометрические циклы в графе // Вісник Запорізького національного університету: Збірник наук. статей. Фіз.-мат. науки. - № 2. - Запоріжжя: ЗНУ. - 2017.- С.116-130.
25. Курош А.Г. Курс высшей алгебры (9-е издание). / А.Г. Курош М.: Наука, 1968 – 431 с.
26. Липский В. Комбинаторика для программистов: Пер. с польс. / В. Липский -М.:Мир.- 1988.- 213 с.
27. Мак-Лейн С. Комбинаторное условие для плоских графов / С. Мак-Лейн //В кн.: Кибернетический сборник. Новая серия. - 1970.-вып. 7.- С.68-77.
28. Мелихов А.Н. Применение графов для проектирования дискретных устройств. / А.Н.Мелихов, Л.С. Берштейн, В.М. Курейчик - М.: Наука. ГРФМЛ.-1974. – 304 с.

29. Михалевич В. С. Последовательные алгоритмы оптимизации и их применение. I, П.— Кибернетика, 1965, № 1, с. 45—56; № 2, с. 85—88.
30. Нечепуренко М.И. Алгоритмы и программы решения задач на графах и сетях / М.И. Нечепуренко, В.К. Попков, С.М. Майнагашев и др. - Новосибирск: Наука. Сиб. отд-ние.- 1990. – 515 с.
31. Пападимитриу Х. Комбинаторная оптимизация. Алгоритмы и сложность.: Пер. с англ. / Х. Пападимитриу, К. Стайглиц - М.: Мир.- 1985. – 512 с.
32. Перепелица В.А. Дискретная математика и моделирование в условиях неопределенности данных./ В.А. Перепелица, Ф.Б. Тебуева. – Издательство «АкадемияЕстествознания», 2007. – 151 с.
33. Пупырев С.Н., Тихонов А.В. Визуализация динамических графов для анализа сложных сетей // Моделирование и анализ информационный систем. – 2010. – № 1.
34. Раппопорт Л.И. Координатно-базисная система представления топологии электронных схем. / Л.И. Раппопорт, Б.Н. Мороговский, С.А. Поливцев // В кн.: Системы и средства автоматизации очистного и проходческого оборудования - М. :Недра. - 1980. - С.56-64.
35. Раппопорт Л.И. Векторная алгебра и проектирование топологии соединений / Л.И. Раппопорт, Б.Н. Мороговский, С.А. Поливцев // В кн.: Вопросы автоматизации проектирования интегральных схем. - Киев:ИК УССР.- 1976.- С.107-124.
36. Раппопорт Л.И. Векторная алгебра пересечений / Л.И. Раппопорт, Б.Н. Мороговский, С.А. Поливцев // В кн.: Многопроцессорные вычислительные структуры.- Таганрог.- 1982. - вып.2(11).- С.53-56.
37. Рейнгольд Э. Комбинаторные алгоритмы, теория и практика. / Э. Рейнгольд, Ю. Нивергельт, Н. Дер - М.: Мир.- 1980. – 480 с.
38. Рингель Г. Теорема о раскраске карт./ Г. Рингель - М.: Мир.- 1977. – 126 с.
39. Свами М. Графы, сети и алгоритмы: Пер. с англ. / М. Свами, К. Тхуласираман - М.: Мир.- 1984. – 455 с.
40. Трохименко Я. К. Метод обобщенных чисел и анализ линейных цепей, М., Изд-во «Советское радио», 1972, - 314 с.
41. Харари Ф. Теория графов. - пер. с англ. Козырева В.П. / под ред. Гаврилова В.Г. / Ф. Харари – М.: Мир. – 1973. – 300 с.
42. Хопкрофт Дж.Е., Тарьян Р.Е. Изоморфизм планарных графов / Дж.Е. Хопкрофт, Р.Е. Тарьян //В кн.: Кибернетический сборник. Новая серия.- 1975.-вып. 12.- С.39-61.
43. Щемелинин В.М. Задача оптимального представления графа электрической схемы / В.М. Щемелинин //В кн.: Микроэлектроника.-1975.-вып. 9.- С.253-261.

44. Di Battista G. Algorithms for Drawing Graphs: an Annotated Bibliography / G. Di Battista, P. Eades, R. Tamassia, I.G. Tollis // Computational Geometry, Theory and Applications. – 1994. – N 4. – P. 235–282.
45. Card S. K. Readings in Information Visualization: Using Vision to Think. / S. K. Card, J. D. Mackinlay, B. Shneiderman – San Francisco: Morgan Kaufmann, 1999.
46. Kavitha T., Liebchen C., Mehlhorn K., Dimitrios M., Romeo Rizzi, Ueckerdt T., Katharina A. Cycle Bases in Graphs Characterization, Algorithms, Complexity, and Applications // Comput. Sci. Rev. 2009. V. 3. P. 199–243.
47. Tamassia Roberto. Handbook of Graph Drawing and Visualization / Roberto Tamassia Chapman and Hall/CRC. – 2013.- 844 p.
48. S. V. Kurapov and A. V. Tolok. Construction of a Topological Drawing of the Most Planar Subgraph of the Non-planar Graph / AUTOMATION AND REMOTE CONTROL Vol. 79. - No. 5. - 2018. – 793-810p.